\pdfoutput=1

\makeatletter
\IfFileExists{./numapde-latex/numapde-packages.sty}{\providecommand*{\input@path}{}\edef\input@path{{./numapde-latex/}\input@path}}{}
\makeatother

\RequirePackage{algpseudocode}
\documentclass{numapde-preprint}

\usepackage{numapde-semantic}
\usepackage{numapde-style}
\RequirePackage{tikz}
\usepackage{numapde-local}

\hypersetup{
	pdftitle={A Manifold of Planar Triangular Meshes with Complete Riemannian Metric},
	pdfauthor={Roland Herzog, Estefanía Loayza-Romero},
	pdfkeywords={triangular meshes, Riemannian manifold, complete Riemannian metric, geodesic curves, mesh deformation}
}

\title{A Manifold of Planar Triangular Meshes with Complete Riemannian Metric}
\subtitle{}
\shorttitle{Planar Triangular Meshes with Complete Riemannian Metric}

\author{Roland Herzog\thanks{Interdisciplinary Center for Scientific Computing, Heidelberg University, 69120 Heidelberg, Germany (\email{roland.herzog@iwr.uni-heidelberg.de}, \url{https://scoop.iwr.uni-heidelberg.de}, \orcid{0000-0003-2164-6575}).}
\and
Estefanía Loayza-Romero\thanks{Institute for Analysis and Numerics,  University of Münster,   48149 Münster,  Germany (\email{estefania.loayza-romero@uni-muenster.de}, \url{https://www.uni-muenster.de/AMM/num/wirth/people/Loayza/index.html}, \orcid{0000-0001-7919-9259}).}}
\shortauthor{R. Herzog and E. Loayza-Romero}

\dedication{}

\IfFileExists{numapde-uncolorizeHyperrefIfChanges.sty}{\RequirePackage{numapde-uncolorizeHyperrefIfChanges}}{}

\begin{document}
\maketitle

\begin{abstract}
Shape spaces are fundamental in a variety of applications including image registration, morphing, matching, interpolation, and shape optimization.
In this work, we consider two-dimensional shapes represented by triangular meshes of a given connectivity.
We show that the collection of admissible configurations representable by such meshes form a smooth manifold.
For this manifold of planar triangular meshes we propose a geodesically complete Riemannian metric.
It is a distinguishing feature of this metric that it preserves the mesh connectivity and prevents the mesh from degrading along geodesic curves.
We detail a symplectic numerical integrator for the geodesic equation in its Hamiltonian formulation.
Numerical experiments show that the proposed metric keeps the cell aspect ratios bounded away from zero and thus avoids mesh degradation along arbitrarily long geodesic curves.\end{abstract}

\begin{keywords}
triangular meshes, Riemannian manifold, complete Riemannian metric, geodesic curves, mesh deformation\end{keywords}

\begin{AMS}
\href{https://mathscinet.ams.org/msc/msc2010.html?t=53Z50}{53Z50}, \href{https://mathscinet.ams.org/msc/msc2010.html?t=57Z25}{57Z25}, \href{https://mathscinet.ams.org/msc/msc2010.html?t=53C22}{53C22}
\end{AMS}

\section{Introduction}
\label{section:introduction}

The notion of shapes is essential for the study of many problems in modeling, computer graphics, geometry processing, computer vision and other fields.
Regardless of the field of application, shapes can be described in many different ways.
We mention sets of points with a connectivity structure \cite{Kendall1984,Bookstein1986}, closed curves describing the boundary of planar shapes \cite{KlassenSrivastavaMioShantanu2004,MichorMumford2006} or surfaces in higher dimension \cite{BauerHarmsMichor2011,BauerHarmsMichor2012}, characteristic functions \cite{DelfourZolesio2011,Zolesio2007}.
For a general overview on shape representations, we direct the reader to \cite{Younes2010}.

In some of these references, the collection of shapes in the respective class is endowed with the structure of a Riemannian manifold.
For instance, the pioneering work of
\makeatletter
\ltx@ifclassloaded{mcom-l}{Kendall~\cite{Kendall1984}}{\cite{Kendall1984}}
\makeatother
defines shapes in terms of a set of points called \emph{landmarks} together with their connectivity.
The author considers configurations modulo translations, rotations, and scaling.
In \cite{MichorMumford2006}, the authors define shapes as embeddings of the unit circle in the plane modulo reparametrizations, which constitutes an infinite-dimensional manifold.
For an easily accessible overview of shape manifolds and their metrics, we refer the reader to \cite{Younes2012}.

In this paper, we are considering two-dimensional shapes represented by triangular meshes, see for instance \Cref{fig:triangulation}.
Mesh representations are frequently used in image registration and image morphing applications; see, \eg, \cite{AlexaCohenOrLevin2000,Alexa2002,BaghaieYuDsouza2014}.
Moreover, meshes frequently appear when partial differential equations are involved, notably in PDE-constrained shape optimization problems; see, \eg, \cite{MohammadiPironneau:2010:1,Paganini:2016:1,SchulzSiebenbornWelker2015:2}.

Although PDE-constrained shape optimization will not play a role in this paper, we draw our motivation from the following two observations.
First, the prevalent approach for the evaluation of shape derivatives is of the optimize--then--discretize category.
Without going into details here, the application of the continuous shape derivative on discrete meshes can cause additional errors to appear; see the discussion in \cite{Berggren2010} for details.
Second, when transitioning from one mesh to another, one frequently allows each mesh vertex to move along a straight line with an assigned velocity (a Euclidean geodesic).
This almost always causes the mesh quality to deteriorate and the mesh to become degenerate in finite time; compare \Cref{fig:euclidean_geodesic}.
Both effects represent severe obstacles in computational shape optimization, and many approaches have been proposed to mitigate mesh deterioration; see for instance the overview provided in \cite{EtlingHerzogLoayzaWachsmuth2018:1}.

In this paper, we are proposing to represent shapes as the collection of all admissible vertex positions of a given reference mesh.
While the question whether a given mesh is admissible or not (see \cref{fig:admissible_and_inadmissible_meshes}) can be answered intuitively, the mathematical description of \emph{all} admissible meshes is non-trivial.
To the best of our knowledge, this set has not been thoroughly studied.
Clearly, it does not bear the structure of a vector space, although it is a subset of the vector space of all possible vertex coordinates in the plane.
Indeed, using the language of simplicial complexes, we prove that this set turns out to be a finite dimensional manifold, which we endow with a suitable Riemannian metric.
The metric on this \emph{manifold of planar triangular meshes} is designed to be \emph{geodesically complete}.
This allows us to achieve mesh deformations along geodesic curves with any given initial velocity and for arbitrarily long times while maintaining non-generate meshes with a given connectivity.
We also observe numerically that cell aspect ratios remain bounded away from zero.
Notice that, in contrast to other shape manifolds, \eg, \cite{Kendall1984,Younes1998}, the manifold of planar triangular meshes does \emph{not} identify meshes which are obtained by rotations, translations, and scalings of each other.
This is motivated by our goal of using this space for PDE-constrained shape optimization problems, where these identifications cannot be made due to space dependent data.
First results in this direction can be found in \cite{HerzogLoayzaRomero:2021:1}.

To the best of our knowledge, the literature in which meshes are considered points on a shape manifold is still scarce.
We are aware of a series of contributions considering meshes representing the surface of three-dimensional objects instead of flat shapes.
In their seminal paper \cite{KilianMitraPottmann2007}, the authors consider two Riemannian metrics on this space based on isometric and rigid deformations.
In \cite{HeerenRumpfWardetzkyWirth2012,HeerenRumpfSchroederWardetzkyWirth2014}, the authors propose a Riemannian metric for the space of discrete thin shells, which takes into account viscous dissipation in terms of membrane and bending energies to reflect physical behavior.
The authors in \cite{YangYangYi-JunPottmannMitra2011} provide a framework to characterize spaces of quad meshes implicitly prescribed by a collection of non-linear constraints.
They also propose ways to explore this shape space using tangent vectors and quadratically parametrized osculant surfaces.
In \cite{AmentaRojas2020}, the authors consider the embedding of a triangular surface mesh into $\R^3$ up to translations, rotations and scale, using its vector of dihedral angles.
Finally, \cite{LiuShiDinovMio2010} construct shape spaces equipped with Riemannian metrics measuring how costly it is to interpolate two shapes through elastic deformations.
Their representation of shapes is based on the discrete exterior derivative, or coboundary operator, of parametrizations over a finite simplicial complex.

In contrast to all of the above, our emphasis is on providing a rigorous differential geometric framework for the admissible vertex configurations of planar triangular meshes and a geodesically complete Riemannian metric.
The latter avoids mesh degradation while still allowing large deformations without jeopardizing mesh quality.
Our work has direct applications in shape optimization and it allows us to formulate and pursue a novel and fully discrete approach to PDE-constrained shape optimization.
This will be discussed elsewhere.

This manuscript is organized as follows.
In \cref{section:motivation_1D} we present a brief motivation by considering the space of point configurations on the real line and equipping it with a complete Riemannian metric.
The formal presentation of the manifold of planar triangular meshes and its principal properties can be found in \Cref{section:discrete_shape_manifold}.
In \cref{section:complete_Riemannian_metric} we propose a Riemannian metric and prove its completeness.
Since closed form solutions of the ensuing geodesic equations are not available, we consider in \cref{section:geodesic_equation} the numerical integration of the geodesic in Hamiltonian formulation by the symplectic Störmer--Verlet scheme.
Finally, we present numerical experiments in \cref{section:numerical_experiments} which illustrate the behavior of mesh transformations along the geodesics associated with the proposed metric.
The conditions which define admissible vertex configurations are formulated in the language of abstract and geometric simplicial complexes.
The reader will find relevant background material in the appendix, along with auxiliary results.
We assume that the reader has a certain familiarity with basic notions of differential geometry.
For an accessible introduction into this topic we refer to \cite{Lee:2018:1, DoCarmo1992}.

\section{Motivation in 1D}
\label{section:motivation_1D}

In order to illustrate and motivate the concepts developed in this manuscript for two-dimensional triangular meshes, we briefly sketch a one-dimensional analogue first.
In this simplified model, we consider a number $N_V \ge 2$ of points arranged along the real line.
The variable point coordinates are denoted by $q^i \in \R$, $i = 1, \ldots, N_V$.
We are interested in all admissible point configurations, given by
\begin{equation}
	\label{eq:configuration_manifold_1D}
	\cM_+
	\coloneqq
	\setDef[auto]{Q = (q^1, \ldots, q^{N_V})^\transp \in \R^{N_V}}{q^1 < q^2 < \cdots < q^{N_V}}
	.
\end{equation}
The only condition for non-degeneracy is that the points maintain their ordering and positive distances from each other.
Moreover, we can characterize the boundary of this set as follows:
\begin{equation}
	\label{eq:boundary_configuration_manifold_1D}
	\partial \cM_+
	=
	\bigcup_{i=1}^{N_V-1} \setDef[auto]{Q = (q^1, \ldots, q^{N_V})^\transp \in \R^{N_V}}{q^1 \le \cdots \le q^i = q^{i+1} \le \cdots \le q^{N_V}}.
\end{equation}
Clearly, $\cM_+$ is an open submanifold of $\R^{N_V}$.
When endowed with the Euclidean Riemannian metric, $\cM_+$ is not geodesically complete.
In fact, geodesics in this metric are the curves $\gamma(t) = Q^0 + t \, V^0$ with initial positions $Q^0 \in \cM_+$ and initial velocity $V^0 \in \R^{N_V}$.
It is easy to see that Euclidean geodesics exist for all $t \in \R$ if and only if the components of the initial velocity are all identical, \ie, $V^0 = \kappa \, (1, 1, \ldots, 1)^\transp$ for some $\kappa \in \R$.

Our interest is to establish a geodesically complete Riemannian metric for $\cM_+$, which ensures that all geodesics exist for all $t \in \R$.
The key component is to find a function $f \colon \cM_+ \to \R$ which is proper, \ie, for which the preimages $f^{-1}(K)$ of compact sets $K \subset \R$ are compact in $\cM_+$.
One particular choice of such a function is
\begin{equation}
	\label{eq:proper_function_1D}
	f(Q)
	\coloneqq
	\sum_{i=1}^{N_V-1} \frac{1}{q^{i+1} - q^i} + \frac{\beta}{2} \abs{Q - \Qref}_2^2
\end{equation}
with $\beta > 0$.
Note that this function goes to infinity when $Q \to \partial \cM_+$.
A theorem of \cite{Gordon1973} then shows that the augmented Riemannian metric represented by
\begin{equation}
	\label{eq:augmented_metric_1D}
	g_{ab}
	=
	\delta_a^b + \frac{\partial f}{\partial q^a} \, \frac{\partial f}{\partial q^b}
\end{equation}
is geodesically complete.
Here $\delta_a^b$ is the Kronecker delta symbol representing the Euclidean metric.
Note that a Riemannian metric is nothing but a smoothly varying inner product in each point of $\cM_+ \subset \R^{N_V}$.
\Cref{eq:augmented_metric_1D} tells us that the matrix representing this inner product is the identity matrix plus a rank-one perturbation by the vector of partial derivatives of $f$.
The first-order derivative of the augmentation function $f$ given in \eqref{eq:proper_function_1D} with respect to $q^i$ is given by the following expression:
\begin{equation*}
	\frac{\partial f}{\partial q^i}
		=
		\frac{1}{q^{i+1}-q^i}- \frac{1}{q^i-q^{i-1}} + \beta \, (q^i-q_{\textup{ref}}^i)
	.
\end{equation*}
Thus, some of the metric coefficients in \eqref{eq:augmented_metric_1D} go to $\infty$ as $q^i \to q^{i+1}$ for some index~$i$.
	This renders points on the boundary $\partial \cM_+$ infinitely far away and resembles the geodesically complete metric of the Poincaré ball model, see, \eg, \cite[Thm.~3.7, p.~62]{Lee:2018:1}.

Unfortunately and, perhaps, surprisingly, even in this one-dimensional setting, the structure of the geodesic equation associated with this Riemannian metric is little intuitive.
To start, let us recall the expression of the geodesic equation (\cf \cite[Eq.~(4.16), p.~103]{Lee:2018:1}):
\begin{equation*}
	\ddot{q}^c + \sum_{a,b = 1}^2 \Gamma_{ab}^c \, \dot{q}^a \, \dot{q}^b = 0, \quad c = 1,2,
\end{equation*}
where $\Gamma_{ab}^c$ are the Christoffel symbols given by
\begin{equation*}
	\Gamma_{ab}^c = \frac12 \sum_{d=1}^2g^{cd}\paren[auto](){%
	\frac{\partial g_{da}}{\partial q^b}
	+
	\frac{\partial g_{db}}{\partial q^a}
	-
	\frac{\partial g_{ab}}{\partial q^d}
	}
	.
\end{equation*}
Note that $g^{cd}$ denotes one of the components of the inverse of the matrix representing the metric.
A formal derivation of the Christoffel symbols can be found for example in \cite[Cor.~5.11, p.~123]{Lee:2018:1}.

Let us consider the simplest case, \ie, $N_V = 2$ and the marginal choice $\beta = 0$.
After some computations involving the first- and second-order derivatives of the function $f$ given in~\eqref{eq:proper_function_1D}, which have been omitted, we find that the Christoffel symbols are given by:
\begin{align*}
	\begin{bmatrix}
		\Gamma^1_{11} & \Gamma^1_{12} \\ \Gamma^1_{21} & \Gamma^1_{22}
	\end{bmatrix}
	&
	=
	\frac{2}{(q^2-q^1)^5 + 2 \, (q^2-q^1)}
	\begin{bmatrix}
		1 & -1 \\ -1 & 1
	\end{bmatrix}
	,
	\\
	\begin{bmatrix}
		\Gamma^2_{11} & \Gamma^2_{12} \\ \Gamma^2_{21} & \Gamma^2_{22}
	\end{bmatrix}
	&
	=
	\frac{2}{(q^2-q^1)^5 + 2 \, (q^2-q^2)}
	\begin{bmatrix}
		-1 & 1 \\ 1 & -1
	\end{bmatrix}
	.
\end{align*}
Therefore the geodesic equation is expressed by the following system of  second-order differential equations
\begin{equation}
	\label{eq:1D_geodesic_standard_coordinates}
	\begin{aligned}
		\ddot{q}^1 + \frac{2}{(q^2-q^1)^5 + 2 \, (q^2-q^1)} (\dot{q}^2-\dot{q}^1)^2
		&
		=
		0
		,
		\\
		\ddot{q}^2 - \frac{2}{(q^2-q^1)^5 + 2 \, (q^2-q^1)} (\dot{q}^2-\dot{q}^1)^2
		&
		=
		0
		.
	\end{aligned}
\end{equation}
One can infer from adding the two components of \eqref{eq:1D_geodesic_standard_coordinates} that the midpoint of $q^1$ and $q^2$ moves with constant speed.
Indeed, utilizing the change of variables $\varepsilon = q^2-q^1$ and $\mu = q^1 + q^2$ allows us to re-state the geodesic equation in more compact form as follows,
\begin{equation*}
	\ddot \mu
	=
	0
	\quad
	\text{and}
	\quad
	\ddot \varepsilon - 2 \, \frac{2}{\varepsilon^5 + 2 \varepsilon } \, [\dot \varepsilon]^2
	=
	0
	.
\end{equation*}

\begin{figure}[htb]
	\begin{subfigure}{0.45\textwidth}
		\centering
		\includegraphics[width=1.0\linewidth]{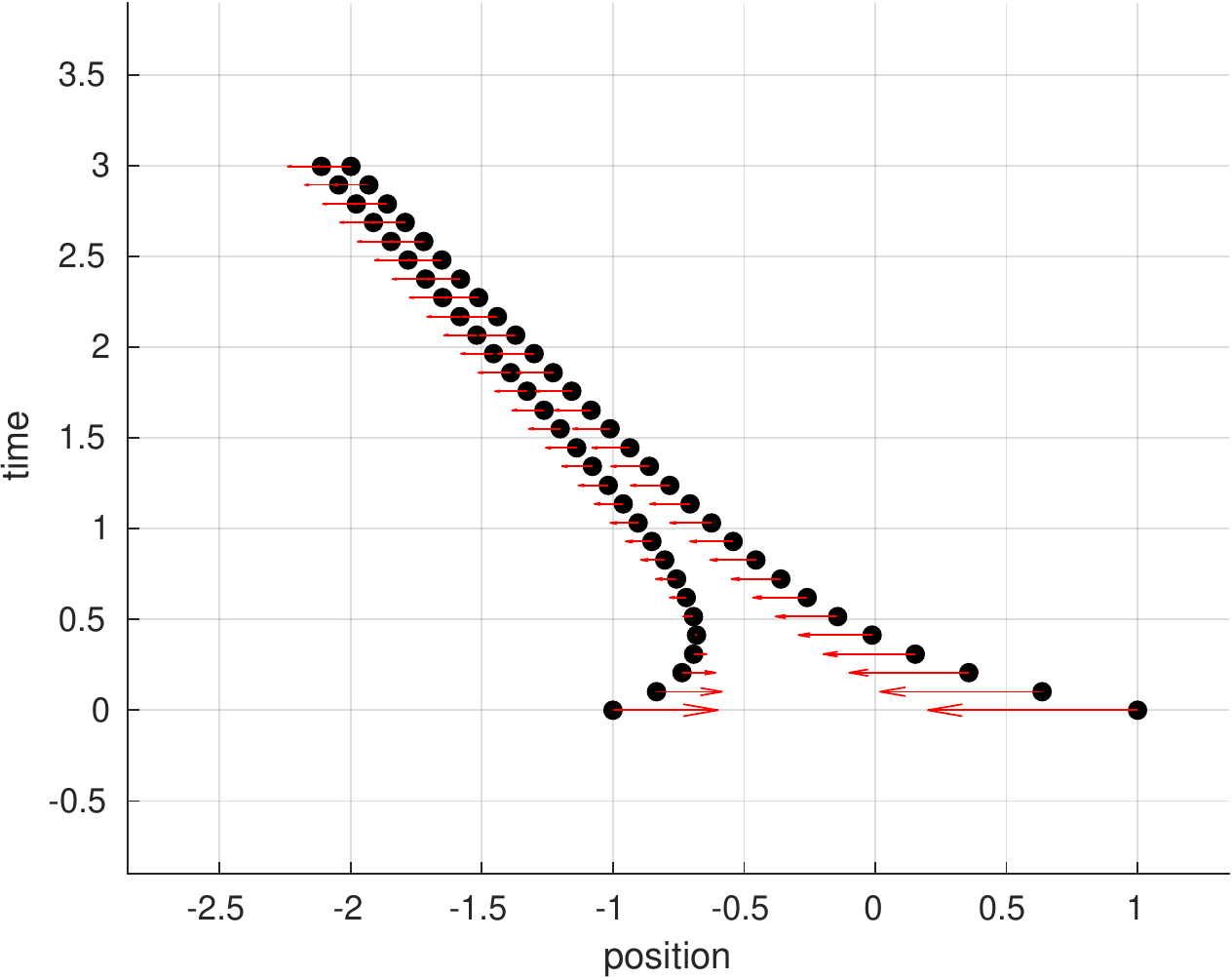}
		\caption{Example geodesic on $\cM_+$ with $N_V = 2$.}
	\end{subfigure}
	\hfill
	\begin{subfigure}{0.45\textwidth}
		\centering
		\includegraphics[width=1.0\linewidth]{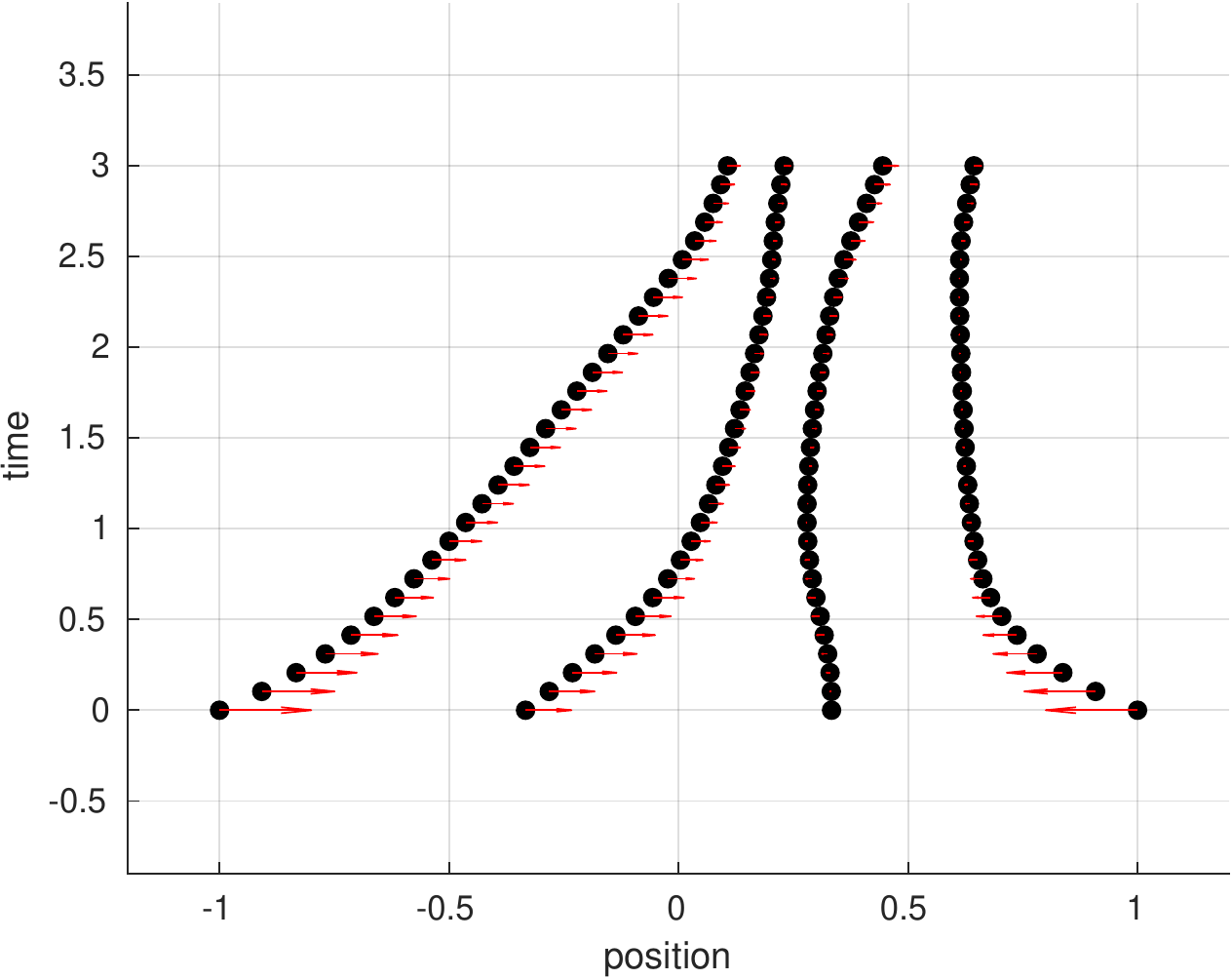}
		\caption{Example geodesic on $\cM_+$ with $N_V = 4$.}
	\end{subfigure}
	\caption{Snapshots of geodesics on the 1D configuration space~$\cM_+$ from \eqref{eq:configuration_manifold_1D} with the metric \eqref{eq:augmented_metric_1D} and $\beta = 1$. Velocities are shown as arrows with their Euclidean lengths.}
	\label{fig:geodesic_1D}
\end{figure}

Some trajectories of the geodesic equation for various initial conditions and numbers of points and with $\beta = 1$, obtained using numerical integration, are shown in \Cref{fig:geodesic_1D}.
It can clearly be seen that the goal of maintaining non-degenerate configurations is achieved for all $t \in \R$.

The transition from the simple one-dimensional setting to two-dimensional meshes is non-trivial.
First of all, the conditions which describe the admissible vertex configurations, and thus the manifold of planar triangular meshes, require some thought.
Second, the design of a proper function~$f$, and the proof of its properness, are significantly more complex.

\section{Manifold of Planar Triangular Meshes}
\label{section:discrete_shape_manifold}

As was already mentioned, we consider shapes represented by triangular meshes in $\R^2$ like those frequently used, \eg, in finite element computations.
Informally, a triangular mesh is a finite collection of non-degenerate triangles satisfying the condition that the intersection of any two triangles is either empty, a common edge, or a common vertex; see, \eg, \cite[Ch.~3]{QuarteroniValli1994:1}.
In the language of simplicial complexes, such meshes are precisely the \emph{pure simplicial $2$-complexes}.

\begin{figure}[htb]
	\centering
	\includegraphics[width=0.25\linewidth]{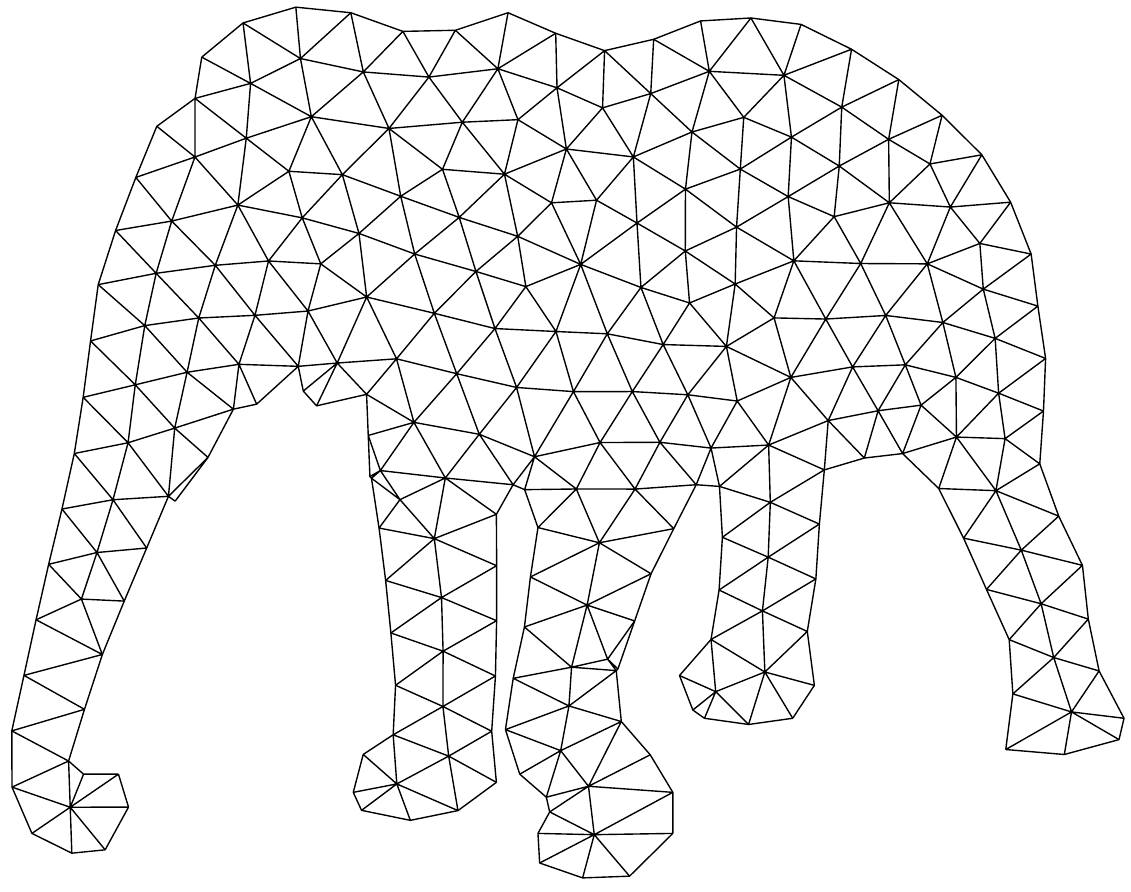}
	\includegraphics[width=0.25\linewidth]{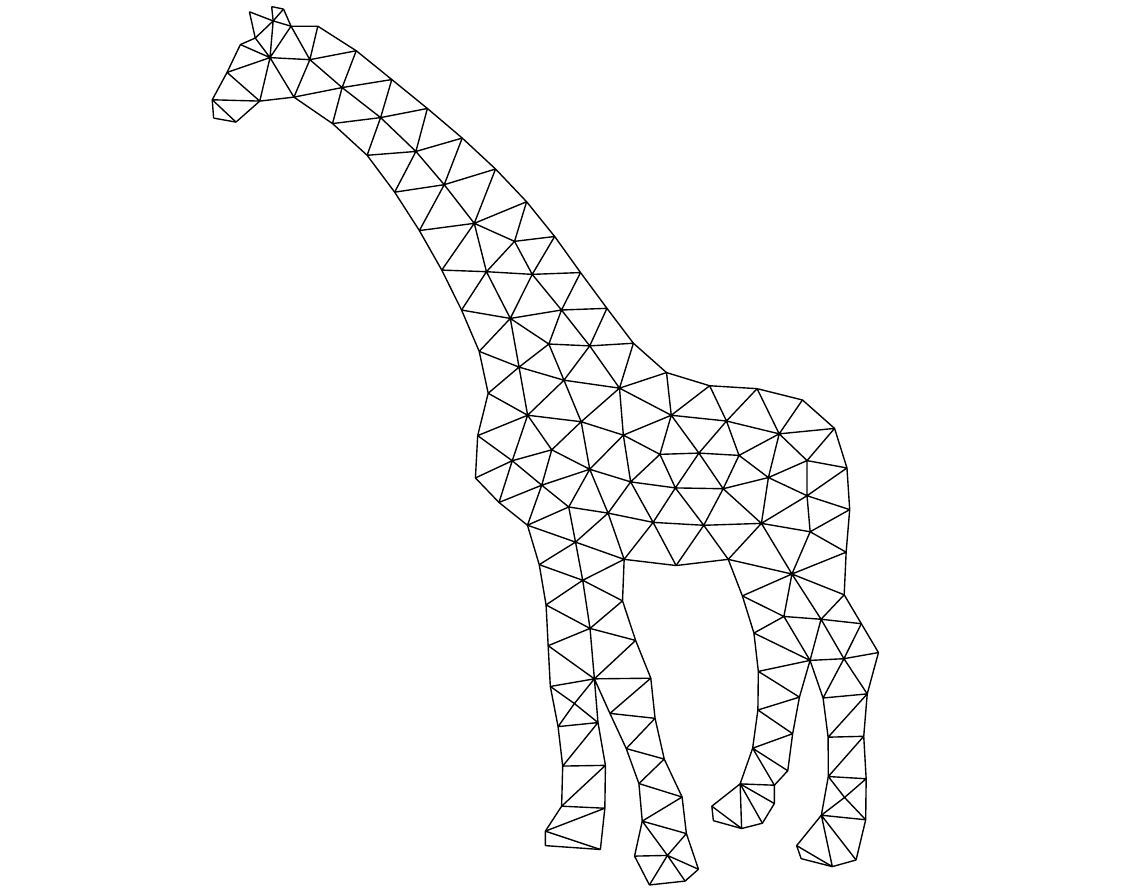}
	\caption{Two examples of triangular meshes, similar to \cite[Fig.~4]{AlexaCohenOrLevin2000}.}
	\label{fig:triangulation}
\end{figure}

In fact, we restrict the discussion to meshes with the property that any two distinct triangles $T_A$, $T_B$ are connected by a finite sequence of triangles $T_A = T_0, T_1, \ldots, T_n = T_B$ such that $T_i \cap T_{i+1}$ is a common edge of both for all $i = 0, \ldots, n-1$.
We refer to this property as \emph{$2$-path connectedness} below.
Clearly, this implies that the union of all triangles forms a connected subset of $\R^2$.
Two examples of such meshes are shown in \Cref{fig:triangulation}.

\begin{figure}[htb]
	\centering
	\includegraphics[width=0.25\linewidth]{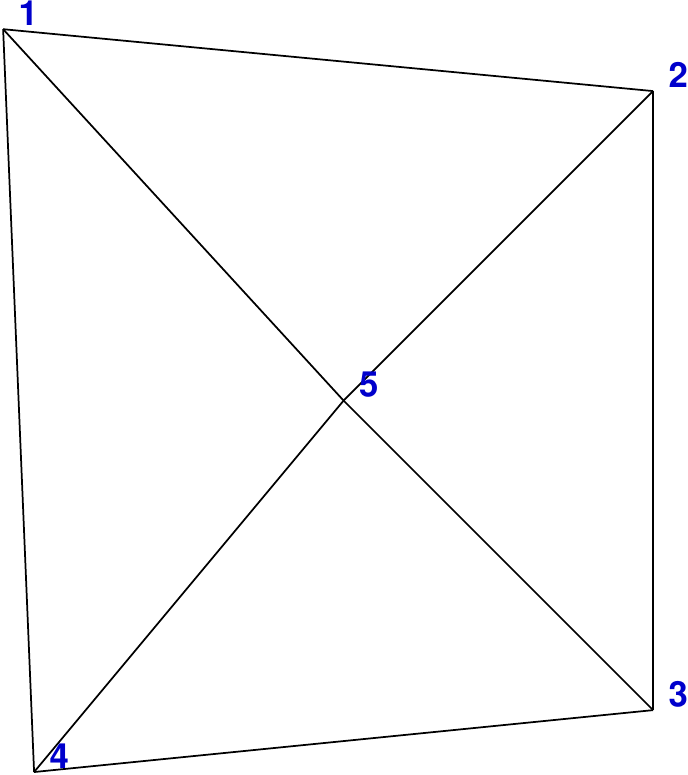}
	\includegraphics[width=0.25\linewidth]{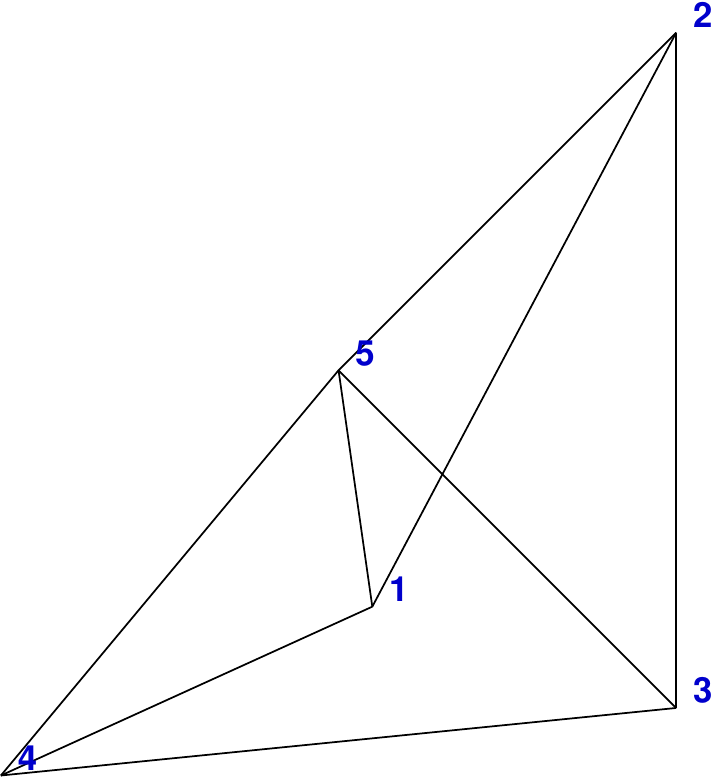}
	\caption{Admissible (left) and inadmissible (right) assignment of vertex coordinates for two meshes sharing the same connectivity.}
	\label{fig:admissible_and_inadmissible_meshes}
\end{figure}

Intuitively, a mesh is described by its connectivity and the positions of its vertices.
We recall that we are interested in all possible configurations of vertex positions a mesh of a given connectivity can attain.
To this end, we need to formulate conditions which avoid triangles becoming degenerate and vertices entering triangles to which they are not incident; see \Cref{fig:admissible_and_inadmissible_meshes}.

\subsection{Admissible Meshes}
\label{subsection:admissible_meshes}

In order to make these ideas formal, we utilize the concept of abstract simplicial complexes as well as geometric simplicial complexes; see \Cref{section:geometric_abstract_simplicial_complexes}.
We emphasize that all geometric simplicial complexes throughout the manuscript have vertices in $\R^2$.
From now on, we will use interchangeably the terms $0$-face or vertex, $1$-face or edge, and $2$-faces or triangle.
For simplicity of notation and without loss of generality, the $0$-faces will be numbered $1, \ldots, N_V$.

We begin with the formalization of \emph{connectivity}.
\begin{definition}
	\label{definition:connectivity_complex}
	Suppose that $\Delta$ is an abstract simplicial $2$-complex with vertex set~$V = \{1, \ldots, N_V\}$.
	We say that $\Delta$ is a \textbf{connectivity complex}, provided that
	\begin{enumerate}
		\item
			\label{item:connectivity_complex_1}
			$\Delta$ is pure,
		\item
			\label{item:connectivity_complex_2}
			$\Delta$ is $2$-path connected,
	\end{enumerate}
\end{definition}

\begin{example}
	\label{example:connectivity_complex}
	The abstract simplicial $2$-complex
	\begin{align}
		\Delta
		&
		=
		\paren[big]\{.{%
			\{1\},
			\{2\},
			\{3\},
			\{4\},
			\{5\},
			\{1,2\},
			\{2,3\},
			\{3,4\},
			\{4,1\},
			\{1,5\},
			\{2,5\},
		}%
		\nonumber
		\\
		&
		\phantom{{}={}\big\{}
		\paren[big].\}{%
			\{3,5\},
			\{4,5\},
			\{1,2,5\},
			\{2,3,5\},
			\{3,4,5\},
			\{4,1,5\}
		}%
		\label{eq:example_connectivity_complex}
	\end{align}
	is a connectivity complex.
	By contrast, the following abstract simplicial $2$-complexes are not connectivity complexes since they violate \cref{item:connectivity_complex_1,item:connectivity_complex_2}, respectively:
	\makeatletter
	\begin{equation*}
		\begin{aligned}
			\Delta
			&
			=
			\paren[big]\{\}{%
				\{1\},
				\{2\},
				\{3\},
				\{4\},
				\{1,2\},
				\{2,3\},
				\{3,1\},
				\{1,4\},
				\{1,2,3\}
			}%
			,
			\\
			\ltx@ifclassloaded{mcom-l}{%
				\Delta
				&
				=
				\paren[big]\{.{%
					\{1\},
					\{2\},
					\{3\},
					\{4\},
					\{5\},
					\{1,2\},
					\{2,3\},
					\{3,1\},
					\{3,4\},
					\{4,5\},
					\{5,3\},
				}
				\\
				&
				\phantom{{}={}\big\{}
					\paren[big].\}{%
					\{1,2,3\},
					\{3,4,5\}
				}%
				.
				}{%
				\Delta
				&
				=
				\paren[big]\{\}{%
					\{1\},
					\{2\},
					\{3\},
					\{4\},
					\{5\},
					\{1,2\},
					\{2,3\},
					\{3,1\},
					\{3,4\},
					\{4,5\},
					\{5,3\},
					\{1,2,3\},
					\{3,4,5\}
				}%
				.
			}
		\end{aligned}
	\end{equation*}
	\makeatother
\end{example}

A connectivity complex~$\Delta$ is a purely combinatorial object, which we can think of as a recipe for constructing meshes.
In order to achieve the latter, we need to assign vertex positions.
These can be summarized in a matrix
\begin{equation}
	\label{eq:vertex_position_matrix}
	Q
	=
	\begin{bmatrix}
		Q_1, \; Q_2, \; \ldots, \; Q_{N_V}
	\end{bmatrix}
	\in
	\R^{2 \times N_V}
	.
\end{equation}
Note that $Q$ with a lower index, such as $Q_1$, is a vector in $\R^2$ and it refers to one of the columns of $Q$ and describes the position of one of the vertices.
As is illustrated in \Cref{fig:admissible_and_inadmissible_meshes}, not all assignments of vertex positions will give rise to an admissible mesh.
In order to distinguish those which do from those which don't, we require further notation.
Given a connectivity complex~$\Delta$ and an assignment~$Q$ of its vertex positions, we define
\begin{equation}
	\label{eq:collection_of_convex_hulls}
	\Sigma_\Delta (Q)
	\coloneqq
	\setDef[auto]{\conv\{Q_{i_0}, \ldots, Q_{i_k}\}}{\{i_0, \ldots, i_k\} \in \Delta}
	\subset
	\powerset{\R^2}
	,
\end{equation}
where $\powerset{\R^2}$ denotes the power set of $\R^2$.
In other words, $\Sigma_\Delta(Q)$ collects the convex hulls of the vertices of all $0$-, $1$- and $2$-faces in $\Delta$.

We can now formalize the set of all admissible meshes with a given connectivity as follows.
\begin{definition}
	\label{definition:zeromanifold}
	Suppose that $\Delta$ is a connectivity complex with vertex set~$V = \{1, \ldots, N_V\}$.
	Then we define the \textbf{set of admissible meshes with connectivity~$\Delta$} as
	\begin{equation}
		\label{eq:zeromanifold}
		\zeromanifold
		\coloneqq
		\setDef[auto]{Q\in \R^{2 \times N_V}}{%
			\begin{aligned}
				&
				\Sigma_\Delta(Q) \text{ is a geometric simplicial complex in $\R^2$}
				\\
				&
				\text{whose associated abstract simplicial complex is $\Delta$}
			\end{aligned}
		}
		.
	\end{equation}
\end{definition}
Notice that $Q \in \zeromanifold$ implies that $\{Q_{i_0}, Q_{i_1}, Q_{i_2}\}$ are affine independent for all $\{i_0,i_1,i_2\} \in \Delta$.

The conditions of $\Sigma_\Delta(Q)$ formalize the idea of an admissible mesh, \ie, that all triangles be non-degenerate and that any two intersecting triangles can only intersect in a common edge or a common vertex.
We remark that we need to insist in definition \eqref{eq:zeromanifold} that the abstract simplicial complex associated with $\Sigma_\Delta(Q)$ agrees with $\Delta$.
As an example, consider
\makeatletter
\begin{equation*}
	\begin{aligned}
		\ltx@ifclassloaded{mcom-l}{%
			\Delta
			&
			=
			\paren[big]\{.{%
				\{1\},
				\{2\},
				\{3\},
				\{4\},
				\{1,2\},
				\{2,3\},
				\{3,1\},
				\{1,4\},
				\{4,2\},
				\{2,1\},
			}%
			\\
			&
			\phantom{{}={}\big\{}
				\paren[big].\}{%
				\{1,2,3\},
				\{2,1,4\}
			}%
			}{%
			\Delta
			&
			=
			\paren[big]\{\}{%
				\{1\},
				\{2\},
				\{3\},
				\{4\},
				\{1,2\},
				\{2,3\},
				\{3,1\},
				\{1,4\},
				\{4,2\},
				\{2,1\},
				\{1,2,3\},
				\{2,1,4\}
			}%
		}
	\end{aligned}
\end{equation*}
\makeatother
and choose coordinates such that $\conv\{Q_1, Q_2, Q_3\}$ and $\conv\{Q_2, Q_1, Q_4\}$ are $2$-simplices but $Q_3 = Q_4$ holds.
Then $\Sigma_\Delta(Q)$ is a geometric simplicial complex but its abstract simplicial complex is smaller than $\Delta$ since the two triangles coincide.

It is easy to see that there exist connectivity complexes~$\Delta$ for which $\zeromanifold$ is empty.
This is the case, for instance, when
\makeatletter
\begin{align}
	\ltx@ifclassloaded{mcom-l}{%
		\Delta
		&
		=
		\paren[big]\{.{%
			\{1\},
			\{2\},
			\{3\},
			\{4\},
			\{5\},
			\{1,2\},
			\{2,3\},
			\{3,1\},
			\{1,4\},
			\{2,4\},
			\{1,5\},
			\{2,5\},
		}%
		\nonumber
		\\
		&
		\phantom{{}={}\big\{}
			\paren[big].\}{%
			\{1,2,3\},
			\{1,2,4\},
			\{1,2,5\}
		}
		.
		}{%
		\Delta
		&
		=
		\paren[big]\{\}{%
			\{1\},
			\{2\},
			\{3\},
			\{4\},
			\{5\},
			\{1,2\},
			\{2,3\},
			\{3,1\},
			\{1,4\},
			\{2,4\},
			\{1,5\},
			\{2,5\},
			\{1,2,3\},
			\{1,2,4\},
			\{1,2,5\}
		}%
		.
	}
	\label{eq:example_connectivity_complex_three_triangles_sharing_an_edge}
\end{align}
\makeatother
Notice that there are three triangles with one common edge, which represents an impossible configuration for a geometric simplicial $2$-complex in $\R^2$.

The possible emptyness of $\zeromanifold$ will not be a cause of concern in what follows.

\begin{figure}[htb]
	\begin{center}
	\makeatletter
	\ltx@ifclassloaded{mcom-l}{%
	\begin{tikzpicture}[scale=0.75, dot/.style={circle,inner sep=1pt,fill,name=#1},
	  extended line/.style={shorten >=-#1,shorten <=-#1},
	  extended line/.default=1cm]%
	}{%
	\begin{tikzpicture}[dot/.style={circle,inner sep=1pt,fill,name=#1},
	  extended line/.style={shorten >=-#1,shorten <=-#1},
	  extended line/.default=1cm]
		}
		\makeatother

	\draw [fill, blue!20] (0,5) -- (1.5,6.5) -- (3,5) -- cycle;
	\draw                 (0,5) -- (1.5,6.5) -- (3,5) -- cycle;
	\draw [fill, gray!20] (0,5) -- (1.5,3.5) -- (3,5) -- cycle;
	\draw                 (0,5) -- (1.5,3.5) -- (3,5) -- cycle;
	\draw [dashed, red]   (-1,5) -- (4,5);
	\draw [red,thick]           (0,5) -- (3,5);

	\draw [fill, blue!20] (5,6) -- (6.5,5) -- (7.5,5) -- cycle;
	\draw                 (5,6) -- (6.5,5) -- (7.5,5) -- cycle;
	\draw [fill, gray!20] (5.5,4) -- (6.5,5) -- (7.5,4.7) -- cycle;
	\draw                 (5.5,4) -- (6.5,5) -- (7.5,4.7) -- cycle;
	\draw [dashed, red]   (4.5,5) -- (8.5,5);
	\draw[red,fill=red] (6.5,5) circle (.4ex);

	\draw [fill, blue!20] (9.5,5) -- (11,6.5) -- (12.5,5) -- cycle;
	\draw                 (9.5,5) -- (11,6.5) -- (12.5,5) -- cycle;
	\draw [red,dashed] (10.5,7) -- (13,4.5);
	\draw [red,thick] (11,6.5) -- (12.5,5);

	\draw [fill, blue!20] (0,2) -- (1.5,3) -- (3,2) -- cycle;
	\draw                 (0,2) -- (1.5,3) -- (3,2) -- cycle;
	\draw                 (1.5,3) -- (4,2.5);
	\draw [red,dashed] (1,3.33) -- (3.5,1.67);
	\draw[red,fill=red] (1.5,3) circle (.4ex);

	\draw [fill, blue!20] (4.5,2) -- (6,3) -- (7.5,2) -- cycle;
	\draw                 (4.5,2) -- (6,3) -- (7.5,2) -- cycle;
	\draw[red, dashed]    (5.5,3.33) -- (8,1.67);
	\draw[red,fill=red] (6,3) circle (.4ex);

	\draw [blue, thick]   (9,2.5) -- (10,3.5);
	\draw [thick]         (10,3.5) -- (11,3);
	\draw [red,dashed]    (8.5,2) -- (10.5,4);
	\draw[red,fill=red]  (10,3.5) circle (.4ex);

	\draw [blue, thick]   (12,2.5) -- (13,3.5);
	\draw [red, dashed]   (11.5,2) -- (13.5,4);
	\draw[red,fill=red]  (13,3.5) circle (.4ex);
	\end{tikzpicture}
\end{center}
	\caption{Illustration for the proof of \cref{proposition:zeromanifold_is_open}, showing (from top left to bottom right) the representative cases $(k,\ell,m) \in \{(2,2,1), (2,2,0), (2,1,1), (2,1,0), (2,0,0), (1,1,0), (1,0,0)\}$. $\sigma(Q)$ is shown in blue, $\sigma'(Q)$ is shown in black, and their intersection $\tau(Q)$ is shown in red. A possible separating hyperplane is displayed as a dashed red line.}
	\label{fig:proof_zeromanifold_is_open}
\end{figure}
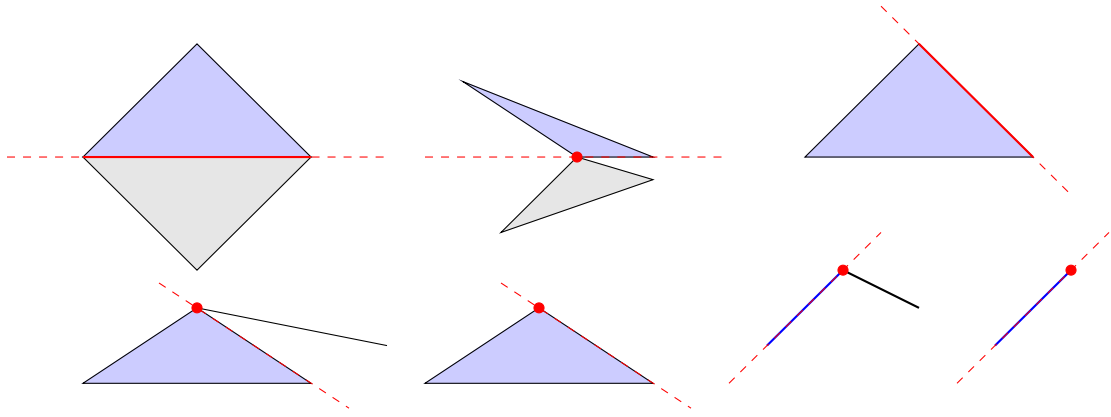

\begin{proposition}
	\label{proposition:zeromanifold_is_open}
	For any given connectivity complex~$\Delta$ with vertex set~$V = \{1, \ldots, N_V\}$, the set $\zeromanifold$ is an open (possibly empty) subset of $\R^{2 \times N_V}$.
\end{proposition}
\begin{proof}
	We can assume that $\zeromanifold$ is non-empty.
	Let $Q \in \zeromanifold$ be arbitrary.
	We need to prove that there exists $\delta > 0$ such that the open ball $B_\delta(Q) \subset \R^{2 \times N_V}$, \eg, in the Frobenius norm, belongs to $\zeromanifold$.
	We proceed in the following steps, selecting a suitable $\delta > 0$ along the way.
	We show that, for all $U \in B_\delta(Q)$,
	\begin{enumerate}
		\item
			\label[statement]{item:zeromanifold_is_open_i}
			$\Sigma_\Delta(U)$ is a geometric simplicial $2$-complex and

		\item
			\label[statement]{item:zeromanifold_is_open_ii}
			$\Sigma_\Delta(U)$ has associated abstract simplicial complex $\Delta$.
	\end{enumerate}
	We begin with \cref{item:zeromanifold_is_open_i}.
	Suppose that $\{i_0,i_1,i_2\}$ is an arbitrary $2$-face in $\Delta$.
	Since $\{Q_{i_0},Q_{i_1},Q_{i_2}\}$ is affine independent, $\det \begin{bmatrix} Q_{i_1} - Q_{i_0}, \; Q_{i_2} - Q_{i_1} \end{bmatrix} \neq 0$.
	By continuity of the determinant function, we can find $\delta > 0$ such that the determinant $\det \begin{bmatrix} U_{i_1} - U_{i_0}, \; U_{i_2} - U_{i_1} \end{bmatrix}$ has the same sign as before for all $U \in B_\delta(Q)$.
	Therefore, $\conv\{U_{i_0},U_{i_1},U_{i_2}\}$ is a $2$-simplex.
	Since the number of $2$-faces in $\Delta$ is finite, a joint value of $\delta > 0$ can be found which is valid for all $2$-faces in $\Delta$.
	Clearly, the same reasoning also applies to the $1$-faces and $0$-faces.
	Consequently, for all $U \in B_\delta(Q)$, $\Sigma_\Delta(U)$ consists of a collection of simplices whose dimension agrees with the dimension of the corresponding face of $\Delta$.

	In the following, let $\sigma$ and $\sigma'$ be any two \emph{distinct} faces in $\Delta$.
	We denote by $\sigma(Q)$ and $\sigma'(Q)$ the corresponding simplices in $\Sigma_\Delta(Q)$.
	For instance, when $\sigma = \{i_0, \ldots, i_k\}$, then $\sigma(Q) = \conv\{Q_{i_0}, \ldots, Q_{i_k}\}$.
	By construction, it is clear that when $\tau \subset \sigma$ holds, then $\tau(U)$ is a face of $\Sigma_\Delta(U)$, for all $U \in B_\delta(Q)$.
	We also know that $\tau(Q) \coloneqq \sigma(Q) \cap \sigma'(Q)$ is either empty or a face of both.
	In order to conclude \cref{item:zeromanifold_is_open_i}, we now show that this property extends to all $U \in B_\delta(Q)$, possibly for a smaller value of $\delta > 0$ than previously chosen.
	We distinguish two cases.
	$\bullet$
	When $\tau(Q) = \emptyset$, then since $\sigma(Q)$ and $\sigma'(Q)$ are compact and convex, there exists an affine linear functional $\varphi$ such that $\varphi < 0$ on $\sigma(Q)$ and $\varphi > 0$ on $\sigma'(Q)$; see, \eg, \cite[Cor.~4.1.3, p.52]{HiriartUrrutyLemarechal2001}.
	Possibly by making $\delta$ smaller, we retain $\varphi < 0$ on $\sigma(U)$ and $\varphi > 0$ on $\sigma'(U)$ for all $U \in B_\delta(Q)$.
	Consequently, $\tau(U) = \sigma(U) \cap \sigma'(U) = \emptyset$.
	$\bullet$
	Suppose that $\tau(Q)$ is a face of both $\sigma(Q)$ and $\sigma'(Q)$, say, $\sigma(Q) = \conv\{Q_{i_0}, \ldots, Q_{i_k}\}$, $\sigma'(Q) = \conv\{Q_{j_0}, \ldots, Q_{j_\ell}\}$ and $\tau(Q) = \conv\{Q_{i_0}, \ldots, Q_{i_m}\} = \conv\{Q_{j_0}, \ldots, Q_{j_m}\}$ with some $0 \le m \le \min\{k,\ell\}$.
	We have already proved that $\sigma(U)$ and $\sigma'(U)$ are simplices of dimensions~$k$ and $\ell$, respectively, for all $U \in B_\delta(Q)$.
	Therefore, the only concern is that $\tau(U)$ is larger than $\conv\{U_{i_0}, \ldots, U_{i_m}\}$.
	In each case, however, we can construct a hyperplane, defined by two vertices of either $\sigma$ or $\sigma'$, which separates $\sigma(U) \setminus \sigma'(U)$ from $\sigma'(U) \setminus \sigma(U)$, for all $U \in B_\delta(Q)$, possibly for a smaller value of $\delta > 0$ than previously chosen.
	See \Cref{fig:proof_zeromanifold_is_open} for an illustration.
	$\bullet$
	Altogether, this confirms that $\tau(U) \coloneqq \sigma(U) \cap \sigma'(U)$ is either empty or a face of both for $U$ in a suitable ball $B_\delta(Q)$.
	While looping over all pairs of distinct faces $\{\sigma,\sigma'\}$, $\delta$ needs to be reduced only finitely many times, therefore we have shown \cref{item:zeromanifold_is_open_i}.

	In order to show \cref{item:zeromanifold_is_open_ii}, recall from \cref{subsection:abstract_simplicial_complexes} that the abstract simplicial complex associated with $\Sigma_\Delta(Q)$ is defined as
	\begin{equation*}
		\Delta(Q)
		=
		\setDef[big]{ \sigma \subset \{1, \ldots, N_V\}}{\conv\{Q_i\}_{i \in \sigma} \in \Sigma_\Delta(Q)}
		.
	\end{equation*}
	Since $Q \in \zeromanifold$, $\Delta(Q) = \Delta$ holds by definition.
	Moreover, by \eqref{eq:collection_of_convex_hulls} we clearly have $\Delta \subset \Delta(U)$ for all $U \subset \R^{2 \times N_V}$ and thus for all $U \in B_\delta(Q)$.
	However, the considerations above show that there can be no additional simplices in $\Sigma_\Delta(U)$ than those coming from \eqref{eq:collection_of_convex_hulls}.
	In other words, $\Delta = \Delta(U)$ holds for all $U \in B_\delta(Q)$, which is \cref{item:zeromanifold_is_open_ii}.
\end{proof}

\begin{figure}[htb]
	\centering
	\includegraphics[width=0.25\linewidth]{Figures/nonoriented_mesh.pdf}
	\qquad
	\includegraphics[width=0.25\linewidth]{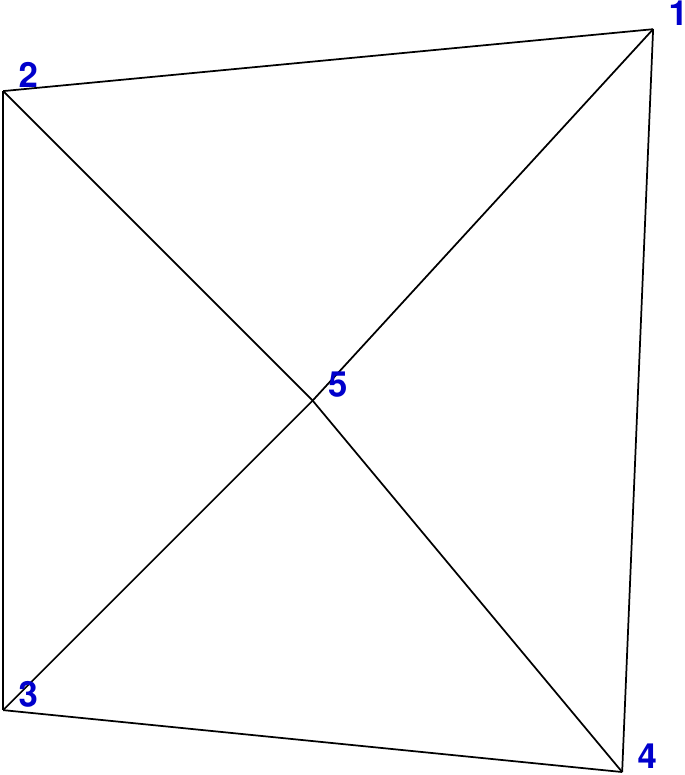}
	\caption{Two admissible meshes with the same connectivity complex~$\Delta$ whose vertex position matrices~$Q$ lie in different connected components of $\zeromanifold$.}
	\label{fig:zeromanifold_is_not_connected}
\end{figure}

\Cref{proposition:zeromanifold_is_open} shows that $\zeromanifold$ is an open submanifold of $\R^{2 \times N_V}$.
It is however easy to see that any non-empty $\zeromanifold$ is not path-connected and, equivalently, not connected.
In fact, $Q \in \zeromanifold$ implies that $-Q$ lies in $\zeromanifold$ as well but in a different connected component.
This is true even if $\Delta$ contains only a single $2$-face.
For example, the two meshes shown in \Cref{fig:zeromanifold_is_not_connected} cannot be joined by a continuous path of vertex positions which remains inside $\zeromanifold$.
In order to resolve this issue, we need to consider orientations.

\subsection{Orientation}
\label{subsection:orientation}

Recall that an orientation of an abstract $2$-face $\sigma = \{i_0,i_1,i_2\}$ is an equivalence class of orderings of $\sigma$, where two orderings $i_0,i_1,i_2$ and $i_{\pi(0)}$, $i_{\pi(1)}$ and $i_{\pi(2)}$ are equivalent if and only if $\pi$ is an even permutation of $0, 1, 2$.
We denote the equivalence class represented by the ordering $i_0,i_1,i_2$ by $[i_0,i_1,i_2]$.
Clearly, each abstract $2$-face has precisely two orientations.
The orientation $[i_0,i_1,i_2]$ of the abstract $2$-face $\sigma = \{i_0,i_1,i_2\}$ induces an orientation on each $1$-face contained in $\sigma$, namely $[i_1, i_2]$, $[i_2, i_0]$, and $[i_0,i_1]$, respectively.

\begin{definition}
	\label{definition:consistently_orientated_abstract_simplicial_2-complex}
	Let $\Delta$ be an abstract simplicial $2$-complex.
	\begin{enumeratelatin}
		\item
			Suppose that each of the $2$-simplices in $\Delta$ has an assigned orientation.
			We say that these orientations are consistent and that $\Delta$ is \textbf{consistently oriented} if and only if the orientations of any two $2$-simplices in $\Delta$ sharing a $1$-face induce opposite orientations on that $1$-face.
		\item
			We say that $\Delta$ is \textbf{orientable} if there exist orientations of all $2$-simplices in $\Delta$ which render $\Delta$ consistently oriented.
	\end{enumeratelatin}
\end{definition}

For example, the abstract simplicial $2$-complex in \eqref{eq:example_connectivity_complex} with orientations $[1, 2, 5]$, $[2, 3, 5]$, $[3, 4, 5]$ and $[4, 1, 5]$ is consistently oriented.
Moreover, the abstract simplicial $2$-complex from \eqref{eq:example_connectivity_complex_three_triangles_sharing_an_edge} is non-orientable.

\begin{remark}
	\label{remark:connectivity_complex_and_connectivity_matrix}
	It is easy to see that a connectivity complex~$\Delta$ can be encoded as a \textbf{connectivity matrix}~$C \in \R^{3 \times N_T}$, where $N_T$ is the number of triangles in~$\Delta$ and each column of~$C$ lists the vertices of one of the triangles.
	For instance, the connectivity complex from \eqref{eq:example_connectivity_complex} could be summarized as
	\begin{equation*}
		C
		=
		\begin{bmatrix}
			2 & 3 & 4 & 1 \\
			1 & 2 & 3 & 4 \\
			5 & 5 & 5 & 5
		\end{bmatrix}
		.
	\end{equation*}
	Its consistent orientation is reflected in the fact that triangles sharing an edge have their two common vertices appear in opposite orderings.
	Consider for instance the third and fourth columns, in which the shared vertices appear in the orders (after an even permutation) $[5,4]$ and $[4,5]$, respectively.
\end{remark}

In the sequel, we will say that \eqq{$[i_0,i_1,i_2]$ is a $2$-face} instead of \eqq{$\{i_0,i_1,i_2\}$ is a $2$-face with orientation $[i_0,i_1,i_2]$}.
Similarly, we write $[i_0,i_1]$ to denote a $1$-face together with its orientation.
For consistency of notation, we will also write $[i_0]$ instead of $\{i_0\}$ for a $0$-face although orientation does not matter there.

\begin{assumption}
	\label{assumption:oriented_connectivity_complex}
	For the remainder of this paper, we suppose that $\Delta$ is a consistently oriented connectivity complex with vertex set~$V = \{1, \ldots, N_V\}$.
	The matrix $Q \in \R^{2 \times N_V}$ denotes an arbitrary assignment of its vertex positions.
	Further assumptions on $Q$ are specified as needed.
\end{assumption}

\begin{definition}
	\label{definition:signed_area}
	Suppose that $[i_0,i_1,i_2]$ is a $2$-face of $\Delta$.
	We define its \textbf{signed area} as
	\begin{equation}
		\label{eq:signed_area}
		\area{Q}[i_0,i_1,i_2]
		\coloneqq
		\frac{1}{2} \det
		\begin{bmatrix}
			Q_{i_1} - Q_{i_0}, \; Q_{i_2} - Q_{i_1}
		\end{bmatrix}
		.
	\end{equation}
\end{definition}
Notice that this definition is independent of the particular ordering of vertices representing the orientation since $\area{Q}[i_0,i_1,i_2] = \area{Q}[i_1,i_2,i_0] = \area{Q}[i_2,i_0,i_1]$ holds.
Specifically, $\area{Q}[i_0,i_1,i_2] > 0$ indicates that the vertex positions $Q_{i_0}$, $Q_{i_1}$, $Q_{i_2}$ are in counterclockwise ordering.
The opposite ordering leads to a change in sign.
Moreover, $\area{Q}[i_0,i_1,i_2] \neq 0$ holds if and only if $\{Q_{i_0}, Q_{i_1}, Q_{i_2}\}$ are affine independent, regardless of the orientation.

\begin{definition}
	\label{definition:plusmanifold}
	We define the \textbf{set of admissible oriented meshes with connectivity~$\Delta$} as
	\begin{equation}
		\label{eq:plusmanifold}
		\plusmanifold
		\coloneqq
		\setDef[big]{Q \in \zeromanifold}{\area{Q}[i_0,i_1,i_2] > 0 \text{ for all $2$-faces $[i_0,i_1,i_2]$ of $\Delta$}}
		.
	\end{equation}
\end{definition}

\begin{figure}[htb]
	\centering
		\includegraphics[width=0.25\linewidth]{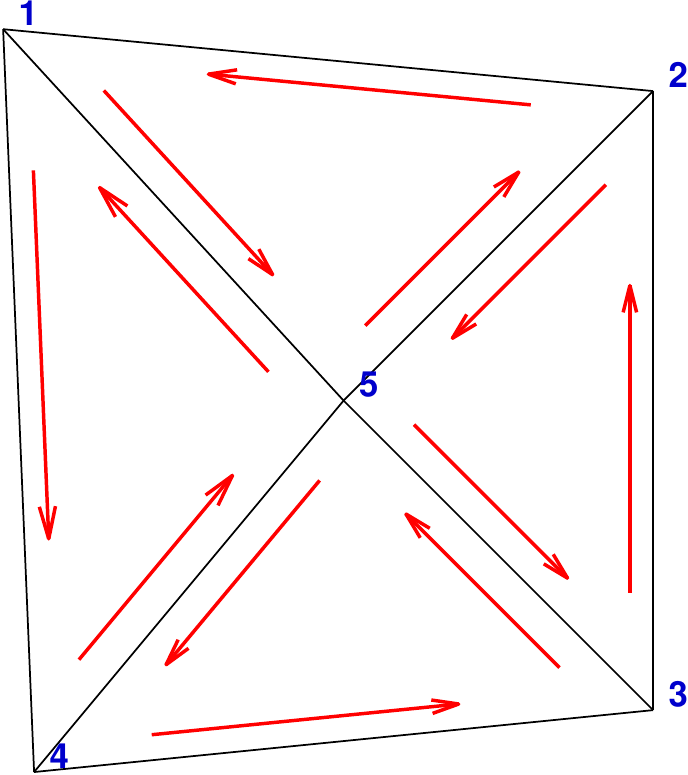}
		\qquad
		\includegraphics[width=0.40\linewidth]{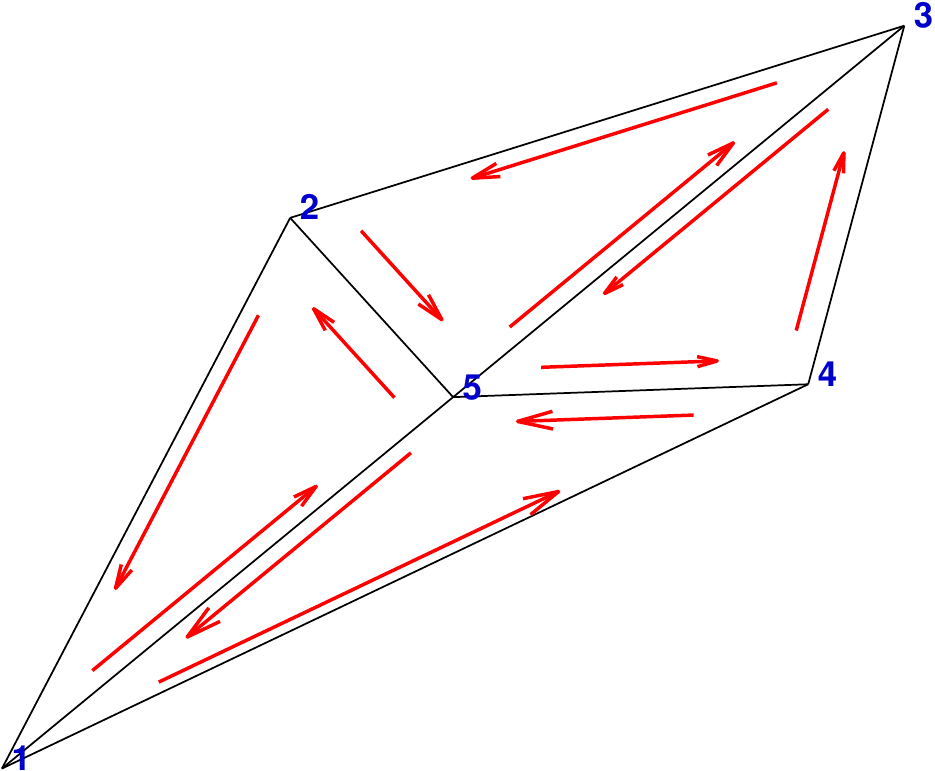}
		\caption{Two admissible oriented meshes, \ie, elements of $\plusmanifold$, with the same consistently oriented connectivity complex~$\Delta$ and different vertex positions~$Q$. The orientations of the $1$-faces are $[2,1,5]$, $[3,2,5]$, $[4,3,5]$, and $[1,4,5]$.}
	\label{fig:points_in_plusmanifold}
\end{figure}

\Cref{fig:points_in_plusmanifold} illustrates two elements of $\plusmanifold$.

\begin{proposition}
	\label{proposition:plusmanifold_is_open}
	The set $\plusmanifold$ is an open (possibly empty) subset of $\R^{2 \times N_V}$.
\end{proposition}
\begin{proof}
	Due to the continuity of the determinant and \cref{proposition:zeromanifold_is_open}, the set described in \eqref{eq:plusmanifold} is a finite intersection of open subsets of $\R^{2 \times N_V}$.
\end{proof}

\subsection{Definition of the Manifold of Planar Triangular Meshes}
\label{subsection:discrete_shape_manifold}

The set of admissible oriented meshes with connectivity $\Delta$,  $\plusmanifold$ in \cref{definition:plusmanifold} answers our initial question, which vertex positions a mesh, or rather an oriented mesh, can attain?
However, in order to endow it with a complete Riemannian metric, we would need this set to be connected, which unfortunately it is not.
This is shown by the following counterexample, provided by one of the reviewers.

\begin{example}
	\label{example:plusmanifoldNotConnected}
	Let us consider the following connectivity complex
	\begin{align*}
		\Delta
		&
		=
		\paren[big]\{.{%
				\{1\},
				\{2\},
				\{3\},
				\{4\},
				\{5\},
				\{6\},
				\{1,2\},
				\{2,3\},
				\{1,3\},
				\{1,5\},
				\{2,5\},
				\{2,6\},
				\{3,6\},
			}%
			\nonumber
			\\
			&
			\phantom{{}={}\big\{}
				\{3,4\},
				\{1,4\},
				\{4,5\},
				\{5,6\},
				\{4,6\},
				\{1,3,4\},
				\{1,2,5\},
				\{2,3,6\},
				\{3,4,6\},
				\\
				&
				\phantom{{}={}\big\{}
					\paren[big].\}{%
					\{1,4,5\},
					\{2,5,6\}
				}%
				,
			\end{align*}
		with orientations of its $2$-faces as follows: $[1,4,3]$, $[1,2,5]$, $[2,3,6]$, $[3,4,6]$, $[1,5,4]$ and $[2,6,5]$.
		Moreover, let us consider
		\begin{align*}
			Q
			&
			=
			\begin{bmatrix}
				0.75 & 1.25 &1   &-0.5 &1 &2.5
				\\
				1.25 & 1.25 &0.75& 0 &2.5 &0
			\end{bmatrix}
			,
			\\
			\widetilde Q
			&
			=
			\begin{bmatrix}
				-0.5 &2.5 &1 &0.75 &1 &1.25
				\\
				0 &0 &2.5 &1.25 &0.75 &1.25
			\end{bmatrix}
			.
		\end{align*}
		It is easy to verify that both $Q$ and $\widetilde Q$ belong to $\plusmanifold$.
		However, there exists no continuous path in $\plusmanifold$ joining~$Q$ and $\widetilde Q$; see \cref{fig:counterExamplePlusManifold} for an illustration.
\end{example}

\begin{figure}
	\begin{subfigure}{0.45\textwidth}
	\centering
	\begin{tikzpicture}[scale=1.4]
\fill[blue!20,opacity = 0.5] (-0.5,0) -- (0.75,1.25) -- (1,2.5) -- cycle;
\draw[black] (-0.5,0) -- (0.75,1.25) -- (1,2.5) -- cycle;
\fill[blue!20,opacity = 0.5] (0.75,1.25) -- (1.25,1.25) -- (1,2.5) -- cycle;
\draw[black] (0.75,1.25) -- (1.25,1.25) -- (1,2.5) -- cycle;
\fill[blue!20,opacity = 0.5] (2.5,0) -- (1.25,1.25) -- (1,2.5) -- cycle;
\draw[black] (2.5,0) -- (1.25,1.25) -- (1,2.5) -- cycle;
\fill[blue!20,opacity = 0.5] (2.5,0) -- (1.25,1.25) -- (1,0.75) -- cycle;
\draw[black](2.5,0) -- (1.25,1.25) -- (1,0.75) -- cycle;
\fill[blue!20,opacity = 0.5] (2.5,0) -- (-0.5,0) -- (1,0.75) -- cycle;
\draw[black](2.5,0) -- (-0.5,0) -- (1,0.75) -- cycle;
\fill[blue!20,opacity = 0.5] (0.75,1.25) -- (-0.5,0) -- (1,0.75) -- cycle;
\draw[black](0.75,1.25) -- (-0.5,0) -- (1,0.75)-- cycle;

\draw (0.75,1.25) node[left]{$Q_1$};
\draw (1.25,1.25) node[right]{$Q_2$};
\draw (1,0.75) node[below]{$Q_3$};
\draw (-0.5,0) node[left]{$Q_4$};
\draw (1,2.5) node[left]{$Q_5$};
\draw (2.5,0) node[right]{$Q_6$};

\draw[blue] (0.75,1.25)  node[circle,fill, black, scale = 0.2]{1};
\draw[blue] (1.25,1.25)  node[circle,fill, black, scale = 0.2]{2};
\draw[blue] (1,0.75)  node[circle,fill, black, scale = 0.2]{3};
\draw[blue] (-0.5,0)  node[circle,fill, black, scale = 0.2]{4};
\draw[blue] (1,2.5)  node[circle,fill, black, scale = 0.2]{5};
\draw[blue] (2.5,0)  node[circle,fill, black, scale = 0.2]{6};
\end{tikzpicture}
	\end{subfigure}
	~
\begin{subfigure}{0.45\textwidth}
\centering
\begin{tikzpicture}[scale=1.4]
\fill[blue!20,opacity = 0.5] (-0.5,0) -- (0.75,1.25) -- (1,2.5) -- cycle;
\draw[black] (-0.5,0) -- (0.75,1.25) -- (1,2.5) -- cycle;
\fill[blue!20,opacity = 0.5] (0.75,1.25) -- (1.25,1.25) -- (1,2.5) -- cycle;
\draw[black] (0.75,1.25) -- (1.25,1.25) -- (1,2.5) -- cycle;
\fill[blue!20,opacity = 0.5] (2.5,0) -- (1.25,1.25) -- (1,2.5) -- cycle;
\draw[black] (2.5,0) -- (1.25,1.25) -- (1,2.5) -- cycle;
\fill[blue!20,opacity = 0.5] (2.5,0) -- (1.25,1.25) -- (1,0.75) -- cycle;
\draw[black](2.5,0) -- (1.25,1.25) -- (1,0.75) -- cycle;
\fill[blue!20,opacity = 0.5] (2.5,0) -- (-0.5,0) -- (1,0.75) -- cycle;
\draw[black](2.5,0) -- (-0.5,0) -- (1,0.75) -- cycle;
\fill[blue!20,opacity = 0.5] (0.75,1.25) -- (-0.5,0) -- (1,0.75) -- cycle;
\draw[black](0.75,1.25) -- (-0.5,0) -- (1,0.75)-- cycle;

\draw (0.75,1.25) node[left]{$\widetilde Q_4$};
\draw (1.25,1.25) node[right]{$\widetilde Q_6$};
\draw (1,0.75) node[below]{$\widetilde Q_5$};
\draw (-0.5,0) node[left]{$\widetilde Q_1$};
\draw (1,2.5) node[left]{$\widetilde Q_3$};
\draw (2.5,0) node[right]{$\widetilde Q_2$};

\draw[blue] (0.75,1.25)  node[circle,fill, black, scale = 0.2]{1};
\draw[blue] (1.25,1.25)  node[circle,fill, black, scale = 0.2]{2};
\draw[blue] (1,0.75)  node[circle,fill, black, scale = 0.2]{3};
\draw[blue] (-0.5,0)  node[circle,fill, black, scale = 0.2]{4};
\draw[blue] (1,2.5)  node[circle,fill, black, scale = 0.2]{5};
\draw[blue] (2.5,0)  node[circle,fill, black, scale = 0.2]{6};
\end{tikzpicture}
\end{subfigure}
	\caption{Illustration of the non-connectedness of $\plusmanifold$; see \cref{example:plusmanifoldNotConnected} for details.}
	\label{fig:counterExamplePlusManifold}
\end{figure}
Since a complete metric requires the underlying manifold to be connected, we cannot use all of $\plusmanifold$.
Therefore, rather than considering all oriented meshes which can be generated from an oriented connectivity complex~$\Delta$, we start from a given oriented reference mesh (represented by vertex coordinates~$\Qref$) and ask, which other oriented meshes with the same connectivity complex can be obtained through continuous deformations within $\plusmanifold$?

This leads us to the following
\begin{definition}
	\label{definition:planarmanifold}
	Suppose that $\Qref \in \plusmanifold$ are the coordinates of a given reference mesh.
	We define the \textbf{manifold of planar triangular meshes} as
	\makeatletter
	\ltx@ifclassloaded{mcom-l}{%
		\begin{equation}
			\label{eq:planarmanifold}
			\planarmanifold
			\coloneqq
			\setDef[auto]{Q \in \plusmanifold}{%
			\begin{aligned}
			&
			\text{there exists a continuous path}
			\\
			&
			\text{in $\plusmanifold$ from $\Qref$ to $Q$.}
			\end{aligned}
			}.
		\end{equation}
		}{%
		\begin{equation}
			\label{eq:planarmanifold}
			\planarmanifold
			\coloneqq
			\setDef[big]{Q \in \plusmanifold}{\text{there exists a continuous path in $\plusmanifold$ from $\Qref$ to $Q$}}
			.
		\end{equation}
	}%
	\makeatother
\end{definition}

The set $\planarmanifold$ is our primary object of interest.
Let us summarize some essential properties.
\begin{theorem}
	\label{theorem:properties_planarmanifold}
	The set $\planarmanifold$ is
	\begin{enumeratelatin}
		\item
			\label[statement]{item:planarmanifold_is_a_path_component}
			a path component of $\plusmanifold$, and thus path-connected,
		\item
			\label[statement]{item:planarmanifold_is_an_open_submanifold}
			an open submanifold of $\R^{2 \times N_V}$.
	\end{enumeratelatin}
\end{theorem}
\begin{proof}
	\Cref{item:planarmanifold_is_a_path_component} is immediate from the definition of $\planarmanifold$.
	Moreover, since $\plusmanifold$ is locally path-connected, its path components are open.
	Since $\plusmanifold$ is itself open in $\R^{2 \times N_V}$ by \cref{proposition:plusmanifold_is_open}, \cref{item:planarmanifold_is_an_open_submanifold} follows.
	We refer the reader to \cite[Ch.~4]{Lee2011} for a background on components and path components of manifolds.
\end{proof}

\begin{remark}
	Notice that as a submanifold, $\planarmanifold$ inherits the Hausdorff and second countability properties of $\R^{2 \times N_V}$.
\end{remark}

\section{A Complete Riemannian Metric for \texorpdfstring{$\planarmanifold$}{the Manifold of Planar Triangular Meshes}}
\label{section:complete_Riemannian_metric}

\begin{figure}[htb]
	\centering
	\includegraphics[width=0.25\linewidth]{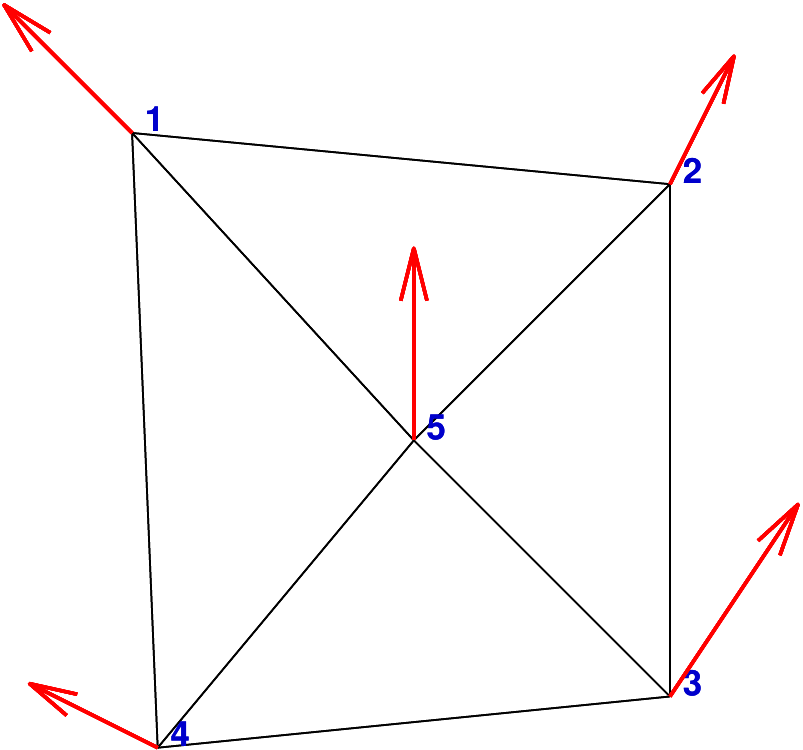}
	\caption{Example of a tangent vector.}
	\label{fig:tangent_vector}
\end{figure}

In this section we propose a Riemannian metric for manifold of planar triangular meshes $\planarmanifold$ and show its completeness.
The tangent space at $Q \in \planarmanifold$ will be denoted by $\tangent{Q}[\planarmanifold]$ and we remark that it is a vector space isomorphic to $\R^{2 \times N_V}$.
In fact, we can visualize a tangent vector~$V \in \R^{2 \times N_V}$ as a collection of vectors in $\R^2$, one per vertex; see \Cref{fig:tangent_vector}.
We will occasionally need the vectorization operation $\vvec \colon \R^{2 \times N_V} \to \R^{2 \, N_V}$ which stacks $Q \in \R^{2 \times N_V}$ column by column.

\begin{figure}[htb]
	\begin{subfigure}{0.3\textwidth}
		\centering
		\includegraphics[width=\linewidth]{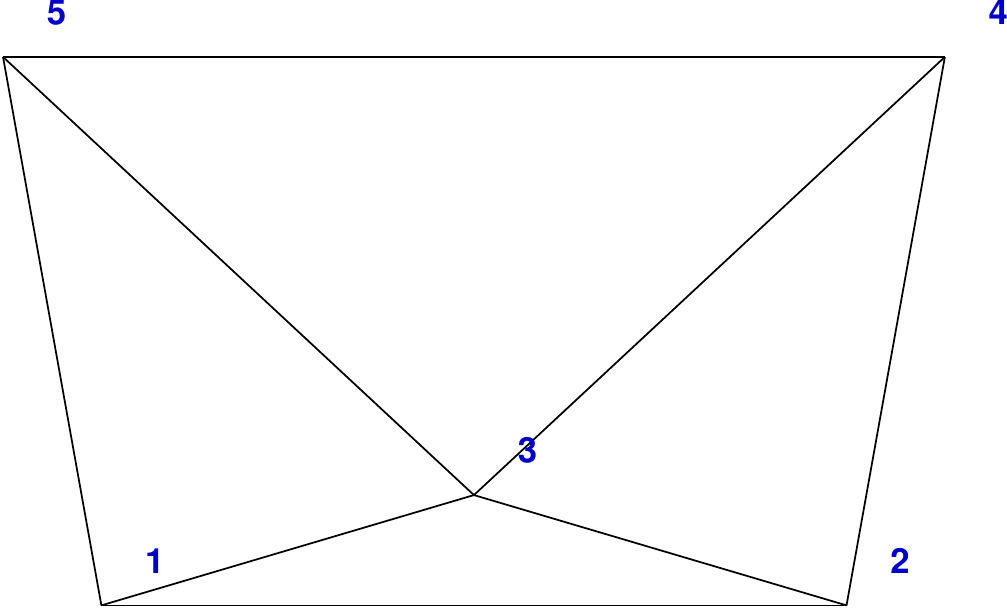}
		\caption{$t=0.3$}
	\end{subfigure}
	\hfill
	\begin{subfigure}{0.3\textwidth}
		\centering
		\raisebox{5mm}{%
		\includegraphics[width=\linewidth]{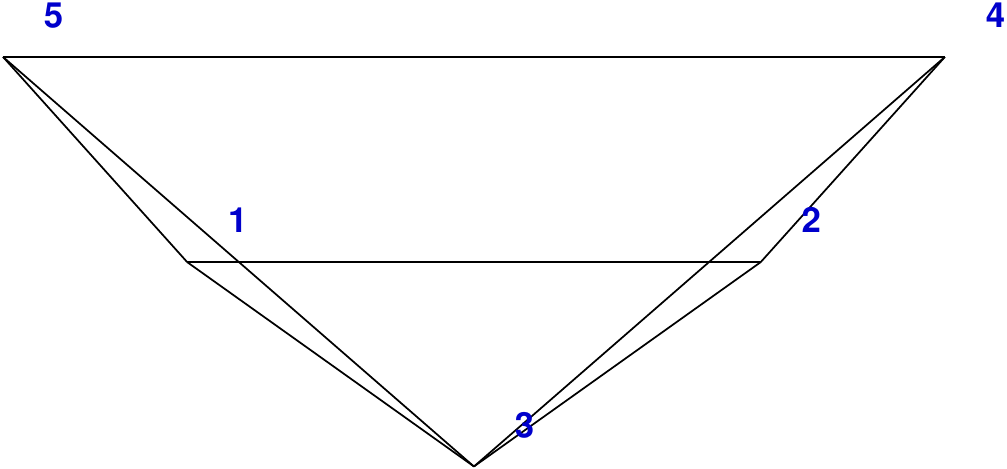}}
		\caption{$t=0.6$}
	\end{subfigure}
	\hfill
	\begin{subfigure}{0.3\textwidth}
		\centering
		\includegraphics[width=\linewidth]{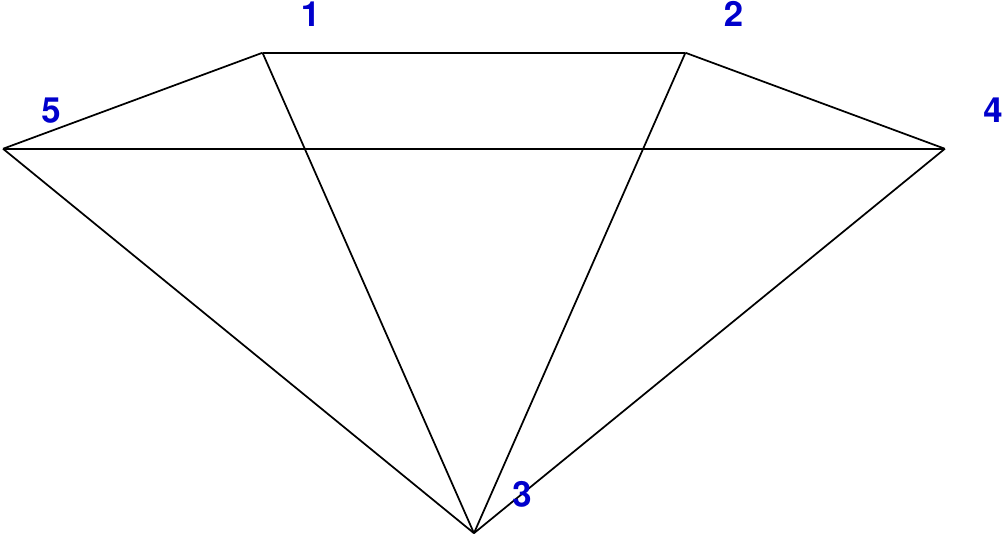}
		\caption{$t=0.9$}
	\end{subfigure}
	\caption{Snapshots of a geodesic on the manifold $\planarmanifold$ of planar triangular meshes with respect to the Euclidean Riemannian metric. Just before $t = 0.6$, the mesh degenerates.}
	\label{fig:euclidean_geodesic}
\end{figure}

The simplest choice of a Riemannian metric for $\planarmanifold$ is the Euclidean one, whose components (with respect to the $\vvec$ chart) are given by
\begin{equation}
	\label{eq:components_Euclidean_metric}
	\widetilde g_{ab}(Q)
	\equiv
	\delta_a^b
	,
	\quad
	a, b = 1, \ldots, 2 \, N_V
	.
\end{equation}
Here $\delta_a^b$ is the Kronecker delta symbol.
Geodesic curves~$\geodesic<s>$ with respect to $\widetilde g$ are trivial to compute.
They simply consist of a collection of straight lines, each emanating from a vertex and traveled along with constant Euclidean velocity.
More precisely, given a point $Q \in \planarmanifold$ and a tangent vector $V \in \R^{2 \times N_V}$, the unique geodesic curve $\geodesic<p>{Q}{V}$ passing through~$Q$ at time zero with velocity~$V$ is
\begin{equation*}
	\geodesic<p>{Q}{V}(t)
	=
	Q + t \, V
	.
\end{equation*}
Since the vertices move without any interference from their neighbors, geodesics \wrt\ $\widetilde g$ generally cease to exist in $\planarmanifold$ after a finite time, terminating in a degenerate mesh.
\Cref{fig:euclidean_geodesic} shows this situation with the Euclidean metric.

\begin{figure}[htb]
	\centering
	\includegraphics[width=0.2\linewidth]{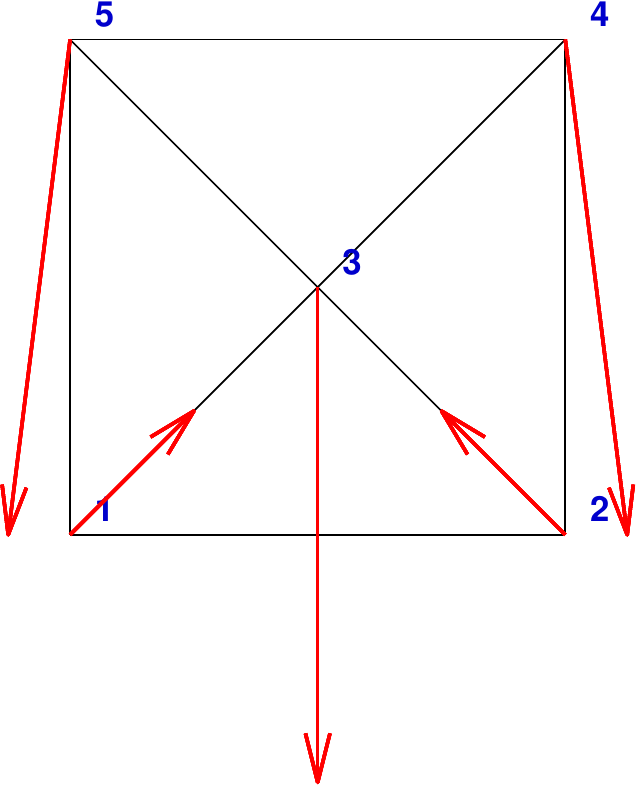}
	\caption{Initial tangent vector and mesh used as an example for the computation of geodesics on the manifold $\planarmanifold$ of planar triangular meshes.}
	\label{fig:euclidean_geodesic_tangent_vector}
\end{figure}

It is the main purpose of the present paper to define a \emph{complete} Riemannian metric instead, which avoids degenerate meshes by design.

\begin{definition}[\protect{\cite[Ch.~7, Def.~2.2]{DoCarmo1992}}]
	A Riemannian manifold $\cM$ is said to be \textbf{geodesically complete} if all geodesic curves $\geodesic<p>{Q}{V}(t)$ exist for all $t \in \R$.
\end{definition}

A recipe to construct complete Riemannian metrics on connected smooth manifolds such as $\planarmanifold$ was presented in \cite{Gordon1973}.
Recall that a function $f: \cM \to \R$ is \textbf{proper} if the preimages $f^{-1}(K)$ of compact sets $K \subset \R$ are compact in $\cM$.
\begin{theorem}[\protect{\cite[Thm.~1]{Gordon1973}}]
	\label{theorem:construction_of_complete_Riemannian_metric}
	Suppose that $\cM$ is a connected manifold of class $C^3$, endowed with a (not necessarily complete) Riemannian metric~$\widetilde g$ with component functions $\widetilde g_{ab}$.
	If $f \colon \cM \to \R$ is any proper function of class~$\cC^3$, then the Riemannian metric~$g$ defined in terms of coordinates $q^a$ by
	\begin{equation}
		\label{eq:augmentation_of_non-complete_Riemannian_metric}
		g_{ab}
		=
		\widetilde g_{ab} + \frac{\partial f}{\partial q^a} \, \frac{\partial f}{\partial q^b}
	\end{equation}
	is geodesically complete.
\end{theorem}
Moreover, \cite[Thm.~2]{Gordon1973} shows that this construction is the only way to obtain complete Riemannian metrics on connected smooth manifolds.

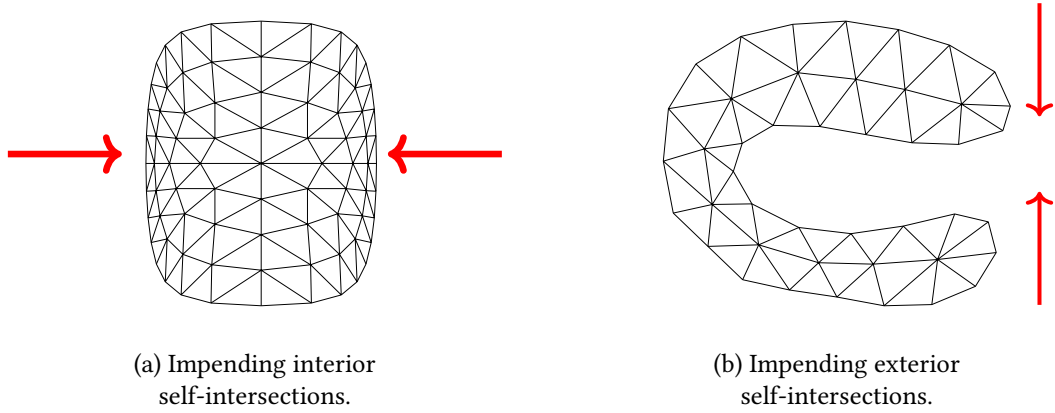
\begin{figure}[htb]
\captionsetup[subfigure]{justification=centering}
	\begin{center}
		\hfill
		\begin{subfigure}{0.455\textwidth}
			\centering
			\makeatletter
			\ltx@ifclassloaded{mcom-l}{
\begin{tikzpicture}[scale = 0.35]
	\draw [->,red,line width = 1.5pt](1,6)--(4,6);
	\draw [<-,red,line width = 1.5pt](11,6)--(14,6);
	}{%
	\begin{tikzpicture}[scale = 0.5]
		\draw [->,red,line width = 2.5pt](1,6)--(4,6);
		\draw [<-,red,line width = 2.5pt](11,6)--(14,6);
	}
	\makeatother
\begin{axis}[%
width=4.521in,
height=3.566in,
at={(0.758in,0.481in)},
scale only axis,
unbounded coords=jump,
xmin=-2.38992324567567,
xmax=2.38992324567567,
ymin=-1.88495559215388,
ymax=1.88495559215388,
axis line style={draw=none},
ticks=none,
legend style={legend cell align=left, align=left, draw=white!15!black}
]
\addplot [color=black]
  table[row sep=crcr]{%
-1.16942282481576	-1.26605553798781e-16\\
-1.26262725567891	-1.92367069372179e-16\\
nan	nan\\
1.29637203307674e-16	-1.17809724509617\\
1.92367069372179e-16	-1.5707963267949\\
nan	nan\\
1.16942282481576	1.76697482303529e-17\\
1.26262725567891	0\\
nan	nan\\
1.29637203307674e-16	1.17809724509617\\
1.92367069372179e-16	1.5707963267949\\
nan	nan\\
-1.03022178401005	-0.833040550904694\\
-1.147810884057	-1.11072073453959\\
nan	nan\\
1.03022178401005	-0.833040550904693\\
1.147810884057	-1.11072073453959\\
nan	nan\\
1.03022178401005	0.833040550904694\\
1.147810884057	1.11072073453959\\
nan	nan\\
-1.03022178401005	0.833040550904694\\
-1.147810884057	1.11072073453959\\
nan	nan\\
-0.665773750028354	4.91747734801386e-18\\
-5.85523948858412e-17	7.06789929214115e-17\\
nan	nan\\
4.17747117866381e-18	-0.392699081698724\\
-5.85523948858412e-17	7.06789929214115e-17\\
nan	nan\\
0.665773750028354	5.30092446910586e-17\\
-5.85523948858412e-17	7.06789929214115e-17\\
nan	nan\\
4.17747117866381e-18	0.392699081698724\\
-5.85523948858412e-17	7.06789929214115e-17\\
nan	nan\\
-0.506949338720432	-0.277680183634898\\
-5.85523948858412e-17	7.06789929214115e-17\\
nan	nan\\
0.506949338720432	-0.277680183634898\\
-5.85523948858412e-17	7.06789929214115e-17\\
nan	nan\\
0.506949338720432	0.277680183634898\\
-5.85523948858412e-17	7.06789929214115e-17\\
nan	nan\\
-0.506949338720432	0.277680183634898\\
-5.85523948858412e-17	7.06789929214115e-17\\
nan	nan\\
-1.15035476592186	-0.300558864942174\\
-1.23899735859896	-0.601117729884347\\
nan	nan\\
-1.10846264380141	-0.578239048577071\\
-1.23899735859896	-0.601117729884347\\
nan	nan\\
-0.541240955222119	-1.11831236973358\\
-0.876973216239766	-1.45122657606972\\
nan	nan\\
-0.857833089058431	-1.00329347166976\\
-0.876973216239766	-1.45122657606972\\
nan	nan\\
0.541240955222119	-1.11831236973358\\
0.876973216239766	-1.45122657606972\\
nan	nan\\
0.857833089058432	-1.00329347166976\\
0.876973216239766	-1.45122657606972\\
nan	nan\\
1.15035476592186	-0.300558864942173\\
1.23899735859896	-0.601117729884346\\
nan	nan\\
1.10846264380141	-0.578239048577071\\
1.23899735859896	-0.601117729884346\\
nan	nan\\
1.15035476592186	0.300558864942173\\
1.23899735859896	0.601117729884346\\
nan	nan\\
1.10846264380141	0.578239048577071\\
1.23899735859896	0.601117729884346\\
nan	nan\\
0.541240955222119	1.11831236973358\\
0.876973216239766	1.45122657606972\\
nan	nan\\
0.857833089058431	1.00329347166976\\
0.876973216239766	1.45122657606972\\
nan	nan\\
-0.541240955222119	1.11831236973358\\
-0.876973216239767	1.45122657606972\\
nan	nan\\
-0.857833089058432	1.00329347166976\\
-0.876973216239767	1.45122657606972\\
nan	nan\\
-1.15035476592186	0.300558864942173\\
-1.23899735859896	0.601117729884346\\
nan	nan\\
-1.10846264380141	0.578239048577071\\
-1.23899735859896	0.601117729884346\\
nan	nan\\
-1.16942282481576	-1.26605553798781e-16\\
-1.00388482185389	-6.08440382253838e-17\\
nan	nan\\
-0.665773750028354	4.91747734801386e-18\\
-1.00388482185389	-6.08440382253838e-17\\
nan	nan\\
-1.15035476592186	-0.300558864942174\\
-1.00388482185389	-6.08440382253838e-17\\
nan	nan\\
-1.15035476592186	0.300558864942173\\
-1.00388482185389	-6.08440382253838e-17\\
nan	nan\\
-0.929958789499237	-0.277680183634898\\
-1.00388482185389	-6.08440382253838e-17\\
nan	nan\\
-0.929958789499237	0.277680183634898\\
-1.00388482185389	-6.08440382253838e-17\\
nan	nan\\
1.29637203307674e-16	-1.17809724509617\\
6.69073372431689e-17	-0.785398163397448\\
nan	nan\\
4.17747117866381e-18	-0.392699081698724\\
6.69073372431689e-17	-0.785398163397448\\
nan	nan\\
-0.541240955222119	-1.11831236973358\\
6.69073372431689e-17	-0.785398163397448\\
nan	nan\\
0.541240955222119	-1.11831236973358\\
6.69073372431689e-17	-0.785398163397448\\
nan	nan\\
-0.506949338720432	-0.670379265333622\\
6.69073372431689e-17	-0.785398163397448\\
nan	nan\\
0.506949338720432	-0.670379265333622\\
6.69073372431689e-17	-0.785398163397448\\
nan	nan\\
1.16942282481576	1.76697482303529e-17\\
1.00388482185389	3.53394964607057e-17\\
nan	nan\\
0.665773750028354	5.30092446910586e-17\\
1.00388482185389	3.53394964607057e-17\\
nan	nan\\
1.15035476592186	-0.300558864942173\\
1.00388482185389	3.53394964607057e-17\\
nan	nan\\
1.15035476592186	0.300558864942173\\
1.00388482185389	3.53394964607057e-17\\
nan	nan\\
0.929958789499237	-0.277680183634898\\
1.00388482185389	3.53394964607057e-17\\
nan	nan\\
0.929958789499237	0.277680183634898\\
1.00388482185389	3.53394964607057e-17\\
nan	nan\\
1.29637203307674e-16	1.17809724509617\\
6.69073372431689e-17	0.785398163397448\\
nan	nan\\
4.17747117866381e-18	0.392699081698724\\
6.69073372431689e-17	0.785398163397448\\
nan	nan\\
0.541240955222119	1.11831236973358\\
6.69073372431689e-17	0.785398163397448\\
nan	nan\\
-0.541240955222119	1.11831236973358\\
6.69073372431689e-17	0.785398163397448\\
nan	nan\\
0.506949338720432	0.670379265333622\\
6.69073372431689e-17	0.785398163397448\\
nan	nan\\
-0.506949338720432	0.670379265333622\\
6.69073372431689e-17	0.785398163397448\\
nan	nan\\
-1.03022178401005	-0.833040550904694\\
-0.837806492185101	-0.555360367269796\\
nan	nan\\
-0.506949338720432	-0.277680183634898\\
-0.837806492185101	-0.555360367269796\\
nan	nan\\
-1.10846264380141	-0.578239048577071\\
-0.837806492185101	-0.555360367269796\\
nan	nan\\
-0.857833089058431	-1.00329347166976\\
-0.837806492185101	-0.555360367269796\\
nan	nan\\
-0.929958789499237	-0.277680183634898\\
-0.837806492185101	-0.555360367269796\\
nan	nan\\
-0.506949338720432	-0.670379265333622\\
-0.837806492185101	-0.555360367269796\\
nan	nan\\
1.03022178401005	-0.833040550904693\\
0.837806492185101	-0.555360367269796\\
nan	nan\\
0.506949338720432	-0.277680183634898\\
0.837806492185101	-0.555360367269796\\
nan	nan\\
0.857833089058432	-1.00329347166976\\
0.837806492185101	-0.555360367269796\\
nan	nan\\
1.10846264380141	-0.578239048577071\\
0.837806492185101	-0.555360367269796\\
nan	nan\\
0.506949338720432	-0.670379265333622\\
0.837806492185101	-0.555360367269796\\
nan	nan\\
0.929958789499237	-0.277680183634898\\
0.837806492185101	-0.555360367269796\\
nan	nan\\
1.03022178401005	0.833040550904694\\
0.837806492185101	0.555360367269796\\
nan	nan\\
0.506949338720432	0.277680183634898\\
0.837806492185101	0.555360367269796\\
nan	nan\\
1.10846264380141	0.578239048577071\\
0.837806492185101	0.555360367269796\\
nan	nan\\
0.857833089058431	1.00329347166976\\
0.837806492185101	0.555360367269796\\
nan	nan\\
0.929958789499237	0.277680183634898\\
0.837806492185101	0.555360367269796\\
nan	nan\\
0.506949338720432	0.670379265333622\\
0.837806492185101	0.555360367269796\\
nan	nan\\
-1.03022178401005	0.833040550904694\\
-0.837806492185101	0.555360367269796\\
nan	nan\\
-0.506949338720432	0.277680183634898\\
-0.837806492185101	0.555360367269796\\
nan	nan\\
-0.857833089058432	1.00329347166976\\
-0.837806492185101	0.555360367269796\\
nan	nan\\
-1.10846264380141	0.578239048577071\\
-0.837806492185101	0.555360367269796\\
nan	nan\\
-0.929958789499237	0.277680183634898\\
-0.837806492185101	0.555360367269796\\
nan	nan\\
-0.506949338720432	0.670379265333622\\
-0.837806492185101	0.555360367269796\\
nan	nan\\
-1.16942282481576	-1.26605553798781e-16\\
-1.25697515925325	-0.306447161216168\\
nan	nan\\
-1.15035476592186	-0.300558864942174\\
-1.25697515925325	-0.306447161216168\\
nan	nan\\
-1.03022178401005	-0.833040550904694\\
-1.05050956842496	-1.30606941284641\\
nan	nan\\
-0.857833089058431	-1.00329347166976\\
-1.05050956842496	-1.30606941284641\\
nan	nan\\
1.29637203307674e-16	-1.17809724509617\\
0.549846710530999	-1.5406139158319\\
nan	nan\\
0.541240955222119	-1.11831236973358\\
0.549846710530999	-1.5406139158319\\
nan	nan\\
1.03022178401005	-0.833040550904693\\
1.205180437008	-0.872687681303776\\
nan	nan\\
1.10846264380141	-0.578239048577071\\
1.205180437008	-0.872687681303776\\
nan	nan\\
1.16942282481576	1.76697482303529e-17\\
1.25697515925325	0.306447161216168\\
nan	nan\\
1.15035476592186	0.300558864942173\\
1.25697515925325	0.306447161216168\\
nan	nan\\
1.03022178401005	0.833040550904694\\
1.05050956842496	1.30606941284641\\
nan	nan\\
0.857833089058431	1.00329347166976\\
1.05050956842496	1.30606941284641\\
nan	nan\\
1.29637203307674e-16	1.17809724509617\\
-0.549846710530999	1.5406139158319\\
nan	nan\\
-0.541240955222119	1.11831236973358\\
-0.549846710530999	1.5406139158319\\
nan	nan\\
-1.03022178401005	0.833040550904694\\
-1.205180437008	0.872687681303775\\
nan	nan\\
-1.10846264380141	0.578239048577071\\
-1.205180437008	0.872687681303775\\
nan	nan\\
-1.03022178401005	-0.833040550904694\\
-1.205180437008	-0.872687681303776\\
nan	nan\\
-1.10846264380141	-0.578239048577071\\
-1.205180437008	-0.872687681303776\\
nan	nan\\
1.29637203307674e-16	-1.17809724509617\\
-0.549846710530997	-1.5406139158319\\
nan	nan\\
-0.541240955222119	-1.11831236973358\\
-0.549846710530997	-1.5406139158319\\
nan	nan\\
1.03022178401005	-0.833040550904693\\
1.05050956842496	-1.30606941284641\\
nan	nan\\
0.857833089058432	-1.00329347166976\\
1.05050956842496	-1.30606941284641\\
nan	nan\\
1.16942282481576	1.76697482303529e-17\\
1.25697515925325	-0.306447161216168\\
nan	nan\\
1.15035476592186	-0.300558864942173\\
1.25697515925325	-0.306447161216168\\
nan	nan\\
1.03022178401005	0.833040550904694\\
1.205180437008	0.872687681303776\\
nan	nan\\
1.10846264380141	0.578239048577071\\
1.205180437008	0.872687681303776\\
nan	nan\\
1.29637203307674e-16	1.17809724509617\\
0.549846710530998	1.5406139158319\\
nan	nan\\
0.541240955222119	1.11831236973358\\
0.549846710530998	1.5406139158319\\
nan	nan\\
-1.03022178401005	0.833040550904694\\
-1.05050956842496	1.30606941284641\\
nan	nan\\
-0.857833089058432	1.00329347166976\\
-1.05050956842496	1.30606941284641\\
nan	nan\\
-1.16942282481576	-1.26605553798781e-16\\
-1.25697515925325	0.306447161216167\\
nan	nan\\
-1.15035476592186	0.300558864942173\\
-1.25697515925325	0.306447161216167\\
nan	nan\\
-1.15035476592186	-0.300558864942174\\
-1.16942282481576	-1.26605553798781e-16\\
nan	nan\\
-1.15035476592186	0.300558864942173\\
-1.16942282481576	-1.26605553798781e-16\\
nan	nan\\
-0.506949338720432	-0.277680183634898\\
-0.665773750028354	4.91747734801386e-18\\
nan	nan\\
-0.929958789499237	-0.277680183634898\\
-0.665773750028354	4.91747734801386e-18\\
nan	nan\\
-0.506949338720432	0.277680183634898\\
-0.665773750028354	4.91747734801386e-18\\
nan	nan\\
-0.929958789499237	0.277680183634898\\
-0.665773750028354	4.91747734801386e-18\\
nan	nan\\
-1.10846264380141	-0.578239048577071\\
-1.15035476592186	-0.300558864942174\\
nan	nan\\
-0.929958789499237	-0.277680183634898\\
-1.15035476592186	-0.300558864942174\\
nan	nan\\
-1.10846264380141	0.578239048577071\\
-1.15035476592186	0.300558864942173\\
nan	nan\\
-0.929958789499237	0.277680183634898\\
-1.15035476592186	0.300558864942173\\
nan	nan\\
-0.541240955222119	-1.11831236973358\\
1.29637203307674e-16	-1.17809724509617\\
nan	nan\\
0.541240955222119	-1.11831236973358\\
1.29637203307674e-16	-1.17809724509617\\
nan	nan\\
-0.506949338720432	-0.277680183634898\\
4.17747117866381e-18	-0.392699081698724\\
nan	nan\\
-0.506949338720432	-0.670379265333622\\
4.17747117866381e-18	-0.392699081698724\\
nan	nan\\
0.506949338720432	-0.277680183634898\\
4.17747117866381e-18	-0.392699081698724\\
nan	nan\\
0.506949338720432	-0.670379265333622\\
4.17747117866381e-18	-0.392699081698724\\
nan	nan\\
-0.857833089058431	-1.00329347166976\\
-0.541240955222119	-1.11831236973358\\
nan	nan\\
-0.506949338720432	-0.670379265333622\\
-0.541240955222119	-1.11831236973358\\
nan	nan\\
0.857833089058432	-1.00329347166976\\
0.541240955222119	-1.11831236973358\\
nan	nan\\
0.506949338720432	-0.670379265333622\\
0.541240955222119	-1.11831236973358\\
nan	nan\\
1.15035476592186	-0.300558864942173\\
1.16942282481576	1.76697482303529e-17\\
nan	nan\\
1.15035476592186	0.300558864942173\\
1.16942282481576	1.76697482303529e-17\\
nan	nan\\
0.506949338720432	-0.277680183634898\\
0.665773750028354	5.30092446910586e-17\\
nan	nan\\
0.929958789499237	-0.277680183634898\\
0.665773750028354	5.30092446910586e-17\\
nan	nan\\
0.506949338720432	0.277680183634898\\
0.665773750028354	5.30092446910586e-17\\
nan	nan\\
0.929958789499237	0.277680183634898\\
0.665773750028354	5.30092446910586e-17\\
nan	nan\\
1.10846264380141	-0.578239048577071\\
1.15035476592186	-0.300558864942173\\
nan	nan\\
0.929958789499237	-0.277680183634898\\
1.15035476592186	-0.300558864942173\\
nan	nan\\
1.10846264380141	0.578239048577071\\
1.15035476592186	0.300558864942173\\
nan	nan\\
0.929958789499237	0.277680183634898\\
1.15035476592186	0.300558864942173\\
nan	nan\\
0.541240955222119	1.11831236973358\\
1.29637203307674e-16	1.17809724509617\\
nan	nan\\
-0.541240955222119	1.11831236973358\\
1.29637203307674e-16	1.17809724509617\\
nan	nan\\
0.506949338720432	0.277680183634898\\
4.17747117866381e-18	0.392699081698724\\
nan	nan\\
0.506949338720432	0.670379265333622\\
4.17747117866381e-18	0.392699081698724\\
nan	nan\\
-0.506949338720432	0.277680183634898\\
4.17747117866381e-18	0.392699081698724\\
nan	nan\\
-0.506949338720432	0.670379265333622\\
4.17747117866381e-18	0.392699081698724\\
nan	nan\\
0.857833089058431	1.00329347166976\\
0.541240955222119	1.11831236973358\\
nan	nan\\
0.506949338720432	0.670379265333622\\
0.541240955222119	1.11831236973358\\
nan	nan\\
-0.857833089058432	1.00329347166976\\
-0.541240955222119	1.11831236973358\\
nan	nan\\
-0.506949338720432	0.670379265333622\\
-0.541240955222119	1.11831236973358\\
nan	nan\\
-1.10846264380141	-0.578239048577071\\
-1.03022178401005	-0.833040550904694\\
nan	nan\\
-0.857833089058431	-1.00329347166976\\
-1.03022178401005	-0.833040550904694\\
nan	nan\\
-0.929958789499237	-0.277680183634898\\
-0.506949338720432	-0.277680183634898\\
nan	nan\\
-0.506949338720432	-0.670379265333622\\
-0.506949338720432	-0.277680183634898\\
nan	nan\\
-0.929958789499237	-0.277680183634898\\
-1.10846264380141	-0.578239048577071\\
nan	nan\\
-0.506949338720432	-0.670379265333622\\
-0.857833089058431	-1.00329347166976\\
nan	nan\\
0.857833089058432	-1.00329347166976\\
1.03022178401005	-0.833040550904693\\
nan	nan\\
1.10846264380141	-0.578239048577071\\
1.03022178401005	-0.833040550904693\\
nan	nan\\
0.506949338720432	-0.670379265333622\\
0.506949338720432	-0.277680183634898\\
nan	nan\\
0.929958789499237	-0.277680183634898\\
0.506949338720432	-0.277680183634898\\
nan	nan\\
0.506949338720432	-0.670379265333622\\
0.857833089058432	-1.00329347166976\\
nan	nan\\
0.929958789499237	-0.277680183634898\\
1.10846264380141	-0.578239048577071\\
nan	nan\\
1.10846264380141	0.578239048577071\\
1.03022178401005	0.833040550904694\\
nan	nan\\
0.857833089058431	1.00329347166976\\
1.03022178401005	0.833040550904694\\
nan	nan\\
0.929958789499237	0.277680183634898\\
0.506949338720432	0.277680183634898\\
nan	nan\\
0.506949338720432	0.670379265333622\\
0.506949338720432	0.277680183634898\\
nan	nan\\
0.929958789499237	0.277680183634898\\
1.10846264380141	0.578239048577071\\
nan	nan\\
0.506949338720432	0.670379265333622\\
0.857833089058431	1.00329347166976\\
nan	nan\\
-0.857833089058432	1.00329347166976\\
-1.03022178401005	0.833040550904694\\
nan	nan\\
-1.10846264380141	0.578239048577071\\
-1.03022178401005	0.833040550904694\\
nan	nan\\
-0.929958789499237	0.277680183634898\\
-0.506949338720432	0.277680183634898\\
nan	nan\\
-0.506949338720432	0.670379265333622\\
-0.506949338720432	0.277680183634898\\
nan	nan\\
-0.506949338720432	0.670379265333622\\
-0.857833089058432	1.00329347166976\\
nan	nan\\
-0.929958789499237	0.277680183634898\\
-1.10846264380141	0.578239048577071\\
nan	nan\\
};

\addplot [color=black]
  table[row sep=crcr]{%
-1.26262725567891	-1.92367069372179e-16\\
-1.25697515925325	-0.306447161216168\\
nan	nan\\
-1.147810884057	-1.11072073453959\\
-1.05050956842496	-1.30606941284641\\
nan	nan\\
1.92367069372179e-16	-1.5707963267949\\
0.549846710530999	-1.5406139158319\\
nan	nan\\
1.147810884057	-1.11072073453959\\
1.205180437008	-0.872687681303776\\
nan	nan\\
1.26262725567891	0\\
1.25697515925325	0.306447161216168\\
nan	nan\\
1.147810884057	1.11072073453959\\
1.05050956842496	1.30606941284641\\
nan	nan\\
1.92367069372179e-16	1.5707963267949\\
-0.549846710530999	1.5406139158319\\
nan	nan\\
-1.147810884057	1.11072073453959\\
-1.205180437008	0.872687681303775\\
nan	nan\\
-1.23899735859896	-0.601117729884347\\
-1.205180437008	-0.872687681303776\\
nan	nan\\
-0.876973216239766	-1.45122657606972\\
-0.549846710530997	-1.5406139158319\\
nan	nan\\
0.876973216239766	-1.45122657606972\\
1.05050956842496	-1.30606941284641\\
nan	nan\\
1.23899735859896	-0.601117729884346\\
1.25697515925325	-0.306447161216168\\
nan	nan\\
1.23899735859896	0.601117729884346\\
1.205180437008	0.872687681303776\\
nan	nan\\
0.876973216239766	1.45122657606972\\
0.549846710530998	1.5406139158319\\
nan	nan\\
-0.876973216239767	1.45122657606972\\
-1.05050956842496	1.30606941284641\\
nan	nan\\
-1.23899735859896	0.601117729884346\\
-1.25697515925325	0.306447161216167\\
nan	nan\\
-1.25697515925325	-0.306447161216168\\
-1.23899735859896	-0.601117729884347\\
nan	nan\\
-1.05050956842496	-1.30606941284641\\
-0.876973216239766	-1.45122657606972\\
nan	nan\\
0.549846710530999	-1.5406139158319\\
0.876973216239766	-1.45122657606972\\
nan	nan\\
1.205180437008	-0.872687681303776\\
1.23899735859896	-0.601117729884346\\
nan	nan\\
1.25697515925325	0.306447161216168\\
1.23899735859896	0.601117729884346\\
nan	nan\\
1.05050956842496	1.30606941284641\\
0.876973216239766	1.45122657606972\\
nan	nan\\
-0.549846710530999	1.5406139158319\\
-0.876973216239767	1.45122657606972\\
nan	nan\\
-1.205180437008	0.872687681303775\\
-1.23899735859896	0.601117729884346\\
nan	nan\\
-1.205180437008	-0.872687681303776\\
-1.147810884057	-1.11072073453959\\
nan	nan\\
-0.549846710530997	-1.5406139158319\\
1.92367069372179e-16	-1.5707963267949\\
nan	nan\\
1.05050956842496	-1.30606941284641\\
1.147810884057	-1.11072073453959\\
nan	nan\\
1.25697515925325	-0.306447161216168\\
1.26262725567891	0\\
nan	nan\\
1.205180437008	0.872687681303776\\
1.147810884057	1.11072073453959\\
nan	nan\\
0.549846710530998	1.5406139158319\\
1.92367069372179e-16	1.5707963267949\\
nan	nan\\
-1.05050956842496	1.30606941284641\\
-1.147810884057	1.11072073453959\\
nan	nan\\
-1.25697515925325	0.306447161216167\\
-1.26262725567891	-1.92367069372179e-16\\
nan	nan\\
};
\end{axis}
\end{tikzpicture}
			\caption{ Impending interior \\ self-intersections.}
		\end{subfigure}
		\hfill
		\begin{subfigure}{0.45\textwidth}
			\centering
			\makeatletter
			\ltx@ifclassloaded{mcom-l}{
\begin{tikzpicture}[scale = 0.35]
	\draw [->,red,line width = 1.5pt](13,10) -- (13,7);
	\draw [->,red,line width = 1.5pt](13, 2) -- (13,5);
	}{%
	\begin{tikzpicture}[scale = 0.5]
		\draw [->,red,line width = 1.5pt](13,10) -- (13,7);
		\draw [->,red,line width = 1.5pt](13, 2) -- (13,5);

	}
	\makeatother
\begin{axis}[%
width=4.521in,
height=3.566in,
at={(0.758in,0.481in)},
scale only axis,
unbounded coords=jump,
xmin=464.22715807673,
xmax=632.423888373357,
ymin=316.857316502618,
ymax=449.515705397861,
axis line style={draw=none},
ticks=none,
legend style={legend cell align=left, align=left, draw=white!15!black}
]
\addplot [color=black]
  table[row sep=crcr]{%
500.302597593591	367.250476481658\\
485.28	363.91\\
nan	nan\\
596.42287122672	406.211767461236\\
609.25	418.74\\
nan	nan\\
596.42287122672	406.211767461236\\
615.21	405.46\\
nan	nan\\
596.42287122672	406.211767461236\\
612.32	395.92\\
nan	nan\\
507.685467823052	408.342358583091\\
523.77	397.97\\
nan	nan\\
533.368157854896	418.41255299619\\
523.77	397.97\\
nan	nan\\
497.240962483107	388.842474184609\\
511.85	390.98\\
nan	nan\\
507.685467823052	408.342358583091\\
511.85	390.98\\
nan	nan\\
500.302597593591	367.250476481658\\
508.44	380.09\\
nan	nan\\
497.240962483107	388.842474184609\\
508.44	380.09\\
nan	nan\\
500.302597593591	367.250476481658\\
515.59	367.49\\
nan	nan\\
518.212064306024	351.577951933301\\
515.59	367.49\\
nan	nan\\
587.161393569587	345.872855717669\\
593.75	363.57\\
nan	nan\\
587.161393569587	345.872855717669\\
606.7	360.5\\
nan	nan\\
587.161393569587	345.872855717669\\
609.76	348.58\\
nan	nan\\
587.161393569587	345.872855717669\\
601.76	335.47\\
nan	nan\\
518.212064306024	351.577951933301\\
512.02	338.03\\
nan	nan\\
574.611496990506	411.895424755854\\
591.753702564334	429.137902407823\\
nan	nan\\
596.42287122672	406.211767461236\\
591.753702564334	429.137902407823\\
nan	nan\\
555.6445693258	415.862346803508\\
572.092218663461	435.171784436712\\
nan	nan\\
574.611496990506	411.895424755854\\
572.092218663461	435.171784436712\\
nan	nan\\
555.6445693258	415.862346803508\\
551.838078686156	438.460839656591\\
nan	nan\\
533.368157854896	418.41255299619\\
551.838078686156	438.460839656591\\
nan	nan\\
533.368157854896	418.41255299619\\
531.344275569228	437.225770968946\\
nan	nan\\
507.685467823052	408.342358583091\\
511.381629732665	432.515616401857\\
nan	nan\\
533.368157854896	418.41255299619\\
511.381629732665	432.515616401857\\
nan	nan\\
507.685467823052	408.342358583091\\
494.010112577455	421.81714701222\\
nan	nan\\
497.240962483107	388.842474184609\\
483.420365055616	404.465148401567\\
nan	nan\\
507.685467823052	408.342358583091\\
483.420365055616	404.465148401567\\
nan	nan\\
500.302597593591	367.250476481658\\
481.441046450087	384.070758625088\\
nan	nan\\
497.240962483107	388.842474184609\\
481.441046450087	384.070758625088\\
nan	nan\\
555.6445693258	415.862346803508\\
541.779775624362	397.527206859741\\
nan	nan\\
533.368157854896	418.41255299619\\
541.779775624362	397.527206859741\\
nan	nan\\
555.6445693258	415.862346803508\\
559.51688246726	394.383281231928\\
nan	nan\\
574.611496990506	411.895424755854\\
559.51688246726	394.383281231928\\
nan	nan\\
574.611496990506	411.895424755854\\
577.257981983122	391.216058236029\\
nan	nan\\
596.42287122672	406.211767461236\\
577.257981983122	391.216058236029\\
nan	nan\\
596.42287122672	406.211767461236\\
595.247610367635	390.580847445549\\
nan	nan\\
500.302597593591	367.250476481658\\
498.65	350.97\\
nan	nan\\
518.212064306024	351.577951933301\\
498.65	350.97\\
nan	nan\\
518.212064306024	351.577951933301\\
530.06592447438	334.038357435121\\
nan	nan\\
541.308968639345	344.626783783843\\
530.06592447438	334.038357435121\\
nan	nan\\
541.308968639345	344.626783783843\\
548.355162874481	331.307403794776\\
nan	nan\\
562.369748667569	344.075811834234\\
548.355162874481	331.307403794776\\
nan	nan\\
587.161393569587	345.872855717669\\
566.524908573656	327.912182243888\\
nan	nan\\
562.369748667569	344.075811834234\\
566.524908573656	327.912182243888\\
nan	nan\\
587.161393569587	345.872855717669\\
584.941112396442	328.374927989866\\
nan	nan\\
587.161393569587	345.872855717669\\
573.942138211584	358.857353551168\\
nan	nan\\
562.369748667569	344.075811834234\\
573.942138211584	358.857353551168\\
nan	nan\\
541.308968639345	344.626783783843\\
553.815085824668	355.916771829901\\
nan	nan\\
562.369748667569	344.075811834234\\
553.815085824668	355.916771829901\\
nan	nan\\
518.212064306024	351.577951933301\\
533.678291548173	358.313076914187\\
nan	nan\\
541.308968639345	344.626783783843\\
533.678291548173	358.313076914187\\
nan	nan\\
562.369748667569	344.075811834234\\
587.161393569587	345.872855717669\\
nan	nan\\
574.611496990506	411.895424755854\\
555.6445693258	415.862346803508\\
nan	nan\\
533.368157854896	418.41255299619\\
555.6445693258	415.862346803508\\
nan	nan\\
596.42287122672	406.211767461236\\
574.611496990506	411.895424755854\\
nan	nan\\
518.212064306024	351.577951933301\\
500.302597593591	367.250476481658\\
nan	nan\\
497.240962483107	388.842474184609\\
500.302597593591	367.250476481658\\
nan	nan\\
541.308968639345	344.626783783843\\
518.212064306024	351.577951933301\\
nan	nan\\
507.685467823052	408.342358583091\\
497.240962483107	388.842474184609\\
nan	nan\\
562.369748667569	344.075811834234\\
541.308968639345	344.626783783843\\
nan	nan\\
533.368157854896	418.41255299619\\
507.685467823052	408.342358583091\\
nan	nan\\
};

\addplot [color=black]
  table[row sep=crcr]{%
523.77	397.97\\
541.779775624362	397.527206859741\\
nan	nan\\
541.779775624362	397.527206859741\\
559.51688246726	394.383281231928\\
nan	nan\\
559.51688246726	394.383281231928\\
577.257981983122	391.216058236029\\
nan	nan\\
577.257981983122	391.216058236029\\
595.247610367635	390.580847445549\\
nan	nan\\
595.247610367635	390.580847445549\\
612.32	395.92\\
nan	nan\\
511.85	390.98\\
523.77	397.97\\
nan	nan\\
508.44	380.09\\
511.85	390.98\\
nan	nan\\
609.25	418.74\\
591.753702564334	429.137902407823\\
nan	nan\\
591.753702564334	429.137902407823\\
572.092218663461	435.171784436712\\
nan	nan\\
572.092218663461	435.171784436712\\
551.838078686156	438.460839656591\\
nan	nan\\
551.838078686156	438.460839656591\\
531.344275569228	437.225770968946\\
nan	nan\\
531.344275569228	437.225770968946\\
511.381629732665	432.515616401857\\
nan	nan\\
511.381629732665	432.515616401857\\
494.010112577455	421.81714701222\\
nan	nan\\
494.010112577455	421.81714701222\\
483.420365055616	404.465148401567\\
nan	nan\\
483.420365055616	404.465148401567\\
481.441046450087	384.070758625088\\
nan	nan\\
481.441046450087	384.070758625088\\
485.28	363.91\\
nan	nan\\
615.21	405.46\\
609.25	418.74\\
nan	nan\\
515.59	367.49\\
508.44	380.09\\
nan	nan\\
612.32	395.92\\
615.21	405.46\\
nan	nan\\
485.28	363.91\\
498.65	350.97\\
nan	nan\\
498.65	350.97\\
512.02	338.03\\
nan	nan\\
512.02	338.03\\
530.06592447438	334.038357435121\\
nan	nan\\
530.06592447438	334.038357435121\\
548.355162874481	331.307403794776\\
nan	nan\\
548.355162874481	331.307403794776\\
566.524908573656	327.912182243888\\
nan	nan\\
566.524908573656	327.912182243888\\
584.941112396442	328.374927989866\\
nan	nan\\
584.941112396442	328.374927989866\\
601.76	335.47\\
nan	nan\\
601.76	335.47\\
609.76	348.58\\
nan	nan\\
609.76	348.58\\
606.7	360.5\\
nan	nan\\
606.7	360.5\\
593.75	363.57\\
nan	nan\\
593.75	363.57\\
573.942138211584	358.857353551168\\
nan	nan\\
573.942138211584	358.857353551168\\
553.815085824668	355.916771829901\\
nan	nan\\
553.815085824668	355.916771829901\\
533.678291548173	358.313076914187\\
nan	nan\\
533.678291548173	358.313076914187\\
515.59	367.49\\
nan	nan\\
};

\end{axis}
\end{tikzpicture}%
			\caption{Impending exterior \\ self-intersections.}
		\end{subfigure}
		\hspace*{\fill}
	\end{center}
	\caption{Different kinds of impending self-intersections.}
	\label{fig:self-intersections}
\end{figure}

As we pointed out in the discussion of the one-dimensional example from \Cref{section:motivation_1D}, intuitively, it is clear that the metric coefficients must be large whenever the mesh is close to a situation of self-intersection.
These self-intersections can be of internal or external nature; see \Cref{fig:self-intersections}.
Internal self-intersections are impending whenever a triangle is about to collapse to a line segment or even a point.
Exterior self-intersections can be recognized by the distance of one boundary face to another boundary face becoming small.
(The definition of boundary and interior faces can be found in the end of \cref{subsection:geometric_simplicial_complexes}.)
These situations, and the metric which prevents them, can be expressed in terms of the heights of the triangle (for internal self-intersections) and the distance of non-incident boundary vertices and edges  (for external self-intersections).

In the following, we are going to use different notions of distances between geometric objects.
The Euclidean distance between two vertices $q$ and $q'$ is going to be denoted by $\norm{q - q'}$.
When we wish to emphasize the vertex numbers in the simplicial complex and the dependence on the matrix~$Q$ of vertex positions, we shall use the alternative notation $\distance{Q}[i_0][i_1] = \norm{Q_{i_0} - Q_{i_1}}$ instead.
The notation $\distance{Q}[\cdot\,][\cdot]$ is also used to denote the Euclidean distance between higher dimensional geometric objects.
To clarify this notion, we consider the following definition.
\begin{definition}
	\label{definition:distanceConvexSets}
	Let $\cC$ denote the collection of all nonempty, convex, compact subsets of $\R^2$.
	A distance on $\cC$ is a mapping $d \colon \cC \times \cC \to \R$ satisfying the following properties for all $A, B \in \cC$:
	\begin{enumeratelatin}
		\item
			\label[condition]{item:distanceIsNonNegative}
			$\distance{}[A][B] \geq 0$.
		\item
			\label[condition]{item:distanceIsSymmetric}
			$\distance{}[A][B] = \distance{}[B][A]$.
		\item
			\label[condition]{item:distanceIsPositive}
			$\distance{}[A][B]>0$ if and only if $A\cap B = \emptyset$.
	\end{enumeratelatin}
\end{definition}
Consider the function 
\begin{equation}
	\label{eq:distanceConvexSets}
	\distance{}[A][B]
		=
		\min
		\setDef[big]{%
		\norm{a-b}}
		{%
	a\in A \text{ and } b\in B},
\end{equation}
where $\norm{\cdot}$ is a norm on $\R^2$.
It is easy to prove that $\distance{}[A][B]$ given in~\eqref{eq:distanceConvexSets} is well-defined and it is a distance on $\cC$ in the sense of \cref{definition:distanceConvexSets}.
This allows us to define
\begin{equation*}
	\distance{Q}[i_0][[j_0,j_1]]
	\coloneqq
	\min\setDef[big]{\norm{Q_{i_0} - x}}{x \in \conv\{Q_{j_0}, Q_{j_1}\}}
	,
\end{equation*}
\ie, the Euclidean distance between vertex~$i_0$ and the edge spanned by vertices~$j_0$ and $j_1$.
Similarly, $\distance{Q}[[i_0,i_1]][[j_0,j_1]]$ denotes the Euclidean distance between the edges $[i_0,i_1]$ and $[j_0,j_1]$.
In the case of vertex-edge distances, we will also use a distance $\Distance{Q}[i_0][[j_0,j_1]]$ different from the Euclidean one; see \cref{definition:distance_vertex_edge}.

\begin{definition}
	\label{definition:edge_lengths_heights}
	Let us denote as $[i_0,i_1,i_2]$ an arbitrary $2$-face of $\Delta$.
	Given $Q \in \plusmanifold$, we denote by
	\begin{equation}
		\label{eq:edge_length}
		\edgelength{Q}{\ell}[i_0,i_1,i_2]
		\coloneqq
		\norm{Q_{i_{\ell \oplus 1}} - Q_{i_{\ell \oplus 2}}}
		=
		\distance{Q}[i_{\ell \oplus 1}][i_{\ell \oplus 2}]
		,
		\quad
		\ell = 0, 1, 2
	\end{equation}
	the \textbf{length of the $\ell$-th edge} (the one opposite the $\ell$-th vertex).
	Here $\oplus$ denotes addition modulo~3.
	Recall from \eqref{eq:signed_area} the definition of the signed area $\area{Q}[i_0,i_1,i_2]$ and let
	\begin{equation}
		\label{eq:triangle_height}
		\height[auto]{Q}{\ell}[i_0,i_1,i_2]
		=
		\frac{2 \, \area{Q}[i_0,i_1,i_2]}{\edgelength{Q}{\ell}[i_0,i_1,i_2]}
		,
		\quad
		\ell = 0, 1, 2
	\end{equation}
	be the \textbf{$\ell$-th height} (the one through the $\ell$-th vertex); see, \eg, \cite[Ch.~2, Rem.~7, p.41]{AgricolaFriedrich2008}.
\end{definition}

Notice that the signed areas $\area{Q}$ are positive for $Q \in \plusmanifold$ and thus all heights in \eqref{eq:triangle_height} are positive as well.

\begin{figure}[htb]
	\centering
	\begin{tikzpicture}[dot/.style={circle,inner sep=1pt,fill,name=#1},
	  extended line/.style={shorten >=-#1,shorten <=-#1},
	  extended line/.default=1cm]

	\draw [thick] (-1,0.5) --  (1,1);
	\node [dot=v1] at (-1,0.5) {};
	\node [dot=v2] at (1,1) {};
	\node [dot=q]   at (0,2.5) {};
	\draw [red,dashed] ($(v1)!(q)!(v2)$) -- (q);

	\node[below] at (-1,0.5)   {$Q_{j_0}$};
	\node[below] at (1,1)     {$Q_{j_1}$};
	\node[left] at (0,2.5)   {$Q_{i_0}$};

	\draw [decorate,decoration={brace,amplitude=5pt},xshift=-4pt,yshift=0pt]
	(0.2,2.5) -- (0.6,0.9)node [black,midway,xshift=40pt] {$\Distance{Q}[i_0][[j_0,j_1]]$};

	\draw [thick] (4.5,0.5) --  (7,1);
	\node [dot=v1] at (4.5,0.5) {};
	\node [dot=v2] at (7,1) {};
	\node [dot=q]   at (8,2) {};
	\draw [red,dashed] ($(v1)!(q)!(v2)$) -- (q);
	\draw [red,dashed] (7,1) -- (8.15,1.23);

	\node[below] at (4.5,0.5) {$Q_{j_0}$};
	\node[above] at (7,1)   {$Q_{j_1}$};
	\node[above] at (8,2)   {$Q_{i_0}$};
	\node at (5.5,1.7) {$\Distance{Q}[i_0][[j_0,j_1]]=d_1+d_2$};

	\draw [decorate,decoration={brace,amplitude=3pt,raise=4pt},yshift=0pt]
	(8,2) -- (8.15,1.23) node [black,midway,xshift=18pt] {$d_1$};
	\draw [decorate,decoration={brace,amplitude=5pt,mirror,raise=4pt},yshift=0pt]
	(7,1) -- (8.15,1.23) node [black,midway,yshift=-18pt] {$d_2$};
\end{tikzpicture}

	\caption{Illustration of the distance \eqref{eq:distance_vertex_edge} of a $0$-face to a $1$-face in an edge oriented coordinate system. The two cases shown are when the projection of the $0$-face onto the infinite line generated by the $1$-face belongs to the $1$-face (left), and when it does not (right).}
	\label{fig:illustration_distance_vertex_edge}
\end{figure}
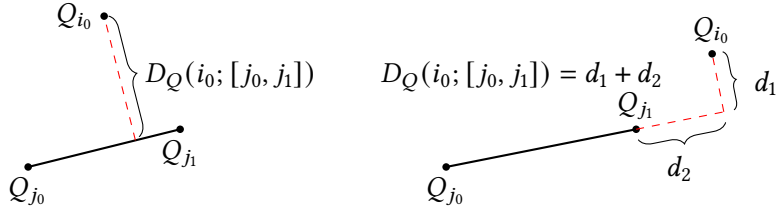

\begin{definition}
	\label{definition:distance_vertex_edge}
	Suppose that $[i_0]$ is a $0$-face and $[j_0,j_1]$ is a $1$-face of $\Delta$.
	We define the \textbf{$1$-norm distance of a $0$-face to a $1$-face}
	\begin{equation*}
		\Distance{Q}[i_0][[j_0,j_1]]
		\coloneqq
		\min\setDef[big]{\norm{Q_{i_0} - x}_1}{x \in \conv\{Q_{j_0}, Q_{j_1}\}}
		,
	\end{equation*}
	where the $1$-norm is based on an edge oriented coordinate system; see \Cref{fig:illustration_distance_vertex_edge} for an illustration.
\end{definition}

It can easily be shown that
\begin{subequations}
	\label{eq:distance_vertex_edge}
	\begin{equation}
		\Distance{Q}[i_0][[j_0,j_1]]
		=
		\norm{Q_{i_0} - Q_{j_0} + t \, (Q_{j_1} - Q_{j_0})}_1
	\end{equation}
	holds, where
	\begin{equation}
		t
		=
		\max\paren[auto]\{\}{0, \min\paren[auto]\{\}{1,\frac{(Q_{i_0} - Q_{j_0})^\transp (Q_{j_1} - Q_{j_0})}{\norm{Q_{j_1} - Q_{j_0}}^2}}}
		.
	\end{equation}
\end{subequations}
Clearly, the well-known equivalence of norms implies
\begin{equation}
	\label{eq:equivalence_of_1-_and_2-norms}
	\distance{Q}[i_0][[j_0,j_1]]
	\le
	\Distance{Q}[i_0][[j_0,j_1]]
	\le
	\sqrt{2} \, \distance{Q}[i_0][[j_0,j_1]]
	.
\end{equation}

We are now in the position to define a preliminary proper function~$f$ which is going to help render the shape manifold $\planarmanifold$ geodesically complete by means of \cref{theorem:construction_of_complete_Riemannian_metric}.
\begin{definition}
	\label{definition:f}
	We denote by $V_\partial$ the set of the boundary $0$-faces and by $E_\partial$ the set of boundary $1$-faces.
	Suppose that the $2$-faces in $\Delta$ are numbered from~$1$ to~$N_T$.
	We define the function $f \colon \plusmanifold \to \R$ by
	\begin{equation}
		\label{eq:f}
		f(Q)
		\coloneqq
		\sum_{k=1}^{N_T} \sum_{\ell=0}^2 \frac{\beta_1}{\height[auto]{Q}{\ell}[i^k_0,i^k_1,i^k_2]}
		+
		\sum_{[j_0,j_1]\in E_\partial}\sum_{\substack{i_0\in V_\partial \\ i_0\neq j_0,j_1}} \frac{\beta_2}{\Distance{Q}[i_0][[j_0,j_1]]}
		+
		\frac{\beta_3}{2} \norm{Q - \Qref}_F^2
		.
	\end{equation}
	Here $\beta_1,\beta_2,\beta_3$ are non-negative parameters, $\height{Q}{\ell}$ are the heights given in \eqref{eq:triangle_height} and $\Distance{Q}$ is the $1$-norm based distance of a vertex to an edge given in \eqref{eq:distance_vertex_edge}.
	Moreover, $\Qref \in \plusmanifold$ serves as a reference configuration and $\norm{\cdot}_F$ denotes the Frobenius norm.
\end{definition}

We remark that the first term in \eqref{eq:f} is designed to avoid interior self-inter\-sec\-tions, which go along with at least one height in a triangle converging to zero.
The second term avoids exterior self-interactions.
The third term penalizes large deviations from the reference mesh.
All three terms are required in order to show the properness of $f$ in \cref{theorem:f_is_proper}.
In preparation, we study several properties of $f$ in what follows.

\begin{lemma}
	\label{lemma:f_is_well-defined_and_continuous}
	For any choice of $\beta_1, \beta_2, \beta_3 \ge 0$, the function~$f$ defined in \eqref{eq:f} is well-defined on $\plusmanifold$ and continuous with values in $(0,\infty)$.
\end{lemma}
\begin{proof}
	From the definition of $\plusmanifold$ it is clear that all areas and the lengths of the $1$-faces are strictly positive and thus the same is true for the heights in \eqref{eq:triangle_height}. In the same way, all the distances from a $0$-face to a $1$-face given in \eqref{eq:distance_vertex_edge} are positive.
	The continuity of $f$ \wrt\ $Q$ on $\plusmanifold$ is obvious.
\end{proof}

The major next step is to prove the properness of $f$ defined in \eqref{eq:f}.
In order to make the proof more readable we present some intermediate results whose proofs are given in \Cref{section:proofs}.
The first result shows important bounds on the heights~$\height{Q}{\ell}$ \eqref{eq:triangle_height}, lengths of $1$-faces~$\edgelength{Q}{\ell}$ \eqref{eq:edge_length}, signed area~$\area{Q}$ \eqref{eq:signed_area}, inradius~$r_Q$ and circumradius~$\circumradius{Q}$ (see \cref{subsection:triangles}) of the $2$-faces of $\Delta$ in terms of the value of~$f$.
Moreover, we prove bounds on the interior angles $\interiorangle{Q}{\ell}$ of the $2$-faces.

\begin{lemma}
	\label{lemma:bounds_of_triangle_properties_in_terms_of_f}
	Let $[i_0,i_1,i_2]$ be an arbitrary $2$-face of $\Delta$.
	Suppose that $\beta_1, \beta_2 \ge 0$ and $\beta_3 > 0$ holds and $Q \in \planarmanifold$.
	Then the following statements hold.
	\begin{enumeratelatin}
		\item \label[statement]{item:bound_height}
			The heights satisfy $\height{Q}{\ell}[i_0,i_1,i_2] \ge \frac{\beta_1}{f(Q)}$, for all $\ell = 0,1,2$.
		\item \label[statement]{item:bound_inradius}
			The inradius satisfies $\inradius{Q}[i_0,i_1,i_2] \ge \frac{\beta_1}{f(Q)}$.
		\item \label[statement]{item:bound_length_edges}
			The lengths of the $1$-faces satisfy $\frac{2 \beta_1}{f(Q)} \le \edgelength{Q}{\ell}[i_0,i_1,i_2] \le 2 \sqrt{\frac{f(Q)}{\beta_3}} + \sqrt{2} \, \norm{\Qref}_F$, for all $\ell = 0,1,2$.
		\item \label[statement]{item:bound_area}
			The area satisfies $\area{Q}[i_0,i_1,i_2]\ge \frac{\pi\beta_1^2}{f^2(Q)}$.
		\item \label[statement]{item:bound_ratio_r_R}
			The inradius and the circumradius satisfy
			\begin{equation*}
				\frac{\inradius{Q}[i_0,i_1,i_2]}{\circumradius{Q}[i_0,i_1,i_2]}
				\ge
				\frac{4\pi\beta_1^3}{\paren[auto](){2 \sqrt{\frac{f(Q)}{\beta_3}} + \sqrt{2} \, \norm{\Qref}_F}^3} \frac{1}{f(Q)^3}
				.
			\end{equation*}
		\item \label[statement]{item:bound_cosines}
			There exists a function $\Psi \colon (0,\infty) \to \R$ which is monotone increasing and takes values in $[0,1)$ such that
			\begin{equation*}
				\abs[big]{\cos\paren[big](){\interiorangle{Q}{\ell}[i_0,i_1,i_2]}}
				\le
				\Psi(f(Q))
				\quad\text{for all } \ell=0,1,2
				.
			\end{equation*}
			In particular, $\abs[big]{\cos\paren[big](){\interiorangle{Q}{\ell}[i_0,i_1,i_2]}} \le \Psi(b) < 1$ holds whenever $f(Q) \le b$.
	\end{enumeratelatin}
\end{lemma}

The next result establishes bounds on the distance from a $0$-face to a $1$-face in terms of the function $f$.
Here we distinguish between boundary and interior $0$- and $1$-faces.
Briefly, a $1$-face is interior if it is shared by two distinct $2$-faces.
Otherwise it belongs to only one $2$-face and is referred to as a boundary $1$-face.
A $0$-face is a called a boundary $0$-face if it belongs to at least one boundary $1$-face.
Otherwise, all of its incident $1$-faces are interior and the $0$-face will be referred to as interior as well.
We refer the reader to \cref{subsection:geometric_simplicial_complexes} for details.

\begin{proposition}
	\label{proposition:bounds_of_distances_in_terms_of_f}
	Suppose that $Q \in \plusmanifold$.
	We consider a $0$-face $[i_0]$ and a $1$-face $[j_0,j_1]$ of $\Delta$ such that $\{Q_{i_0}\} \cap \conv\{Q_{j_0},Q_{j_1}\} = \emptyset$.
	Suppose $\beta_1, \beta_2, \beta_3 \ge 0$.
	Then the following statements hold.
	\begin{enumeratelatin}
		\item \label[statement]{item:distance_boundary_vertex_boundary_edge}
			If $[i_0]$ and $[j_0,j_1]$ are boundary, then
			\begin{equation}
				\label{eq:distance_boundary_vertex_boundary_edge}
				\Distance{Q}[i_0][[j_0,j_1]]
				\ge
				\frac{\beta_2}{f(Q)}
				.
			\end{equation}

		\item \label[statement]{item:distance_interior_vertex_arbitrary_edge}
			If $[i_0]$ is interior and $[j_0,j_1]$ is interior or boundary, then
			\begin{equation}
				\label{eq:distance_interior_vertex_arbitrary_edge}
				\Distance{Q}[i_0][[j_0,j_1]]
				\ge
				\frac{\beta_1}{f(Q)}
				.
			\end{equation}

		\item \label[statement]{item:distance_boundary_vertex_interior_edge}
			If $[i_0]$ is boundary and $[j_0,j_1]$ is interior, then
			\begin{equation}
				\label{eq:distance_boundary_vertex_interior_edge}
				\Distance{Q}[i_0][[j_0,j_1]]
				\ge
				\frac{\min\{\beta_1,\beta_2\}}{\sqrt{2}f(Q)} \min_{\theta \in \Theta} \paren[auto]\{\}{\sqrt{1-\cos^2(\theta)}}
				,
			\end{equation}
			where $\Theta$ is the set of four angles formed by the edge~$[j_0,j_1]$ and the adjacent $1$-faces belonging to $\closedSt{[j_0,j_1]}$; see \Cref{fig:construction_interior_edge_lemma_C1}.
	\end{enumeratelatin}
\end{proposition}

We are now in a position to present the main theorem of this section.
\begin{theorem}
	\label{theorem:f_is_proper}
	Suppose that $\beta_1, \beta_2, \beta_3 > 0$ holds.
	Then the restriction of~$f$ to $\planarmanifold$ is proper.
\end{theorem}
\begin{proof}
	Suppose that $K \subset \R$ is an arbitrary, compact set.
	We can suppose that $f^{-1}(K)$ is non-empty.
	Consequently, we have $K \subset [a,b]$ for some $b > 0$.
	Since $\planarmanifold$ carries the metric subspace topology of $\R^{2\times N_V}$, the compactness of $f^{-1}(K)$ agree with its sequential compactness.
	To verify the latter, suppose that $(Q^n) \subset f^{-1}(K) \subset \planarmanifold$ is an arbitrary sequence.
	We will show that it contains a subsequence which converges in $\planarmanifold$.

	The definition of~$f$ implies
	\begin{equation*}
		\frac{\beta_3}{2} \norm{Q^n - \Qref}_F^2
		\le
		f(Q^n)
		\le
		b
	\end{equation*}
	for all $n \in \N$ and thus $(Q^n)$ is bounded.
	Consequently, there exists a subsequence (which we do not relabel) converging to some $Q^*$ in $\R^{2 \times N_V}$.
	It remains to prove that $Q^*$ belongs to $\planarmanifold$ and that $Q^* \in f^{-1}(K)$ holds.
	We proceed by proving the following results about the limit configuration~$Q^*$:
	\begin{enumerate}
		\item
			\label[statement]{item:limit_has_positive_areas}
			The signed area $\area[auto]{Q^*}[i^k_0,i^k_1,i^k_2]$ of each triangle $k = 1, \ldots, N_T$ is strictly positive.
			In particular, the points $\{Q^*_{i^k_0}, Q^*_{i^k_1}, Q^*_{i^k_2}\}$ are affine independent.

		\item
			\label[statement]{item:limit_is_a_simplicial_complex}
			$\Sigma_\Delta(Q^*)$ is a simplicial $2$-complex whose associated abstract simplicial complex is $\Delta$.
	\end{enumerate}
	\Cref{item:limit_has_positive_areas,item:limit_is_a_simplicial_complex} together prove that $Q^*$ belongs to $\plusmanifold$.
	From there we will proceed to show that
	\begin{enumerate}[resume]
		\item
			\label[statement]{item:limit_has_a_path_to_Qref}
			$Q^*$ belongs to $\planarmanifold$ and $Q^* \in f^{-1}(K)$ holds.
	\end{enumerate}

	To show \cref{item:limit_has_positive_areas}, fix an arbitrary oriented $2$-face $[i_0,i_1,i_2]$ of $\Delta$.
	Thanks to \cref{lemma:bounds_of_triangle_properties_in_terms_of_f}, \cref{item:bound_area} we know
	\begin{equation*}
		\area{Q^n}[i_0,i_1,i_2]
		\ge
		\dfrac{\pi \beta_1^2}{f^2(Q^n)}\ge \frac{\pi \beta_1^2}{b^2}
		.
	\end{equation*}
	Since $\area{Q}[i_0,i_1,i_2]$ depends continuously on $Q$, we can pass to the limit and obtain \cref{item:limit_has_positive_areas}.

	The proof of \cref{item:limit_is_a_simplicial_complex} is broken down into the following steps, according to the definition of $\Sigma_\Delta(Q^*)$, see \eqref{eq:collection_of_convex_hulls}, and the definition of simplicial complexes; see \cref{section:geometric_abstract_simplicial_complexes}.
	\begin{enumerategreek}
		\item
			\label[statement]{item:limit_is_a_simplicial_complex:1}
			$\Sigma_\Delta(Q^*)$ is a non-empty, finite collection of simplices in $\R^2$.

		\item
			\label[statement]{item:limit_is_a_simplicial_complex:2}
			Every face of a simplex in $\Sigma_\Delta(Q^*)$ also belongs to $\Sigma_\Delta(Q^*)$.

		\item
			\label[statement]{item:limit_is_a_simplicial_complex:3}
			The non-empty intersection of any two simplices $\sigma_*, \sigma_*'$ in $\Sigma_\Delta(Q^*)$ is a face of both $\sigma_*$ and $\sigma_*'$.

		\item
			\label[statement]{item:limit_is_a_simplicial_complex:4}
			The abstract simplicial complex underlying $\Sigma_\Delta(Q^*)$ is $\Delta$.
	\end{enumerategreek}
	We proved in \cref{item:limit_has_positive_areas} that $\area{Q^*}[i_0,i_1,i_2] > 0$ holds for all $2$-faces in $\Delta$.
	Therefore, the vertex positions $\paren[big]\{\}{Q^*_{i_0},Q^*_{i_1},Q^*_{i_2}}$ are affine independent so that their convex hulls are $2$-simplices in $\R^2$.
	Since $\Delta$ is an abstract simplicial $2$-complex, all other sets in $\Sigma_\Delta(Q^*)$ are the convex hulls of some subset of the vertices of a triangle, and these vertices are clearly affine independent as well.
	This shows \cref{item:limit_is_a_simplicial_complex:1,item:limit_is_a_simplicial_complex:2}.

	The proof of \cref{item:limit_is_a_simplicial_complex:3,item:limit_is_a_simplicial_complex:4} is the most difficult part.
	In practical terms, we have to show that the structure of the simplicial complex describing the mesh does not change when passing to the limit $n \to \infty$ in the vertex positions.

	The dimensions of the individual simplices in $\Sigma_\Delta(Q^n)$ are easily seen to be stable under this limit.
	Indeed, notice that $1$-faces will not collapse at the limit thanks to $\distance{Q}[i_0][i_1] = \edgelength{Q}{\ell}[i_0,i_1,i_2]$ and \cref{lemma:bounds_of_triangle_properties_in_terms_of_f}, \cref{item:bound_length_edges}.
	In the same way, $2$-faces do not collapse at the limit thanks to the bound on the heights given by \cref{lemma:bounds_of_triangle_properties_in_terms_of_f}, \cref{item:bound_height}.
	Therefore, the only concern is that unwanted intersections might appear.
	For instance, a vertex might converge to meet an edge which it is not supposed to intersect.
	We need to show that the properties of~$f$ prevent this from happening.

	In the following we are considering two arbitrary, \emph{distinct} faces $\sigma$ and $\sigma'$ of $\Delta$.
	(The case of $\Delta$ consisting of a single $2$-simplex is trivial.)
	According to the dimension of $\sigma$ and the vertices involved, we denote it by $[i_0]$ if it is a $0$-face, by $[i_0,i_1]$ if it is a $1$-face, and by $[i_0,i_1,i_2]$ in case of a $2$-face.
	The corresponding vertex indices for $\sigma'$ are $j_0$, $j_1$ and $j_2$.
	We denote the corresponding faces in $\Sigma_\Delta(Q^n)$ by $\sigma_n$ and $\sigma'_n$, and those in $\Sigma_\Delta(Q^*)$ by $\sigma_*$ and $\sigma_*'$, respectively.
	For instance, when $\sigma = [i_0,i_1]$, then $\sigma_n = \conv\paren\{\}{Q^n_{i_0},Q^n_{i_1}}$ and $\sigma_* = \conv\paren\{\}{Q^*_{i_0},Q^*_{i_1}}$.

	We proceed by distinguishing cases, according to the dimensions of $\sigma$ and $\sigma'$.
	Notice that $\dim \sigma = \dim \sigma_n = \dim \sigma_*$ and $\dim \sigma' = \dim \sigma'_n = \dim \sigma_*'$, which follows from the proof of \cref{item:limit_is_a_simplicial_complex:1,item:limit_is_a_simplicial_complex:2} above.
	In each case, we need to verify that the intersection $\sigma_* \cap \sigma_*'$ is of the same type (empty set, a vertex, etc.) as $\sigma \cap \sigma'$ and $\sigma_n \cap \sigma'_n$.
	To this end, we argue that non-intersecting faces maintain a positive distance $\distance{Q^*}$ or $\Distance{Q^*}$ also in the limit.

	In addition to the dimensions of $\sigma$ and $\sigma'$, we need to distinguish whether they are interior or boundary faces.
	To keep the proof more concise, we restrict the following discussion to the case where both $\sigma$ and $\sigma'$ are \emph{boundary} faces.
	Nonetheless, the proof is involving estimates pertaining to interior faces as well.
	The modifications in the remaining cases use very similar arguments and the details are left to the reader.

	\begin{description}
		\item[$\sigma,\sigma'$ are boundary $0$-faces:]
			Since $\sigma = [i_0]$ and $\sigma' = [j_0]$ are distinct vertices, $\sigma \cap \sigma' = \emptyset$ holds.
			We have two cases to consider.
			$\bullet$
			First, if $[i_0,j_0]$ is a face in $\Delta$, then there exists another vertex~$[k]$ such that $[i_0,j_0,k]$ is a $2$-simplex in $\Delta$.
			This is since $\Delta$ is pure.
			Therefore, $\norm{Q^n_{i_0} - Q^n_{j_0}}$ agrees with the length $\edgelength{Q^n}{\ell}[i_0,j_0,k]$ of an edge for some $\ell \in \{0,1,2\}$.
			From \cref{lemma:bounds_of_triangle_properties_in_terms_of_f}, \cref{item:bound_length_edges}, we know $\edgelength{Q^n}{\ell}[i_0,j_0,k] \ge (2 \beta_1) / f(Q^n) \ge (2 \beta_1) / b$ and thus $\norm{Q^*_{i_0} - Q^*_{j_0}} \ge (2 \beta_1) / b$, \ie, $\sigma_* = \{Q^*_{i_0}\}$ and $\sigma_*' = \{Q^*_{j_0}\}$ do not intersect.
			$\bullet$
			Second, if $[i_0,j_0]$ is not a face in $\Delta$, then there exists $[j_1]$ such that $[j_0,j_1]$ is a boundary $1$-face.
			Using \cref{proposition:bounds_of_distances_in_terms_of_f}, \cref{item:distance_boundary_vertex_boundary_edge} we know $\Distance[auto]{Q^n}[i_0][[j_0,j_1]] \ge \beta_2/b$ and thus $\Distance[auto]{Q^*}[i_0][[j_0,j_1]] \ge \beta_2/b$.
			In particular, $\norm{Q^*_{i_0} - Q^*_{j_0}} \ge \beta_2/b$ holds.
			Again, $\sigma_* \cap \sigma_*' = \emptyset$.

		\item[$\sigma$ is a boundary $0$-face and $\sigma'$ is a boundary $1$-face:]
			$\bullet$
			First assume that $\sigma = [i_0]$ and $\sigma' = [j_0,j_1]$ do not intersect.
			Using \cref{proposition:bounds_of_distances_in_terms_of_f}, \cref{item:distance_boundary_vertex_boundary_edge} we obtain $\Distance[auto]{Q^n}[i_0][[j_0,j_1]] \ge \beta_2/b$ and thus $\Distance[auto]{Q^*}[i_0][[j_0,j_1]] \ge \beta_2/b$, \ie, $\sigma_* \cap \sigma_*' = \emptyset$.
			$\bullet$
			Second, if $\sigma$ and $\sigma'$ intersect, then necessarily $\sigma \cap \sigma' = \sigma$.
			Without loss of generality, suppose $[i_0] = [j_0] \neq [j_1]$.
			Moreover, since $\Delta$ is pure, there exists a vertex~$[k]$ such that $[j_0,j_1,k]$ is a $2$-simplex in $\Delta$.
			Therefore, $\norm{Q^n_{i_0} - Q^n_{j_1}} = \norm{Q^n_{j_0} - Q^n_{j_1}}$ agrees with the length $\edgelength{Q^n}{\ell}[j_0,j_1,k]$ of an edge for some $\ell \in \{0,1,2\}$.
			Thanks to \cref{lemma:bounds_of_triangle_properties_in_terms_of_f}, \cref{item:bound_length_edges}, we have $\edgelength{Q^n}{\ell}[j_0,j_1,k] \ge (2 \beta_1) / f(Q^n) \ge (2 \beta_1)/b$ and thus $\norm{Q^*_{i_0} - Q^*_{j_1}} \ge (2 \beta_1)/b$.
			Therefore, $\sigma_* \cap \sigma_*' = \sigma_*$ as desired.

		\item[$\sigma$ is a boundary $0$-face and $\sigma'$ is a boundary $2$-face:]
			We begin by noticing that for any $Q \in \zeromanifold$ and in particular for any $Q \in \planarmanifold$, the distance between a $0$-face and a $2$-face satisfies
			\begin{align}
				\distance{Q}[i_0][[j_0,j_1,j_2]]
				&
				=
				\min \paren[auto]\{\}{%
					\distance{Q}[i_0][[j_0,j_1]],
					\distance{Q}[i_0][[j_1,j_2]],
					\distance{Q}[i_0][[j_2,j_0]]
				}
				\nonumber
				\\
				&
				\ge
				\frac{1}{\sqrt{2}}
				\min \paren[auto]\{\}{%
					\Distance{Q}[i_0][[j_0,j_1]],
					\Distance{Q}[i_0][[j_1,j_2]],
					\Distance{Q}[i_0][[j_0,j_2]]
				}
				.
				\label{eq:distance_vertex_triangle}
			\end{align}

			$\bullet$
			Now, we suppose that $\sigma = [i_0]$ and $\sigma' = [j_0,j_1,j_2]$ do not intersect.
			Since $\sigma'$ is a boundary face, at least one of its $1$-faces is a boundary face.
			Without loss of generality, let $[j_0,j_1]$ be a boundary $1$-face.
			As in the previous case, \cref{proposition:bounds_of_distances_in_terms_of_f}, \cref{item:distance_boundary_vertex_boundary_edge} gives a bound on the first term of \eqref{eq:distance_vertex_triangle}.
			Moreover, \wolog we assume $[j_1,j_2]$ is and interior edge. Then, thanks to \cref{proposition:bounds_of_distances_in_terms_of_f}, \cref{item:distance_boundary_vertex_interior_edge} the second term of \eqref{eq:distance_vertex_triangle} is bounded.
			In the same way, if $[j_0,j_2]$ is either interior or boundary can be bounded using one of the previous arguments.
			In summary, we get
			\begin{equation*}
				\distance{Q^n}[i_0][[j_0,j_1,j_2]]
				\ge
				\frac{\min\{\beta_1,\beta_2\}}{2b}
				\min_{\theta\in\Theta}
				\paren[auto]\{\}{%
					\sqrt{1-\cos^2(\theta)}
				}
			\end{equation*}
			Notice, that thanks to \cref{lemma:bounds_of_triangle_properties_in_terms_of_f}, \cref{item:bound_cosines}, $\abs{\cos(\theta)} \le \Psi(b) < 1$ holds for all $\theta\in\Theta$, which implies
			\begin{equation*}
				\distance{Q^n}[i_0][[j_0,j_1,j_2]]
				\ge
				\frac{\min\{\beta_1,\beta_2\}}{2b}
				\paren[auto](){\sqrt{1-\Psi^2(b)}}
			\end{equation*}
			and thus the same bound is valid for $\distance{Q^*}[i_0][[j_0,j_1,j_2]]$.
			This shows that $\sigma_* \cap \sigma_*' = \emptyset$ holds.

			$\bullet$
			Second, if $\sigma$ and $\sigma'$ intersect, then necessarily $\sigma \cap \sigma' = \sigma$.
			Without loss of generality, suppose $[i_0] = [j_0]$.
			Therefore, $\norm{Q^n_{i_0} - Q^n_{j_1}} = \norm{Q^n_{j_0} - Q^n_{j_1}}$ agrees with the length $\edgelength{Q^n}{\ell}[j_0,j_1,j_2]$ of an edge for some $\ell \in \{0,1,2\}$.
			In the same way, the distance $\distance{Q^n}[i_0][[j_1,j_2]]$ agrees with $\height{Q^n}{\ell}[j_0,j_1,j_2]$ for some $\ell = 0,1,2$.
			Thus, \cref{lemma:bounds_of_triangle_properties_in_terms_of_f}, \cref{item:bound_length_edges,item:bound_height} show $\sigma_* \cap \sigma_*' = \sigma_*$ as desired.

		\item[$\sigma$ and $\sigma'$ are boundary $1$-faces:]
			As in the previous case, we begin by noticing that for any $Q \in \zeromanifold$ and in particular for any $Q \in \planarmanifold$, the distance $\distance{Q}[[i_0,i_1]][[j_0,j_1]]$ of two $1$-faces satisfies
			\makeatletter
			\begin{align}
				\ltx@ifclassloaded{mcom-l}{%
					\MoveEqLeft
					\distance{Q}[[i_0,i_1]][[j_0,j_1]]
					\nonumber
					\\
					&
					=
					\min\paren[auto]\{\}{%
						\distance{Q}[i_0][[j_0,j_1]],
						\distance{Q}[i_1][[j_0,j_1]],
						\distance{Q}[j_0][[i_0,i_1]],
						\distance{Q}[j_1][[i_0,i_1]]
					}
					\nonumber
					\\
					&
					\ge
					\frac{1}{\sqrt{2}}
					\min\paren[auto]\{.{%
						\Distance{Q}[i_0][[j_0,j_1]],
						\Distance{Q}[i_1][[j_0,j_1]],
						\Distance{Q}[j_0][[i_0,i_1]],
					}
					\nonumber
					\\
					&
					\phantom{{}\ge{}\frac{1}{\sqrt{2}}\min\{}
					\paren[auto].\}{%
						\Distance{Q}[j_1][[i_0,i_1]]
					}
					.
					}{%
					\distance{Q}[[i_0,i_1]][[j_0,j_1]]
					&
					=
					\min\paren[auto]\{\}{%
						\distance{Q}[i_0][[j_0,j_1]],
						\distance{Q}[i_1][[j_0,j_1]],
						\distance{Q}[j_0][[i_0,i_1]],
						\distance{Q}[j_1][[i_0,i_1]]
					}
					\nonumber
					\\
					&
					\ge
					\frac{1}{\sqrt{2}}
					\min\paren[auto]\{\}{%
						\Distance{Q}[i_0][[j_0,j_1]],
						\Distance{Q}[i_1][[j_0,j_1]],
						\Distance{Q}[j_0][[i_0,i_1]],
						\Distance{Q}[j_1][[i_0,i_1]]
					}
					.
				}
				\label{eq:distance_edge_edge}
			\end{align}
			\makeatother
			$\bullet$
			First suppose that $\sigma = [i_0,i_1]$ and $\sigma' = [j_0,j_1]$ do not intersect.
			Since both $\sigma$ and $\sigma'$ are boundary $1$-faces, their vertices are boundary vertices.
			Therefore, by \cref{proposition:bounds_of_distances_in_terms_of_f}, \cref{item:distance_boundary_vertex_boundary_edge}, all four terms on the right hand side of \eqref{eq:distance_edge_edge} with $Q$ replaced by $Q^n$ are bounded below by $\beta_2/b$.
			Consequently, we obtain $\distance{Q^*}[[i_0,i_1]][[j_0,j_1]] \ge \beta_2/(\sqrt{2}b)$, \ie, $\sigma_* \cap \sigma_*' = \emptyset$.
			$\bullet$
			Second, when $\sigma = [i_0,i_1]$ and $\sigma' = [j_0,j_1]$ intersect, they intersect in a vertex.
			Without loss of generality, $[i_0] = [j_0]$ holds.
			Thanks to \cref{proposition:bounds_of_distances_in_terms_of_f}, \cref{item:distance_boundary_vertex_boundary_edge}, we find that $\Distance[auto]{Q^n}[i_1][[j_0,j_1]] \ge \beta_2/b$ holds and thus $\Distance[auto]{Q^*}[i_1][[j_0,j_1]] \ge \beta_2/b$ as well.
			The same argument also shows $\Distance{Q^*}[j_1][[i_0,i_1]] \ge \beta_2/b$.
			As mentioned previously, the edge lengths $\norm{Q^*_{i_0} - Q^*_{i_1}}$ and $\norm{Q^*_{j_0} - Q^*_{j_1}}$ remain positive by \cref{lemma:bounds_of_triangle_properties_in_terms_of_f}, \cref{item:bound_length_edges}.
			This implies $\sigma_* \cap \sigma_*' = \{Q^*_{i_0}\} = \{Q^*_{j_0}\}$.

		\item[$\sigma$ is a boundary $1$-face and $\sigma'$ is a boundary $2$-face:]
			We begin by noticing that for any $Q \in \zeromanifold$ and in particular for any $Q \in \planarmanifold$, the distance $\distance{Q}[[i_0,i_1]][[j_0,j_1,j_2]]$ of any $1$-face and any $2$-face satisfies
			\makeatletter
			\begin{align}
				\ltx@ifclassloaded{mcom-l}{%
					\MoveEqLeft
					\distance{Q}[[i_0,i_1]][[j_0,j_1,j_2]]
					\nonumber
					\\
					&
					=
					\min\paren[auto]\{\}{%
						\distance{Q}[[i_0,i_1]][[j_0,j_1]],
						\distance{Q}[[i_0,i_1]][[j_1,j_2]],
						\distance{Q}[[i_0,i_1]][[j_2,j_0]]
					}
					\nonumber
					\\
					&
					\ge
					\frac{1}{\sqrt{2}}
					\min\big\{%
					\Distance{Q}[i_0][[j_0,j_1]],
					\Distance{Q}[i_0][[j_1,j_2]],
					\Distance{Q}[i_0][[j_2,j_0]],
					\nonumber
					\\
					&
					\phantom{{}=\frac{1}{\sqrt{2}}\min\big\{}
					\Distance{Q}[i_1][[j_0,j_1]],
					\Distance{Q}[i_1][[j_1,j_2]],
					\Distance{Q}[i_1][[j_2,j_0]],
					\nonumber
					\\
					&
					\phantom{{}=\frac{1}{\sqrt{2}}\min\big\{}
					\Distance{Q}[j_0][[i_0,i_1]],
					\Distance{Q}[j_1][[i_0,i_1]],
					\Distance{Q}[j_2][[i_0,i_1]]
					\big\}
					.
				}{%
					\distance{Q}[[i_0,i_1]][[j_0,j_1,j_2]]
					&
					=
					\min\paren[auto]\{\}{%
					\distance{Q}[[i_0,i_1]][[j_0,j_1]],
					\distance{Q}[[i_0,i_1]][[j_1,j_2]],
					\distance{Q}[[i_0,i_1]][[j_2,j_0]]
					}
					\nonumber
					\\
					&
					\ge
					\frac{1}{\sqrt{2}}
					\min\big\{%
					\Distance{Q}[i_0][[j_0,j_1]],
					\Distance{Q}[i_0][[j_1,j_2]],
					\Distance{Q}[i_0][[j_2,j_0]],
					\nonumber
					\\
					&
					\phantom{{}=\min\big\{}
					\Distance{Q}[i_1][[j_0,j_1]],
					\Distance{Q}[i_1][[j_1,j_2]],
					\Distance{Q}[i_1][[j_2,j_0]],
					\nonumber
					\\
					&
					\phantom{{}=\min\big\{}
					\Distance{Q}[j_0][[i_0,i_1]],
					\Distance{Q}[j_1][[i_0,i_1]],
					\Distance{Q}[j_2][[i_0,i_1]]
					\big\}
					.
				}
				\label{eq:distance_edge_triangle}
			\end{align}
			\makeatother
			The estimate follows from \eqref{eq:distance_edge_edge}.
			$\bullet$
			First suppose that $\sigma = [i_0,i_1]$ and $\sigma' = [j_0,j_1,j_2]$ do not intersect.
			Some of the terms on the right hand side of \eqref{eq:distance_edge_triangle} with $Q$ replaced by $Q^n$ are distances between boundary vertices and boundary edges, for which \cref{proposition:bounds_of_distances_in_terms_of_f}, \cref{item:distance_boundary_vertex_boundary_edge} provides a lower bound $\beta_2/b$.
			The remaining terms are distances between interior vertices and interior edges, or boundary vertices and interior edges, which can be bounded below by \cref{proposition:bounds_of_distances_in_terms_of_f}, \cref{item:distance_interior_vertex_arbitrary_edge} or \cref{item:distance_boundary_vertex_interior_edge}, respectively.
			Altogether, we obtain
			\begin{equation*}
				\distance{Q^n}[[i_0,i_1]][[j_0,j_1,j_2]]
				\ge
				\frac{\min\{\beta_1,\beta_2\}}{2b}
				\sqrt{1-\Psi^2(b)}
				>
				0
				.
			\end{equation*}
			$\bullet$
			Second, assume that $\sigma$ and $\sigma'$ intersect in a vertex.
			Without loss of generality, $[i_0] = [j_0]$ holds.
			We recall \eqref{eq:distance_vertex_triangle}:
			\makeatletter
			\ltx@ifclassloaded{mcom-l}{%
				\begin{align*}
					\MoveEqLeft
					\distance{Q}[i_1][[j_0,j_1,j_2]]
					\\
					&
					\ge
					\frac{1}{\sqrt{2}}
					\min\paren[auto]\{\}{%
						\Distance{Q}[i_1][[j_0,j_1]],
						\Distance{Q}[i_1][[j_1,j_2]],
						\Distance{Q}[i_1][[j_2,j_0]]
					}
					.
				\end{align*}
				}{%
				\begin{equation*}
					\distance{Q}[i_1][[j_0,j_1,j_2]]
					\ge
					\frac{1}{\sqrt{2}}
					\min\paren[auto]\{\}{%
						\Distance{Q}[i_1][[j_0,j_1]],
						\Distance{Q}[i_1][[j_1,j_2]],
						\Distance{Q}[i_1][[j_0,j_2]]
					}
					.
				\end{equation*}
			}
			\makeatother
			As above, some of the terms on the right hand side of \eqref{eq:distance_vertex_triangle} with $Q$ replaced by $Q^n$ are distances between the boundary vertex~$[i_1]$ and boundary edges, and the remaining terms are distances between $[i_1]$ and interior edges.
			An application of \cref{proposition:bounds_of_distances_in_terms_of_f}, \cref{item:distance_boundary_vertex_boundary_edge} and \cref{item:distance_boundary_vertex_interior_edge} yields a uniformly positive lower bound for $\distance{Q^n}[i_1][[j_0,j_1,j_2]]$ and thus for $\distance{Q^*}[i_1][[j_0,j_1,j_2]]$.
			This implies $\sigma_* \cap \sigma_*' = \{Q^*_{i_0}\} = \{Q^*_{j_0}\}$.
			$\bullet$
			Third, assume that $\sigma \cap \sigma' = \sigma$.
			Without loss of generality, suppose $[i_0] = [j_0]$ and $[i_1] = [j_1]$.
			Then we can use $\distance[auto]{Q^n}[j_2][[i_0,i_1]] \ge \height{Q^n}{\ell}[j_0,j_1,j_2]$ for $\ell = 2$, and thus \cref{lemma:bounds_of_triangle_properties_in_terms_of_f}, \cref{item:bound_height} yields a lower bound of $\beta_1/b$.
			As mentioned previously, the edge length $\norm{Q_{i_0} - Q_{i_1}} = \norm{Q_{j_0} - Q_{j_1}}$ remains positive by \cref{lemma:bounds_of_triangle_properties_in_terms_of_f}, \cref{item:bound_length_edges}.
			From here we conclude $\sigma_* \cap \sigma_*' = \sigma_*$.

		\item[$\sigma$ and $\sigma'$ are boundary $2$-faces:]
			We first notice that
			\begin{equation}
				\label{eq:distance_triangle_triangle}
				\distance{Q}[[i_0,i_1,i_2]][[j_0,j_1,j_2]]
				\ge
				\frac{1}{\sqrt{2}}
				\min_{\substack{\ell = 0,1,2 \\ \hat{\ell} = 0,1,2}
				}
				\paren[auto]\{\}{\Distance{Q}[i_\ell][[j_{\hat{\ell}},j_{\hat{\ell}\oplus 1}]], \Distance{Q}[j_{\hat{\ell}}][[i_{\ell},i_{\ell\oplus 1}]]}
			\end{equation}
			holds.
			$\bullet$
			First suppose that $\sigma = [i_0,i_1,i_2]$ and $\sigma' = [j_0,j_1,j_2]$ do not intersect.
			Each term on the right hand side of \eqref{eq:distance_triangle_triangle} with $Q$ replaced by $Q^n$ can be estimated below by \cref{proposition:bounds_of_distances_in_terms_of_f}, \cref{item:distance_boundary_vertex_boundary_edge,item:distance_interior_vertex_arbitrary_edge,item:distance_boundary_vertex_interior_edge}.
			Thanks to \cref{lemma:bounds_of_triangle_properties_in_terms_of_f}, \cref{item:bound_cosines}, we obtain the uniform lower bound
			\begin{equation*}
					\distance{Q^n}[[i_0,i_1,i_2]][[j_0,j_1,j_2]]
					\ge
					\frac{\min\{\beta_1,\beta_2\}}{2b}
					\sqrt{1-\Psi^2(b)}
					>
					0
			\end{equation*}
			and thus $\sigma_* \cap \sigma_*' = \emptyset$ holds.

			$\bullet$
			Second, suppose that $\sigma$ and $\sigma'$ intersect in a vertex.
			Without loss of generality, suppose $[i_0] = [j_0]$.
			We need to consider
			\makeatletter
			\ltx@ifclassloaded{mcom-l}{%
				\begin{align*}
					\MoveEqLeft
					\distance{Q}[i_1][[j_0,j_1,j_2]]
					\\
					&
					\ge
					\frac{1}{\sqrt{2}}
					\min\paren[auto]\{\}{%
						\Distance{Q}[i_1][[j_0,j_1]],
						\Distance{Q}[i_1][[j_1,j_2]],
						\Distance{Q}[i_1][[j_2,j_0]]
					}
				\end{align*}
				}{%
				\begin{equation*}
					\distance{Q}[i_1][[j_0,j_1,j_2]]
					\ge
					\frac{1}{\sqrt{2}}
					\min\paren[auto]\{\}{%
						\Distance{Q}[i_1][[j_0,j_1]],
						\Distance{Q}[i_1][[j_1,j_2]],
						\Distance{Q}[i_1][[j_2,j_0]]
					}
				\end{equation*}
			}
			\makeatother
			and the same with $i_1$ replaced by $i_2$ and show that these expressions are bounded away from zero for $Q = Q^n$.
			Regardless of whether $[i_1]$ and $[i_2]$, and $[j_0,j_1]$, $[j_1,j_2]$, $[j_2,j_0]$ are interior or boundary faces, in each case, one of \cref{proposition:bounds_of_distances_in_terms_of_f}, \cref{item:distance_boundary_vertex_boundary_edge,item:distance_interior_vertex_arbitrary_edge,item:distance_boundary_vertex_interior_edge} applies and provides this lower bound.
			The same argument applies with the roles of $\sigma$ and $\sigma'$ reversed.
			We can conclude that $\sigma_* \cap \sigma_*' = \{Q^*_{i_0}\} = \{Q^*_{j_0}\}$ as desired.
			$\bullet$
			Third, suppose that $\sigma$ and $\sigma'$ intersect in a common edge, which is necessarily an interior $1$-face.
			Without loss of generality, suppose that $[i_0] = [j_0]$ and $[i_1] = [j_2]$.
			We need to estimate only the distances from $[i_2]$ to $[j_0,j_1,j_2]$ and from $[j_1]$ to $[i_0,i_1,i_2]$.
			To this end, we use
			\begin{equation*}
				\begin{aligned}
					\distance{Q^n}[i_2][[j_0,j_1,j_2]]
					&
					\ge
					\frac{\min\{\beta_1,\beta_2\}}{2b}
					\sqrt{1-\Psi^2(b)}
					&
					>0
					\\
					\distance{Q^n}[j_1][[i_0,i_1,i_2]]
					&
					\ge
					\frac{\min\{\beta_1,\beta_2\}}{2b}
					\sqrt{1-\Psi^2(b)}
					&
					>
					0
				\end{aligned}
			\end{equation*}
			which follow \cref{proposition:bounds_of_distances_in_terms_of_f}, \cref{item:distance_boundary_vertex_boundary_edge,item:distance_interior_vertex_arbitrary_edge,item:distance_boundary_vertex_interior_edge}.
			Therefore, $\sigma_* \cap \sigma_*' = \conv\paren\{\}{Q^*_{i_0},Q^*_{i_1}} = \conv\paren\{\}{Q^*_{j_2},Q^*_{j_0}}$ as desired.
			Moreover, $Q^*_{i_0}$ and $Q^*_{i_1}$ are going to remain distinct points since $\norm{Q^*_{i_0} - Q^*_{i_1}} = \edgelength{Q^*}{2}[i_0,i_1,i_2]$, which remains bounded away from zero by \cref{lemma:bounds_of_triangle_properties_in_terms_of_f}, \cref{item:bound_length_edges}.
	\end{description}

	 As we have already mentioned, similar arguments as above can be used in the case when $\sigma$ or $\sigma'$ are interior faces.
	 To summarize, we conclude that the limiting simplices $\sigma_*$, $\sigma_*'$ intersect in the same way as $\sigma_n$ and $\sigma'_n$, and the dimension of this intersection is in turn dictated by the underlying abstract simplicial complex.
	 We have thus shown that $\Sigma_\Delta(Q^*)$ is a simplicial $2$-complex whose associated abstract simplicial complex is $\Delta$, which concludes the proof of \cref{item:limit_is_a_simplicial_complex}.

	\Cref{item:limit_has_positive_areas,item:limit_is_a_simplicial_complex} together show that $Q^*$ belongs to $\plusmanifold$.
	It remains to confirm \cref{item:limit_has_a_path_to_Qref}.
	Since $\plusmanifold$ is locally connected, its connected components and path components agree.
	Since connected components are closed, and $\planarmanifold$ is by definition a path component of $\plusmanifold$, $\planarmanifold$ is closed in $\plusmanifold$.
	(We refer the reader to \cite[Ch.~4]{Lee2011} for more details.)

	Having shown that $Q^n$ converges to $Q^*$ in $\plusmanifold$ and $(Q^n) \subset \planarmanifold$, we can conclude that $Q^*$ belongs to $\planarmanifold$ as well.
	Finally, using the continuity of $f$ on $\plusmanifold$, see \cref{lemma:f_is_well-defined_and_continuous}, and thus on $\planarmanifold$, we infer that $Q^* \in f^{-1}(K)$ holds.
	This confirms the sequential compactness of $f^{-1}(K)$ in $\planarmanifold$ and concludes the proof.
\end{proof}

\begin{remark}
	For the proof of \cref{theorem:f_is_proper}, we used that $\Delta$ is a connectivity complex according to \cref{definition:connectivity_complex}, \ie, $\Delta$ is a pure, abstract simplicial $2$-complex, which is $2$-path connected.
	The purity was used whenever we embedded a lower-dimensional face into a $2$-face.
	The $2$-path connectedness enters through \cref{proposition:bounds_of_distances_in_terms_of_f}.
	We encourage the reader to check the proofs in \cref{section:geometric_abstract_simplicial_complexes}, especially \cref{lemma:link_polygonal_chain,lemma:link_closed_polygonal_chain}.
\end{remark}

Although \cref{theorem:f_is_proper} shows the properness of $f$ on $\planarmanifold$, unfortunately we cannot directly use it in \cref{theorem:construction_of_complete_Riemannian_metric} since $f$ lacks differentiability due to the occurrence of the vertex-edge distance~$Q \mapsto \Distance{Q}$ in \eqref{eq:f}, which is only Lipschitz.
Therefore, we replace $\Distance{Q}$ by a regularized function $Q \mapsto \regularizedDistance{Q}$ of class~$\cC^3$.
This then gives rise to the following regularized function $f^\mu \colon \plusmanifold \to \R$:
\begin{equation}
	\label{eq:f_mu}
	f^\mu(Q)
	\coloneqq
	\sum_{k=1}^{N_T} \sum_{\ell=0}^2 \frac{\beta_1}{\height[auto]{Q}{\ell}[i^k_0,i^k_1,i^k_2]}
	+
	\sum_{[j_0,j_1]\in E_\partial}\sum_{\substack{i_0\in V_\partial \\ i_0\neq j_0,j_1}} \frac{\beta_2}{\regularizedDistance{Q}[i_0][[j_0,j_1]]}
	+
	\frac{\beta_3}{2} \norm{Q - \Qref}_F^2
	,
\end{equation}
which we can use in place of \eqref{eq:f}.
Provided that we choose $f^\mu \ge f$, \cref{proposition:f_mu_is_proper} below implies that $f^\mu$ is proper as well.
Indeed, we are providing in \cref{section:example_regularization_vertex-edge_distance} an entire family of regularized functions~$f^\mu$ with parameter~$\mu$ which, in addition to satisfying $f^\mu \ge f$, approximate~$f$ arbitrarily well.

\begin{proposition}
	\label{proposition:f_mu_is_proper}
	Let $X$ be a metric space and consider a function $f\colon X \rightarrow \R$ which is proper.
	Suppose that $\widehat f \colon X \to \R$ is continuous and it satisfies $0 \le f \le \widehat f$.
	Then $\widehat f$ is proper.
\end{proposition}
\begin{proof}
	Let $K \subset \R$ be compact.
	We need to prove that $(\widehat f)^{-1}(K)$ is compact.
	In case $(\widehat f)^{-1}(K) = \emptyset$, nothing is to be shown.
	Otherwise, since $\widehat f$ is non-negative, we can suppose that $K \subset [a,b]$ with $a \ge 0$.
	Since $(\widehat f)^{-1}(K) \subset X$ and $X$ is a metric space, the compactness of $(\widehat f)^{-1}(K)$ is equivalent to its sequential compactness.
	To show the latter, consider a sequence $(x_n) \subset (\widehat f)^{-1}(K)$, \ie, $\widehat f(x_n) \in K$.
	Using the assumption $0 \le f \le \widehat f$, we obtain $0 \le f(x_n) \le \widehat f(x_n) \in K \subset [a,b]$ and therefore $f(x_n) \in [0,b]$ for all $n \in \N$.
	Since $f$ is proper and $[0,b]$ is compact, there exists a subsequence, still labeled $(x_n)$, such that $x_n \to x^*$ in $[0,b]$.
	Thanks to the continuity of $\widehat f$, we have $\widehat f(x_n) \to \widehat f(x^*)$.
	Since all $\widehat f(x_n) \in K$, the limit $\widehat f(x^*)$ belongs to $K$ as well.
\end{proof}

We summarize these results in the following theorem.
\begin{theorem}
	\label{theorem:f_mu_is_proper_and_Riemannian_metric_is_complete}
	Suppose that $\beta_1, \beta_2, \beta_3 > 0$ holds.
	Consider a continuous function~$Q \mapsto \regularizedDistance{Q}$ on $\plusmanifold$ which satisfies $0 < \regularizedDistance{Q} \le \Distance{Q}$.
	Then the following statements hold.
	\begin{enumeratelatin}
		\item
			\label[statement]{item:f_mu_is_proper}
			The restriction of $f^\mu$ to $\planarmanifold$ defined in \eqref{eq:f_mu} is proper.

		\item
			\label[statement]{item:Riemannian_metric_is_complete}
			Suppose in addition that $Q \mapsto \regularizedDistance{Q}$ is of class $\cC^3$ on $\planarmanifold$.
			Then $\planarmanifold$, endowed with the Riemannian metric whose components (with respect to the $\vvec$ chart) are given by
			\begin{equation}
				\label{eq:complete_metric_for_planar_triangular_meshes}
				g_{ab}
				=
				\delta_a^b + \dfrac{\partial f^\mu}{\partial (\vvec Q)^a} \dfrac{\partial f^\mu}{\partial (\vvec Q)^b}
				,
				\quad
				a, b = 1, \ldots, 2 \, N_V
				,
			\end{equation}
			is geodesically complete.
	\end{enumeratelatin}
\end{theorem}
\begin{proof}
	The definition of $f$ in \eqref{eq:f}, the definition of $f^\mu$ in \eqref{eq:f_mu} and the assumption $0 < \regularizedDistance{Q} \le \Distance{Q}$ imply $0 < f \le f^\mu$ on
	$\plusmanifold$.
	\Cref{item:f_mu_is_proper} now follows from \cref{proposition:f_mu_is_proper}.
	\Cref{item:Riemannian_metric_is_complete} follows immediately from \Cref{theorem:construction_of_complete_Riemannian_metric}.
\end{proof}

\begin{remark}
	\label{remark:f_mu_with_threshold}
	\cref{theorem:f_mu_is_proper_and_Riemannian_metric_is_complete} remains valid when \eqref{eq:f_mu} is replaced by the slightly more general function
	\makeatletter
	\ltx@ifclassloaded{mcom-l}{%
		\begin{multline}
			\label{eq:f_mu_with_cutoffs}
			\widetilde f^\mu(Q)
			\coloneqq
			\sum_{k=1}^{N_T} \sum_{\ell=0}^2 \chi_1 \paren[auto](){\frac{\beta_1}{{\height[auto]{Q}{\ell}[i^k_0,i^k_1,i^k_2]}}}
			+
			\sum_{[j_0,j_1]\in E_\partial}\sum_{\substack{i_0\in V_\partial \\ i_0\neq j_0,j_1}} \chi_2 \paren[auto](){\frac{\beta_2}{\regularizedDistance{Q}[i_0][[j_0,j_1]]}}
			\\
			+
			\frac{\beta_3}{2} \norm{Q - \Qref}_F^2
			.
		\end{multline}
		}{%
		\begin{equation}
			\label{eq:f_mu_with_cutoffs}
			\widetilde f^\mu(Q)
			\coloneqq
			\sum_{k=1}^{N_T} \sum_{\ell=0}^2 \chi_1 \paren[auto](){\frac{\beta_1}{{\height[auto]{Q}{\ell}[i^k_0,i^k_1,i^k_2]}}}
			+
			\sum_{[j_0,j_1]\in E_\partial}\sum_{\substack{i_0\in V_\partial \\ i_0\neq j_0,j_1}} \chi_2 \paren[auto](){\frac{\beta_2}{\regularizedDistance{Q}[i_0][[j_0,j_1]]}}
			+
			\frac{\beta_3}{2} \norm{Q - \Qref}_F^2
			.
		\end{equation}
	}
	\makeatother
	Here $\chi_1$ is a cut-off function of class~$C^3$ which satisfies $\chi_1(s) = 0$ on some interval $[0,\underline s]$ and $\chi_1(s) = s$ for $s \ge \overline s$.
	The same holds for $\chi_2$.
	In other words, the first and second terms, which were seen to be responsible to avoid interior and exterior self-intersections, respectively, can safely be turned off when the heights, or the distances of boundary vertices to non-incident boundary edges, respectively, are larger than a threshold~$1/\underline s$.
	We will exploit this in our numerical experiments.
\end{remark}

The geodesic equation associated with the metric \eqref{eq:complete_metric_for_planar_triangular_meshes} is quite involved.
	We describe in \cref{section:geodesic_equation} a numerical approach to its solution, and present numerical results in \cref{section:numerical_experiments}.
In the remainder of this section, we study an invariance property of the proposed metric \eqref{eq:complete_metric_for_planar_triangular_meshes}.
Moreover, we show that this metric agrees with the Euclidean metric for tangent vectors representing translations in case $\beta_3 = 0$.

\begin{proposition}[Invariance of \eqref{eq:complete_metric_for_planar_triangular_meshes} with respect to mesh translations and rotations]
	\label{proposition:invariance_of_metric_wrt_translations_rotations}
	Suppose that $T \colon \R^2 \to \R^2$ is defined by $T(x) = R \, x + b$ with $R \in \text{SO}(2)$ and $b \in \R^2$.
	We extend $R$ and $T$ to $\R^{2 \times N_V}$, operating column by column.
	We denote by $g_Q \colon \tangent{Q}[\planarmanifold] \times \tangent{Q}[\planarmanifold] \to \R$ the metric \eqref{eq:complete_metric_for_planar_triangular_meshes} at an arbitrary point $Q \in \planarmanifold$.
	Moreover, we denote by $\overline g_Q$ the metric similar to \eqref{eq:complete_metric_for_planar_triangular_meshes} obtained by replacing $\Qref$ by $T(\Qref)$ in \eqref{eq:f_mu}.
	Then
	\begin{equation}
		\label{eq:invariance_of_metric_wrt_translations_rotations}
		g_Q(V,W)
		=
		\overline g_{T(Q)}(RV,RW)
	\end{equation}
	holds for all $V, W \in \tangent{Q}[\planarmanifold]$.
\end{proposition}
\begin{proof}
	Since $T$ is a rotation and translation, $Q \in \planarmanifold$ implies $T(Q) \in \planarmanifold$.
	Moreover, the heights in the first term of \eqref{eq:f_mu} depend only on the relative positions of the vertices to each other, \ie,
	\begin{equation*}
		\height[big]{Q}{\ell}[i_0,i_1,i_2]
		=
		\height[big]{T(Q)}{\ell}[i_0,i_1,i_2]
	\end{equation*}
	holds for all $2$-faces $[i_0,i_1,i_2]$.
	Similarly, the vertex-edge distances are also invariant with respect to rotation/translation, \ie,
	\begin{equation*}
		\regularizedDistance{Q}[i_0][[j_0,j_1]]
		=
		\regularizedDistance{T(Q)}[i_0][[j_0,j_1]]
	\end{equation*}
	holds for all boundary vertices~$i_0$ and non-incident boundary edges~$[j_0,j_1]$.
	This shows that the second sum in \eqref{eq:f_mu} is also invariant with respect to $T$.
	Finally, we have $\norm{Q - \Qref}_F = \norm{T(Q - \Qref)}_F$.
	This shows that $f^\mu(Q;\Qref) = f^\mu(T(Q);T(\Qref))$ holds, where the right hand side term uses $T(\Qref)$ in place of $\Qref$ in \eqref{eq:f_mu}.
	\Cref{eq:invariance_of_metric_wrt_translations_rotations} now follows easily from the chain rule.
\end{proof}

\begin{proposition}[Agreement of \eqref{eq:complete_metric_for_planar_triangular_meshes} with the Euclidean metric for translations]
	\label{proposition:complete_metric_agrees_with_euclidean_for_translations}
	We denote by $g_Q$ the metric \eqref{eq:complete_metric_for_planar_triangular_meshes} at an arbitrary point $Q \in \planarmanifold$.
	Suppose that $V \in \tangent{Q}[\planarmanifold] \cong \R^{2 \times N_V}$ is a tangent vector satisfying $V = \begin{bmatrix} V_0, V_0, \ldots, V_0 \end{bmatrix}$ for some $V_0 \in \R^2$, representing a translation.
	Moreover, $W = \begin{bmatrix} W_1, W_2, \ldots, W_{N_V} \end{bmatrix}$ is an arbitrary vector $\in \tangent{Q}[\planarmanifold]$.
	If $\beta_3 = 0$ holds, then
	\begin{equation}
		\label{eq:complete_metric_agrees_with_euclidean_for_translations}
		g_Q(V,W)
		=
		\trace(V^\transp W)
		=
		\sum_{j=1}^{N_V} V_0^\transp W_j
	\end{equation}
	holds, \ie, the action of \eqref{eq:complete_metric_for_planar_triangular_meshes} on $(V,W)$ agrees with the action of the Euclidean metric.
\end{proposition}
\begin{proof}
	We consider a representative term for the first and second sum in \eqref{eq:f_mu}.
	Suppose that $[i_0,i_1,i_2]$ is a $2$-face of $\Delta$.
	Since heights and distances are translation invariant, we have $\height[auto]{Q}{\ell}[i_0,i_1,i_2] = \height[auto]{Q + t V}{\ell}[i_0,i_1,i_2]$ as well as $\regularizedDistance{Q}[i_0][[j_0,j_1]] = \regularizedDistance{Q + t V}[i_0][[j_0,j_1]]$ for all $t \in \R$.
	Consequently, the directional derivative of $f^\mu$ at $Q$, in the direction of $V$, is equal to zero.
	Taking into account the definition \eqref{eq:complete_metric_for_planar_triangular_meshes} of the metric, the claim follows.
\end{proof}

An immediate consequence of \cref{proposition:complete_metric_agrees_with_euclidean_for_translations} is that geodesics with respect to the metric \eqref{eq:complete_metric_for_planar_triangular_meshes}, whose initial tangent vectors~$V$ represent a translation, will be identical to Euclidean geodesics, \ie, $\gamma(t) = Q + t \, V$ holds, provided that $\beta_3 = 0$ holds.

\section{Geodesic Equation}
\label{section:geodesic_equation}

Geodesics are acceleration-free curves $t \mapsto \gamma(t)$ on a Riemannian manifold.
They are uniquely defined by a choice of an initial point $\gamma(0)$ and initial velocity $\dot \gamma(0)$.
In case of the Euclidean metric, geodesics are straight lines with constant velocity.
As illustrated in \Cref{fig:euclidean_geodesic}, they need not, in general, exist for all $t$ on $\planarmanifold$.
Therefore, we introduced the metric \eqref{eq:complete_metric_for_planar_triangular_meshes}, which renders the manifold of planar triangular meshes $\planarmanifold$ geodesically complete.

Suppose that $q^a$, $a = 1, \ldots, n$ are local coordinates on a Riemannian manifold~$\cM$ of dimension~$n$, and let $\gamma^a \coloneqq q^a \circ \gamma$ denote the coordinates of a curve $\gamma$.
Then $\gamma$ is a geodesic if and only if its coordinate curves solve the following system of nonlinear second-order ordinary differential equations
\begin{equation}
	\label{geodesic_equation}
	\ddot{\gamma}^c + \Gamma_{ab}^c \, \dot{\gamma}^a \, \dot{\gamma}^b
	=
	0
	,
	\quad
	c = 1, \ldots, n
	.
\end{equation}
In \eqref{geodesic_equation}, $\Gamma_{ab}^c$ is evaluated at $\gamma(t)$.
Here and it the following we use Einstein's summation convention.
Furthermore, $\Gamma^c_{ab}$ are the Christoffel symbols defined by
\begin{equation}
	\label{eq:christoffel_symbols}
	\Gamma^c_{ab}
	=
	\dfrac{1}{2} g^{cd} \paren[auto](){\dfrac{\partial g_{da}}{\partial q^b}+\dfrac{\partial g_{db}}{\partial q^a}-\dfrac{\partial g_{ab}}{\partial q^d}}
\end{equation}
where $g_{ab}$ and $g^{ab}$ are the components of the metric tensor and of its inverse, respectively, evaluated at $\gamma(t)$.

The unique solution of \eqref{geodesic_equation} with initial conditions $\gamma(0) = Q \in \cM$ and $\dot \gamma(0) = V \in \tangent{Q}[\cM]$ is denoted by $\geodesic<p>{Q}{V}(t)$.
It exists for some open interval containing~$0$.
\begin{definition}
	\label{definition:exponential_map}
	Let $Q \in \cM$ and $V \in \tangent{Q}[\cM]$.
	The \textbf{exponential map} at $Q$ is defined as
	\begin{equation*}
		V
		\mapsto
		\exponential{Q}{V}
		\coloneqq
		\geodesic<p>{Q}{V}(1)
		.
	\end{equation*}
\end{definition}

Generally, the exponential map is defined only for $V$ in a neighborhood of zero in $\tangent{Q}[\cM]$.
However, on geodesically complete manifolds, $\exponential{Q}$ is defined on all of $\tangent{Q}[\cM]$ thanks to the Hopf--Rinow theorem (see, \eg, \cite[Ch.~7, Thm.~2.8, p.146]{DoCarmo1992}) and $\exponential{Q}$ is a global diffeomorphism for all $Q \in \cM$.

\subsection{Hamiltonian Formulation for Geodesics}
\label{subsection:hamiltonian_approach}

Unfortunately, the geodesic equation \eqref{geodesic_equation} for the manifold of planar triangular meshes $\planarmanifold$ \wrt\ the Riemannian metric \eqref{eq:complete_metric_for_planar_triangular_meshes} does not have a closed form solution.
Therefore, we need to resort to a numerical integrator.
Our method of choice here is the Störmer--Verlet scheme, which is based on the Hamiltonian formulation of the geodesic equation and which preserves the kinetic energy of the initial velocity.

In the following we denote by $q^a$, $a = 1, \ldots, 2 \, N_V$, the coordinates of a point $Q \in \planarmanifold$ with respect the chart $\vvec \colon \R^{2 \times N_V} \to \R^{2 \, N_V}$, which stacks $Q \in \R^{2 \times N_V}$ column by column.
Similarly, we denote by $v^a$ the components of a tangent vector $V \in \tangent{Q}[\planarmanifold]$ in the chart induced basis of the tangent space.
Finally, the components of a cotangent vector $P \in \cotangent{Q}[\planarmanifold]$ are denoted by $p_a$.

The geodesic equation \eqref{geodesic_equation} can be understood as the Euler--Lagrange equation
\begin{equation}
	\label{eq:Euler-Lagrange_equation}
	\dfrac{\d}{\d t}\left(\dfrac{\partial \cL(Q,V)}{\partial v^a}\right)- \dfrac{\partial \cL(Q,V)}{\partial q^a}
	=
	0
\end{equation}
associated with the \textbf{Lagrangian}
\begin{equation}
	\label{eq:Lagrangian}
	\cL(Q,V)
	\coloneqq
	\dfrac{1}{2} g_{ab}(Q) \, v^a \, v^b
\end{equation}
defined on the tangent bundle $\tangentBundle[\planarmanifold]$; see for instance \cite[Ch.~5, Thm.~31.2.1, p.321]{DubrovinFomenkoNovikov1992}.

The derivation of the Hamiltonian formulation of the geodesic can be found, for instance, in \cite[§33]{DubrovinFomenkoNovikov1992}.
The corresponding \textbf{Hamiltonian} is formulated in terms of position~$Q$ and \textbf{momentum}~$p_a = g_{ab}(Q) \, v^b$.
It is defined as
\begin{equation}
	\label{eq:Hamiltonian}
	\cH(Q,P)
	\coloneqq
	PV - \cL(Q,V)
	=
	\frac{1}{2}
	g^{ab}(Q) \, p_a \, p_b
\end{equation}
on the cotangent bundle $\cotangentBundle[\planarmanifold]$.
Therefore, the Hamiltonian system is given by
\begin{subequations}
	\label{eq:Hamiltonian_system}
	\begin{align}
		\dot{p_a}
		&
		=
		- \dfrac{\partial \cH(Q,P)}{\partial q^a}
		=
		\dfrac{1}{2}\dfrac{\partial g_{bd}(Q)}{\partial q^a} \, g^{bc}(Q) \, p_c \, g^{de}(Q) \, p_e
		\label{eq:Hamiltonian_system_1}
		\\
		\dot{q^a}
		&
		=
		+ \dfrac{\partial \cH(Q,P)}{\partial p_a}
		=
		g^{ab}(Q) \, p_b
		.
		\label{eq:Hamiltonian_system_2}
	\end{align}
\end{subequations}
In \eqref{eq:Hamiltonian_system_1} we used the rule of the derivative of the inverse function.

An advantage of the Hamiltonian structure \eqref{eq:Hamiltonian_system} is that we can use energy preserving integrators, such as the Störmer--Verlet scheme.

\subsection{Störmer--Verlet Scheme for Hamiltonian Systems}
\label{subsection:stoermer_verlet}

The flow of a Hamiltonian system as in \eqref{eq:Hamiltonian_system} on a Riemannian manifold~$\cM$ is defined as $\varphi_t \colon \cotangentBundle[\cM] \to \cotangentBundle[\cM]$, mapping the initial condition $(Q^0,P^0)$ to the solution at time~$t$.
An important geometric property of a Hamiltonian systems is that its flow $\varphi_t$ is \emph{symplectic}, \ie, the derivative $\varphi_t' = \frac{\partial \varphi_t}{\partial(Q,P)}$ of the flow satisfies
\begin{equation*}
	\paren[auto](){\varphi_t'}^\transp J \, \varphi_t'
	=
	J
	\quad
	\text{with }
	J
	=
	\begin{bmatrix}
		0 & \id \\ -\id & 0
	\end{bmatrix}
	,
\end{equation*}
where $\id$ denotes the identity matrix of the dimension of $\cM$.
A numerical integrator which preserves this property is said to be \emph{symplectic}.
A prominent example is the Störmer--Verlet scheme, see for instance \cite[eq.(2.10)]{HairerLubichWanner2003}.
Its application to the Hamiltonian system \eqref{eq:Hamiltonian_system} is given in \Cref{algorithm:stoermer-verlet}.

\begin{algorithm}[Störmer--Verlet scheme for the geodesic equation \eqref{geodesic_equation} in Hamiltonian form \eqref{eq:Hamiltonian_system}] \hfill
  \label{algorithm:stoermer-verlet}
	\begin{algorithmic}[1]
		\Require abstract simplicial complex $\Delta$ with $N_V$ vertices
		\Require matrix $\Qref \in \plusmanifold \subset \R^{2 \times N_V}$ of reference vertex positions
		\Require matrix $Q^0 \in \planarmanifold \subset \R^{2 \times N_V}$ of initial vertex positions
		\Require initial tangent vector $V^0 \in \R^{2 \times N_V}$
		\Require final time~$T$; number of time steps~$N$; time step size $\Delta t \coloneqq \frac{T}{N}$
		\Ensure approximate solution of the geodesic equation \eqref{geodesic_equation} with initial conditions $\gamma(0) = Q^0$ and $\dot \gamma(0) = V^0$ on $\planarmanifold$ at times $t_n = n \, \Delta t$, $n = 0, \ldots, N$
		\State  set initial momentum $p_a^0 = g_{ab}(Q^0) \, v^{0,b}$\;
		\For{$n \gets 0$ \textbf{to} $N-1$}
		\State solve $p_a^{n+\nicefrac12} = p_a^n - \dfrac{\Delta t}{4}\paren[auto]\{\}{\dfrac{\partial g_{bd}(Q^n)}{\partial q^a} \, g^{bc}(Q^n) \, p_c^{n+\nicefrac12} \, g^{de}(Q^n) \, p_e^{n+\nicefrac12}}$
		\newline
		\hspace*{\fill} for $a = 1, \ldots, 2 \, N_V$
		\label{line:implicit_equation_for_P}
		\State solve $q^{n+1,a} = q^{n,a} + \dfrac{\Delta t}{2} \paren[auto]\{\}{g^{ab}(Q^n) \, p_a^{n+\nicefrac12} + g^{ab}(Q^{n+1}) \, p_a^{n+\nicefrac12}}$
		\newline
		\hspace*{\fill} for $a = 1, \ldots, 2 \, N_V$
		\label{line:implicit_equation_for_Q}
		\State set $p_a^{n+1} = p_a^{n+\nicefrac12} - \dfrac{\Delta t}{4}\paren[auto]\{\}{\dfrac{\partial g_{bd}(Q^{n+1})}{\partial q^a} \, g^{bc}(Q^{n+1}) \, p_c^{n+\nicefrac12} \, g^{de}(Q^{n+1}) \, p_e^{n+\nicefrac12}}$
		\newline
		\hspace*{\fill} for $a = 1, \ldots, 2 \, N_V$
		\label{line:explicit_equation_for_P}
		\State solve $v^{n+1,a} = g^{ab}(Q^{n+1}) \, p_b^{n+1}$ for $a = 1, \ldots, 2 \, N_V$
		\label{line:solve_V}
		\EndFor
		\State \Return $Q^0, \ldots, Q^N$, approximating $\gamma$ at $t_0, \ldots, t_N$
		\State \Return $V^0, \ldots, V^N$, approximating $\dot \gamma$ at $t_0, \ldots, t_N$
		\State \Return $P^0, \ldots, P^N$, approximating the momentum $p_a = g_{ab}(Q) \, v^b$ at $t_0, \ldots, t_N$
	\end{algorithmic}
\end{algorithm}

In our implementation, we use a fixed-point iteration to solve the implicit equations for $P^{n+\nicefrac12}$ and $Q^{n+1}$ in \cref{line:implicit_equation_for_P,line:implicit_equation_for_Q}.

More details on an efficient implementation which takes advantage of the specific structure of the Riemannian metric tensor \eqref{eq:complete_metric_for_planar_triangular_meshes} and its inverse is beyond the scope of this paper.

\section{Numerical Experiments}
\label{section:numerical_experiments}

The purpose of this section is to numerically investigate how planar triangular meshes deform under the complete Riemannian metric \eqref{eq:complete_metric_for_planar_triangular_meshes}, and to compare it with the Euclidean metric.
To this end, we implemented the Störmer--Verlet scheme (\Cref{algorithm:stoermer-verlet}) to integrate the geodesic equation numerically.

Our experiments are structured as follows.
In \cref{subsection:elementary_transformations}, we investigate elementary transformations (translation, shearing, scaling, and rotation) of a square mesh with crossed diagonals.
In \cref{subsection:comparison_complete_incomplete_metric} we revisit the example depicted in \Cref{fig:euclidean_geodesic}, whose initial tangent vector gives rise to a more complex transformation dynamic.
The aforementioned experiments confirm that the proposed metric successfully avoids self-intersections of the mesh due to its completeness.
To study this also quantitatively, we conduct in \cref{subsection:mesh_quality} an experiment using a slightly more complex mesh, where we evaluate a mesh quality measure along the geodesics.

To be precise, all of our experiments are based on the function
\begin{equation}
	\label{eq:f_mu_without_exterior_term}
	f^\mu(Q)
	\coloneqq
	\sum_{k=1}^{N_T} \sum_{\ell=0}^2 \frac{\beta_1 }{{\height[auto]{Q}{\ell}[i^k_0,i^k_1,i^k_2]}}
	+
	\frac{\beta_3}{2} \norm{Q - \Qref}_F^2
\end{equation}
which is used to construct the metric
\begin{equation}
	\label{eq:complete_metric_for_planar_triangular_meshes_without_exterior_term}
	g_{ab} = \delta_a^b + \dfrac{\partial f^\mu}{\partial q^a} \dfrac{\partial f^\mu}{\partial q^b}
	,
	\quad
	a, b = 1, \ldots, 2 \, N_V
\end{equation}
as in \eqref{eq:complete_metric_for_planar_triangular_meshes}.
The omission of the $\beta_2$-term in \eqref{eq:f_mu_without_exterior_term} preventing exterior self-intersections is justified by \cref{remark:f_mu_with_threshold}.
For the experiments shown in this paper, exterior self-intersections are never a factor, which we verified a~posteriori, and thus we can choose a cut-off function~$\chi_2$ which effectively removes the second term in \eqref{eq:f_mu}.
The choice of parameters $\beta_1, \beta_3 > 0$, which control the relative importance of each term in $f^\mu$ in relation to the Euclidean base metric in \eqref{eq:complete_metric_for_planar_triangular_meshes_without_exterior_term}, is described below for each experiment individually.

\subsection{Elementary Mesh Transformations}
\label{subsection:elementary_transformations}

In this section we showcase the deformation of a simple mesh under a number of elementary transformations.
To be precise, we consider initial tangent vectors which would produce a translation, shearing, scaling, or rotation, respectively, of the mesh in the Euclidean setting ($\beta_1 = \beta_3 = 0$ in \eqref{eq:f_mu_without_exterior_term}).
In particular, we study numerically the influence of the parameters $\beta_1$ and $\beta_3$.

\begin{figure}[htb]
	\begin{center}
		\begin{subfigure}{0.3\textwidth}
			\centering
			\includegraphics[width=0.75\linewidth]{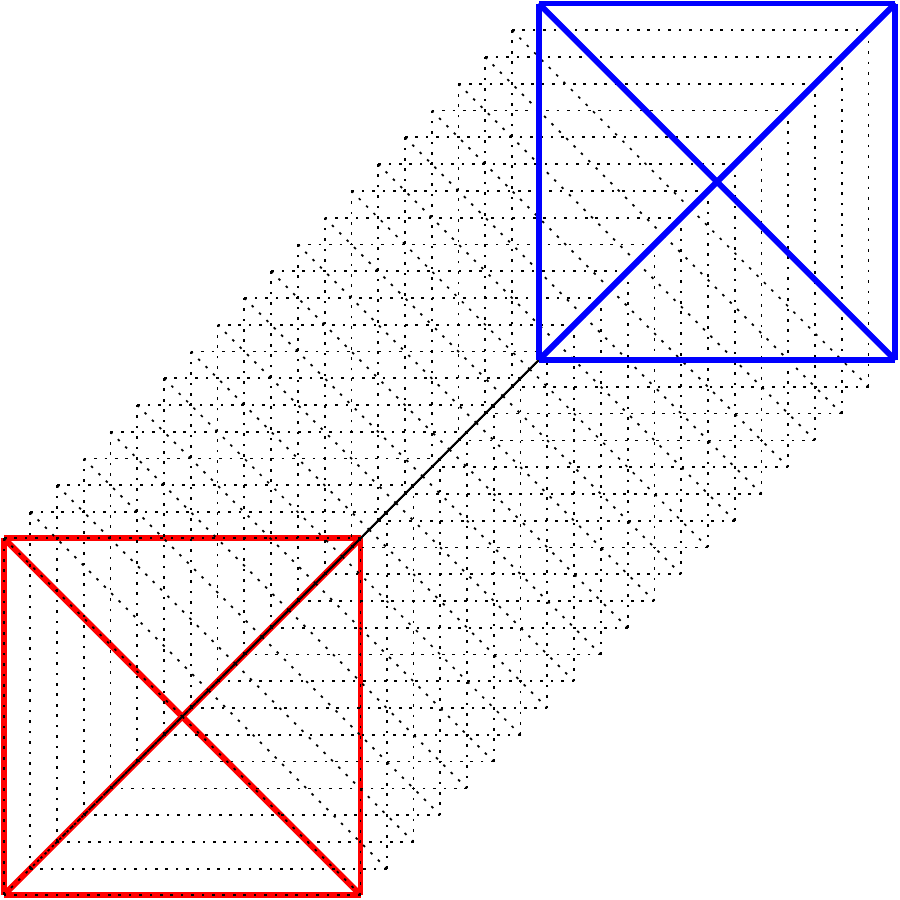}
			\caption{$\beta_1=0$, $\beta_3=0$}
			\label{subfig:translation_beta1_0_beta2_0}
		\end{subfigure}
		\hfill
		\begin{subfigure}{0.3\textwidth}
			\centering
			\includegraphics[width=0.75\linewidth]{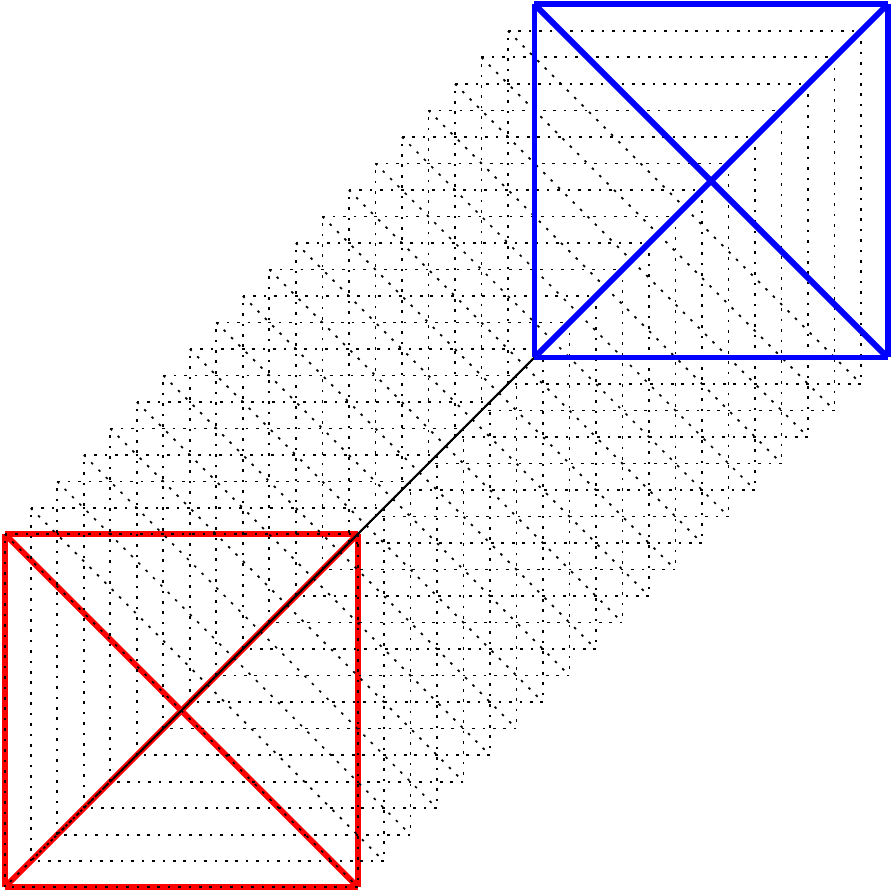}
			\caption{$\beta_1=0.5$, $\beta_3=0$}
			\label{subfig:translation_beta1_5_beta2_0}
		\end{subfigure}
		\hfill
		\begin{subfigure}{0.3\textwidth}
			\centering
			\includegraphics[width=0.75\linewidth]{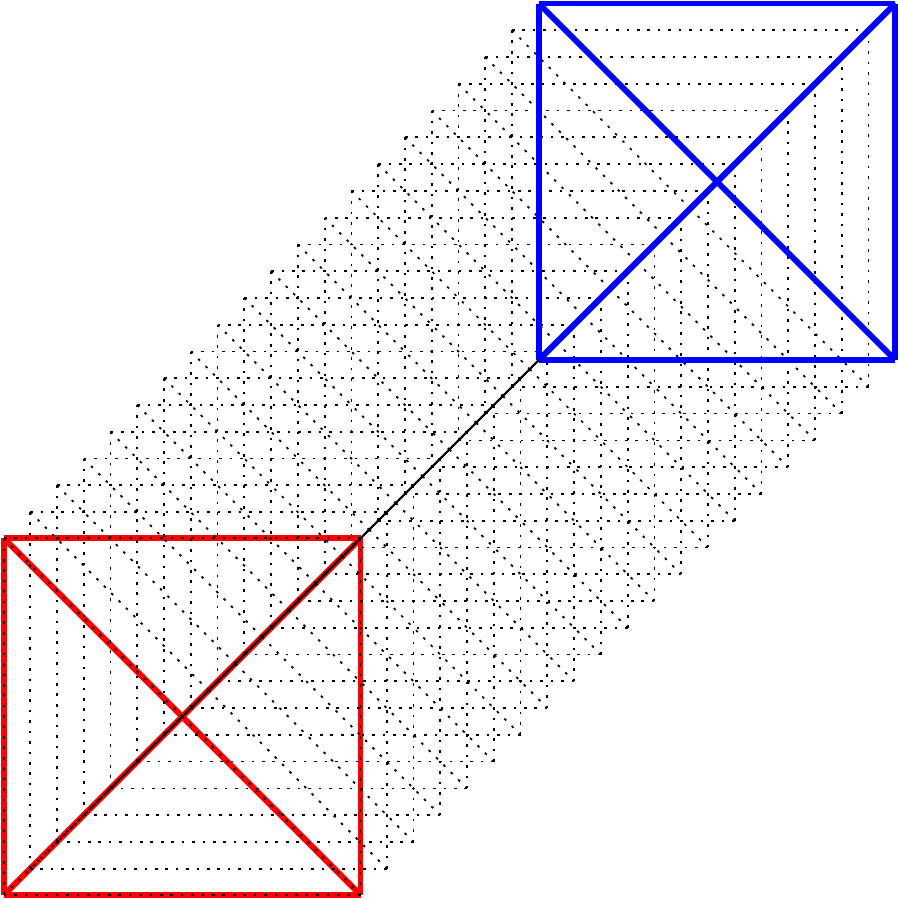}
			\caption{$\beta_1=1.0$, $\beta_3=0$}
			\label{subfig:translation_beta1_10_beta2_0}
		\end{subfigure}
		\\
		\vspace*{0.1cm}
		\begin{subfigure}{0.3\textwidth}
			\centering
			\includegraphics[width=0.75\linewidth]{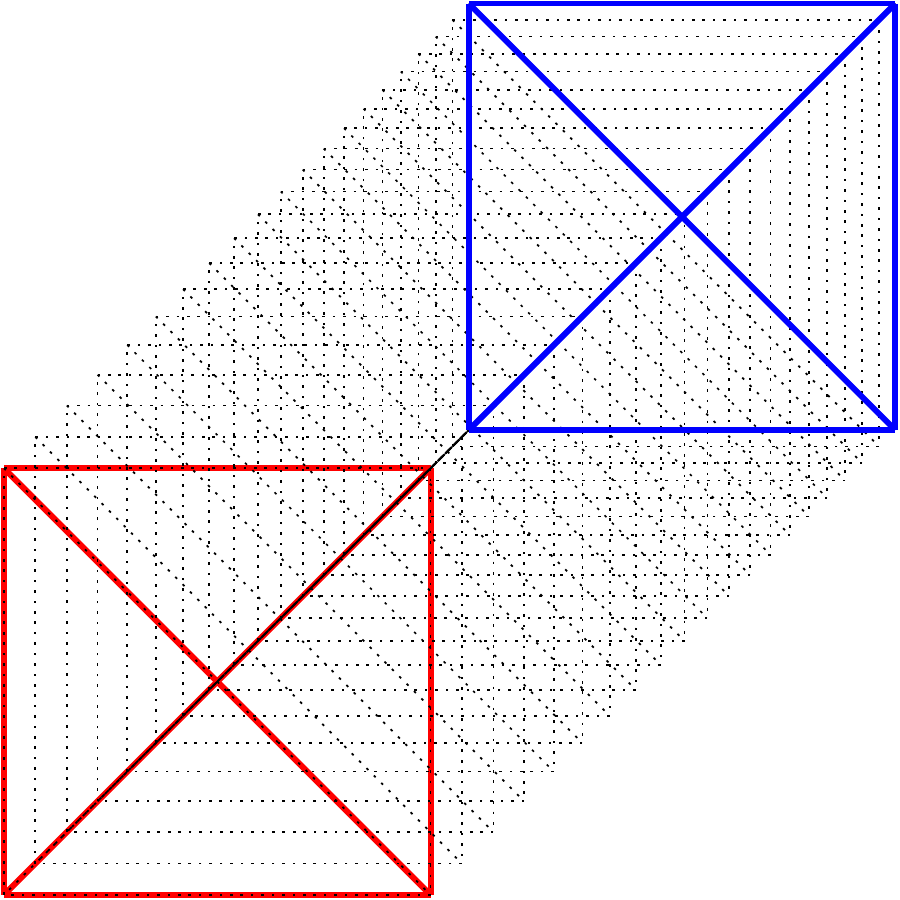}
			\caption{$\beta_1=0$, $\beta_3=0.5$}
			\label{subfig:translation_beta1_0_beta2_5}
		\end{subfigure}
		\hfill
		\begin{subfigure}{0.3\textwidth}
			\centering
			\includegraphics[width=0.75\linewidth]{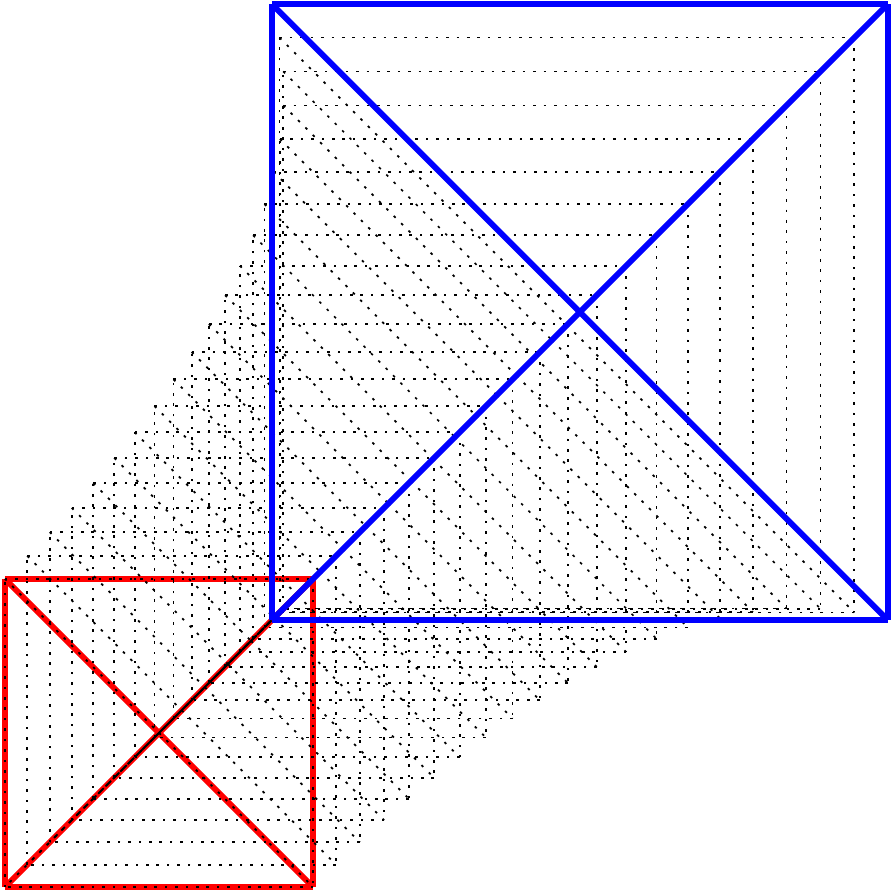}
			\caption{$\beta_1=0.5$, $\beta_3=0.5$}
			\label{subfig:translation_beta1_5_beta2_5}
		\end{subfigure}
		\hfill
		\begin{subfigure}{0.3\textwidth}
			\centering
			\includegraphics[width=0.75\linewidth]{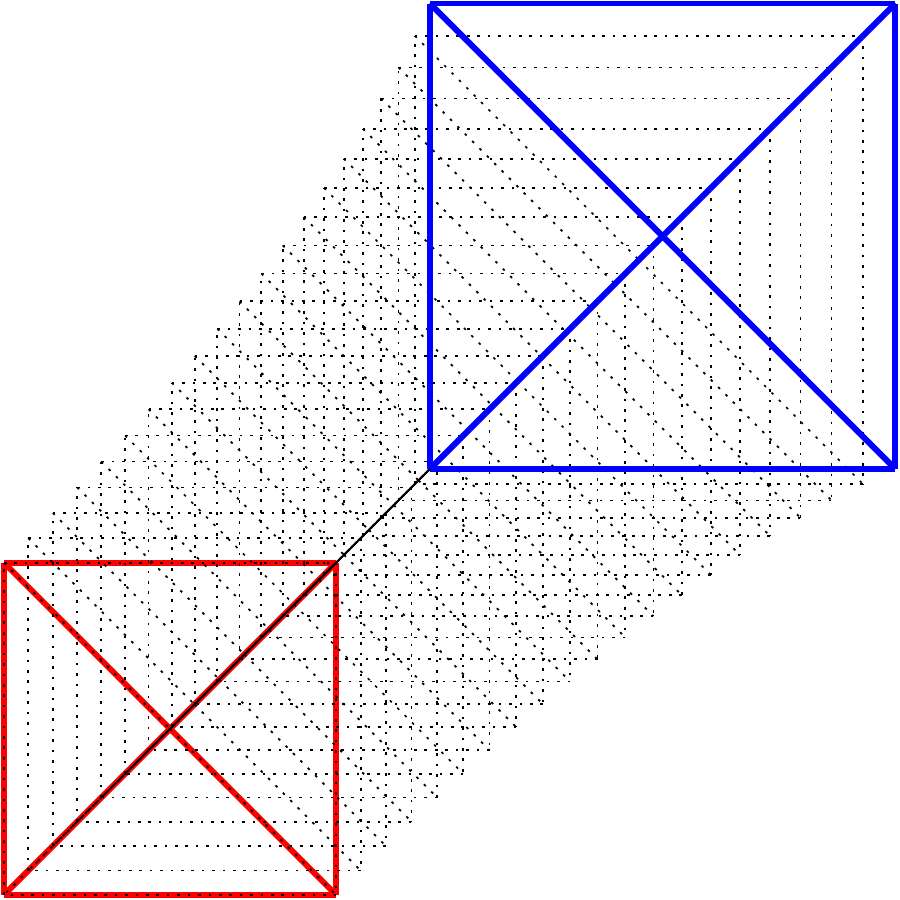}
			\caption{$\beta_1=1.0$, $\beta_3=0.5$}
			\label{subfig:translation_beta1_10_beta2_5}
		\end{subfigure}
		\\
		\vspace*{0.1cm}
		\begin{subfigure}{0.3\textwidth}
			\centering
			\includegraphics[width=0.75\linewidth]{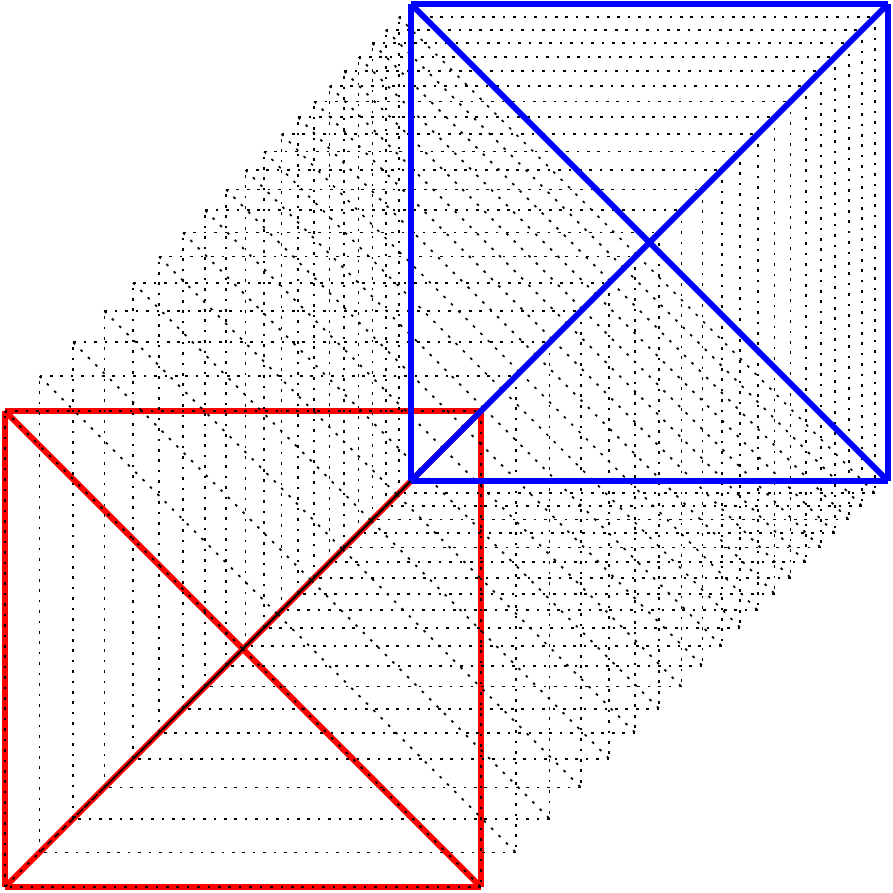}
			\caption{$\beta_1=0$, $\beta_3=1.0$}
			\label{subfig:translation_beta1_0_beta2_10}
		\end{subfigure}
		\hfill
		\begin{subfigure}{0.3\textwidth}
			\centering
			\includegraphics[width=0.75\linewidth]{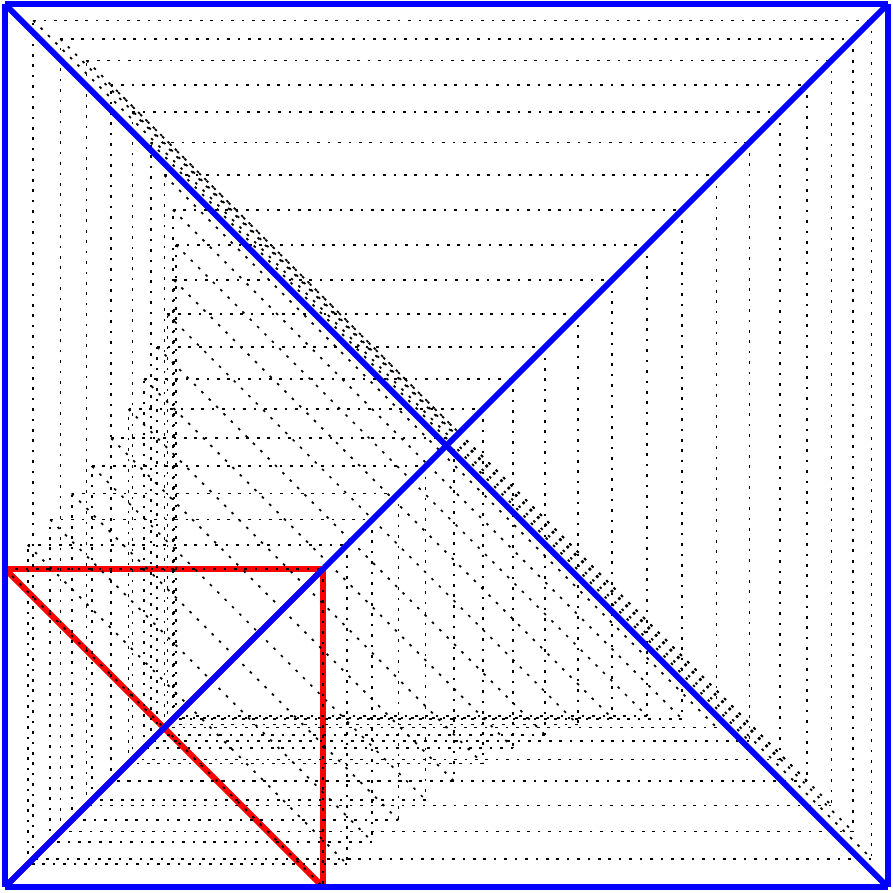}
			\caption{$\beta_1=0.5$, $\beta_3=1.0$}
			\label{subfig:translation_beta1_5_beta2_10}
		\end{subfigure}
		\hfill
		\begin{subfigure}{0.3\textwidth}
			\centering
			\includegraphics[width=0.75\linewidth]{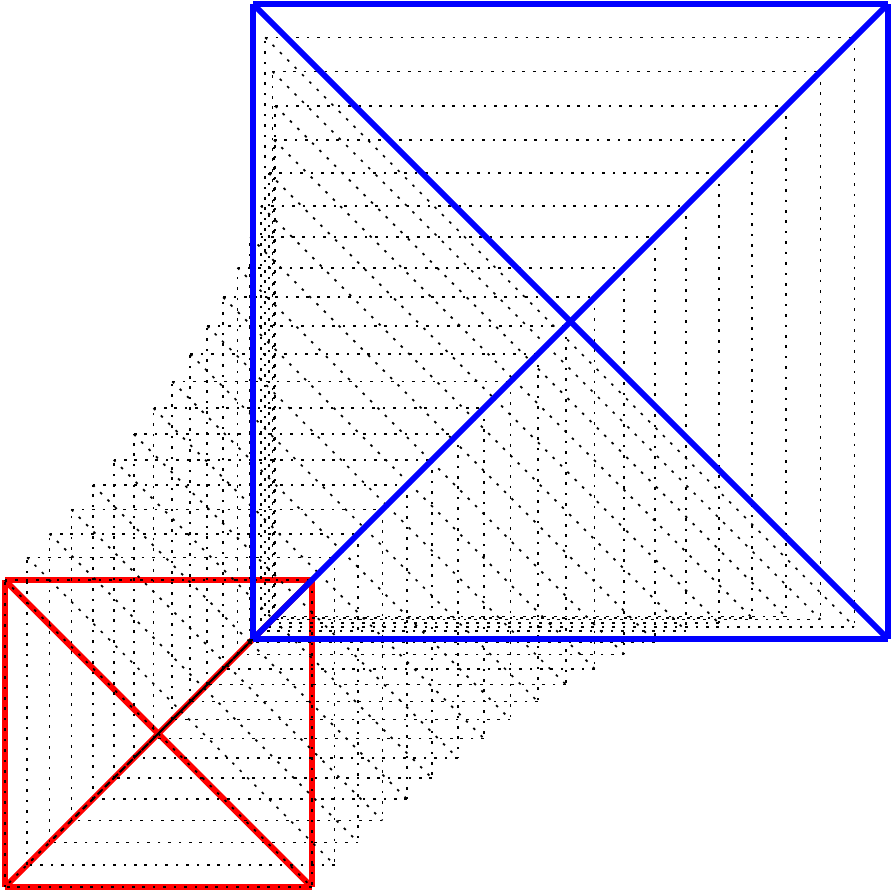}
			\caption{$\beta_1=1.0$, $\beta_3=1.0$}
			\label{subfig:translation_beta1_10_beta2_10}
		\end{subfigure}
	\end{center}
	\caption{20~snapshots of geodesics for different values of $\beta_1$, $\beta_3$, starting from the same initial mesh (shown in red) and produced by the same initial tangent vector, which induces a translation; see \cref{subsection:elementary_transformations}. The final mesh is shown in blue.}
	\label{fig:translations}
\end{figure}

\begin{figure}[htb]
	\begin{center}
		\begin{subfigure}[b]{0.3\textwidth}
			\centering
			\includegraphics[width=0.85\linewidth]{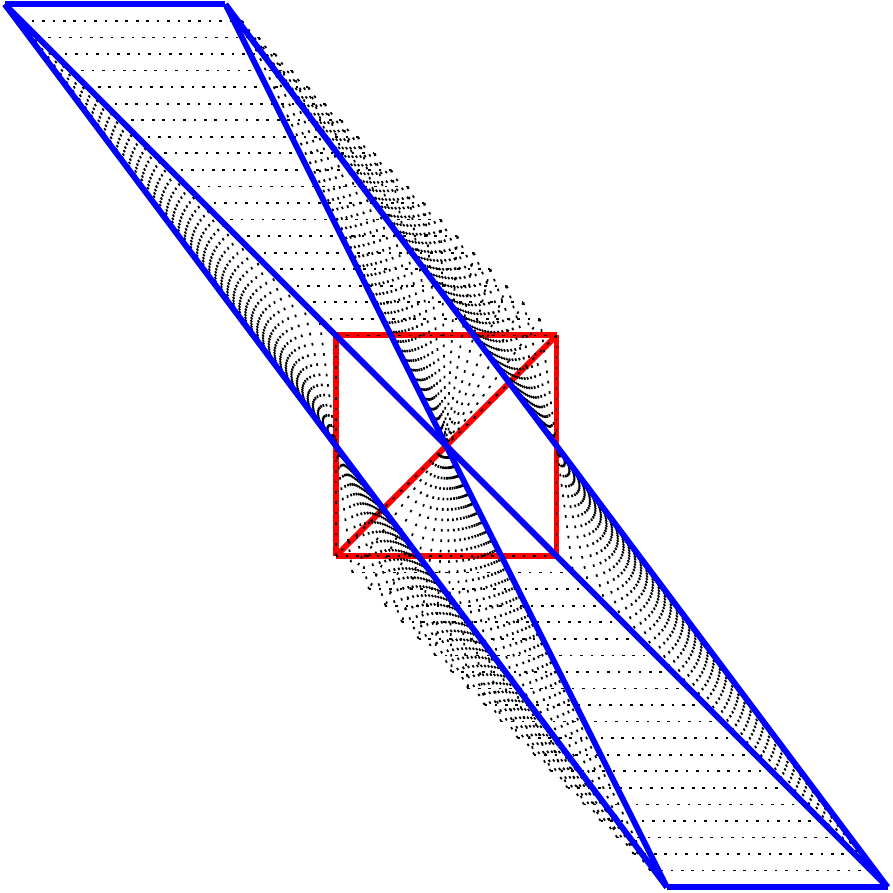}
			\caption{$\beta_1=0$, $\beta_3=0$}
			\label{subfig:shearing_beta1_0_beta2_0}
		\end{subfigure}
		\hfill
		\begin{subfigure}[b]{0.3\textwidth}
			\centering
			\includegraphics[width=0.85\linewidth]{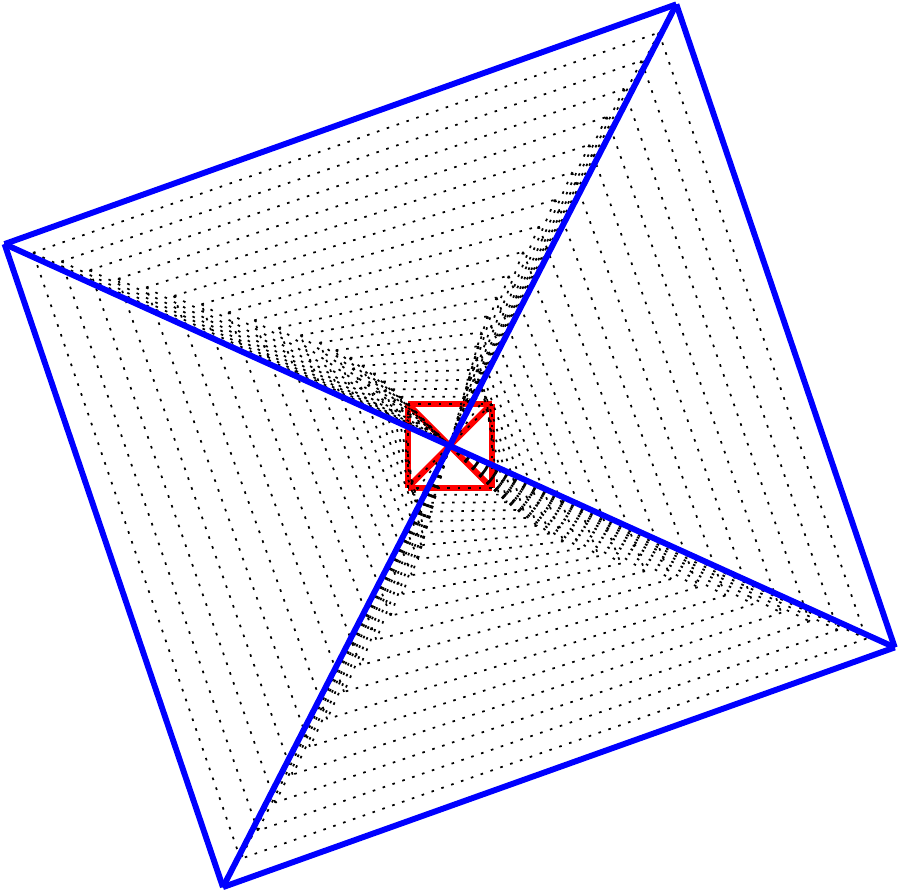}
			\caption{$\beta_1=0.5$, $\beta_3=0$}
			\label{subfig:shearing_beta1_5_beta2_0}
		\end{subfigure}
		\hfill
		\begin{subfigure}[b]{0.3\textwidth}
			\centering
			\includegraphics[width=0.85\linewidth]{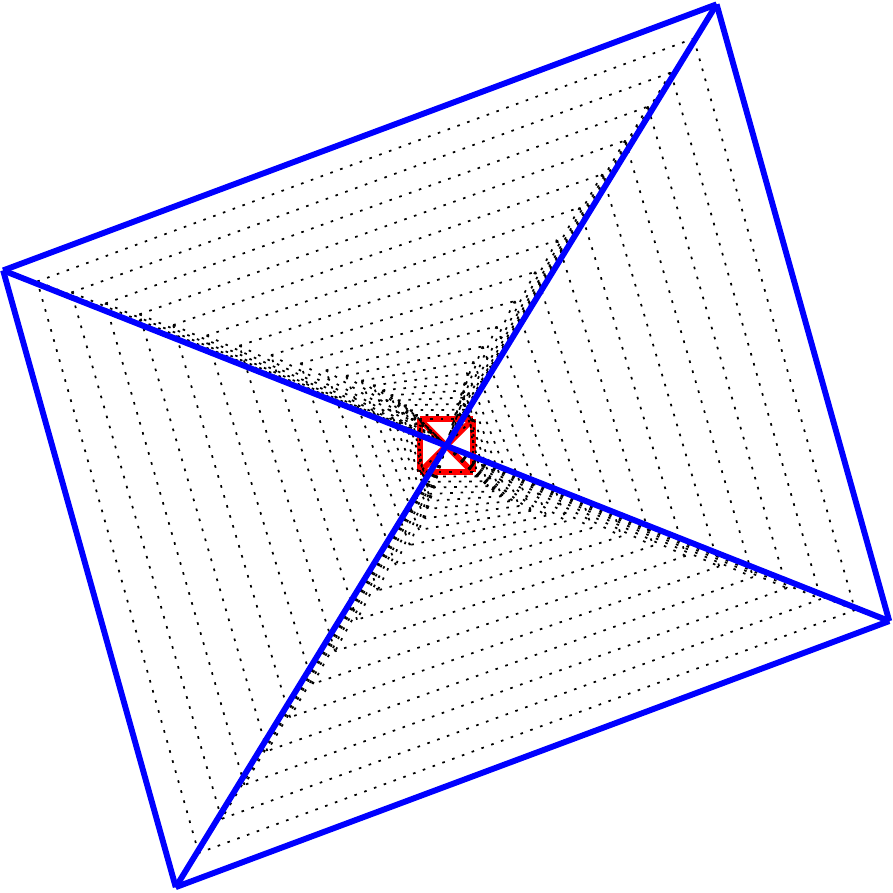}
			\caption{$\beta_1=1.0$, $\beta_3=0$}
			\label{subfig:shearing_beta1_10_beta2_0}
		\end{subfigure}
		\\
		\begin{subfigure}[b]{0.3\textwidth}
			\centering
			\includegraphics[width=0.85\linewidth]{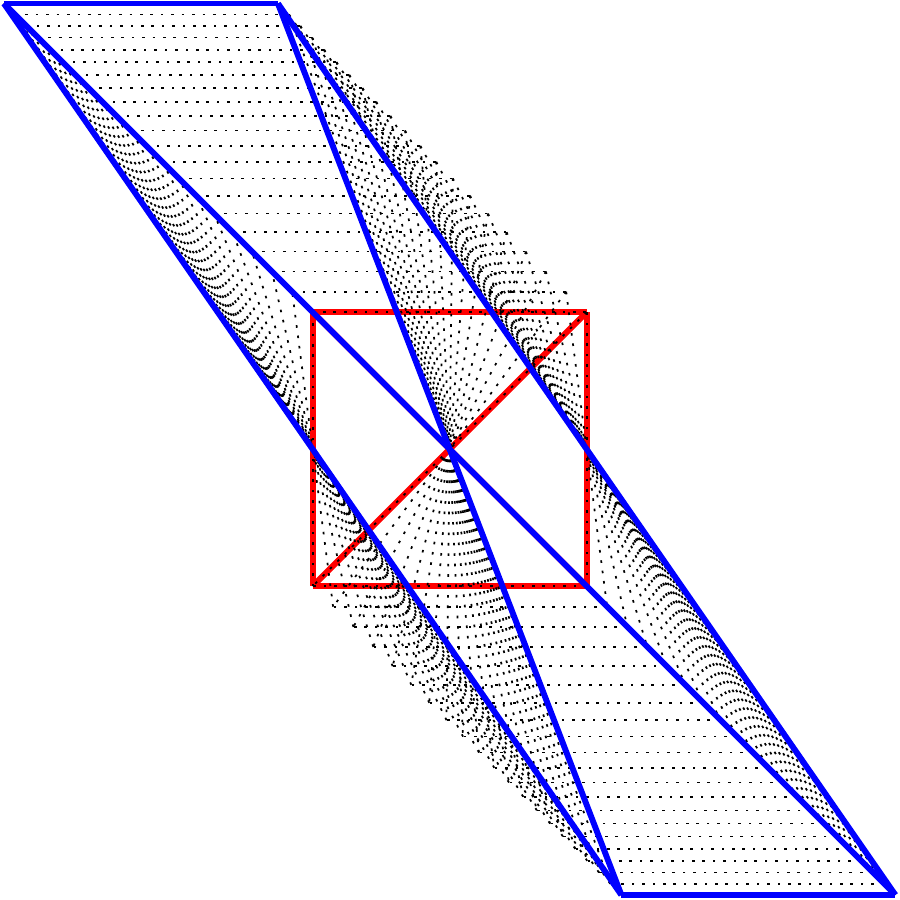}
			\caption{$\beta_1=0$, $\beta_3=0.5$}
			\label{subfig:shearing_beta1_0_beta2_5}
		\end{subfigure}
		\hfill
		\begin{subfigure}[b]{0.3\textwidth}
			\centering
			\raisebox{10mm}{%
				\includegraphics[width=0.85\linewidth]{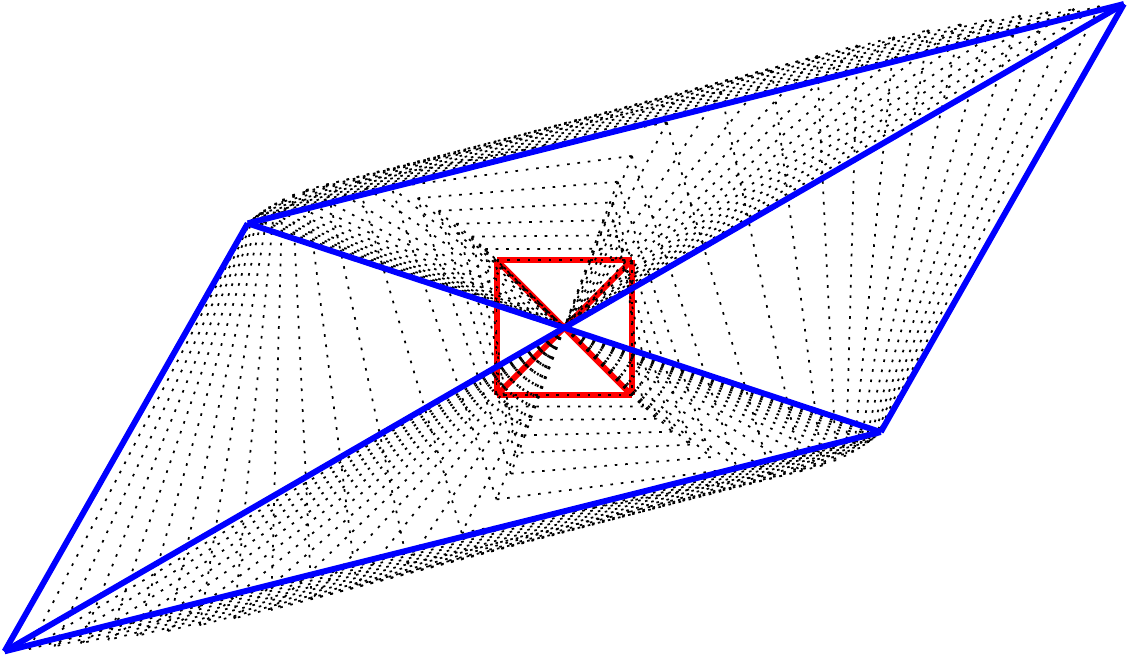}
			}
			\caption{$\beta_1=0.5$, $\beta_3=0.5$}
			\label{subfig:shearing_beta1_5_beta2_5}
		\end{subfigure}
		\hfill
		\begin{subfigure}[b]{0.3\textwidth}
			\centering
			\raisebox{10mm}{%
				\includegraphics[width=0.85\linewidth]{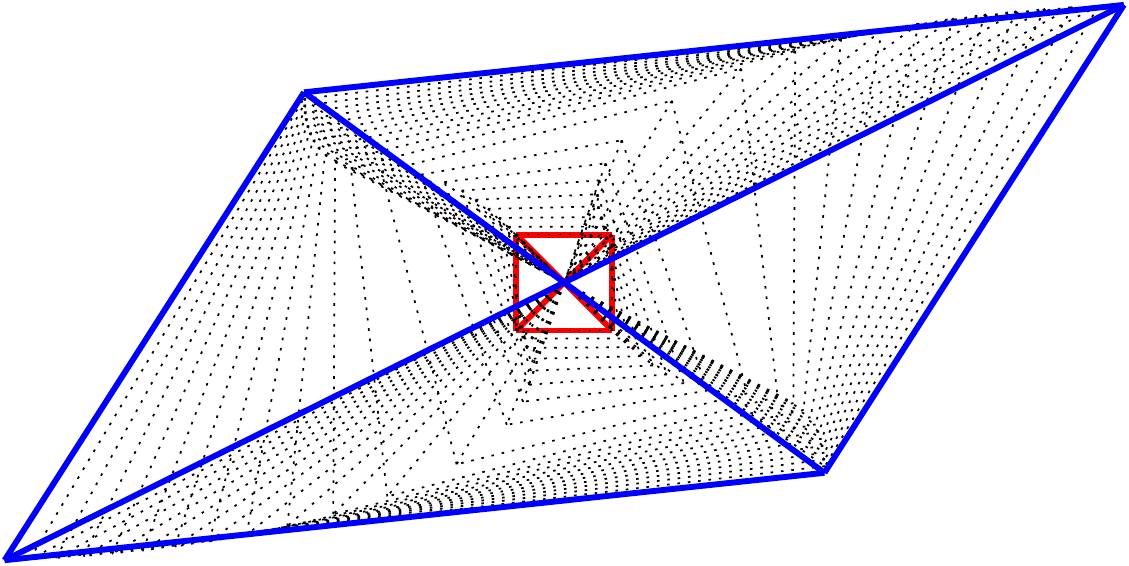}
			}
			\caption{$\beta_1=1.0$, $\beta_3=0.5$}
			\label{subfig:shearing_beta1_10_beta2_5}
		\end{subfigure}
		\\
		\begin{subfigure}[b]{0.3\textwidth}
			\centering
			\includegraphics[width=0.85\linewidth]{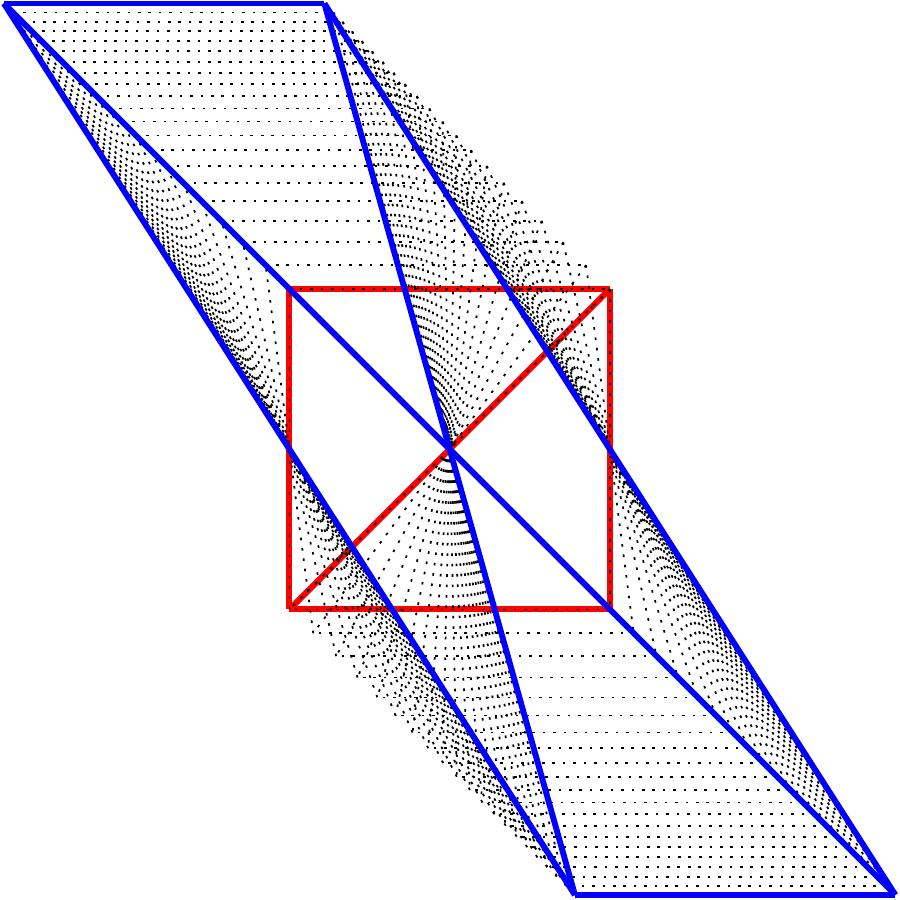}
			\caption{$\beta_1=0$, $\beta_3=1.0$}
			\label{subfig:shearing_beta1_0_beta2_10}
		\end{subfigure}
		\hfill
		\begin{subfigure}[b]{0.3\textwidth}
			\centering
			\raisebox{12mm}{%
				\includegraphics[width=0.85\linewidth]{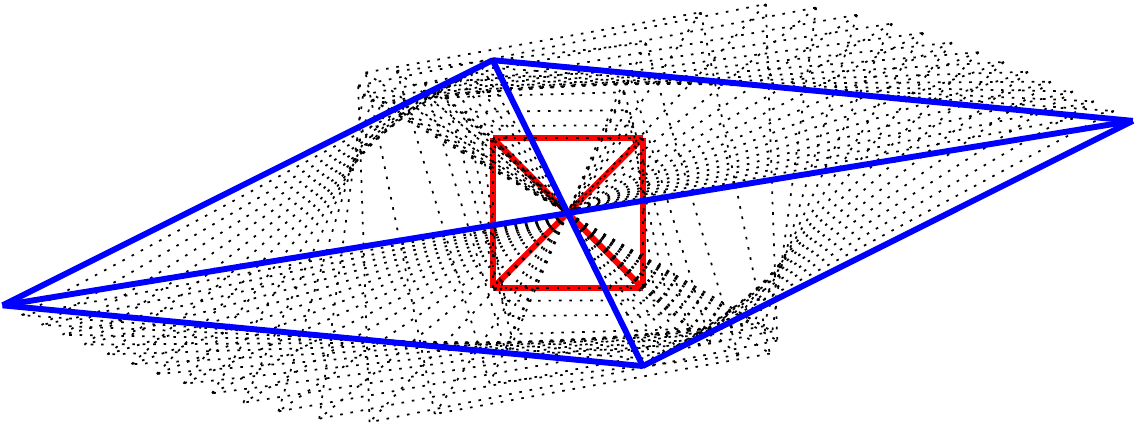}
			}
			\caption{$\beta_1=0.5$, $\beta_3=1.0$}
			\label{subfig:shearing_beta1_5_beta2_10}
		\end{subfigure}
		\hfill
		\begin{subfigure}[b]{0.3\textwidth}
			\centering
			\raisebox{8mm}{%
				\includegraphics[width=0.85\linewidth]{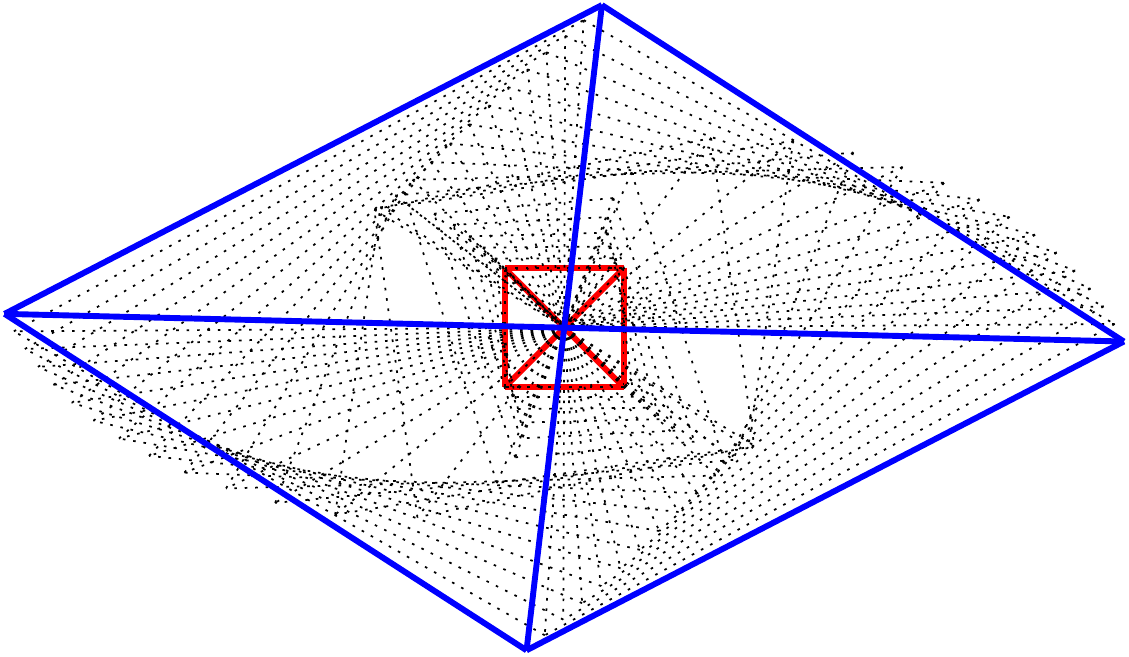}
			}
			\caption{$\beta_1=1.0$, $\beta_3=1.0$}
			\label{subfig:shearing_beta1_10_beta2_10}
		\end{subfigure}
	\end{center}
	\caption{20~snapshots of geodesics for different values of $\beta_1$, $\beta_3$, starting from the same initial mesh (shown in red) and produced by the same initial tangent vector, which induces shearing; see \cref{subsection:elementary_transformations}. The final mesh is shown in blue.}
	\label{fig:shearing}
\end{figure}

\begin{figure}[htb]
	\begin{center}
		\begin{subfigure}{0.3\textwidth}
			\centering
			\includegraphics[width=0.85\linewidth]{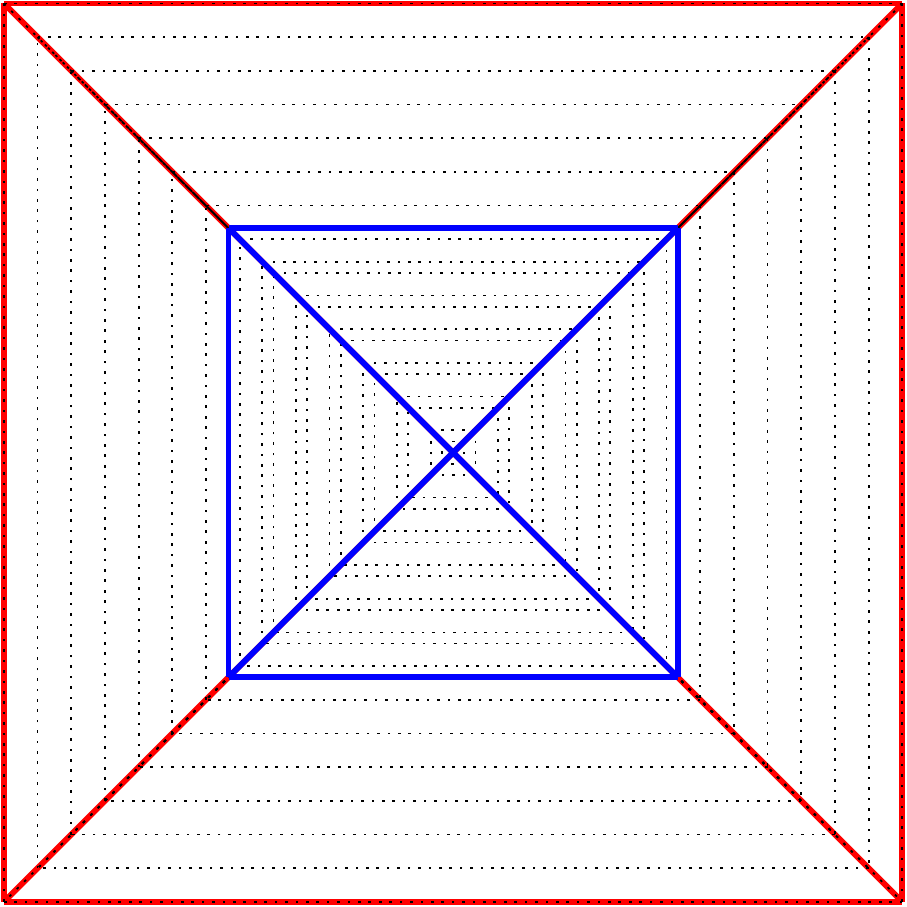}
			\caption{$\beta_1=0$, $\beta_3=0$}
			\label{subfig:scaling_beta1_0_beta2_0}
		\end{subfigure}
		\hfill
		\begin{subfigure}{0.3\textwidth}
			\centering
			\includegraphics[width=0.85\linewidth]{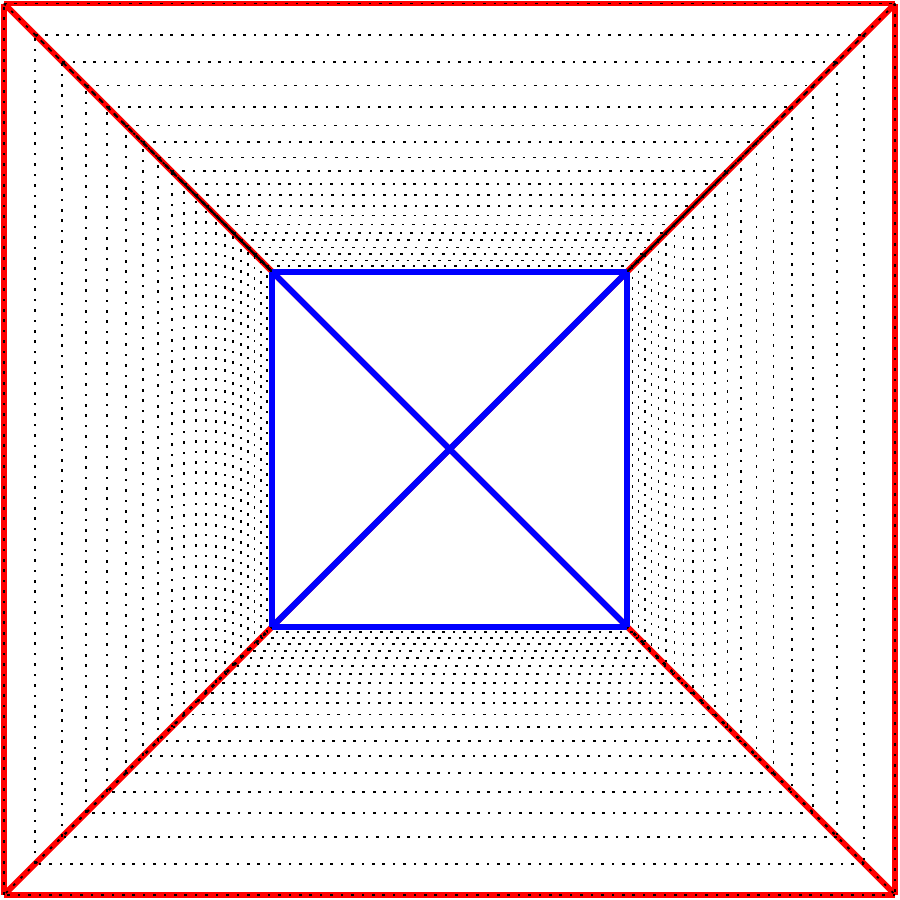}
			\caption{$\beta_1=0.5$, $\beta_3=0$}
			\label{subfig:scaling_beta1_5_beta2_0}
		\end{subfigure}
		\hfill
		\begin{subfigure}{0.3\textwidth}
			\centering
			\includegraphics[width=0.85\linewidth]{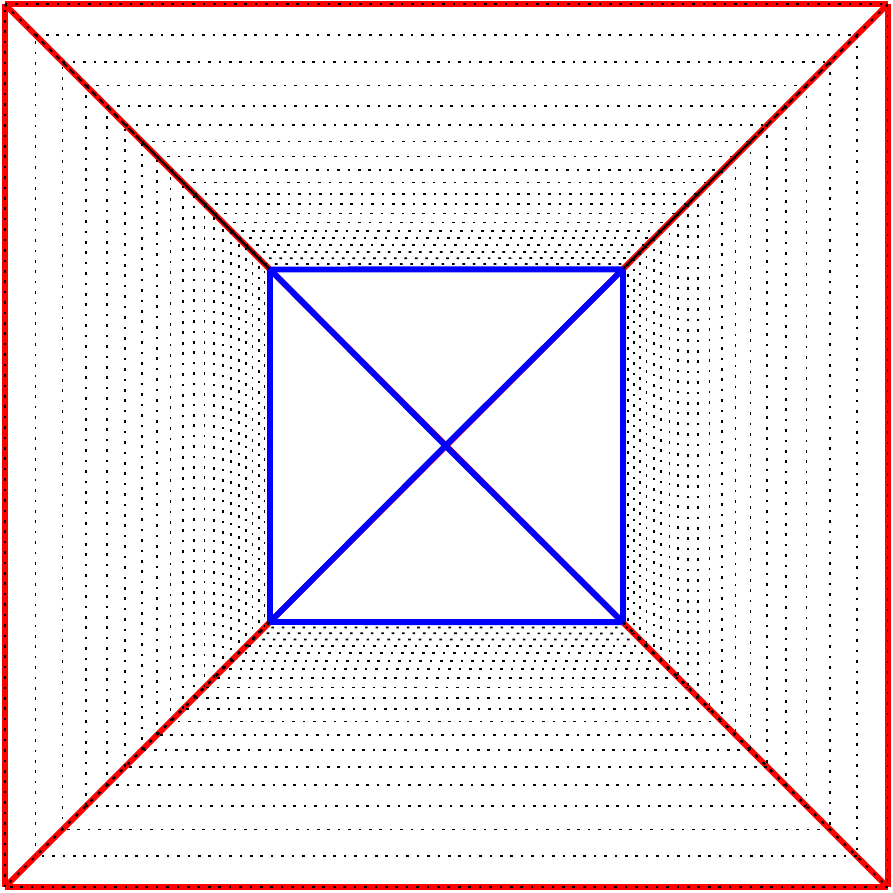}
			\caption{$\beta_1=1.0$, $\beta_3=0$}
			\label{subfig:scaling_beta1_10_beta2_0}
		\end{subfigure}
		\\
		\begin{subfigure}{0.3\textwidth}
			\centering
			\includegraphics[width=0.85\linewidth]{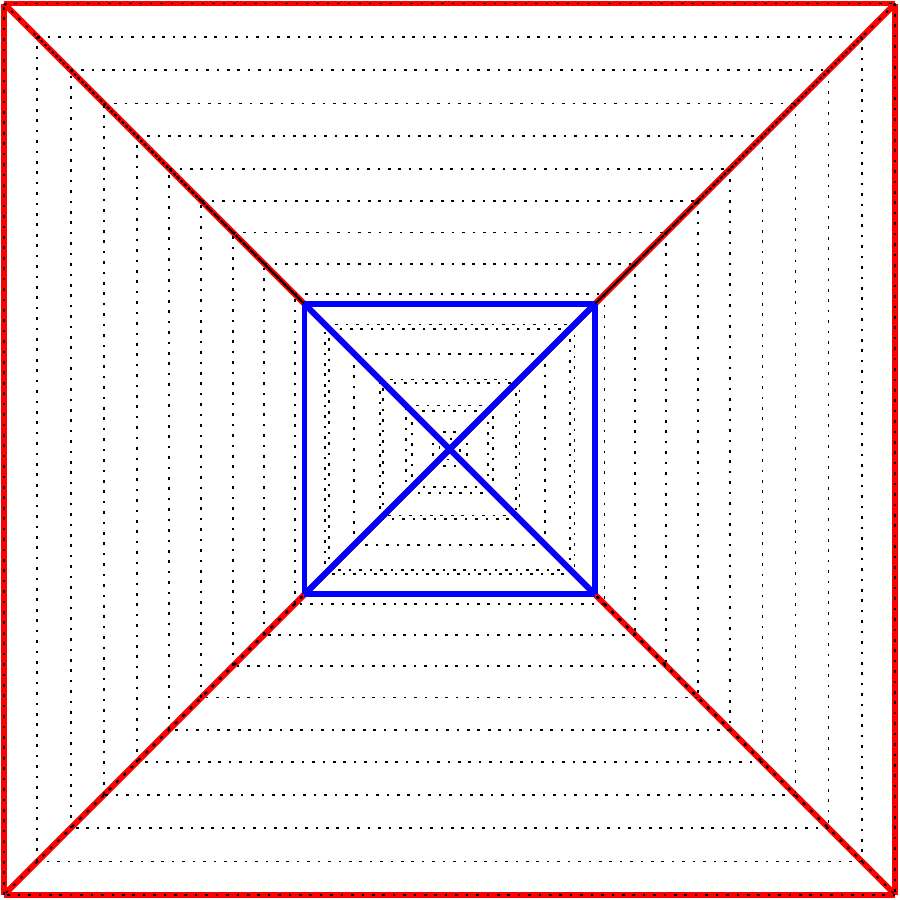}
			\caption{$\beta_1=0$, $\beta_3=0.5$}
			\label{subfig:scaling_beta1_0_beta2_5}
		\end{subfigure}
		\hfill
		\begin{subfigure}{0.3\textwidth}
			\centering
			\includegraphics[width=0.85\linewidth]{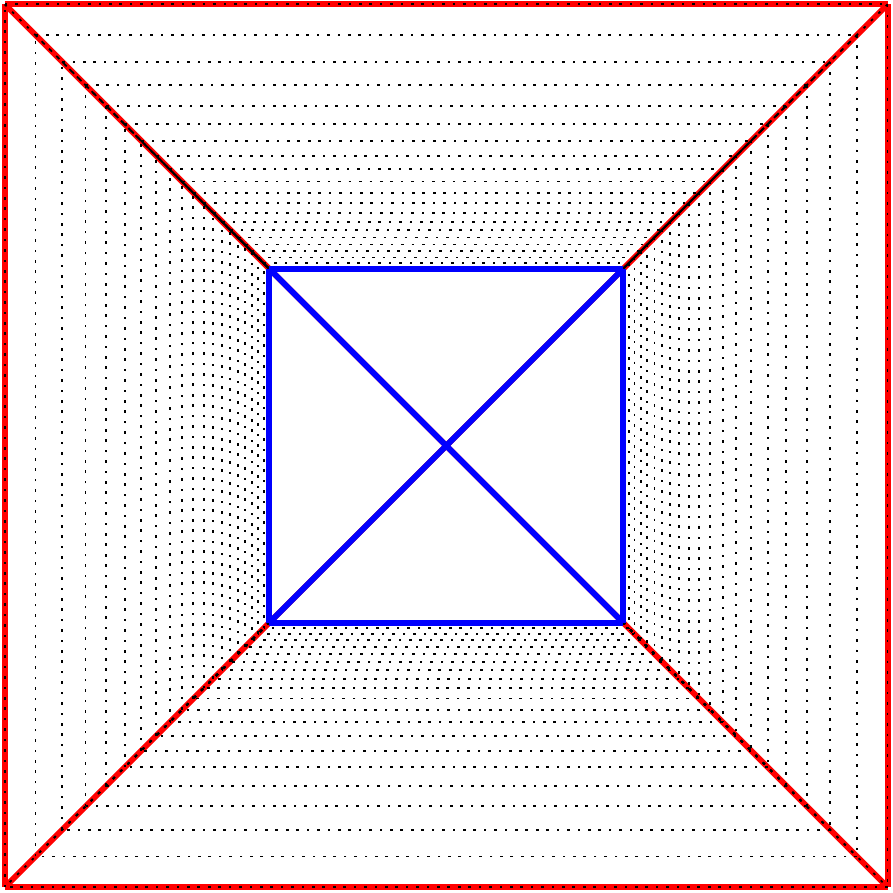}
			\caption{$\beta_1=0.5$, $\beta_3=0.5$}
			\label{subfig:scaling_beta1_5_beta2_5}
		\end{subfigure}
		\hfill
		\begin{subfigure}{0.3\textwidth}
			\centering
			\includegraphics[width=0.85\linewidth]{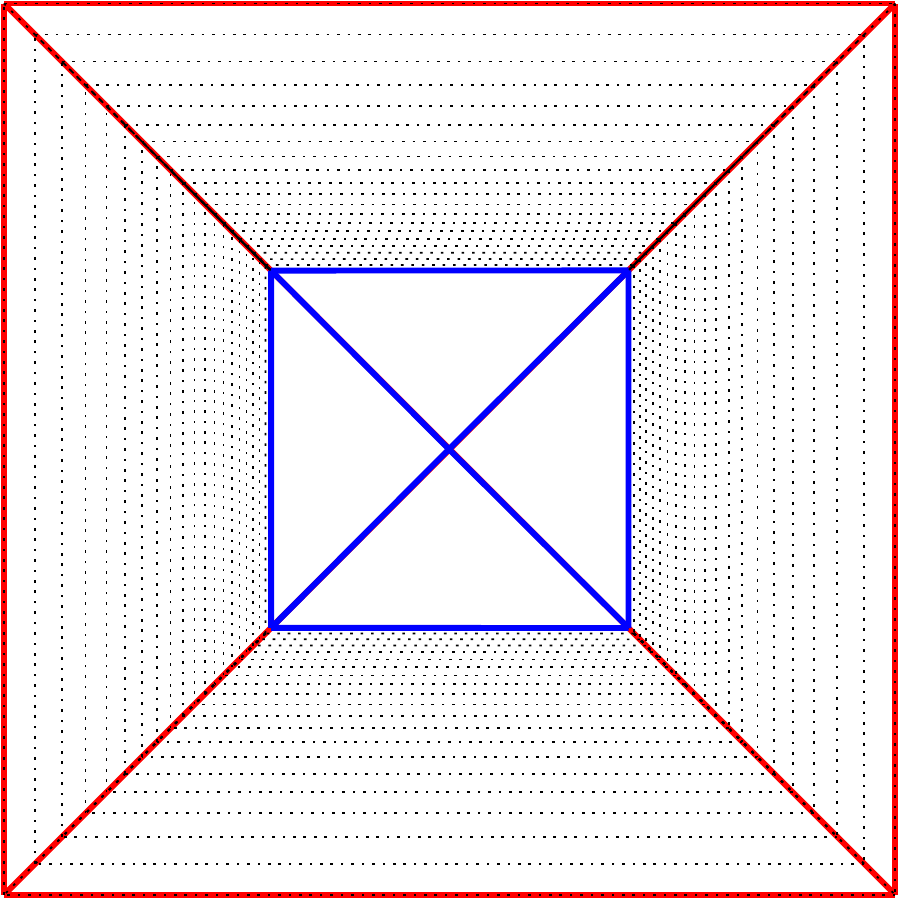}
			\caption{$\beta_1=1.0$, $\beta_3=0.5$}
			\label{subfig:scaling_beta1_10_beta2_5}
		\end{subfigure}
		\\
		\begin{subfigure}{0.3\textwidth}
			\centering
			\includegraphics[width=0.85\linewidth]{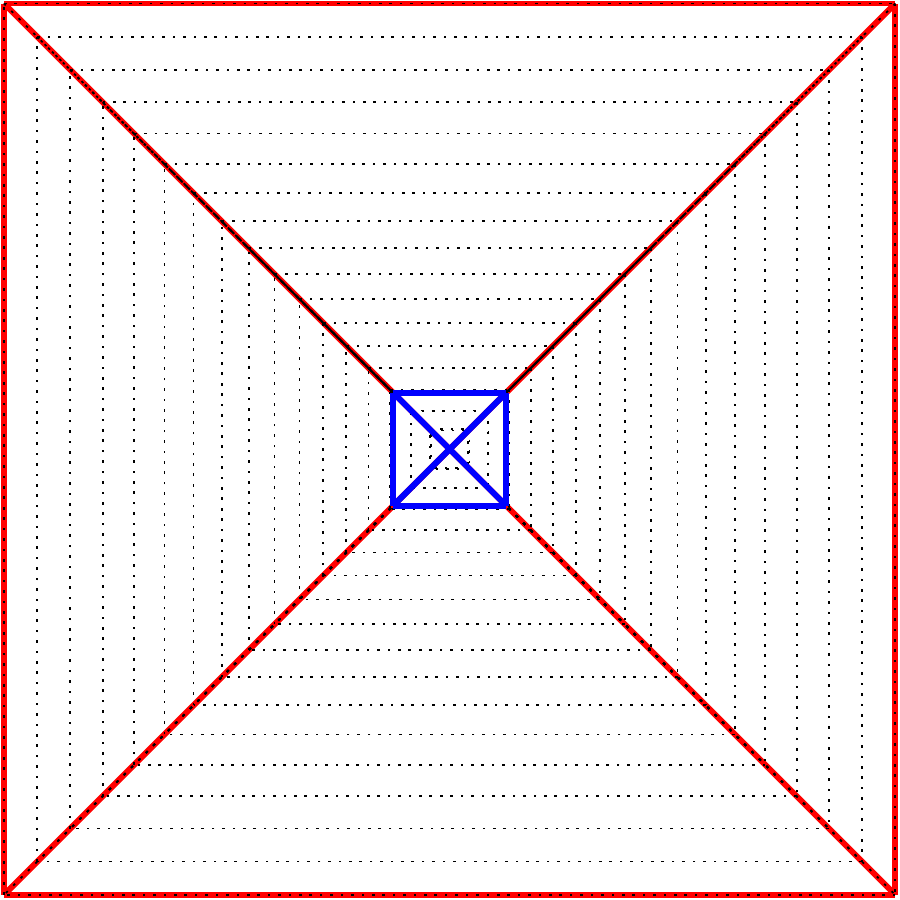}
			\caption{$\beta_1=0$, $\beta_3=1.0$}
			\label{subfig:scaling_beta1_0_beta2_10}
		\end{subfigure}
		\hfill
		\begin{subfigure}{0.3\textwidth}
			\centering
			\includegraphics[width=0.85\linewidth]{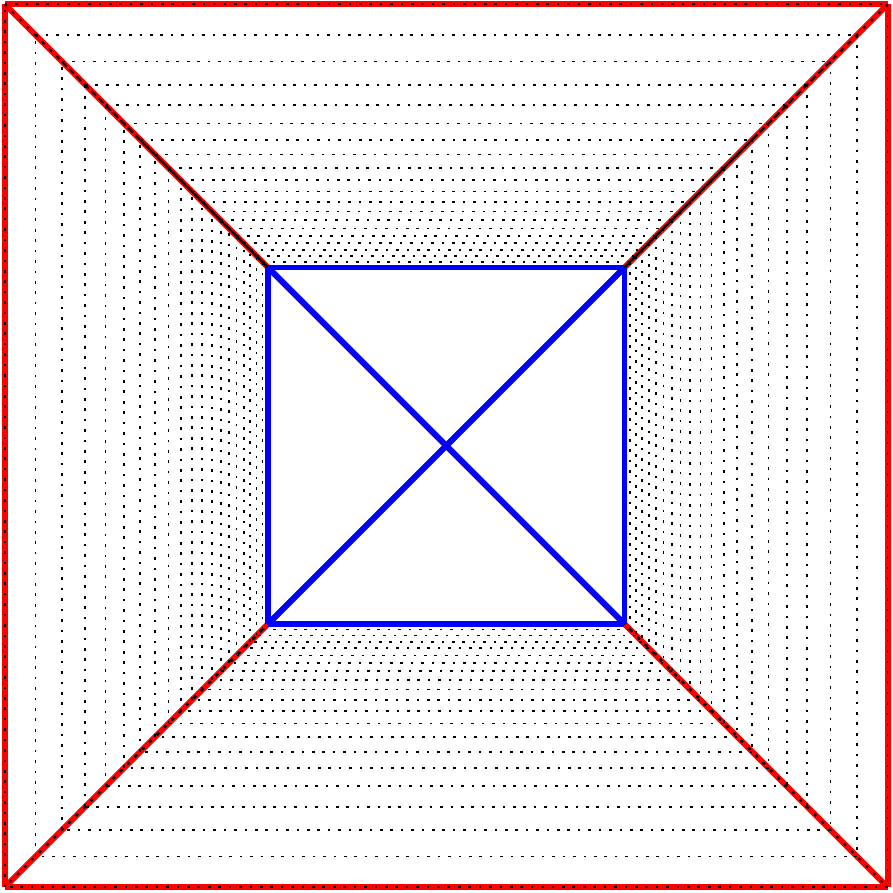}
			\caption{$\beta_1=0.5$, $\beta_3=1.0$}
			\label{subfig:scaling_beta1_5_beta2_10}
		\end{subfigure}
		\hfill
		\begin{subfigure}{0.3\textwidth}
			\centering
			\includegraphics[width=0.85\linewidth]{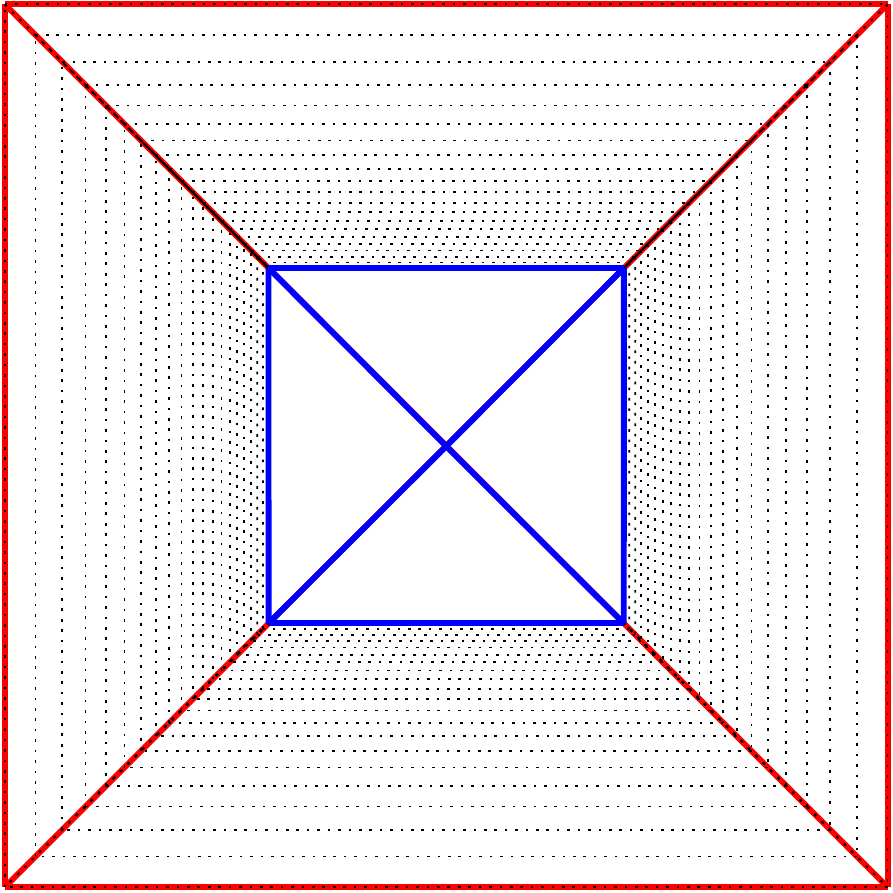}
			\caption{$\beta_1=1.0$, $\beta_3=1.0$}
			\label{subfig:scaling_beta1_10_beta2_10}
		\end{subfigure}
	\end{center}
	\caption{20~snapshots of geodesics for different values of $\beta_1$, $\beta_3$, starting from the same initial mesh (shown in red) and produced by the same initial tangent vector, which induces scaling; see \cref{subsection:elementary_transformations}. The final mesh is shown in blue.}
	\label{fig:scaling}
\end{figure}

\begin{figure}[htb]
	\begin{center}
		\begin{subfigure}{0.3\textwidth}
			\centering
			\includegraphics[width=0.85\linewidth]{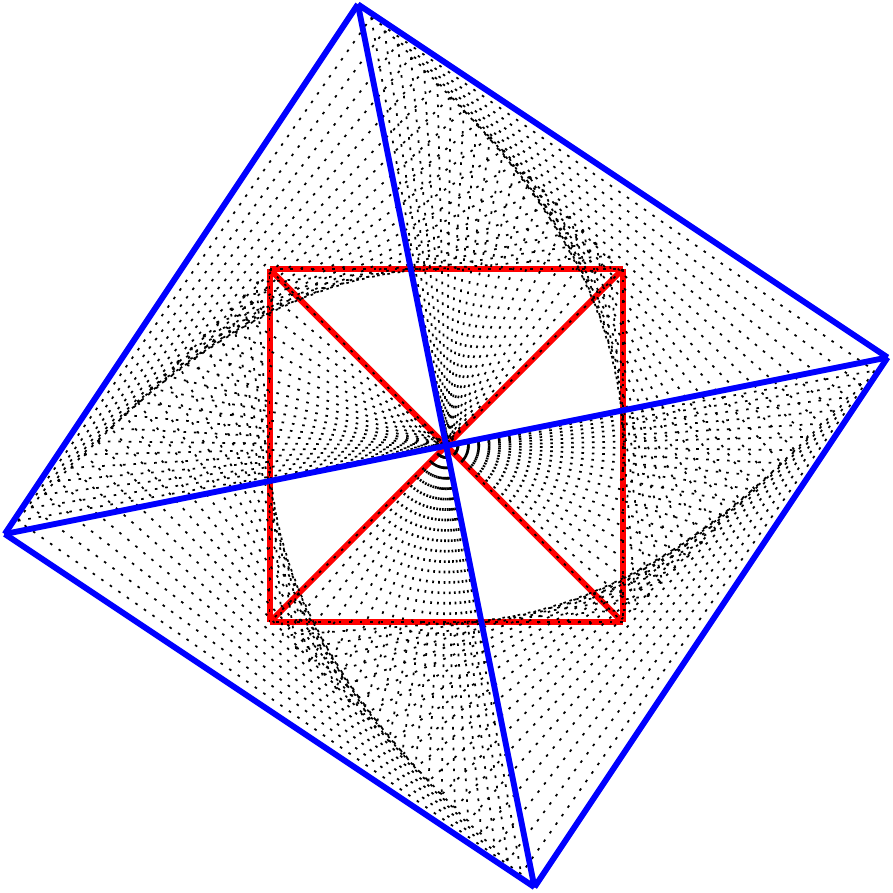}
			\caption{$\beta_1=0$, $\beta_3=0$}
			\label{subfig:rotation_beta1_0_beta2_0}
		\end{subfigure}
		\hfill
		\begin{subfigure}{0.3\textwidth}
			\centering
			\includegraphics[width=0.85\linewidth]{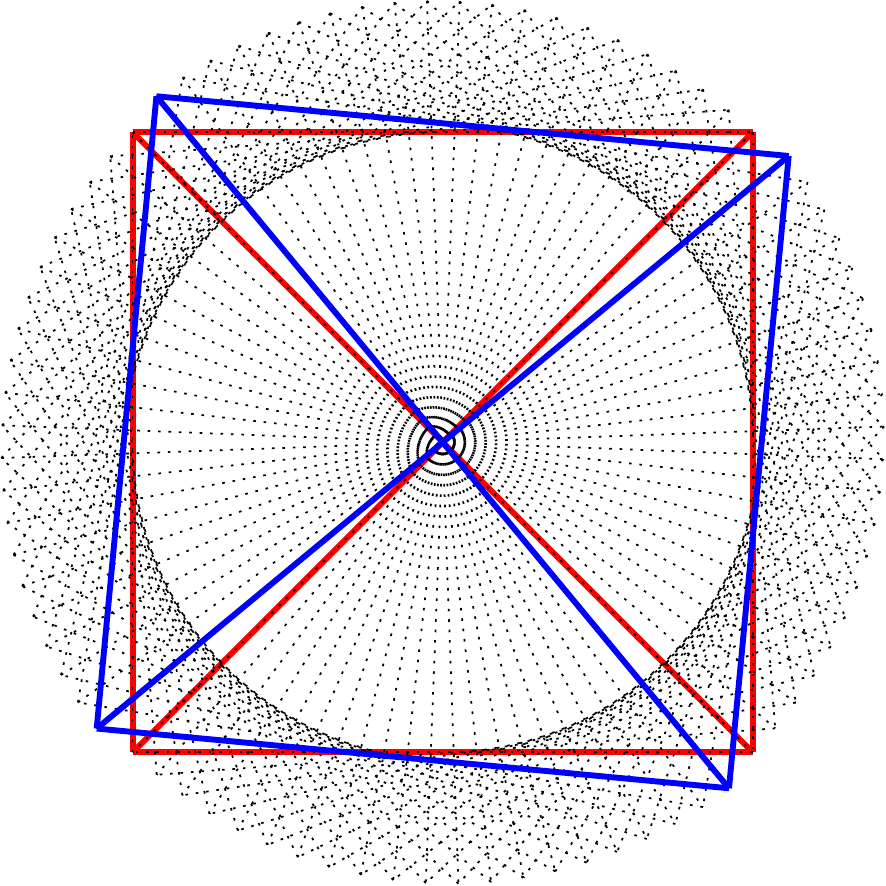}
			\caption{$\beta_1=0.5$, $\beta_3=0$}
			\label{subfig:rotation_beta1_5_beta2_0}
		\end{subfigure}
		\hfill
		\begin{subfigure}{0.3\textwidth}
			\centering
			\includegraphics[width=0.85\linewidth]{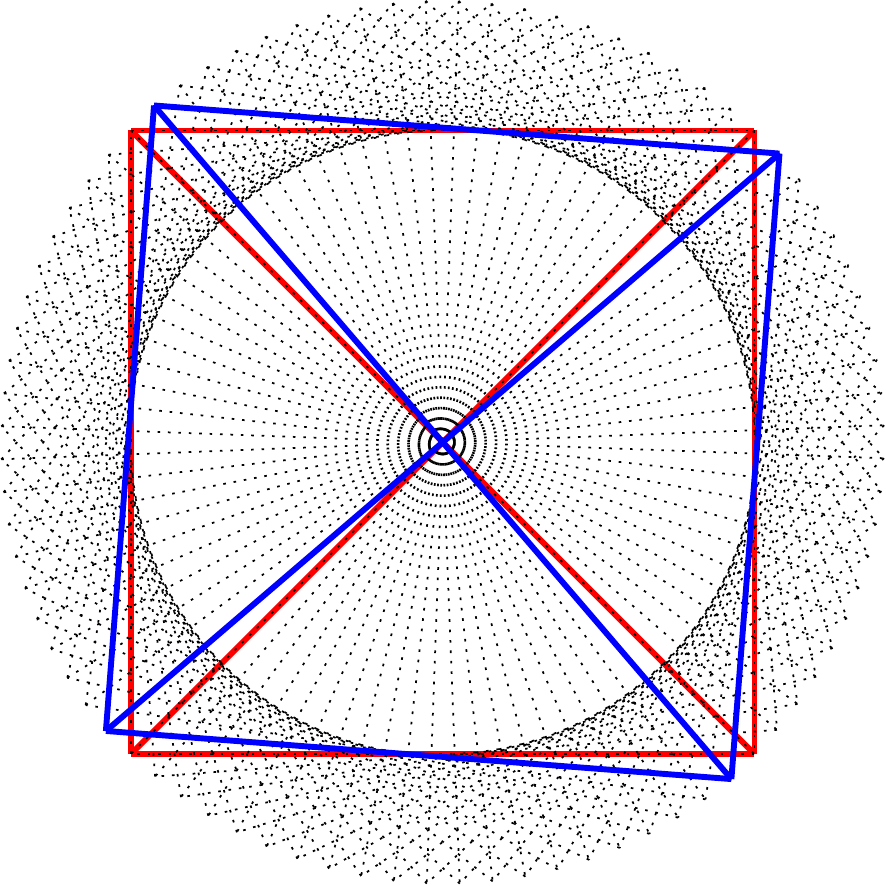}
			\caption{$\beta_1=1.0$, $\beta_3=0$}
			\label{subfig:rotation_beta1_10_beta2_0}
		\end{subfigure}
		\\
		\begin{subfigure}{0.3\textwidth}
			\centering
			\includegraphics[width=0.85\linewidth]{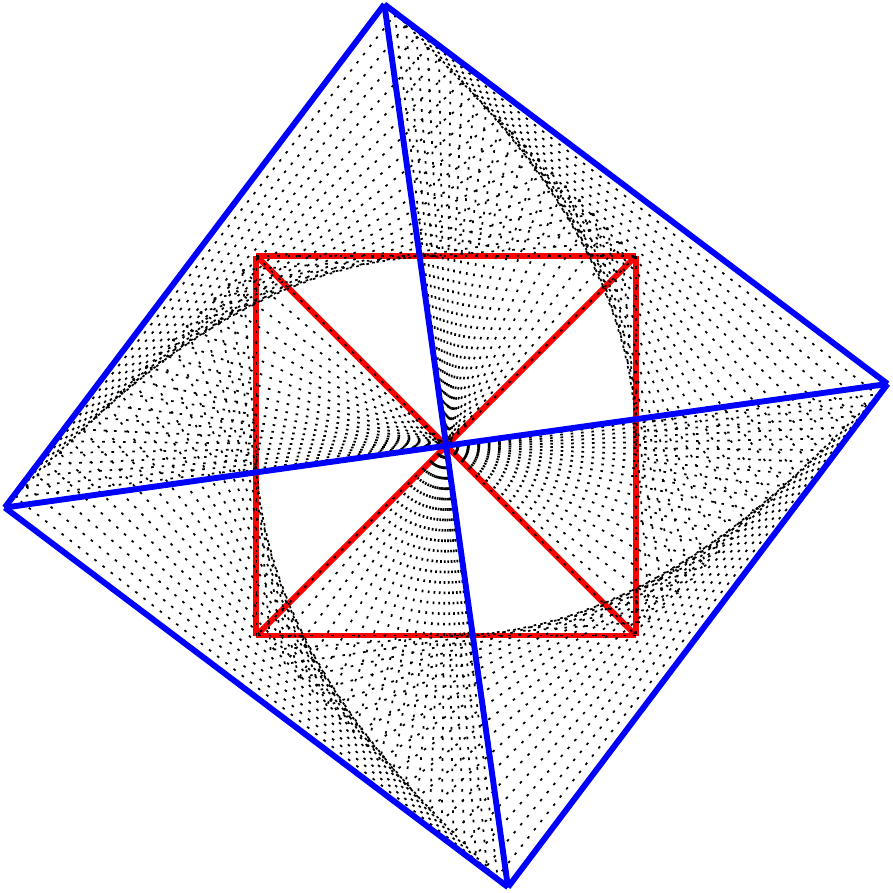}
			\caption{$\beta_1=0$, $\beta_3=0.5$}
			\label{subfig:rotation_beta1_0_beta2_5}
		\end{subfigure}
		\hfill
		\begin{subfigure}{0.3\textwidth}
			\centering
			\includegraphics[width=0.85\linewidth]{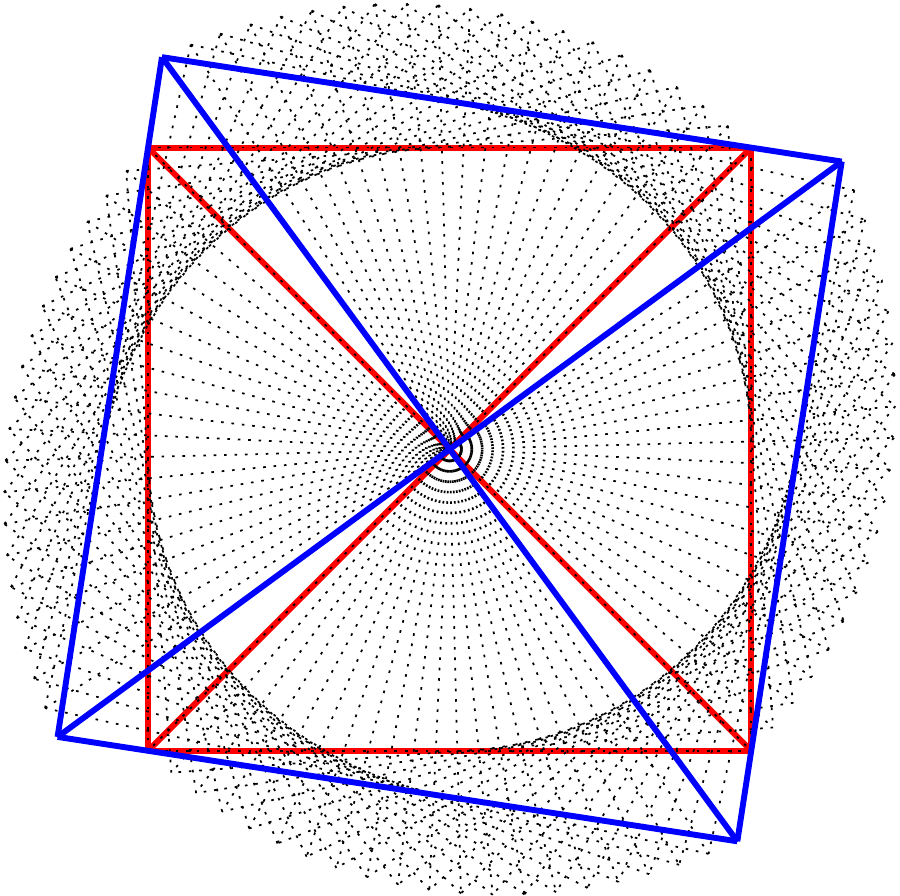}
			\caption{$\beta_1=0.5$, $\beta_3=0.5$}
			\label{subfig:rotation_beta1_5_beta2_5}
		\end{subfigure}
		\hfill
		\begin{subfigure}{0.3\textwidth}
			\centering
			\includegraphics[width=0.85\linewidth]{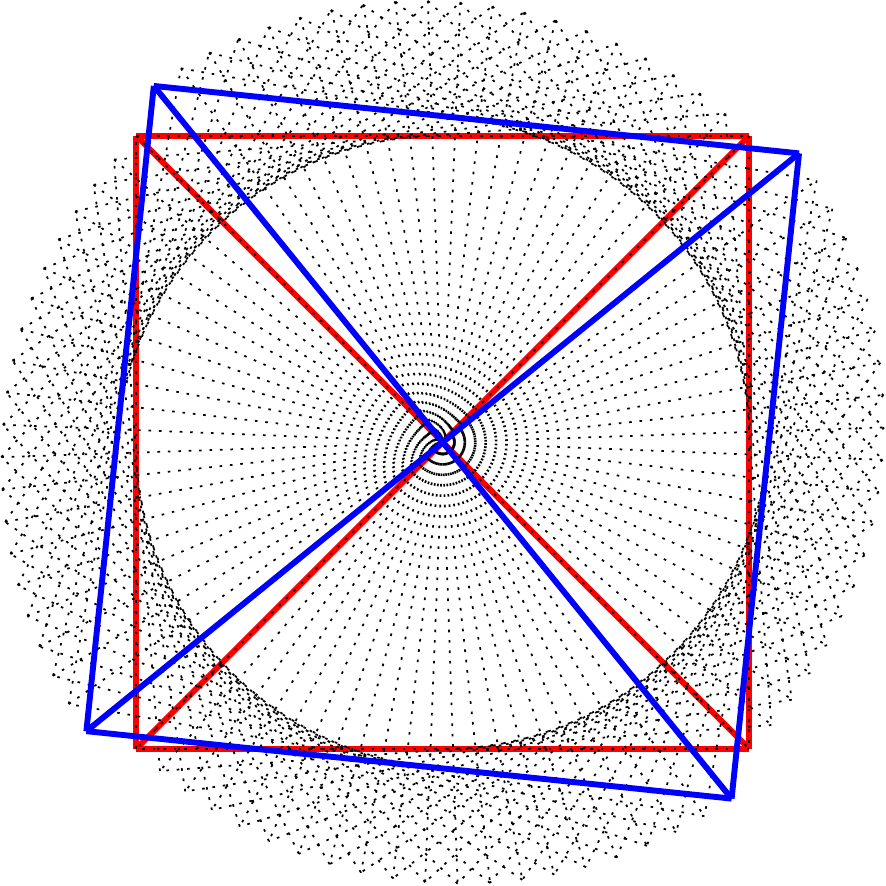}
			\caption{$\beta_1=1.0$, $\beta_3=0.5$}
			\label{subfig:rotation_beta1_10_beta2_5}
		\end{subfigure}
		\\
		\begin{subfigure}{0.3\textwidth}
			\centering
			\includegraphics[width=0.85\linewidth]{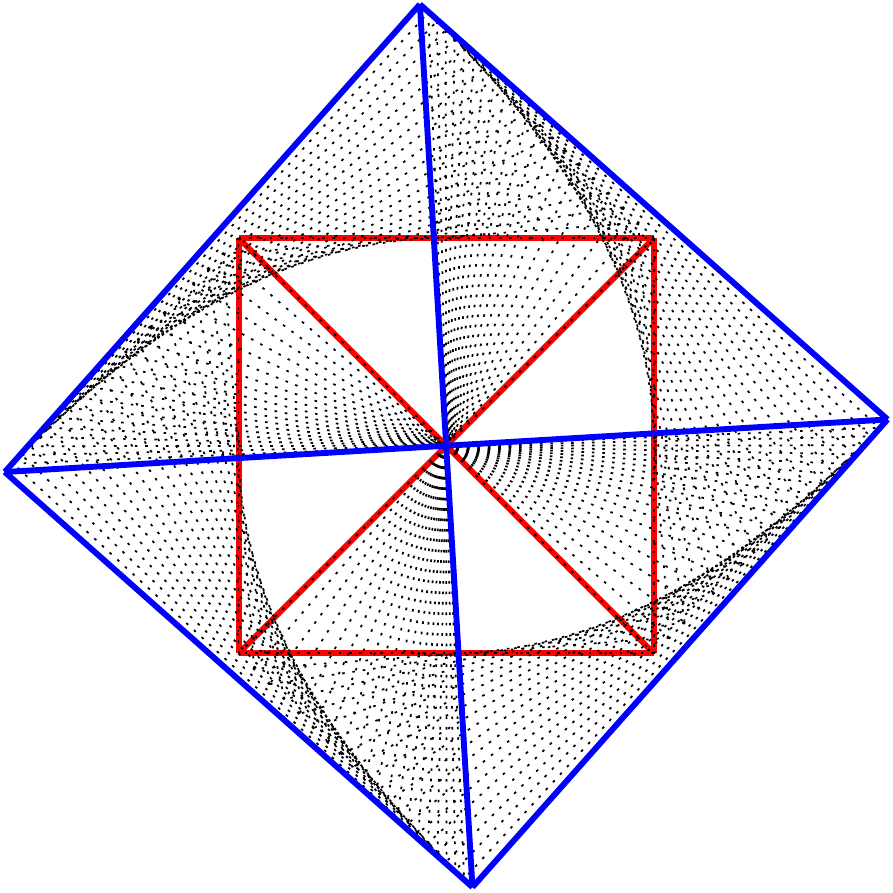}
			\caption{$\beta_1=0$, $\beta_3=1.0$}
			\label{subfig:rotation_beta1_0_beta2_10}
		\end{subfigure}
		\hfill
		\begin{subfigure}{0.3\textwidth}
			\centering
			\includegraphics[width=0.85\linewidth]{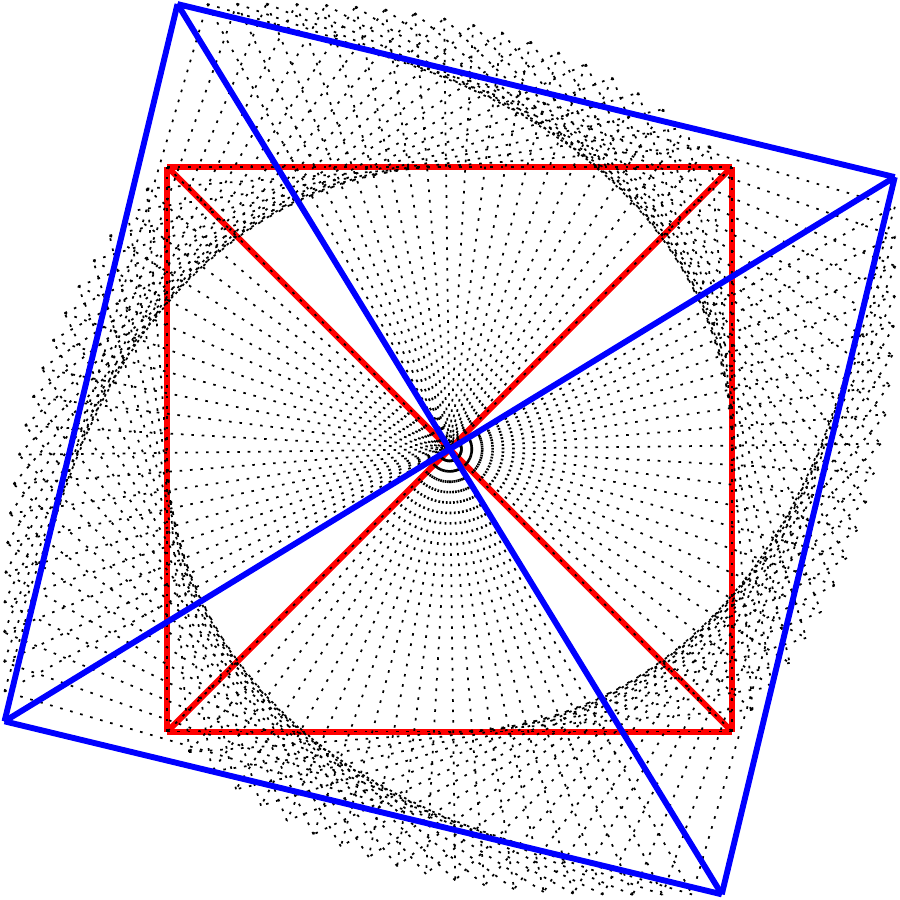}
			\caption{$\beta_1=0.5$, $\beta_3=1.0$}
			\label{subfig:rotation_beta1_5_beta2_10}
		\end{subfigure}
		\hfill
		\begin{subfigure}{0.3\textwidth}
			\centering
			\includegraphics[width=0.85\linewidth]{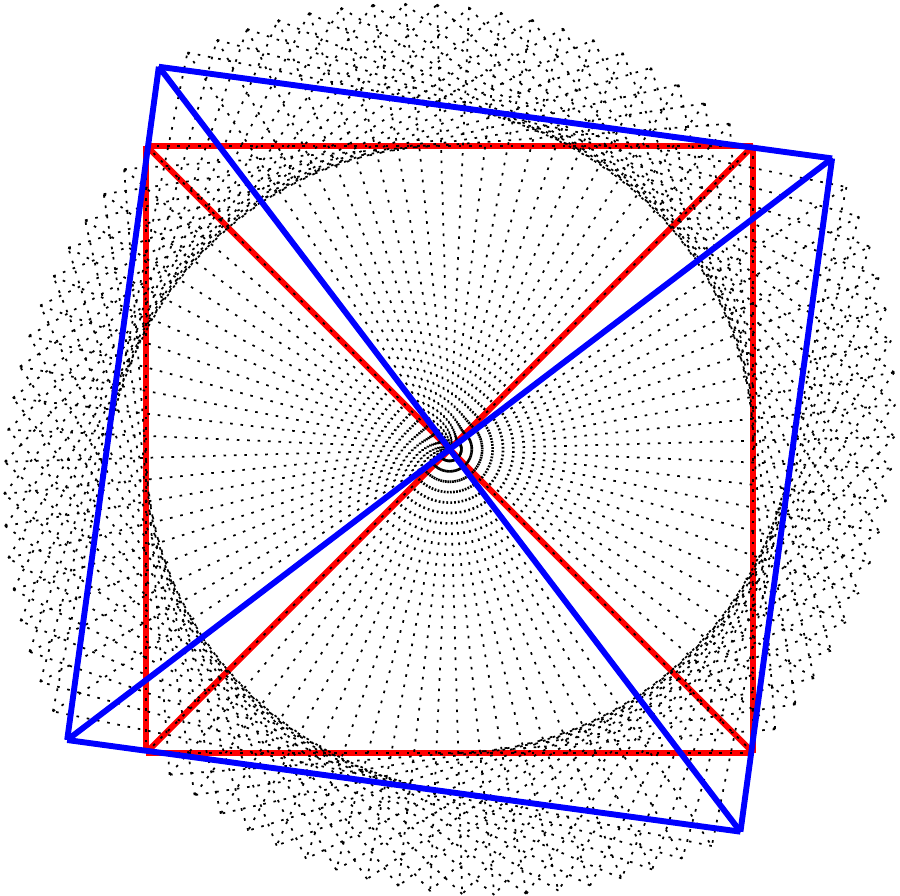}
			\caption{$\beta_1=1.0$, $\beta_3=1.0$}
			\label{subfig:rotation_beta1_10_beta2_10}
		\end{subfigure}
	\end{center}
	\caption{20~snapshots of geodesics for different values of $\beta_1$, $\beta_3$, starting from the same initial mesh (shown in red) and produced by the same initial tangent vector, which induces a rotation; see \cref{subsection:elementary_transformations}. The final mesh is shown in blue.}
	\label{fig:rotations}
\end{figure}

Each of the \Cref{fig:translations,fig:shearing,fig:scaling,fig:rotations} shows 20~snapshots of a geodesic on the interval $[0,3]$, obtained using values $\beta_1, \beta_3 \in \{0,0.5,1.0\}$.
The initial tangent vector is the same for all plots in a figure.
Although the initial tangent vectors are not shown, they can be easily recognized by the displacements they induce in the first time step.
Notice that the scaling of the plots within a figure may vary to make better use of the available space.
In each case, $\Qref$ is chosen to be the initial mesh.

In \Cref{fig:translations}, the initial vector represents a translation of the mesh.
As predicted by \cref{proposition:complete_metric_agrees_with_euclidean_for_translations}, when $\beta_3 = 0$ holds, geodesics \wrt \eqref{eq:complete_metric_for_planar_triangular_meshes_without_exterior_term} coincide with Euclidan geodesics, \ie, they remain translations, as can be seen from the first row in \Cref{fig:translations}.
In case $\beta_1 = 0$, the mesh is merely translated as well, albeit with a speed (in the Euclidean sense) depending on $\beta_3$; see the first column of \Cref{fig:translations}.

In \Cref{fig:shearing} we consider an initial tangent vector which induces a shearing motion.
In this case, there is a pronounced difference between the mesh evolution along the Euclidean geodesic ($\beta_1 = \beta_3 = 0$; see \Cref{subfig:shearing_beta1_0_beta2_0}) and those with $\beta_1 > 0$.
In the latter case, we can clearly see how the term involving the heights in \eqref{eq:f_mu_without_exterior_term} counteracts the impending mesh degeneracy observed along the Euclidean geodesic and helps to maintain a favorable cell aspect ratio.

\Cref{fig:scaling} shows the mesh deformation when the initial tangent vector induces a scaling of the mesh, \ie, the tangent vectors at the four corner vertices are pointing inwards.
In the Euclidean case ($\beta_1 = \beta_3 = 0$), this quickly leads to a non-admissible (flipped) mesh shown in \Cref{subfig:scaling_beta1_0_beta2_0} (top left).
The same is true for the other experiments with $\beta_1 = 0$ (first column).
However, with positive values of $\beta_1$ and $\beta_3$ we see that the completeness of the metric prevents the mesh from shrinking too much and it remains admissible for all times.

In \Cref{fig:rotations} we show some geodesics when the initial tangent vector induces a rotation.
Here the Euclidean geodesic ($\beta_1 = \beta_3 = 0$) does not cause the mesh to become degenerate.
Naturally, the Euclidean geodesic resembles a rotation only for small times since all vertices move along straight lines.
Interestingly, however, in the case $\beta_1 > 0$ and $\beta_3 = 0$ the geodesics appear to produce exact rotations.

The number of time steps used to integrate the geodesics in the inverval~$[0,3]$ using the Störmer--Verlet scheme (\Cref{algorithm:stoermer-verlet}) was $N = 10^3$ (\Cref{fig:translations}), $N = 10^4$ (\Cref{fig:shearing}), $N = 10^5$ (\Cref{fig:scaling}) and $N = 10^5$ (\Cref{fig:rotations}), respectively.
	The reason for using different numbers of time steps is that different initial velocities were observed to require different sizes of the time step~$\Delta t$ in order for the fixed-point iterations in \cref{line:implicit_equation_for_P,line:implicit_equation_for_Q} of the Störmer--Verlet scheme (\cref{algorithm:stoermer-verlet}) to converge.

\subsection{More Complex Initial Tangent Vector}
\label{subsection:comparison_complete_incomplete_metric}
\label{subsection:complete_vs_incomplete}

\begin{figure}
	\captionsetup[subfigure]{labelformat=empty}
	\begin{subfigure}{0.22\textwidth}
		\centering
		\includegraphics[width=\linewidth]{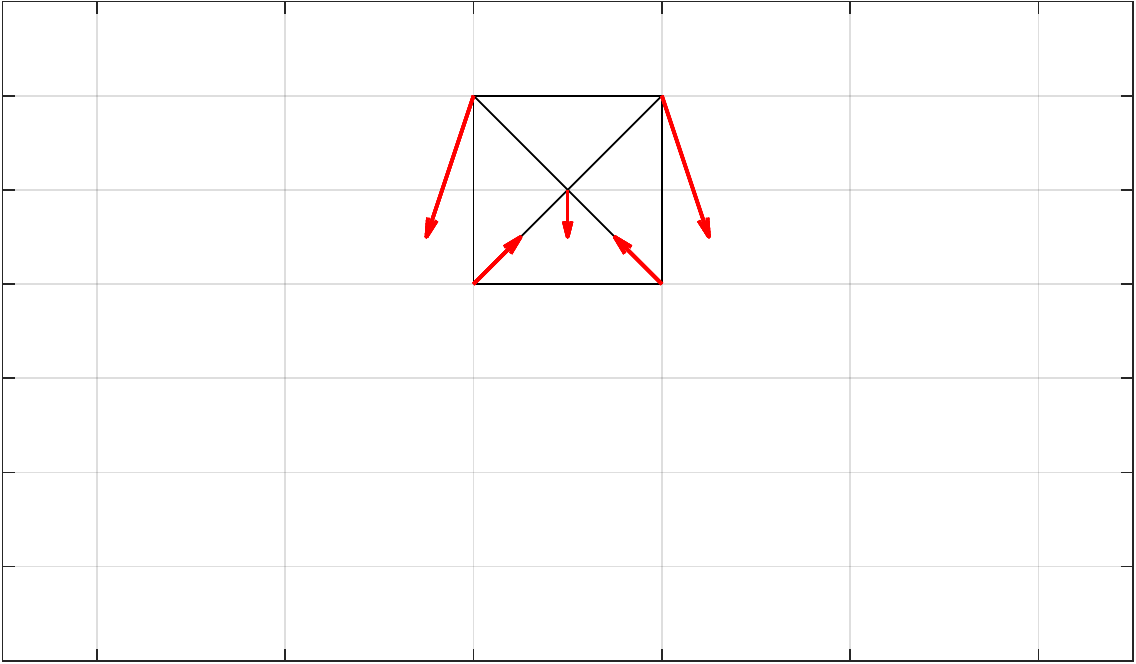}
		\caption{Complete, $t=0.0$}
	\end{subfigure}
	\hfill
	\begin{subfigure}{0.22\textwidth}
		\centering
		\includegraphics[width=\linewidth]{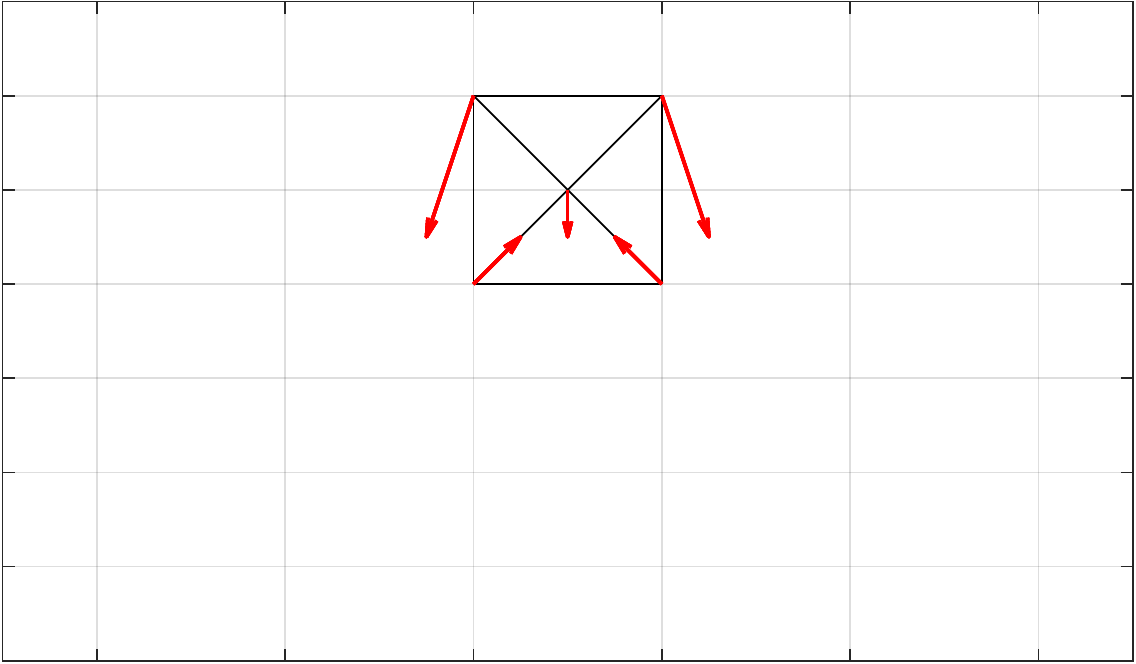}
		\caption{Euclidean, $t=0.0$}
	\end{subfigure}
	\hfill
	\begin{subfigure}{0.22\textwidth}
		\centering
		\includegraphics[width=\linewidth]{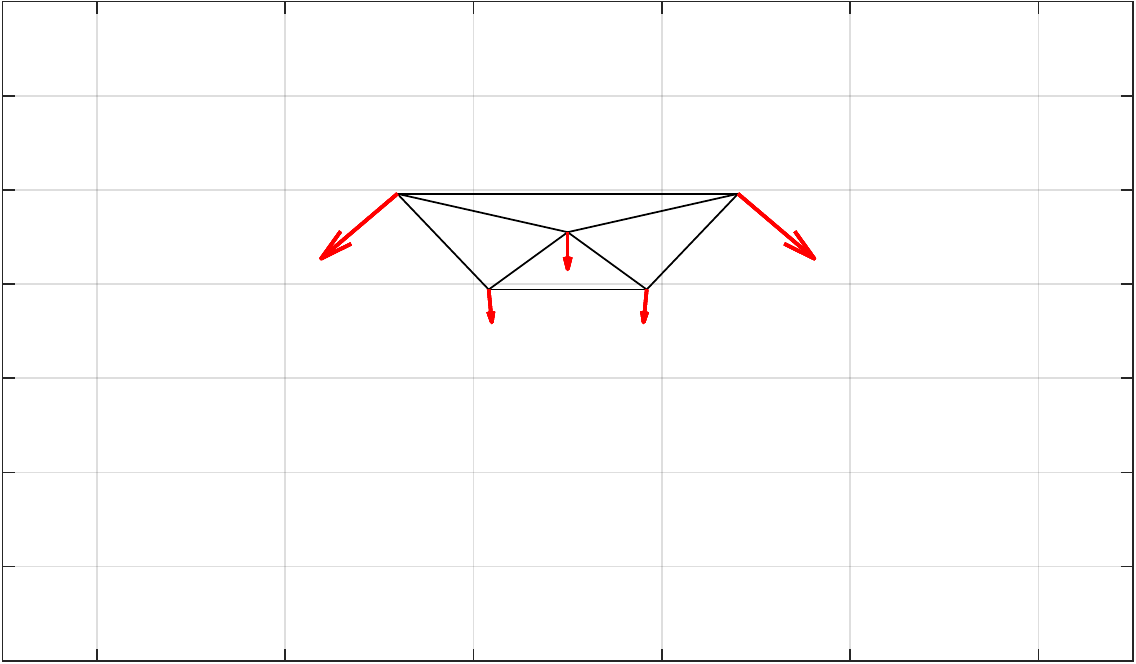}
		\caption{Complete, $t=1.07$}
	\end{subfigure}
	\hfill
	\begin{subfigure}{0.22\textwidth}
		\centering
		\includegraphics[width=\linewidth]{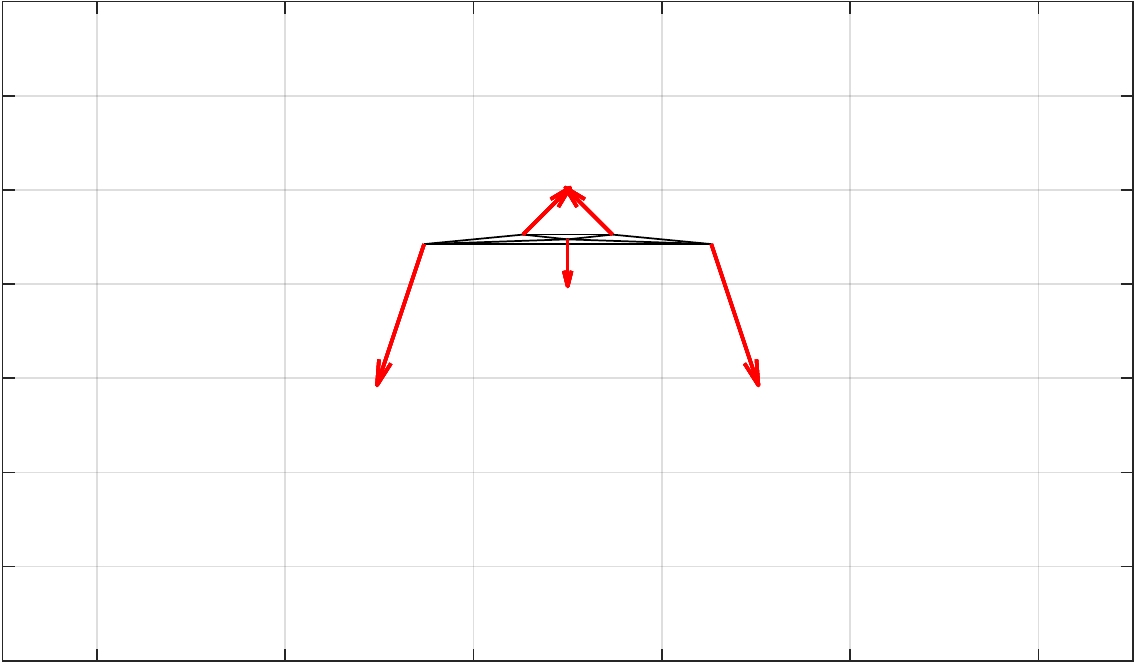}
		\caption{Euclidean, $t=1.07$}
	\end{subfigure}
	\\
	\begin{subfigure}{0.22\textwidth}
		\centering
		\includegraphics[width=\linewidth]{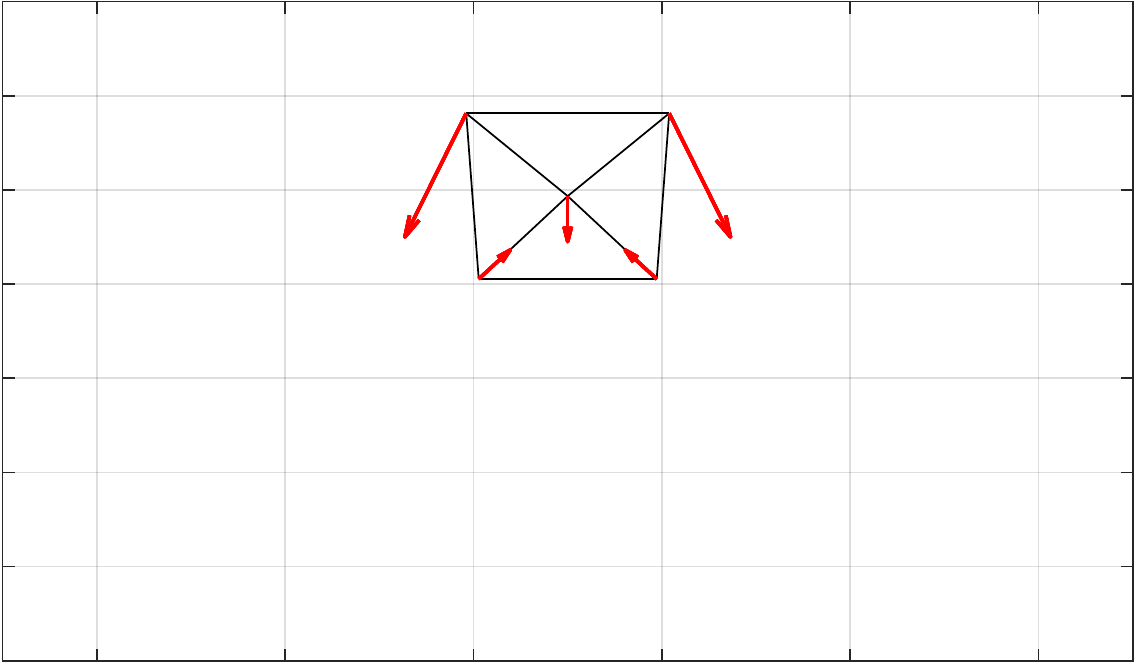}
		\caption{Complete, $t=0.13$}
	\end{subfigure}
	\hfill
	\begin{subfigure}{0.22\textwidth}
		\centering
		\includegraphics[width=\linewidth]{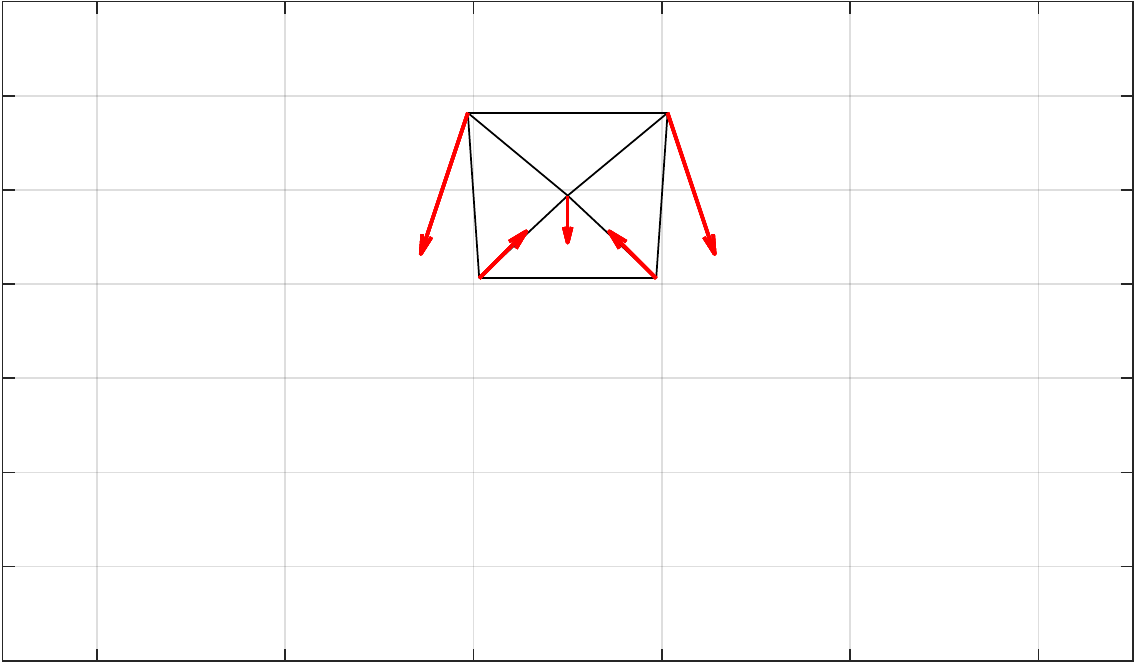}
		\caption{Euclidean, $t=0.13$}
	\end{subfigure}
	\hfill
	\begin{subfigure}{0.22\textwidth}
		\centering
		\includegraphics[width=\linewidth]{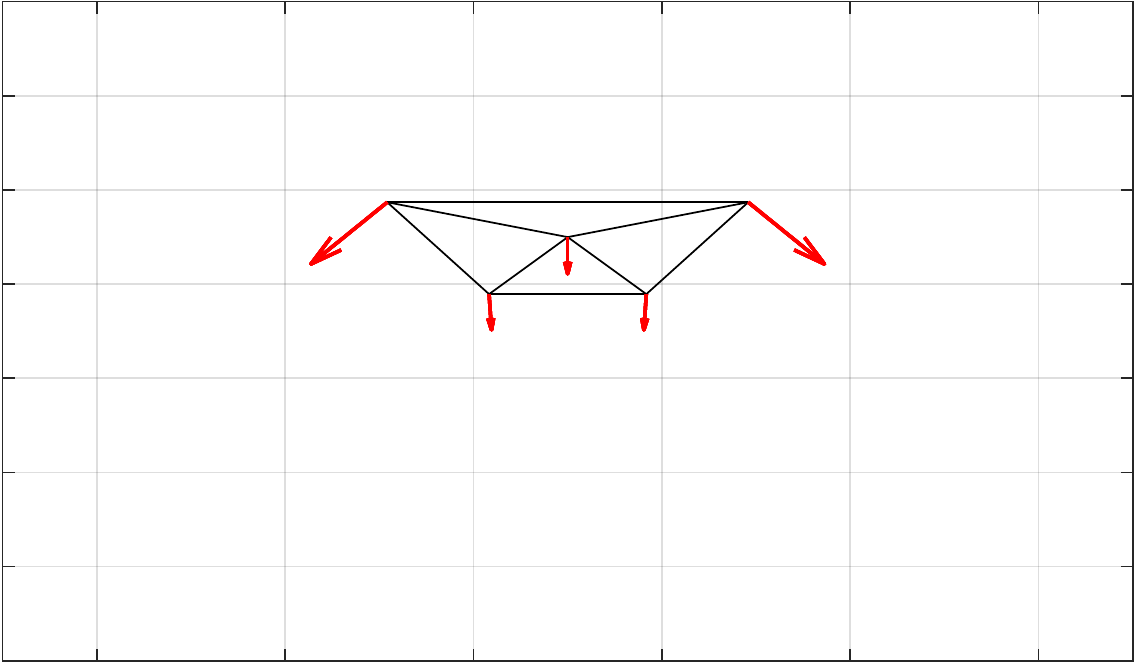}
		\caption{Complete, $t=1.2$}
	\end{subfigure}
	\hfill
	\begin{subfigure}{0.22\textwidth}
		\centering
		\includegraphics[width=\linewidth]{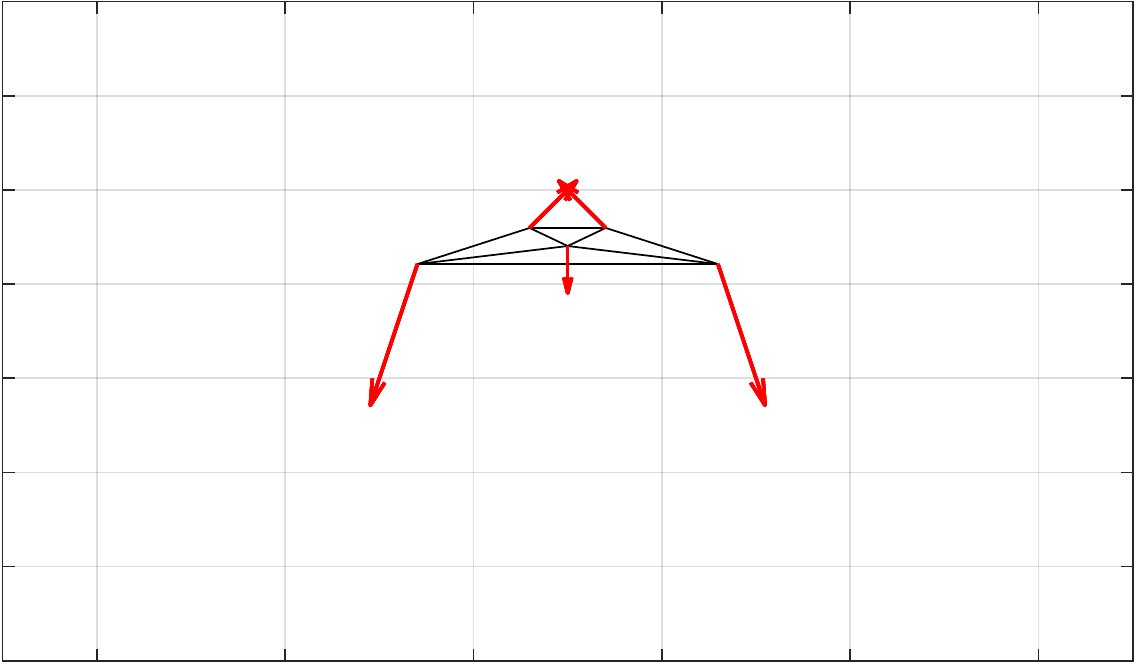}
		\caption{Euclidean, $t=1.2$}
	\end{subfigure}
	\\
	\begin{subfigure}{0.22\textwidth}
		\centering
		\includegraphics[width=\linewidth]{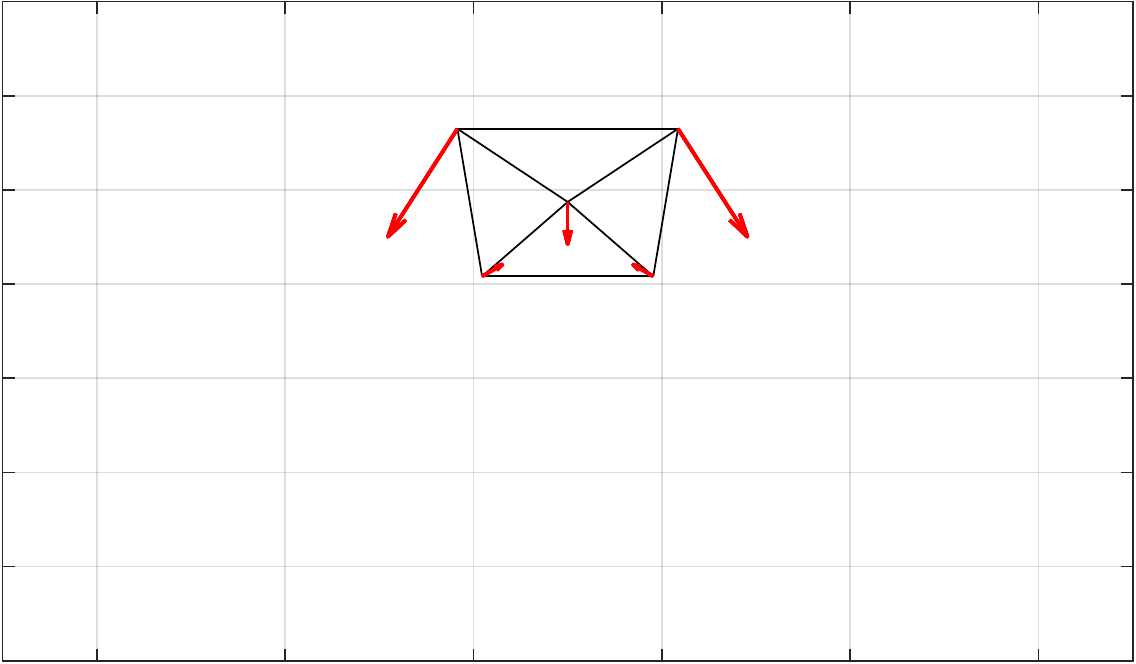}
		\caption{Complete, $t=0.27$}
	\end{subfigure}
	\hfill
	\begin{subfigure}{0.22\textwidth}
		\centering
		\includegraphics[width=\linewidth]{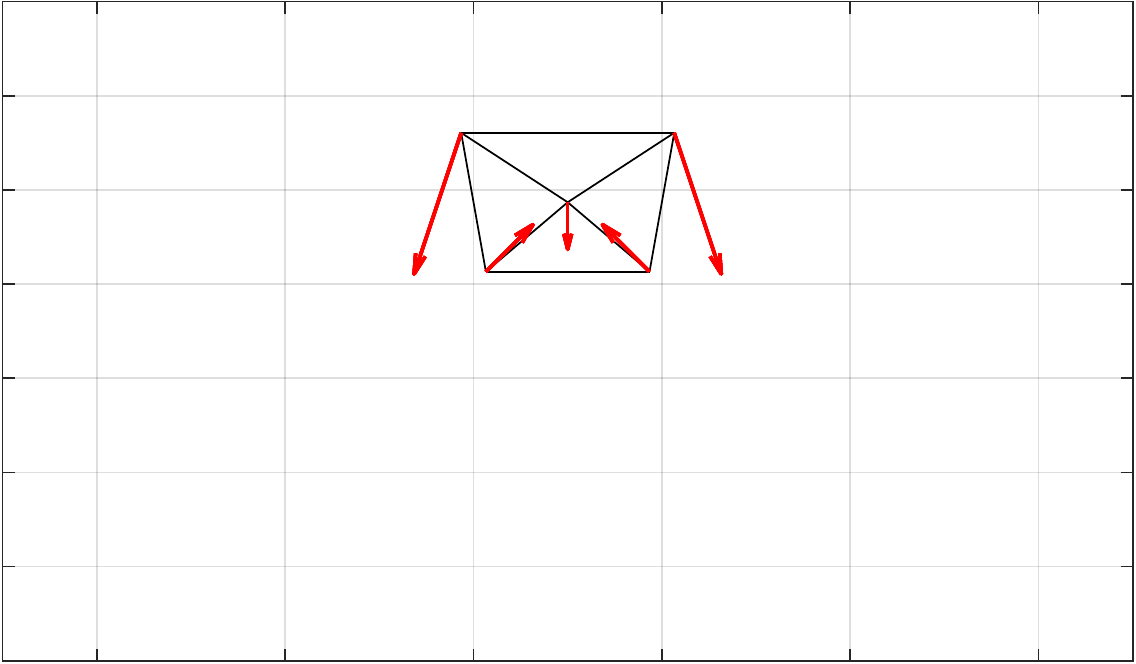}
		\caption{Euclidean, $t=0.27$}
	\end{subfigure}
	\hfill
	\begin{subfigure}{0.22\textwidth}
		\centering
		\includegraphics[width=\linewidth]{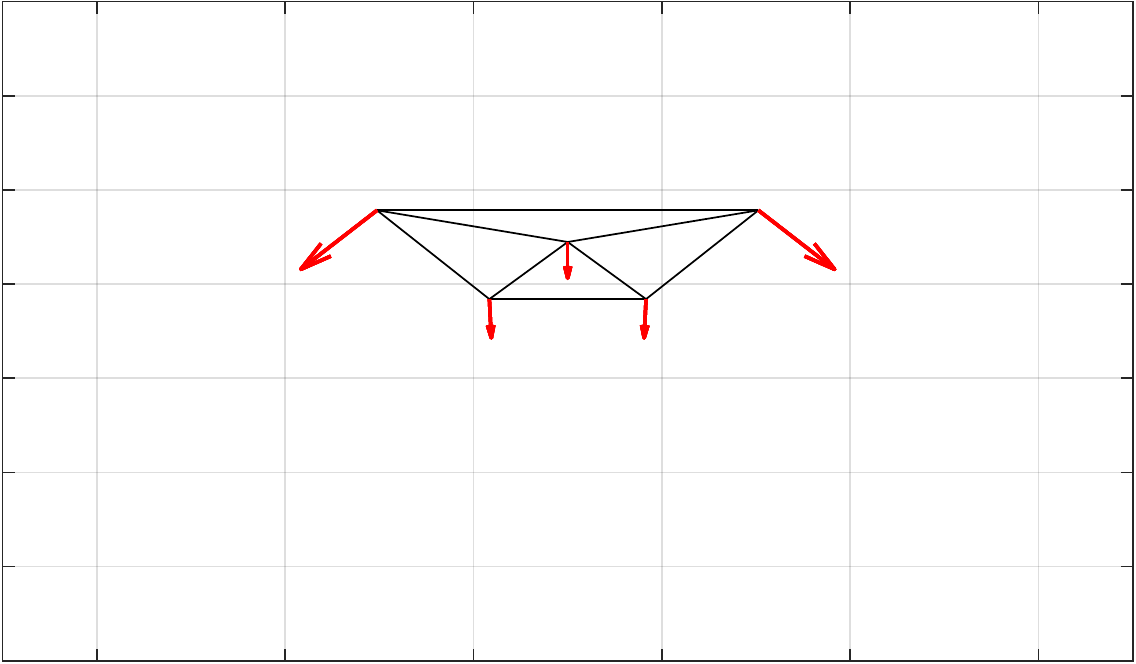}
		\caption{Complete, $t=1.33$}
	\end{subfigure}
	\hfill
	\begin{subfigure}{0.22\textwidth}
		\centering
		\includegraphics[width=\linewidth]{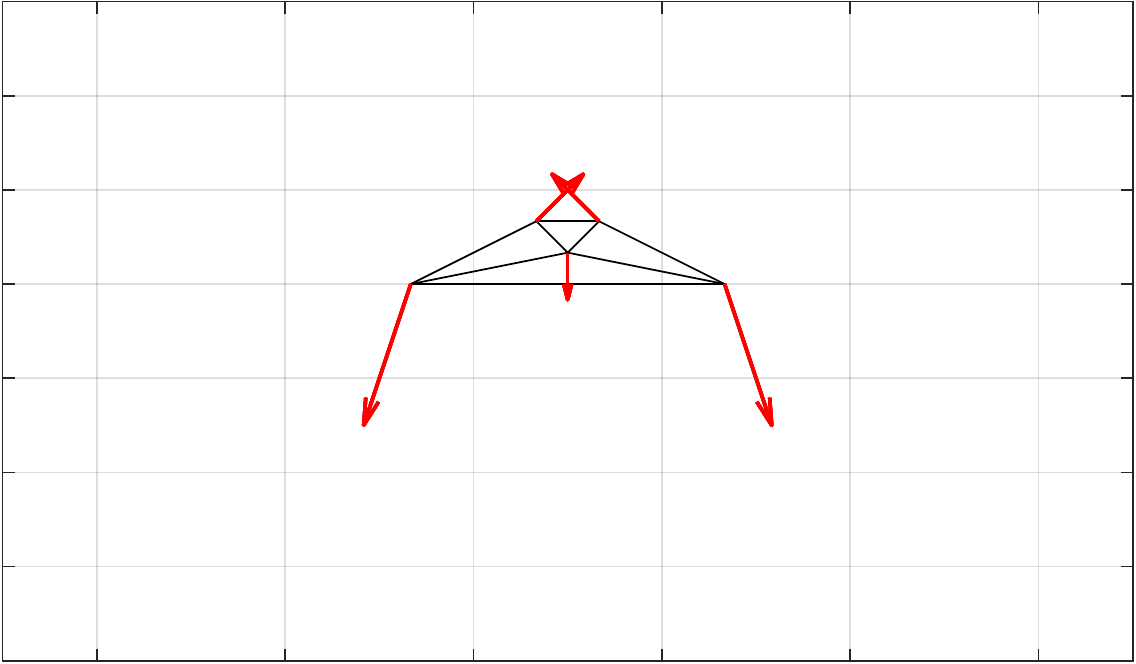}
		\caption{Euclidean, $t=1.33$}
	\end{subfigure}
	\\
	\begin{subfigure}{0.22\textwidth}
		\centering
		\includegraphics[width=\linewidth]{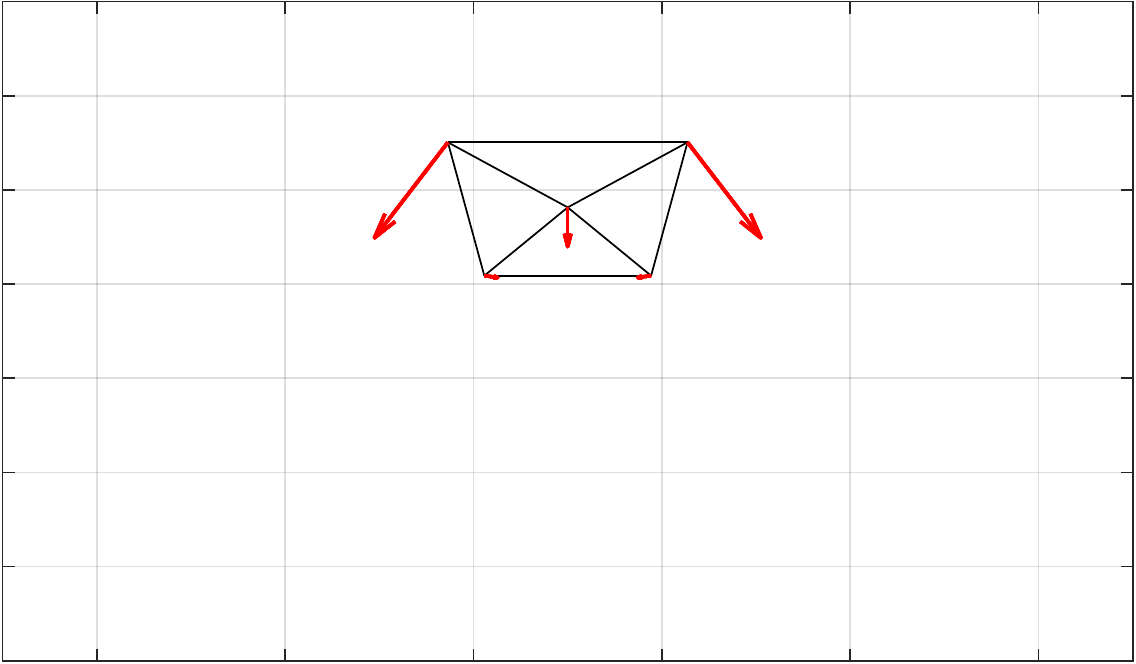}
		\caption{Complete, $t=0.40$}
	\end{subfigure}
	\hfill
	\begin{subfigure}{0.22\textwidth}
		\centering
		\includegraphics[width=\linewidth]{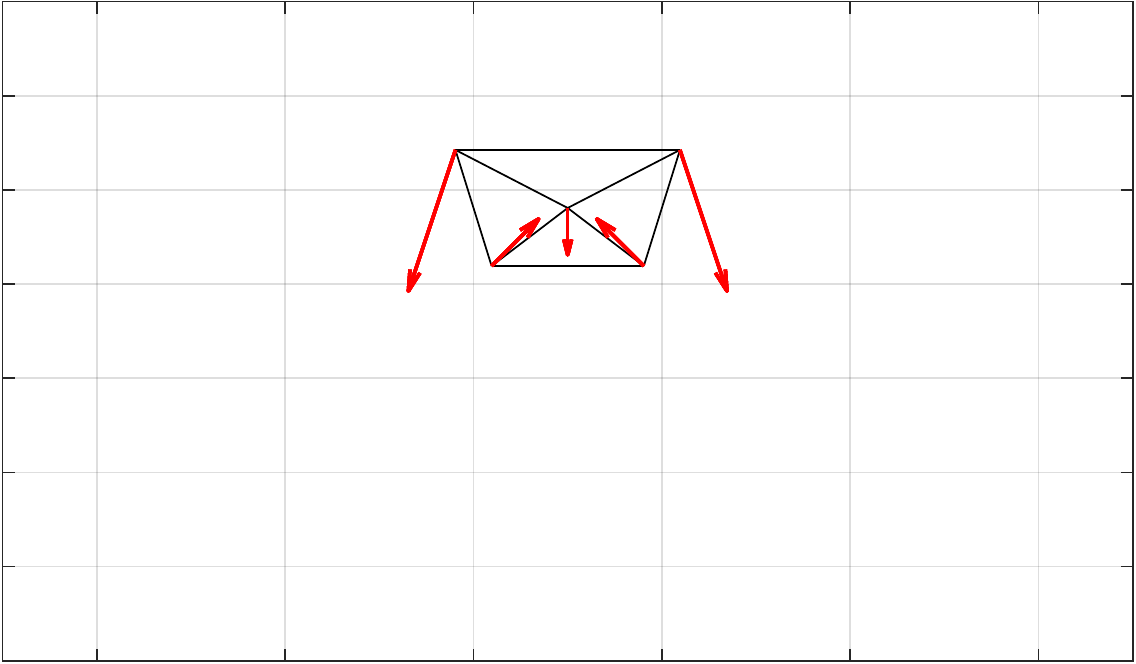}
		\caption{Euclidean, $t=0.40$}
	\end{subfigure}
	\hfill
	\begin{subfigure}{0.22\textwidth}
		\centering
		\includegraphics[width=\linewidth]{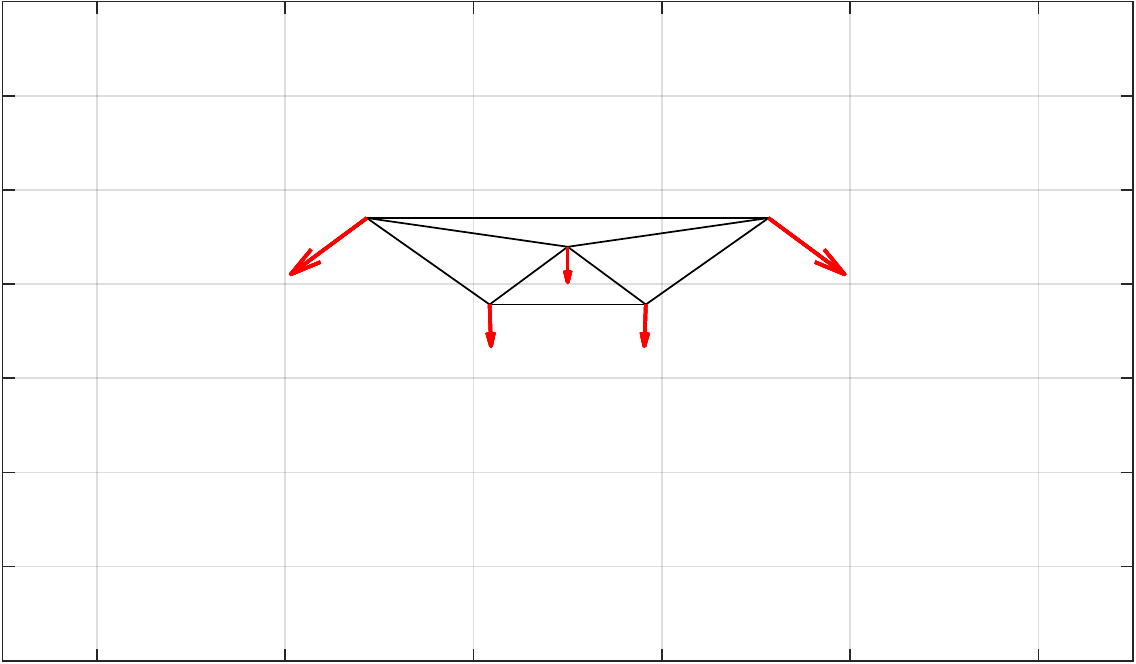}
		\caption{Complete, $t=1.47$}
	\end{subfigure}
	\hfill
	\begin{subfigure}{0.22\textwidth}
		\centering
		\includegraphics[width=\linewidth]{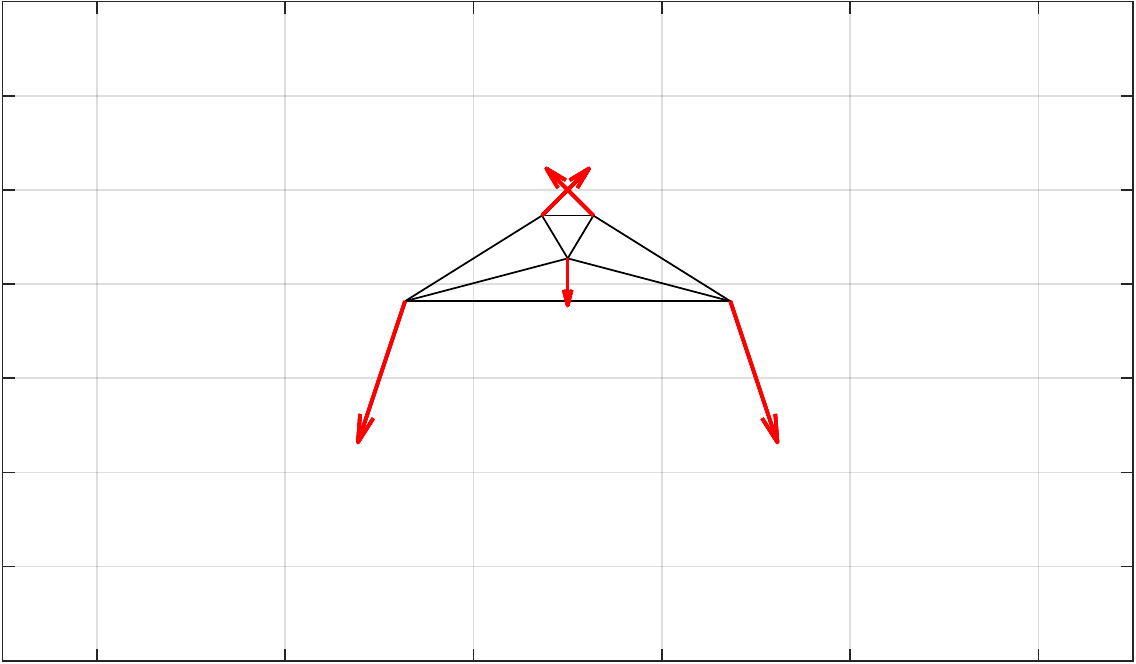}
		\caption{Euclidean, $t=1.47$}
	\end{subfigure}
	\\
	\begin{subfigure}{0.22\textwidth}
		\centering
		\includegraphics[width=\linewidth]{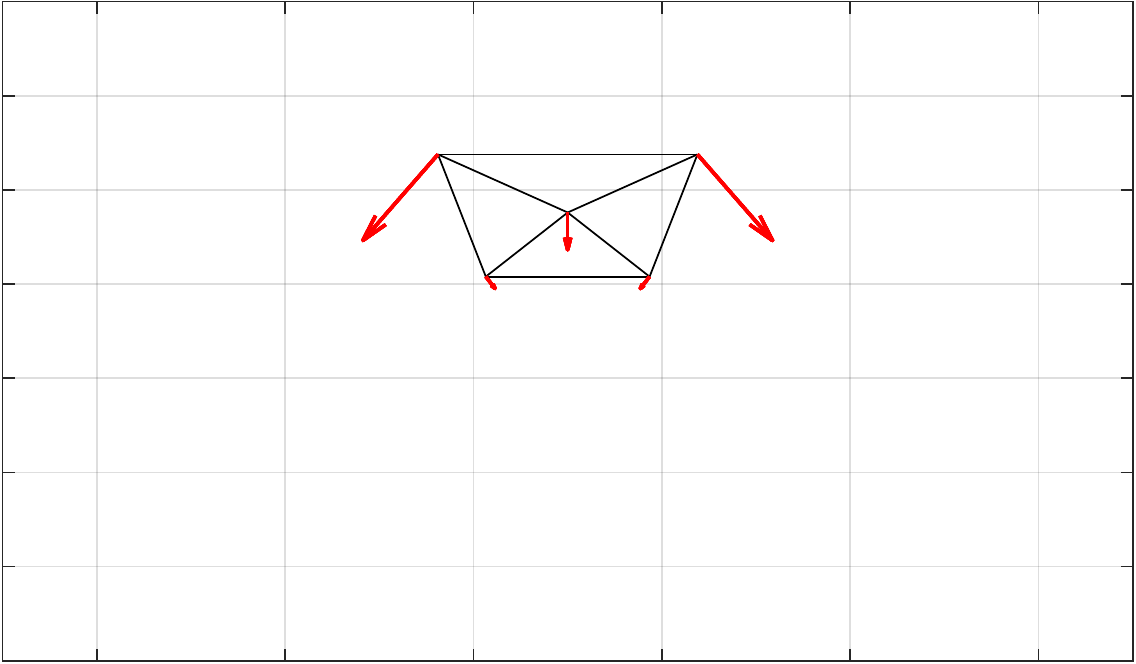}
		\caption{Complete, $t=0.53$}
	\end{subfigure}
	\hfill
	\begin{subfigure}{0.22\textwidth}
		\centering
		\includegraphics[width=\linewidth]{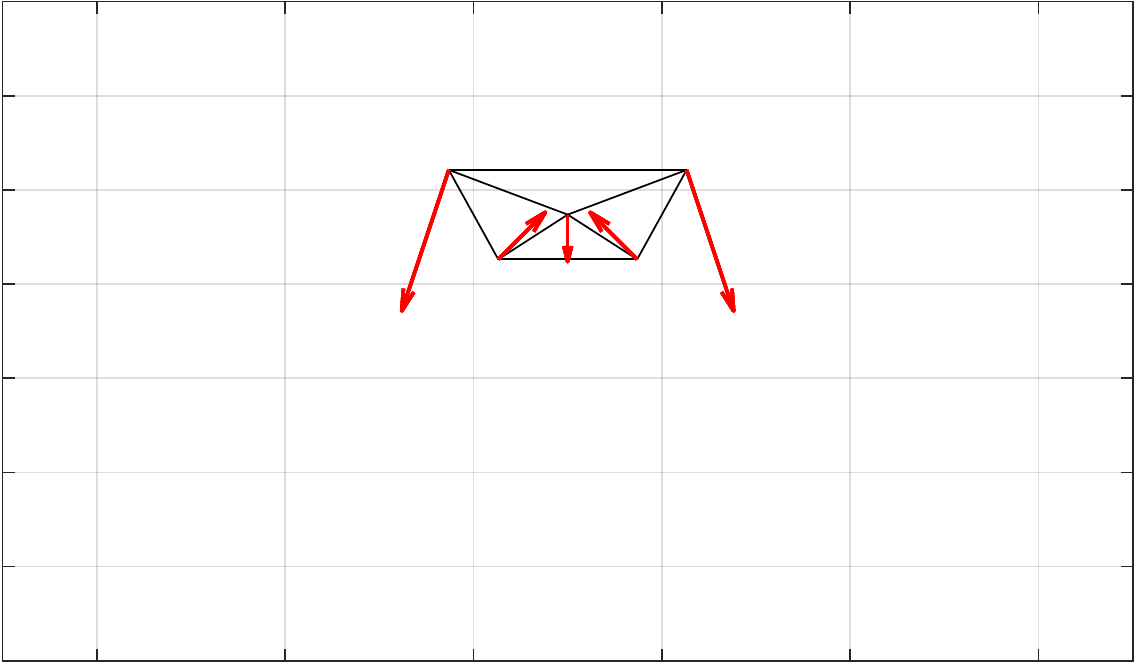}
		\caption{Euclidean, $t=0.53$}
	\end{subfigure}
	\hfill
	\begin{subfigure}{0.22\textwidth}
		\centering
		\includegraphics[width=\linewidth]{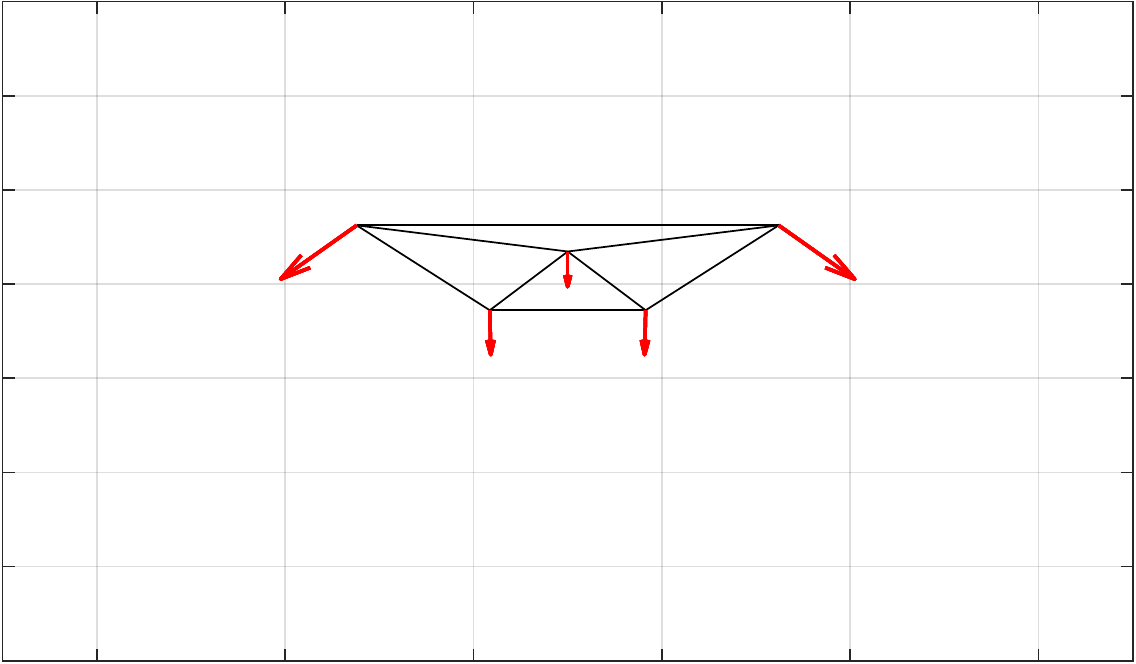}
		\caption{Complete, $t=1.60$}
	\end{subfigure}
	\hfill
	\begin{subfigure}{0.22\textwidth}
		\centering
		\includegraphics[width=\linewidth]{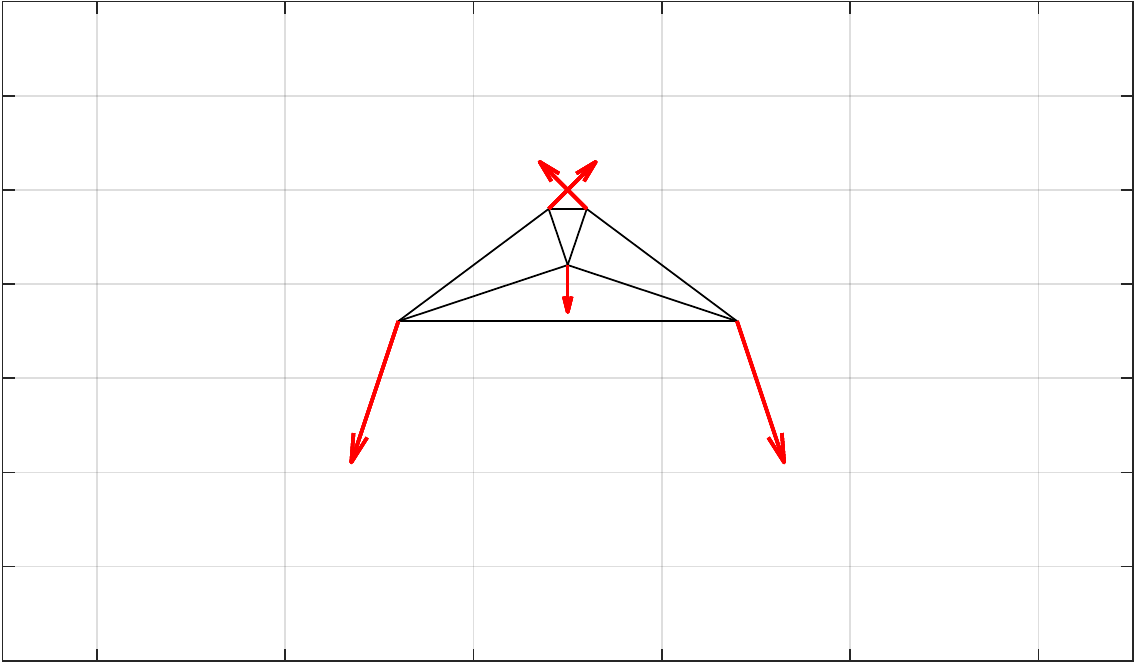}
		\caption{Euclidean, $t=1.60$}
	\end{subfigure}
	\hfill
	\begin{subfigure}{0.22\textwidth}
		\centering
		\includegraphics[width=\linewidth]{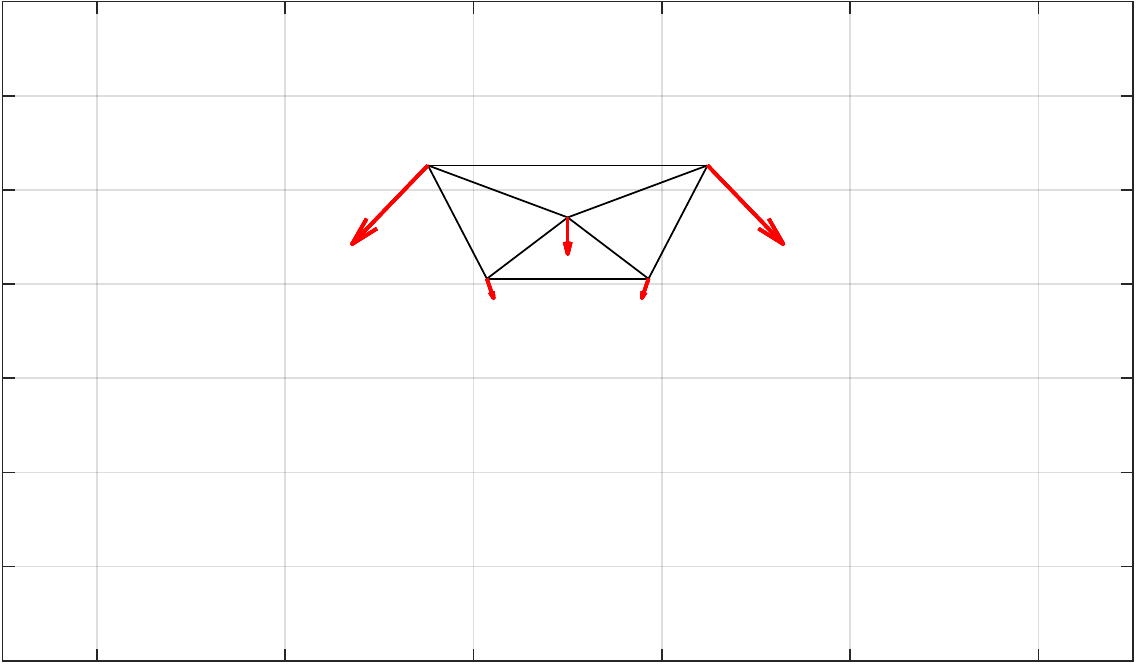}
		\caption{Complete, $t=0.67$}
	\end{subfigure}
	\hfill
	\begin{subfigure}{0.22\textwidth}
		\centering
		\includegraphics[width=\linewidth]{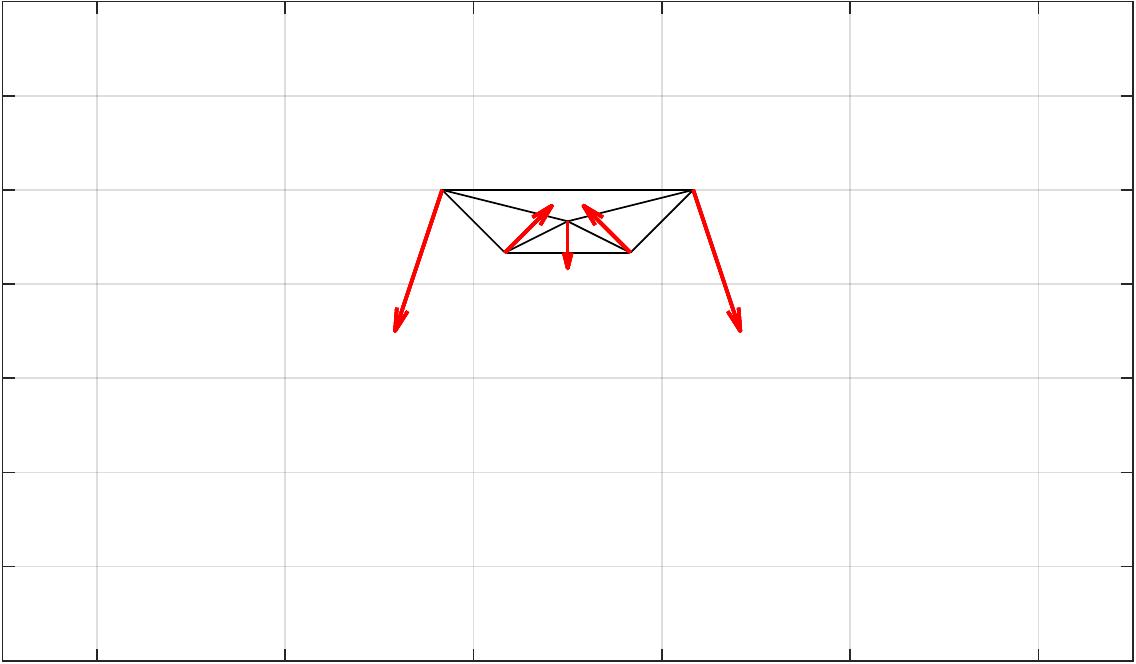}
		\caption{Euclidean, $t=0.67$}
	\end{subfigure}
	\hfill
	\begin{subfigure}{0.22\textwidth}
		\centering
		\includegraphics[width=\linewidth]{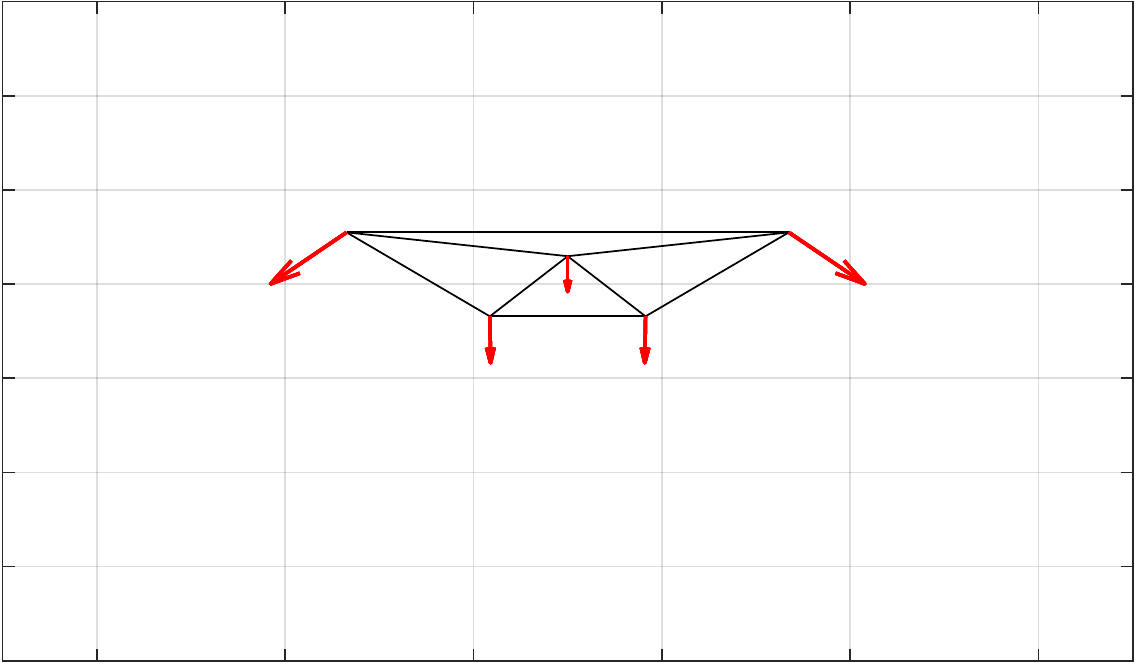}
		\caption{Complete, $t=1.73$}
	\end{subfigure}
	\hfill
	\begin{subfigure}{0.22\textwidth}
		\centering
		\includegraphics[width=\linewidth]{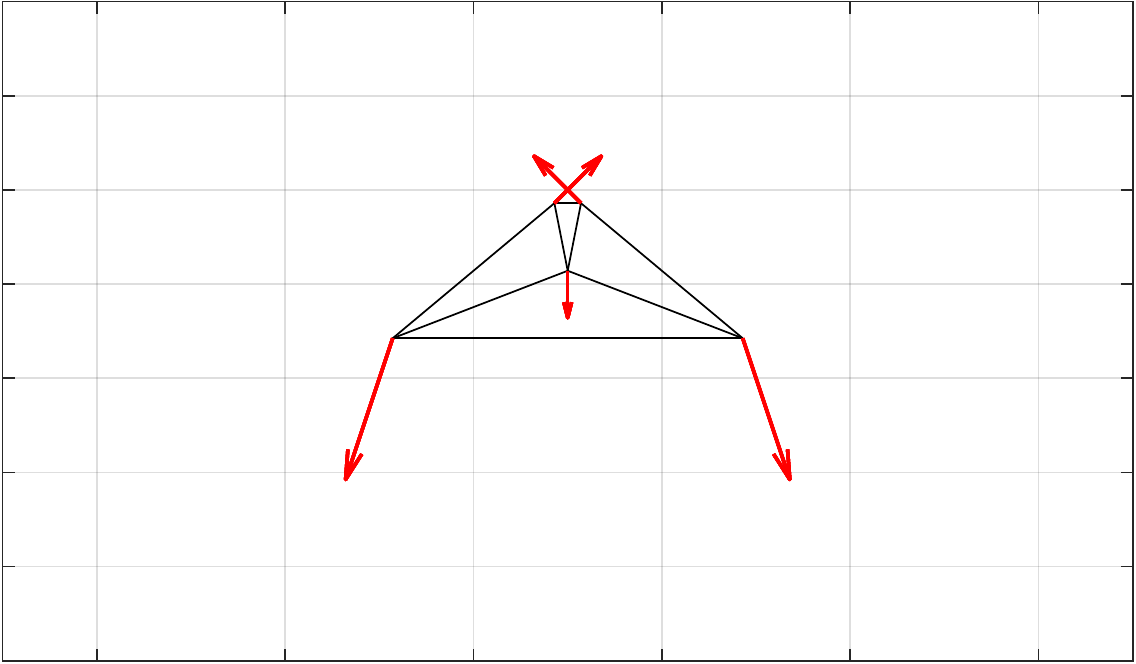}
		\caption{Euclidean, $t=1.73$}
	\end{subfigure}
	\hfill
	\begin{subfigure}{0.22\textwidth}
		\centering
		\includegraphics[width=\linewidth]{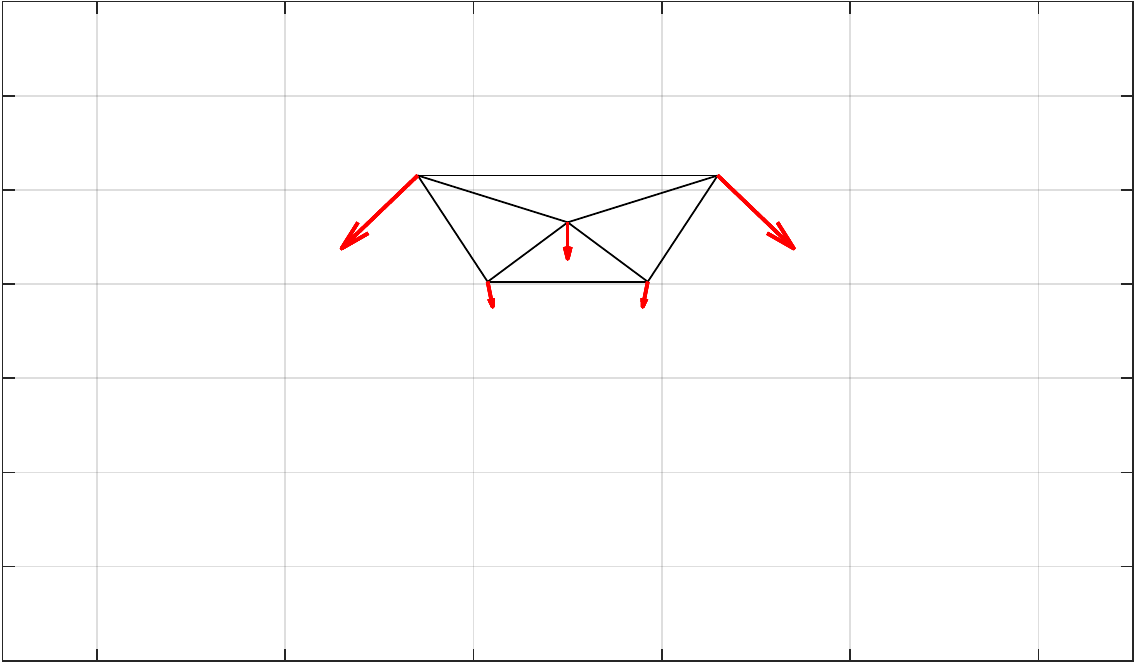}
		\caption{Complete, $t=0.80$}
	\end{subfigure}
	\hfill
	\begin{subfigure}{0.22\textwidth}
		\centering
		\includegraphics[width=\linewidth]{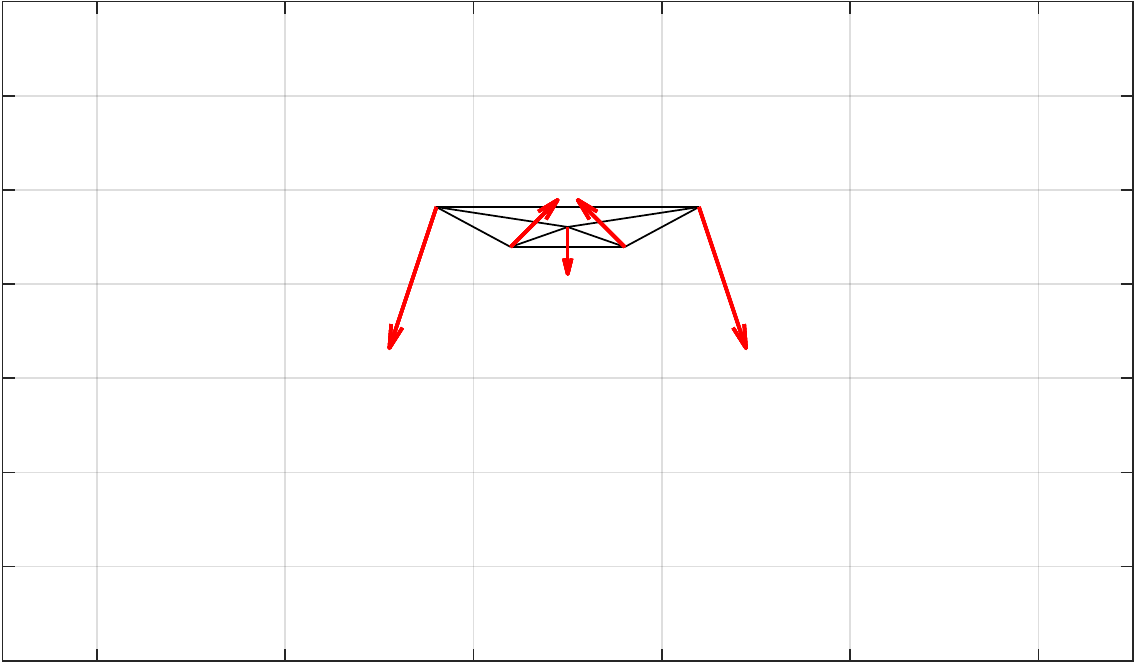}
		\caption{Euclidean, $t=0.80$}
	\end{subfigure}
	\hfill
	\begin{subfigure}{0.22\textwidth}
		\centering
		\includegraphics[width=\linewidth]{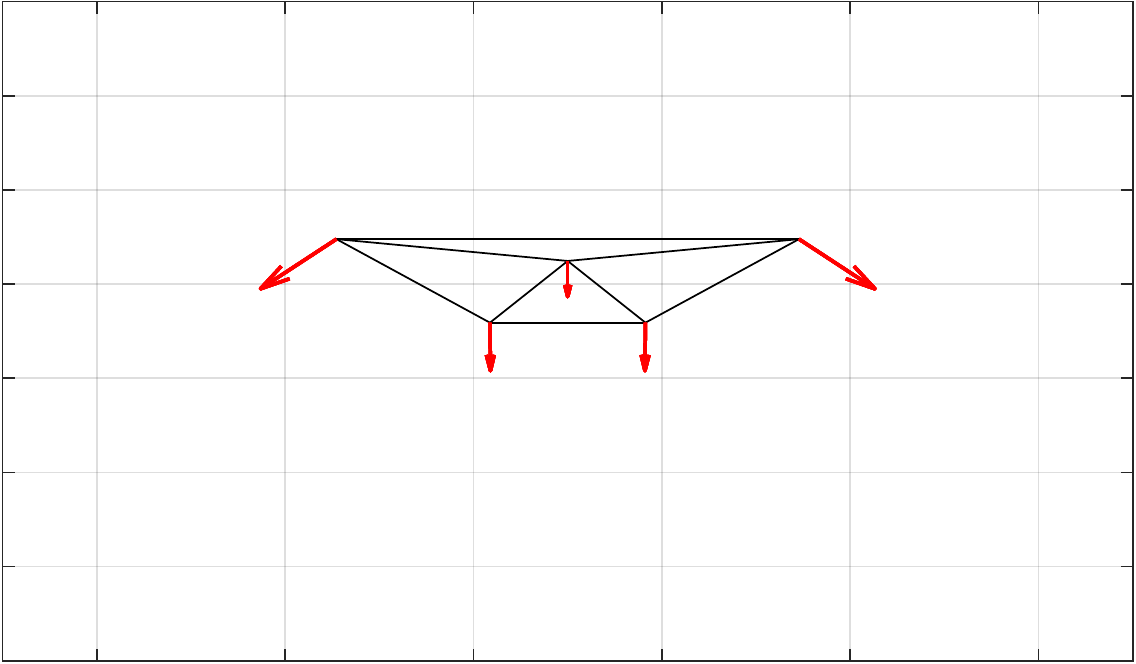}
		\caption{Complete, $t=1.87$}
	\end{subfigure}
	\hfill
	\begin{subfigure}{0.22\textwidth}
		\centering
		\includegraphics[width=\linewidth]{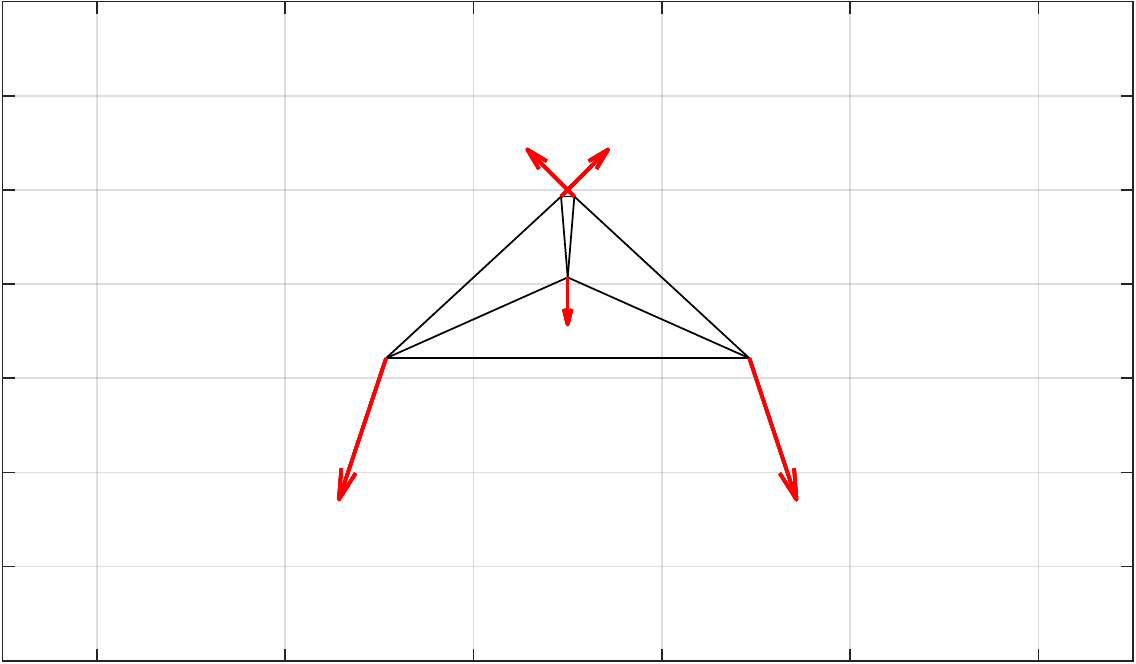}
		\caption{Euclidean, $t=1.87$}
	\end{subfigure}
	\hfill
	\begin{subfigure}{0.22\textwidth}
		\centering
		\includegraphics[width=\linewidth]{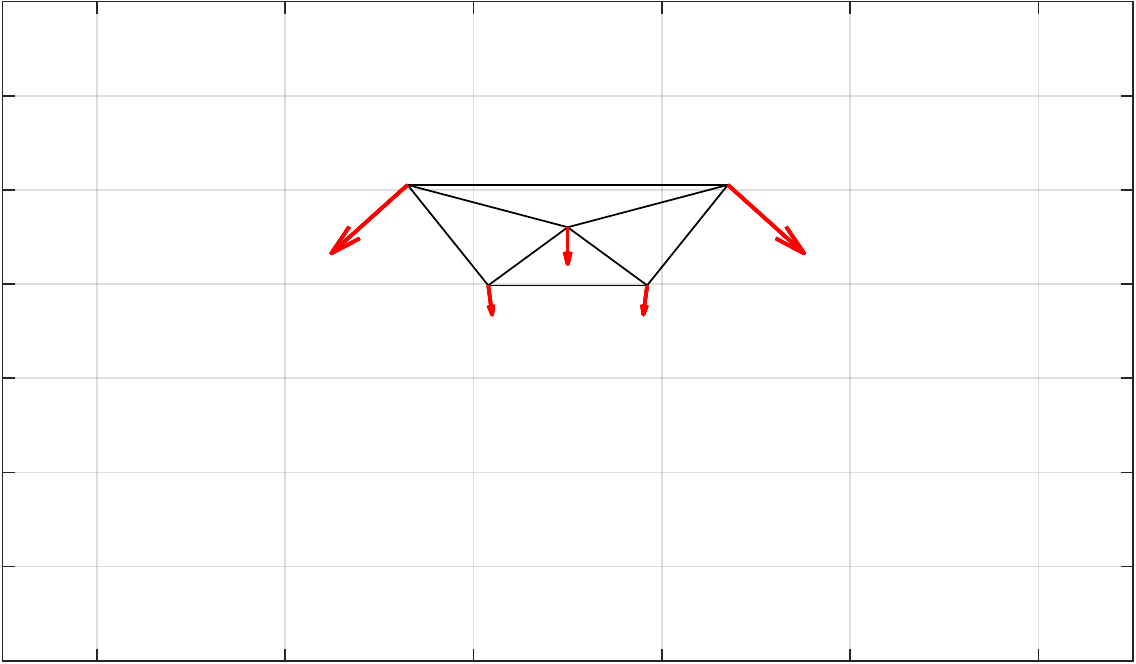}
		\caption{Complete, $t=0.93$}
	\end{subfigure}
	\hfill
	\begin{subfigure}{0.22\textwidth}
		\centering
		\includegraphics[width=\linewidth]{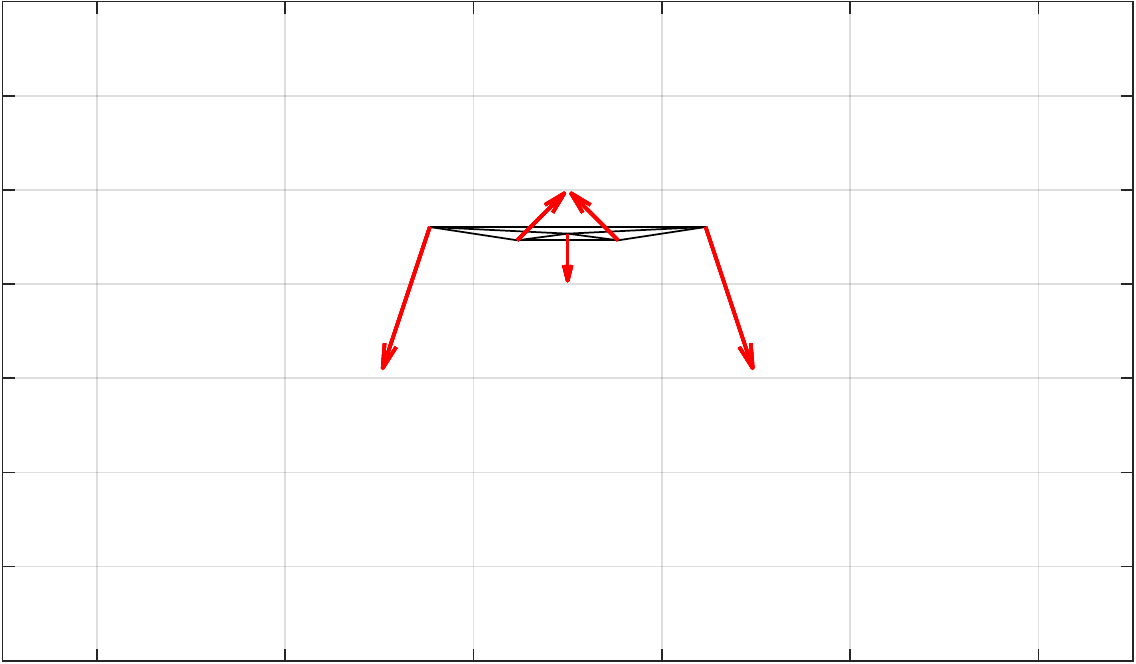}
		\caption{Euclidean, $t=0.93$}
	\end{subfigure}
	\hfill
	\begin{subfigure}{0.22\textwidth}
		\centering
		\includegraphics[width=\linewidth]{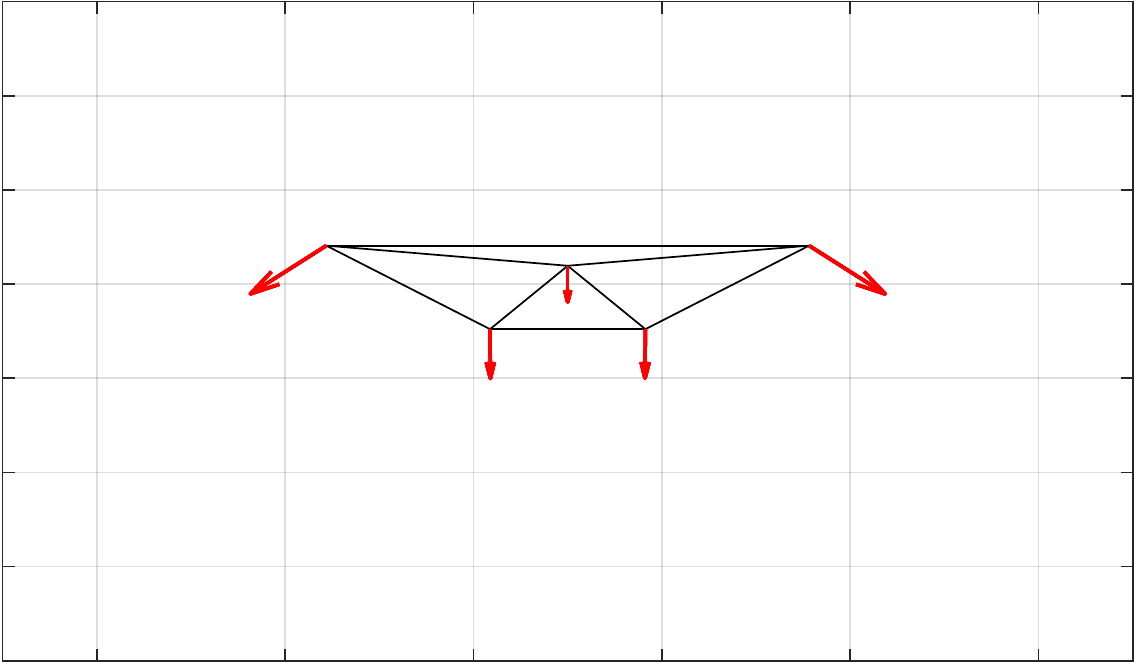}
		\caption{Complete, $t=2.0$}
	\end{subfigure}
	\hfill
	\begin{subfigure}{0.22\textwidth}
		\centering
		\includegraphics[width=\linewidth]{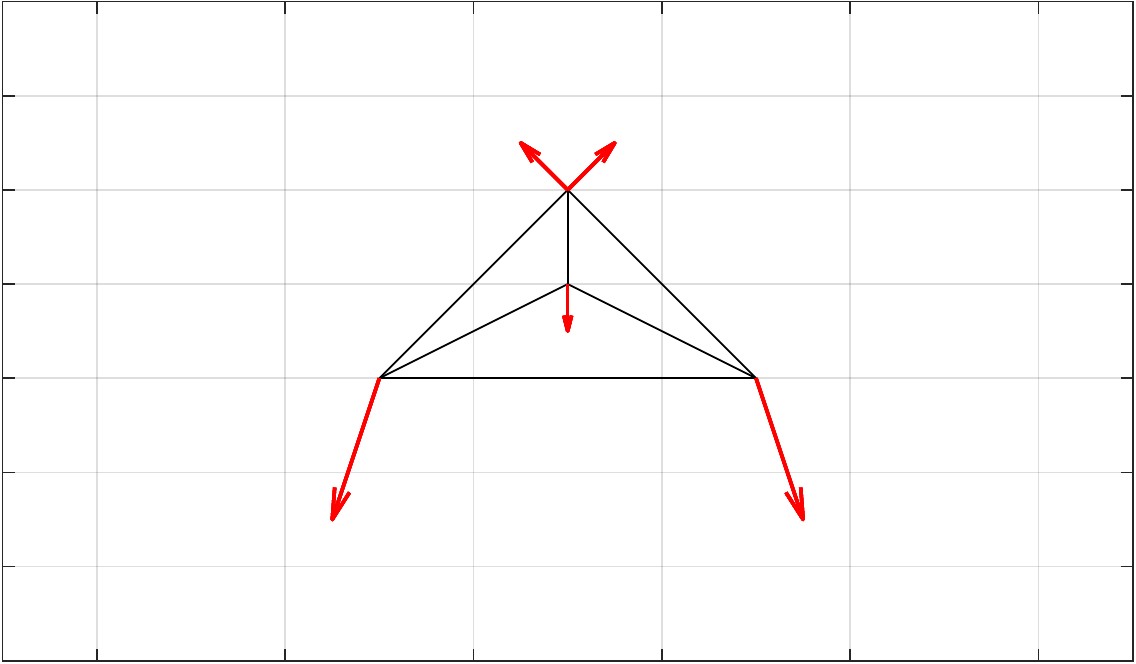}
		\caption{Euclidean, $t=2.0$}
	\end{subfigure}
	\caption{Snapshots of the geodesics described in \Cref{subsection:complete_vs_incomplete}, comparing the complete metric with $\beta_1 = \beta_3 = 1$ (first and third columns) and the Euclidean metric with $\beta_1 = \beta_3 = 0$ (second and fourth columns).}
	\label{fig:complete_geodesic}
\end{figure}

In this example, we revisit the setup depicted in \Cref{fig:euclidean_geodesic}, \ie, we retain a very simple initial mesh but consider the geodesic in the direction of an unfavorable initial tangent vector (\Cref{fig:euclidean_geodesic_tangent_vector}).
We compare the Euclidean geodesic ($\beta_1 = \beta_3 = 0$) with the proposed metric (for values $\beta_1 = \beta_3 = 1$).
In each case, $\Qref$ is chosen to be the initial mesh.
\Cref{fig:complete_geodesic} shows the respective mesh evolution at 16~snapshots within the interval $[0,2]$.
The scaling of the axes is the same in each snapshot.
The integration of the Euclidean geodesics is achieved with $N = 100$ time steps while we use $N = 5 \cdot 10^5$ in the case $\beta_1 = \beta_3 = 1$.

In the Euclidean case, the mesh degenerates very quickly and becomes non-admissible around $t = 1$.
By contrast, the completeness of the proposed metric ensures the mesh to be admissible for all times.
Notice that the inward pointing initial tangent vectors at the two bottom vertices are repelled and reverse direction around $t = 0.53$ as a consequence of the term involving $\beta_1$ in \eqref{eq:f_mu_without_exterior_term}, which avoids interior self-intersections.

\subsection{Mesh Quality Experiment}
\label{subsection:mesh_quality}

In this experiment, we consider a slightly more complex mesh consisting of $N_V = 25$~vertices and $N_T = 32$~triangles.
The initial mesh is a discretized version of the unit circle.
The initial tangent vector considered acts only on the boundary of the mesh, \ie, its components pertaining to interior mesh vertices are zero.
Such a situation occurs frequently in shape optimization, where the Hadamard structure theorem (\cite[Thm.~2.27]{SokolowskiZolesio1992}) provides an expression for the shape derivative which is supported only the boundary of the current mesh.
One would then usually apply an extension technique to obtain a displacement field (tangent vector) supported in all vertices.
A typical example is to achieve this through the solution of an elasticity equation.
With our approach, the displacement of all vertices, including the interior ones, is achieved automatically by following the geodesic with respect to the metric \eqref{eq:complete_metric_for_planar_triangular_meshes_without_exterior_term}.

As a quality measure for the mesh, we consider the mesh aspect ratio
\begin{equation}
	\label{eq:aspect_ratio}
	AR(\Delta;Q)
	\coloneqq
	\min_{[i_0,i_1,i_2]\in \Delta} \frac{2 \, \inradius{Q}[i_0,i_1,i_2]}{\circumradius{Q}[i_0,i_1,i_2]}
	.
\end{equation}
The aspect ratio takes values between $0$ and $1$, where the latter is achieved precisely for unilateral triangles.

\begin{figure}[htbp]
	\begin{center}
		\begin{subfigure}{0.3\textwidth}
			\includegraphics[width=\linewidth]{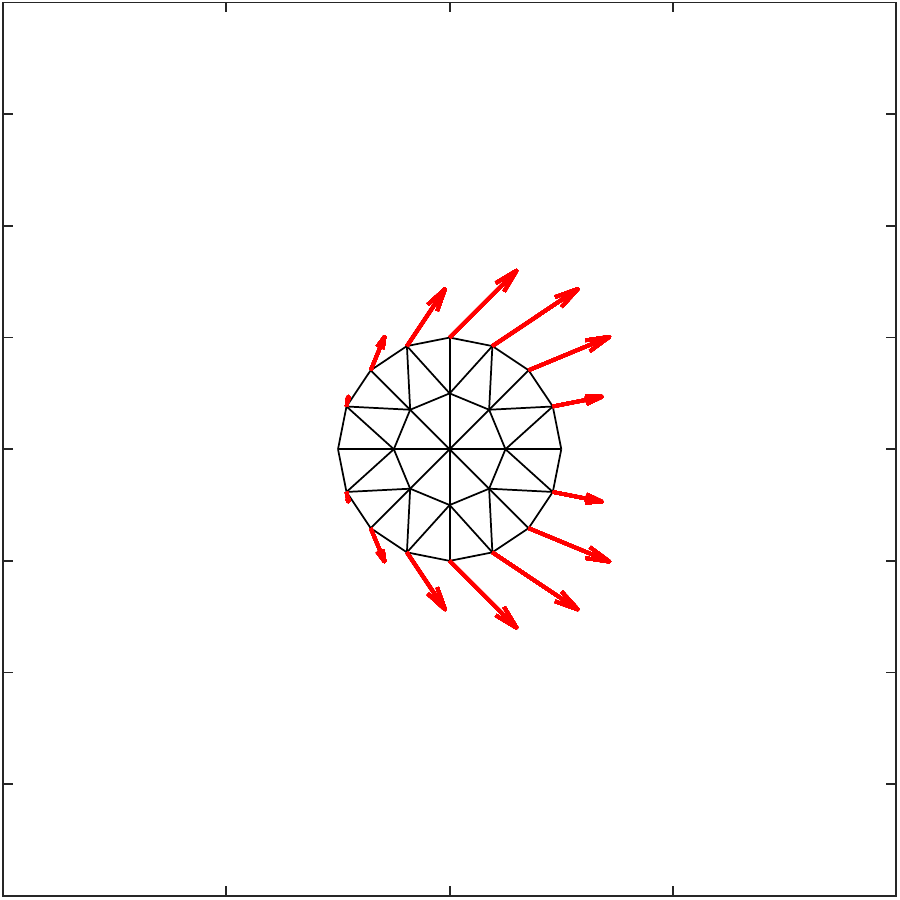}
			\caption{$\beta = 0.0$, $t=0$}
		\end{subfigure}
		\hfill
		\begin{subfigure}{0.3\textwidth}
			\includegraphics[width=\linewidth]{Figures/mesh_quality/initial_tangent_vector_and_mesh_MeshQuality.pdf}
			\caption{$\beta = 0.05$, $t=0$}
		\end{subfigure}
		\hfill
		\begin{subfigure}{0.3\textwidth}
			\includegraphics[width=\linewidth]{Figures/mesh_quality/initial_tangent_vector_and_mesh_MeshQuality.pdf}
			\caption{$\beta = 0.15$, $t=0$}
		\end{subfigure}
		\\
		\begin{subfigure}{0.3\textwidth}
			\includegraphics[width=\linewidth]{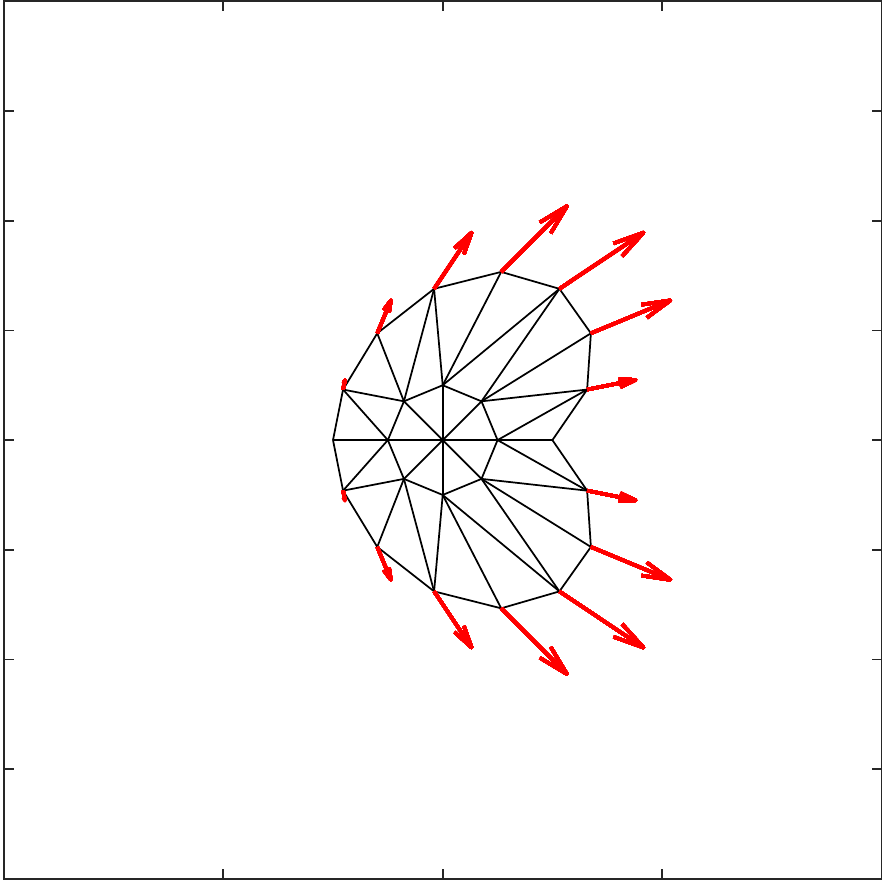}
			\caption{$\beta = 0.0$, $t=1$}
		\end{subfigure}
		\hfill
		\begin{subfigure}{0.3\textwidth}
			\includegraphics[width=\linewidth]{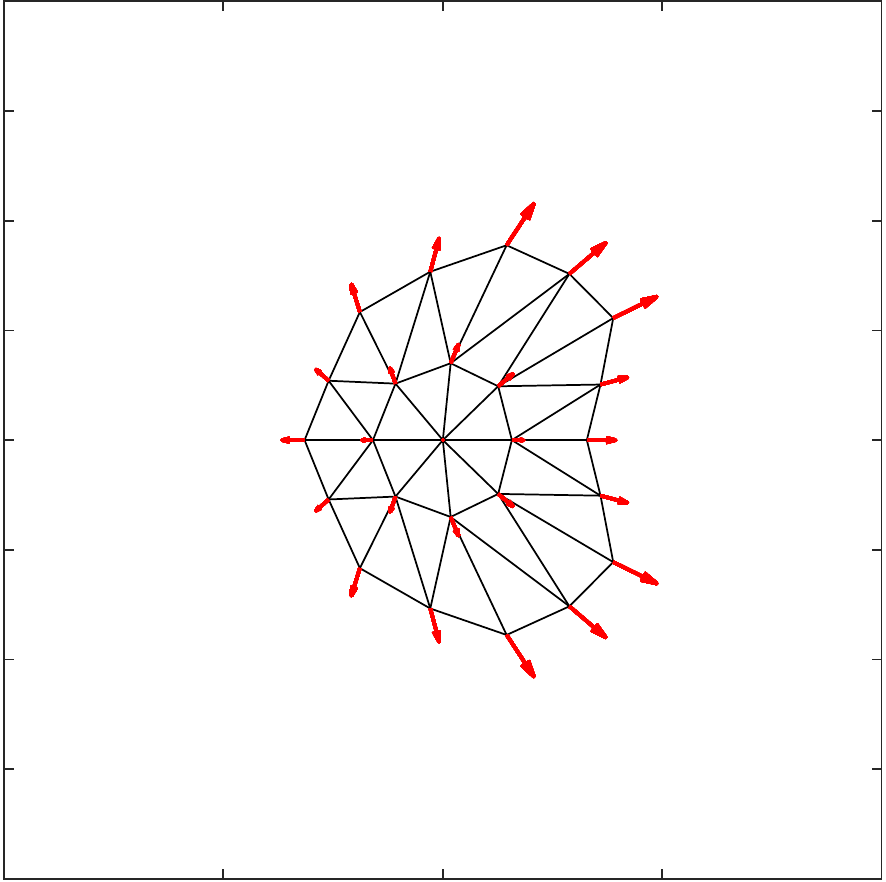}
			\caption{$\beta = 0.05$, $t=1$}
		\end{subfigure}
		\hfill
		\begin{subfigure}{0.3\textwidth}
			\includegraphics[width=\linewidth]{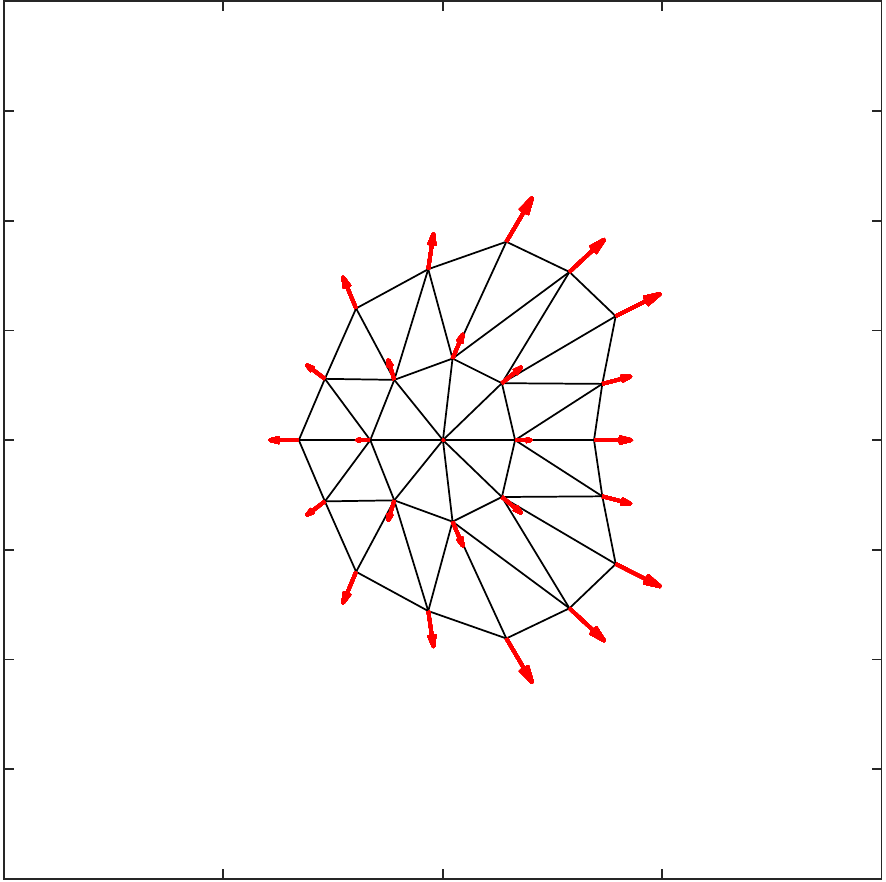}
			\caption{$\beta = 0.15$, $t=1$}
		\end{subfigure}
		\\
		\begin{subfigure}{0.3\textwidth}
			\includegraphics[width=\linewidth]{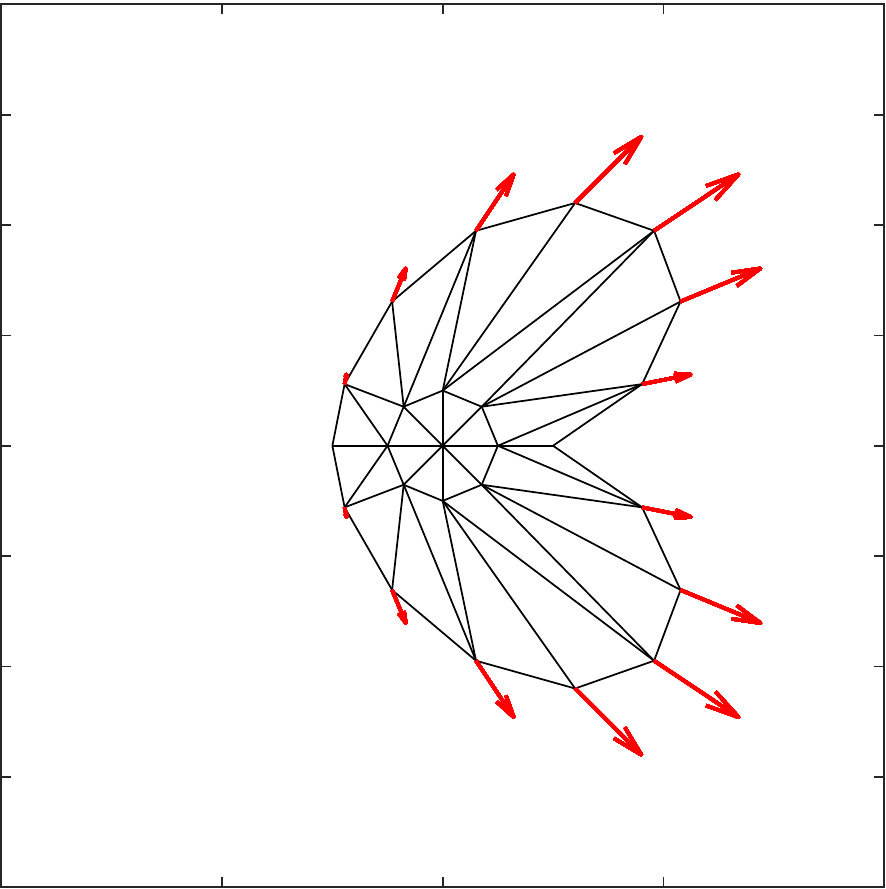}
			\caption{$\beta = 0.0$, $t=2$}
		\end{subfigure}
		\hfill
		\begin{subfigure}{0.3\textwidth}
			\includegraphics[width=\linewidth]{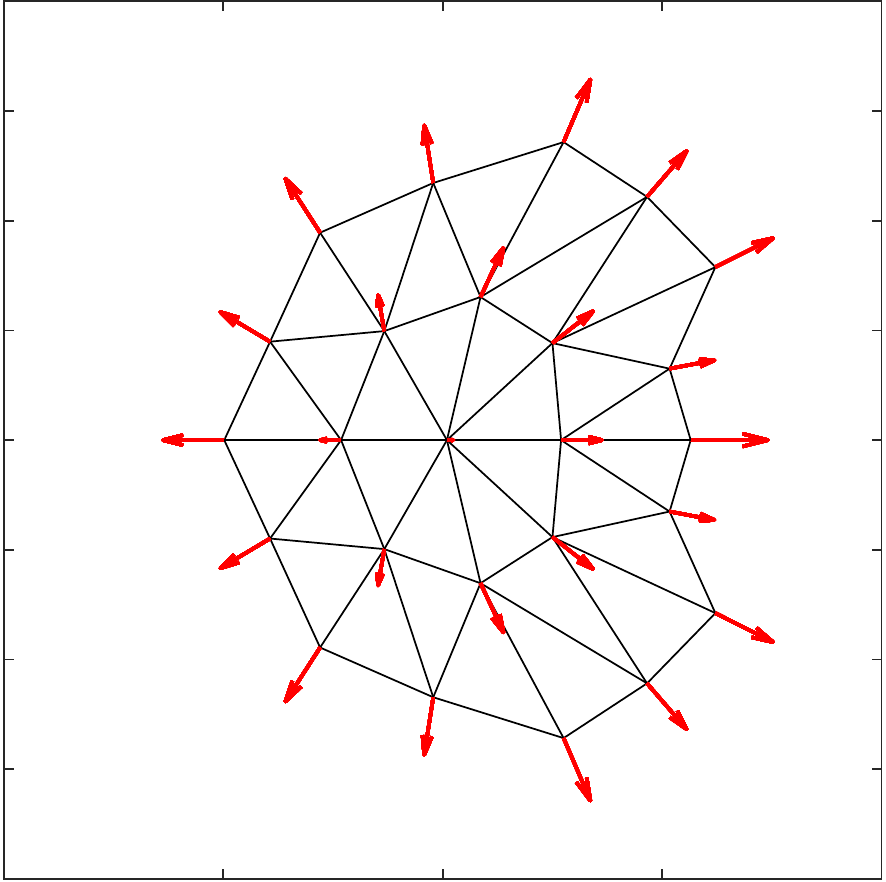}
			\caption{$\beta = 0.05$, $t=2$}
		\end{subfigure}
		\hfill
		\begin{subfigure}{0.3\textwidth}
			\includegraphics[width=\linewidth]{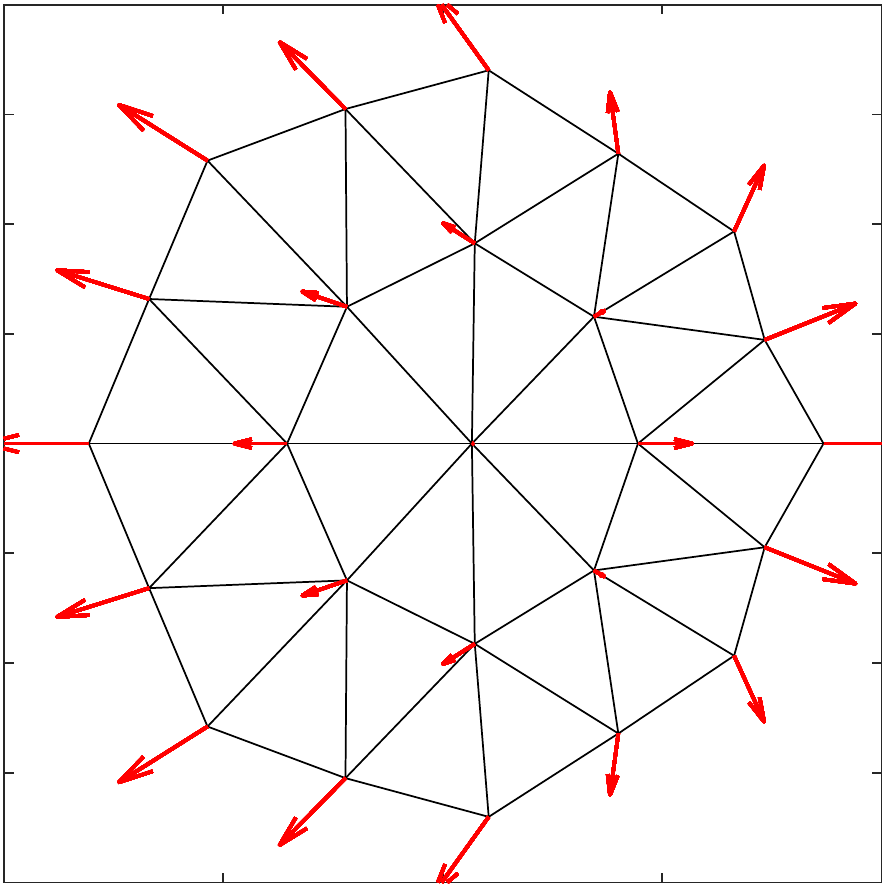}
			\caption{$\beta = 0.15$, $t=2$}
		\end{subfigure}
	\end{center}
	\caption{Snapshots of the geodesics described in \Cref{subsection:mesh_quality} for different values of $\beta = \beta_1 = \beta_3$. The tangent vectors are also shown. They have been scaled to improve the visualization.}
	\label{fig:MeshQuality_snapshots}
\end{figure}

In \Cref{fig:MeshQuality_snapshots} we show three snapshots of three geodesics on the interval $[0,2]$.
We compare the Euclidean geodesic ($\beta_1 = \beta_3 = 0$) with the proposed metric for values $\beta_1 = \beta_3 = 0.05$ and $\beta_1 = \beta_3 = 0.15$.
As in the previous experiments, $\Qref$ is chosen to be the initial mesh.
From here on, we refer to the joint value of $\beta_1$ and $\beta_3$ as $\beta \coloneqq \beta_1 = \beta_3$.

Note that the influence of the augmentation terms in the metric \eqref{eq:complete_metric_for_planar_triangular_meshes_without_exterior_term} becomes more pronounced as the values of $\beta$ grows.
	Increasing values of $\beta$ thus result in the geodesics to depart more and more from the Euclidean geodesics.
	As a consequence, smaller time step sizes~$\Delta t$ are required in order for the fixed-point iterations in \cref{line:implicit_equation_for_P,line:implicit_equation_for_Q} of the Störmer--Verlet scheme (\cref{algorithm:stoermer-verlet}) to converge.
	Consequently, we typically require larger numbers of time steps to sucessfully integrate a geodesic for larger values of $\beta$ on the same time interval.
	In the present experiment, we used $N = 10$ for $\beta = 0$, $N = 5 \cdot 10^2$ for $\beta = 0.05$, $N = 5 \cdot 10^3$ for $\beta = 0.15$.
	We also ran experiments with $\beta = 0.10$ (with $N = 2.5 \cdot 10^3$ time steps), $\beta = 0.2$ (with $N = 7.5 \cdot 10^3$), $\beta = 0.25$ (with $N = 10^4$) and $\beta=0.5$ (with $N = 5 \cdot 10^4$).

\begin{figure}[htb]
	\begin{center}
		\begin{tikzpicture}[scale=0.75]
			\begin{axis}[
				width=4.9in,
				height=4.2in,
				title = Mesh aspect ratios along geodesics,
				xlabel = time,
				ylabel = mesh aspect ratio \eqref{eq:aspect_ratio},
				xmin = 0,
				xmax = 2,
				legend style={draw=none},
				legend pos = south west
				]
				\addplot[solid,              ultra thick, color=TolVibrantBlue]    table [x=time,y=mesh_quality, col sep=comma]{src/MeshQuality_beta_0.txt};
				\addlegendentry{$\beta=0.00$}
				\addplot[dashdotdotted,      ultra thick, color=purple]            table [x=time,y=mesh_quality, col sep=comma]{src/MeshQuality_beta_5.txt};
				\addlegendentry{$\beta=0.05$}
				\addplot[dotted,             ultra thick, color=TolVibrantTeal]    table [x=time,y=mesh_quality, col sep=comma]{src/MeshQuality_beta_10.txt};
				\addlegendentry{$\beta=0.10$}
				\addplot[dashdotted,         ultra thick, color=red]               table [x=time,y=mesh_quality, col sep=comma]{src/MeshQuality_beta_15.txt};
				\addlegendentry{$\beta=0.15$}
				\addplot[loosely dotted,     ultra thick, color=olive]             table [x=time,y=mesh_quality, col sep=comma]{src/MeshQuality_beta_20.txt};
				\addlegendentry{$\beta=0.20$}
				\addplot[densely dashed,     ultra thick, color=TolVibrantMagenta] table [x=time,y=mesh_quality, col sep=comma]{src/MeshQuality_beta_25.txt};
				\addlegendentry{$\beta=0.25$}
				\addplot[densely dashdotted, ultra thick, color=TolBrightBlue]     table [x=time,y=mesh_quality, col sep=comma]{src/MeshQuality_beta_50.txt};
				\addlegendentry{$\beta=0.50$}
			\end{axis}
		\end{tikzpicture}
	\end{center}
	\caption{Mesh quality plot over time along geodesics for different values of $\beta = \beta_1 = \beta_3$; see \cref{subsection:mesh_quality}.}
	\label{fig:MeshQuality}
\end{figure}
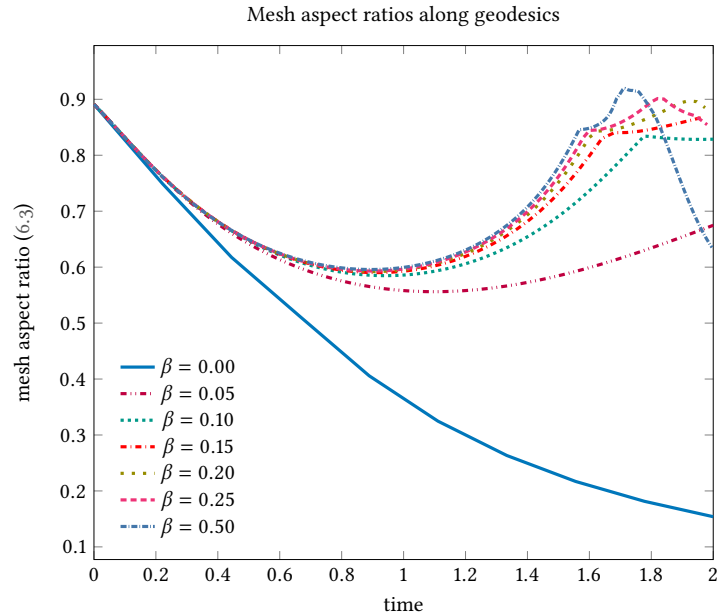

As we can observe from these snapshots, but also from the mesh quality plot in \Cref{fig:MeshQuality}, the mesh aspect ratio over time degenerates for the Euclidean metric, where only the boundary values are displaced.
By contrast, the mesh aspect ratio is observed to be bounded away from zero for the metric \eqref{eq:complete_metric_for_planar_triangular_meshes_without_exterior_term} when $\beta > 0$.
Even small, positive values of $\beta$ were observed to be effective in improving the mesh quality.
Interestingly, the quality was almost identical independently of the choice of $\beta \in \{ 0.15, 0.20, 0.25\}$.

\subsection{Run Times}
\label{subsection:run_times}

We end this section by showing the run times to compute the geodesics depicted in \cref{fig:MeshQuality_snapshots} associated with strictly positive values of $\beta_1$ and $\beta_3$, \ie, $\beta_1 = \beta_3 = 0.05$ and $\beta_1 = \beta_3 = 0.15$.
We also consider the run times for the corresponding geodesics on uniform refinements of the initial, coarse mesh.

In order to define these corresponding geodesics on finer meshes, we aim to generate the same value of the augmentation function for all meshes.
To this end, we have to adapt the values of $\beta_1$ and $\beta_3$.
More precisely, we can find from \eqref{eq:f_mu_without_exterior_term} that we need $2^{\meshLevel-1} N_T \beta_1$ and $N_V \beta_3$ to be constant across all mesh levels.
Here the parameter $\meshLevel$ denotes the mesh level $\meshLevel=1,2,3$, and $N_T$, $N_V$ are the number of triangles and vertices on the respective mesh level.
For instance, in case of the first experiment in \cref{tab:timesOfExecution} we have $2^{\meshLevel-1} N_T \beta_1 = 1.6 $, and $N_V \beta_3 = 1.25$.

The geodesics are computed with $N = 5 \cdot 10^2$ time steps on the time interval $[0,2]$ across all mesh levels for the first experiment and with $N = 5 \cdot 10^3$ for the second.
\Cref{tab:timesOfExecution} collects the CPU times for the computation of the geodesics.
All experiments were performed on a laptop with an Intel Core~i7-1065G7 CPU with 1.3~GHz and 16~GiB of RAM.
\begin{table}
	\begin{center}
		\sisetup{round-mode = places, round-precision = 0}
		\begin{filecontents*}{src/timesOfExecution_meshLevel.txt}
N,meshLevel,Nt,Nq,Ne,beta1,beta3,time
500,1,32,25,16,0.05,0.05,4.8735
500,2,128,81,32,0.00625,0.015432,9.6428
500,3,512,289,64,0.00078125,0.0043253,18.322
5000,1,32,25,16,0.15,0.15,43.703
5000,2,128,81,32,0.01875,0.046296,76.019
5000,3,512,289,64,0.0023437,0.012976,195.94
\end{filecontents*}
\pgfplotstabletypeset[%
		col sep = comma,
		every head row/.style={before row=\toprule,after row=\midrule},
		every last row/.style={after row=\bottomrule},
		every nth row={3}{before row=\midrule},
		columns/N/.style={
		assign cell content/.code={%
			\ifnum\pgfplotstablerow=0
				\pgfkeyssetvalue{/pgfplots/table/@cell content}%
				{\multirow{3}{*}{\pgfmathprintnumber[precision = 0, sci, sci zerofill,retain unit mantissa=false]{##1}}}%
			\else
				\ifnum\pgfplotstablerow=3
					\pgfkeyssetvalue{/pgfplots/table/@cell content}%
					{\multirow{3}{*}{\pgfmathprintnumber[precision = 0, sci, sci zerofill,retain unit mantissa=false]{##1}}}
				\else
					\ifnum\pgfplotstablerow=6
						\pgfkeyssetvalue{/pgfplots/table/@cell content}%
						{\multirow{3}{*}{\pgfmathprintnumber[precision = 0, sci, sci zerofill,retain unit mantissa=false]{##1}}}%
					\else
						\pgfkeyssetvalue{/pgfplots/table/@cell content}{}%
					\fi
				\fi
			\fi
		},
		column name=$N$},
		columns/meshLevel/.style={
		int detect,
		column name={mesh level}},
		columns/Nt/.style={
		int detect,
		column name=$N_T$},
		columns/Nq/.style={
		int detect,
		column name=$N_V$},
		columns/Ne/.style={
		int detect,
		column name=$N_E$},
		columns/beta1/.style={
		fixed, fixed zerofill,
		precision=4,
		column name=$\beta_1$},
		columns/beta3/.style={
		fixed, fixed zerofill,
		precision=4,
		column name=$\beta_3$},
		columns/time/.style={std, precision=2,column name={time},
		assign cell content/.code={%
		\pgfkeyssetvalue{/pgfplots/table/@cell content}%
		{\SI{##1}{\second}}%
		}}
		]{src/timesOfExecution_meshLevel.txt}
	\end{center}
	\caption{CPU times for the computation of the geodesics depicted in \Cref{fig:MeshQuality_snapshots} (center and right column) for different mesh levels, with parameters $\beta_1$ and $\beta_3$ adjusted to the mesh level; see \Cref{subsection:run_times}.}
	\label{tab:timesOfExecution}
\end{table}

\section*{Conclusions and Outlook}

In this paper we studied the set $\zeromanifold$ of all possible configurations (vertex coordinates) for planar triangular meshes of a given connectivity.
These were described using the language of simplicial complexes.
We proved that $\zeromanifold$ is an open submanifold of $\R^{2 \times N_V}$.
Our primary object of interest, however, was $\planarmanifold$, a path component of $\zeromanifold$ with one of the two possible orientations.
Our main result was the construction of a complete Riemannian metric for $\planarmanifold$.
To this end, we devised a proper function and used a theorem of \cite{Gordon1973}.

The main motivation for this study originates in computational shape optimization, where triangular meshes are frequently used to represent planar shapes, in particular when the problem involves a partial differential equation.
It is common practice in computational shape optimization to deform meshes along Euclidean geodesics, which is computationally trivial but inevitably leads to a restriction in step sizes in order to counteract mesh degeneration.
The metric we propose avoids this.

As a drawback, we obtained geodesic equations which can only be integrated numerically.
Our method of choice was the Störmer--Verlet scheme.
Although we put some effort into an efficient implementation which will be described elsewhere, the solution of the geodesic equation remains computationally involved and may easily become the most expensive step in a shape optimization loop.
Therefore, it appears reasonable to investigate compromises between following exact geodesics with respect to the complete metric, and following the straight lines of Euclidean geometry.
This is a topic for future research and partially addressed in \cite{HerzogLoayzaRomero:2021:1}.

We hope that the ideas presented here may be of interest for applications other than shape optimization as well, including mesh morphing, mesh interpolation, and image registration.
In this context, it becomes important to be able to evaluate, besides the exponential map, also the logarithmic map as well as the parallel transport of tangent vectors.
We are planning to investigate these topics in a separate publication.

\appendix

\section{Geometric and Abstract Simplicial Complexes}
\label{section:geometric_abstract_simplicial_complexes}

This section collects basic definitions concerning geometric and abstract simplicial complexes.
For a more thorough introduction, we refer the reader to, \eg, \cite{HorakJost2013,EdelsbrunnerHarer2010,Misztal2010}.

\subsection{Geometric Simplicial Complexes}
\label{subsection:geometric_simplicial_complexes}

We begin with the notion of geometric simplicial complexes.
Throughout, let $Z$ be a finite dimensional vector space.
(In our applications, we will have $Z = \R^2$.)
A \textbf{simplex}~$\sigma$ of dimension~$k \in \N_0$ (or \textbf{$k$-simplex}) in~$Z$ is the convex hull of $k+1$ affine independent points (the \textbf{vertices} of $\sigma$) in $Z$.
A \textbf{face} of \textbf{dimension~$m$} ($0 \le m \le k$) (an \textbf{$m$-face}) of~$\sigma$ is the convex hull of a subset of $m + 1$ of its vertices.
$0$-faces are vertices, $1$-faces are \textbf{edges} and $2$-faces are \textbf{triangles}.

A (finite) \textbf{geometric simplicial complex}~$\Sigma$ in~$Z$ is a non-empty, finite set of simplices in $Z$ satisfying the following conditions:
every face of a simplex~$\sigma \in \Sigma$ also belongs to $\Sigma$, and the non-empty intersection of any two simplices~$\sigma, \sigma'$ in $\Sigma$ is a face of both $\sigma$ and $\sigma'$.
We say that a geometric simplicial complex~$\Sigma$ is of \textbf{dimension} $k \in \N_0$ (or a \textbf{geometric simplicial $k$-complex}) if $k$ is the largest dimension of any simplex in $\Sigma$.
The \textbf{vertex set} of a geometric simplicial complex~$\Sigma$ is the union of the vertices of all of its faces.

A geometric simplicial $k$-complex~$\Sigma$ is \textbf{pure} if all maximal elements of $\Sigma$ (\wrt\ the partial order of set inclusion) have dimension~$k$.
In other words, a geometric simplicial $k$-complex is pure if and only if every simplex in $\Sigma$ is the face of some $k$-simplex in $\Sigma$.
A geometric simplicial $k$-complex~$\Sigma$ in $Z$ is said to be \textbf{$m$-path connected} ($1 \le m \le k$) if for any two distinct $m$-faces $\sigma, \sigma'$ in $\Delta$, there exists a finite sequence of $m$-faces, starting in $\sigma_0 = \sigma$ and ending in $\sigma_n = \sigma'$, such that $\sigma_i \cap \sigma_{i+1}$ is an $(m-1)$-face for $i = 0, \ldots, n-1$.

\begin{figure}[htb]
	\begin{center}
		\begin{subfigure}{0.3\textwidth}
			\centering
			\makeatletter
			\ltx@ifclassloaded{mcom-l}{
			\begin{tikzpicture}[scale = 0.4, line cap=round,line join=round,>=triangle 45,x=1.0cm,y=1.0cm]
				}{%
				\begin{tikzpicture}[scale = 0.5, line cap=round,line join=round,>=triangle 45,x=1.0cm,y=1.0cm]
				}
				\makeatother
\clip(-3,-2.5) rectangle (6,7.5);
\draw [color = gray, line width = 0.2pt, domain=-6.230245364561282:0.0] plot(\x,{(-0--2.6*\x)/-1.5});
\draw [color = gray, line width = 0.2pt, domain=-1.5000000000000004:17.823872468411377] plot(\x,{(-0-2.6*\x)/1.5});
\draw [color = gray, line width = 0.2pt,domain=-6.230245364561282:3.0] plot(\x,{(-15.59--5.2*\x)/-3});
\draw [color = gray, line width = 0.2pt,domain=1.1102230246251565E-15:17.823872468411377] plot(\x,{(--15.59-5.2*\x)/3});
\draw [color = gray, line width = 0.2pt,domain=-6.230245364561282:4.5] plot(\x,{(-15.59--2.6*\x)/-1.5});
\draw [color = gray, line width = 0.2pt,domain=3.0:17.823872468411377] plot(\x,{(--15.59-2.6*\x)/1.5});
\draw [color = gray, line width = 0.2pt,domain=-6.230245364561282:4.5] plot(\x,{(-15.59-0*\x)/-6});
\draw [color = gray, line width = 0.2pt,domain=-1.5000000000000004:17.823872468411377] plot(\x,{(--15.59-0*\x)/6});
\draw [color = gray, line width = 0.2pt,domain=0.0:17.823872468411377] plot(\x,{(-0--5.2*\x)/3});
\draw [color = gray, line width = 0.2pt,domain=-6.230245364561282:3.0] plot(\x,{(-0-5.2*\x)/-3});
\draw [color = gray, line width = 0.2pt,domain=-1.5000000000000004:17.823872468411377] plot(\x,{(--7.79--2.6*\x)/1.5});
\draw [color = gray, line width = 0.2pt,domain=-6.230245364561282:1.1102230246251565E-15] plot(\x,{(-7.79-2.6*\x)/-1.5});
\draw [color = gray, line width = 0.2pt,domain=-6.230245364561282:4.5] plot(\x,{(--7.79-2.6*\x)/-1.5});
\draw [color = gray, line width = 0.2pt,domain=3.0:17.823872468411377] plot(\x,{(-7.79--2.6*\x)/1.5});
\draw [color = gray, line width = 0.2pt,domain=0.0:17.823872468411377] plot(\x,{(-0-0*\x)/3});
\draw [color = gray, line width = 0.2pt,domain=-6.230245364561282:3.0] plot(\x,{(-0-0*\x)/-6});
\draw [color = gray, line width = 0.2pt,domain=1.1102230246251565E-15:17.823872468411377] plot(\x,{(--15.59-0*\x)/3});
\draw [color = gray, line width = 0.2pt,domain=-6.230245364561282:3.0] plot(\x,{(-15.59-0*\x)/-3});

\fill[fill=red,fill opacity=0.2] (0,0) -- (3,0) -- (4.5,2.6) -- (3,5.2) -- (0,5.2) -- (-1.5,2.6) -- cycle;
\draw [color = black, line width=1pt] (0,0)-- (3,0) -- (4.5,2.6) -- (3,5.2) -- (0,5.2) -- (-1.5,2.6) -- cycle;
\fill [color=blue] (1.5,2.6) circle (6pt);
\draw[color=blue] (1.5,2.4) node[below] {$\sigma$};
\end{tikzpicture}
			\caption{$\St(\sigma)$}
			\label{subfigure:star_of_0-face}
		\end{subfigure}
		\hfill
		\begin{subfigure}{0.3\textwidth}
			\makeatletter
			\ltx@ifclassloaded{mcom-l}{
			\begin{tikzpicture}[scale = 0.4, line cap=round,line join=round,>=triangle 45,x=1.0cm,y=1.0cm]
				}{%
				\begin{tikzpicture}[scale = 0.5, line cap=round,line join=round,>=triangle 45,x=1.0cm,y=1.0cm]
				}
				\makeatother
\clip(-3,-2.5) rectangle (6,7.5);
\draw [color = gray, line width = 0.2pt, domain=-6.230245364561282:0.0] plot(\x,{(-0--2.6*\x)/-1.5});
\draw [color = gray, line width = 0.2pt, domain=-1.5000000000000004:17.823872468411377] plot(\x,{(-0-2.6*\x)/1.5});
\draw [color = gray, line width = 0.2pt, domain=-6.230245364561282:3.0] plot(\x,{(-15.59--5.2*\x)/-3});
\draw [color = gray, line width = 0.2pt, domain=1.1102230246251565E-15:17.823872468411377] plot(\x,{(--15.59-5.2*\x)/3});
\draw [color = gray, line width = 0.2pt, domain=-6.230245364561282:4.5] plot(\x,{(-15.59--2.6*\x)/-1.5});
\draw [color = gray, line width = 0.2pt, domain=3.0:17.823872468411377] plot(\x,{(--15.59-2.6*\x)/1.5});
\draw [color = gray, line width = 0.2pt, domain=-6.230245364561282:4.5] plot(\x,{(-15.59-0*\x)/-6});
\draw [color = gray, line width = 0.2pt, domain=-1.5000000000000004:17.823872468411377] plot(\x,{(--15.59-0*\x)/6});
\draw [color = gray, line width = 0.2pt, domain=0.0:17.823872468411377] plot(\x,{(-0--5.2*\x)/3});
\draw [color = gray, line width = 0.2pt, domain=-6.230245364561282:3.0] plot(\x,{(-0-5.2*\x)/-3});
\draw [color = gray, line width = 0.2pt, domain=-1.5000000000000004:17.823872468411377] plot(\x,{(--7.79--2.6*\x)/1.5});
\draw [color = gray, line width = 0.2pt, domain=-6.230245364561282:1.1102230246251565E-15] plot(\x,{(-7.79-2.6*\x)/-1.5});
\draw [color = gray, line width = 0.2pt, domain=-6.230245364561282:4.5] plot(\x,{(--7.79-2.6*\x)/-1.5});
\draw [color = gray, line width = 0.2pt, domain=3.0:17.823872468411377] plot(\x,{(-7.79--2.6*\x)/1.5});
\draw [color = gray, line width = 0.2pt, domain=0.0:17.823872468411377] plot(\x,{(-0-0*\x)/3});
\draw [color = gray, line width = 0.2pt, domain=-6.230245364561282:3.0] plot(\x,{(-0-0*\x)/-6});
\draw [color = gray, line width = 0.2pt, domain=1.1102230246251565E-15:17.823872468411377] plot(\x,{(--15.59-0*\x)/3});
\draw [color = gray, line width = 0.2pt, domain=-6.230245364561282:3.0] plot(\x,{(-15.59-0*\x)/-3});
\fill[fill=red ,fill opacity=0.2] (0,0) -- (3,0) -- (4.5,2.6) -- (3,5.2) -- (0,5.2) -- (-1.5,2.6) -- cycle;
\draw [line width=1pt,color=red!70] (0,0)-- (3,0)-- (4.5,2.6)-- (3,5.2)-- (0,5.2)-- (-1.5,2.6) -- cycle;
\draw [line width=1pt,color=red!70](0,5.2)-- (3,0);
\draw [line width=1pt,color=red!70](3,5.2)-- (0,0);
\draw [line width=1pt,color=red!70](-1.5,2.6)-- (4.5,2.6);
\fill [color=blue] (1.5,2.6) circle (6pt);
\draw[color=blue] (1.5,2.4) node[below] {$\sigma$};
\end{tikzpicture}
			\caption{$\closedSt{\sigma}$}
			\label{subfigure:closed_star_of_0-face}
		\end{subfigure}
		\hfill
		\begin{subfigure}{0.3\textwidth}
			\makeatletter
			\ltx@ifclassloaded{mcom-l}{
			\begin{tikzpicture}[scale = 0.4, line cap=round,line join=round,>=triangle 45,x=1.0cm,y=1.0cm]
				}{%
				\begin{tikzpicture}[scale = 0.5, line cap=round,line join=round,>=triangle 45,x=1.0cm,y=1.0cm]
				}
				\makeatother
\clip(-3,-2.5) rectangle (6,7.5);
\draw [color = gray, line width = 0.2pt, domain=-6.230245364561282:0.0] plot(\x,{(-0--2.6*\x)/-1.5});
\draw [color = gray, line width = 0.2pt, domain=-1.5000000000000004:17.823872468411377] plot(\x,{(-0-2.6*\x)/1.5});
\draw [color = gray, line width = 0.2pt, domain=-6.230245364561282:3.0] plot(\x,{(-15.59--5.2*\x)/-3});
\draw [color = gray, line width = 0.2pt, domain=1.1102230246251565E-15:17.823872468411377] plot(\x,{(--15.59-5.2*\x)/3});
\draw [color = gray, line width = 0.2pt, domain=-6.230245364561282:4.5] plot(\x,{(-15.59--2.6*\x)/-1.5});
\draw [color = gray, line width = 0.2pt, domain=3.0:17.823872468411377] plot(\x,{(--15.59-2.6*\x)/1.5});
\draw [color = gray, line width = 0.2pt, domain=-6.230245364561282:4.5] plot(\x,{(-15.59-0*\x)/-6});
\draw [color = gray, line width = 0.2pt, domain=-1.5000000000000004:17.823872468411377] plot(\x,{(--15.59-0*\x)/6});
\draw [color = gray, line width = 0.2pt, domain=0.0:17.823872468411377] plot(\x,{(-0--5.2*\x)/3});
\draw [color = gray, line width = 0.2pt, domain=-6.230245364561282:3.0] plot(\x,{(-0-5.2*\x)/-3});
\draw [color = gray, line width = 0.2pt, domain=-1.5000000000000004:17.823872468411377] plot(\x,{(--7.79--2.6*\x)/1.5});
\draw [color = gray, line width = 0.2pt, domain=-6.230245364561282:1.1102230246251565E-15] plot(\x,{(-7.79-2.6*\x)/-1.5});
\draw [color = gray, line width = 0.2pt, domain=-6.230245364561282:4.5] plot(\x,{(--7.79-2.6*\x)/-1.5});
\draw [color = gray, line width = 0.2pt, domain=3.0:17.823872468411377] plot(\x,{(-7.79--2.6*\x)/1.5});
\draw [color = gray, line width = 0.2pt, domain=0.0:17.823872468411377] plot(\x,{(-0-0*\x)/3});
\draw [color = gray, line width = 0.2pt, domain=-6.230245364561282:3.0] plot(\x,{(-0-0*\x)/-6});
\draw [color = gray, line width = 0.2pt, domain=1.1102230246251565E-15:17.823872468411377] plot(\x,{(--15.59-0*\x)/3});
\draw [color = gray, line width = 0.2pt, domain=-6.230245364561282:3.0] plot(\x,{(-15.59-0*\x)/-3});
\draw [line width=1pt,color=red!70] (0,0)-- (3,0)-- (4.5,2.6) -- (3,5.2) -- (0,5.2) -- (-1.5,2.6) -- cycle;
\draw [line width=1pt,color=black] (0,5.2)-- (3,0);
\draw [line width=1pt,color=black](3,5.2)-- (0,0);
\draw [line width=1pt,color=black] (-1.5,2.6)-- (4.5,2.6);
\fill [color=blue] (1.5,2.6) circle (6pt);
\draw[color=blue] (1.5,2.4) node[below] {$\sigma$};
\end{tikzpicture}
			\caption{$\lk(\sigma)$}
			\label{subfigure:link_of_0-face}
		\end{subfigure}
		\\
		\begin{subfigure}{0.3\textwidth}
			\makeatletter
			\ltx@ifclassloaded{mcom-l}{
			\begin{tikzpicture}[scale = 0.4, line cap=round,line join=round,>=triangle 45,x=1.0cm,y=1.0cm]
				}{%
				\begin{tikzpicture}[scale = 0.5, line cap=round,line join=round,>=triangle 45,x=1.0cm,y=1.0cm]
				}
				\makeatother
				\clip(-3,-2.5) rectangle (6,7.5);
\draw [color = gray, line width = 0.2pt,domain=-6.230245364561282:0.0] plot(\x,{(-0--2.6*\x)/-1.5});
\draw [color = gray, line width = 0.2pt,domain=-1.5000000000000004:17.823872468411377] plot(\x,{(-0-2.6*\x)/1.5});
\draw [color = gray, line width = 0.2pt,domain=-6.230245364561282:3.0] plot(\x,{(-15.59--5.2*\x)/-3});
\draw [color = gray, line width = 0.2pt,domain=1.1102230246251565E-15:17.823872468411377] plot(\x,{(--15.59-5.2*\x)/3});
\draw [color = gray, line width = 0.2pt,domain=-6.230245364561282:4.5] plot(\x,{(-15.59--2.6*\x)/-1.5});
\draw [color = gray, line width = 0.2pt,domain=3.0:17.823872468411377] plot(\x,{(--15.59-2.6*\x)/1.5});
\draw [color = gray, line width = 0.2pt,domain=-6.230245364561282:4.5] plot(\x,{(-15.59-0*\x)/-6});
\draw [color = gray, line width = 0.2pt,domain=-1.5000000000000004:17.823872468411377] plot(\x,{(--15.59-0*\x)/6});
\draw [color = gray, line width = 0.2pt,domain=0.0:17.823872468411377] plot(\x,{(-0--5.2*\x)/3});
\draw [color = gray, line width = 0.2pt,domain=-6.230245364561282:3.0] plot(\x,{(-0-5.2*\x)/-3});
\draw [color = gray, line width = 0.2pt,domain=-1.5000000000000004:17.823872468411377] plot(\x,{(--7.79--2.6*\x)/1.5});
\draw [color = gray, line width = 0.2pt,domain=-6.230245364561282:1.1102230246251565E-15] plot(\x,{(-7.79-2.6*\x)/-1.5});
\draw [color = gray, line width = 0.2pt,domain=-6.230245364561282:4.5] plot(\x,{(--7.79-2.6*\x)/-1.5});
\draw [color = gray, line width = 0.2pt,domain=3.0:17.823872468411377] plot(\x,{(-7.79--2.6*\x)/1.5});
\draw [color = gray, line width = 0.2pt,domain=0.0:17.823872468411377] plot(\x,{(-0-0*\x)/3});
\draw [color = gray, line width = 0.2pt,domain=-6.230245364561282:3.0] plot(\x,{(-0-0*\x)/-6});
\draw [color = gray, line width = 0.2pt,domain=1.1102230246251565E-15:17.823872468411377] plot(\x,{(--15.59-0*\x)/3});
\draw [color = gray, line width = 0.2pt,domain=-6.230245364561282:3.0] plot(\x,{(-15.59-0*\x)/-3});
\fill[color=black,fill=red,fill opacity=0.2] (3,5.2) -- (1.5,2.6) -- (4.5,2.6) -- cycle;
\fill[color=black,fill=red,fill opacity=0.2] (1.5,2.6) -- (4.5,2.6) -- (3,0) -- cycle;
\draw [line width=1pt,color=blue] (1.5,2.6)-- (4.5,2.6);
\draw [line width=1pt,color=black] (3,5.2)-- (4.5,2.6)-- (3,0) -- (1.5,2.6) -- cycle;
\draw[color=blue] (3,2.4) node[below] {$\sigma$};
\end{tikzpicture}
			\caption{$\St(\sigma)$}
			\label{subfigure:star_of_1-face}
		\end{subfigure}
		\hfill
		\begin{subfigure}{0.3\textwidth}
			\makeatletter
			\ltx@ifclassloaded{mcom-l}{
			\begin{tikzpicture}[scale = 0.4, line cap=round,line join=round,>=triangle 45,x=1.0cm,y=1.0cm]
				}{%
				\begin{tikzpicture}[scale = 0.5, line cap=round,line join=round,>=triangle 45,x=1.0cm,y=1.0cm]
				}
				\makeatother
				\clip(-3,-2.5) rectangle (6,7.5);
\draw [color = gray, line width = 0.2pt,domain=-6.230245364561282:0.0] plot(\x,{(-0--2.6*\x)/-1.5});
\draw [color = gray, line width = 0.2pt,domain=-1.5000000000000004:17.823872468411377] plot(\x,{(-0-2.6*\x)/1.5});
\draw [color = gray, line width = 0.2pt,domain=-6.230245364561282:3.0] plot(\x,{(-15.59--5.2*\x)/-3});
\draw [color = gray, line width = 0.2pt,domain=1.1102230246251565E-15:17.823872468411377] plot(\x,{(--15.59-5.2*\x)/3});
\draw [color = gray, line width = 0.2pt,domain=-6.230245364561282:4.5] plot(\x,{(-15.59--2.6*\x)/-1.5});
\draw [color = gray, line width = 0.2pt,domain=3.0:17.823872468411377] plot(\x,{(--15.59-2.6*\x)/1.5});
\draw [color = gray, line width = 0.2pt,domain=-6.230245364561282:4.5] plot(\x,{(-15.59-0*\x)/-6});
\draw [color = gray, line width = 0.2pt,domain=-1.5000000000000004:17.823872468411377] plot(\x,{(--15.59-0*\x)/6});
\draw [color = gray, line width = 0.2pt,domain=0.0:17.823872468411377] plot(\x,{(-0--5.2*\x)/3});
\draw [color = gray, line width = 0.2pt,domain=-6.230245364561282:3.0] plot(\x,{(-0-5.2*\x)/-3});
\draw [color = gray, line width = 0.2pt,domain=-1.5000000000000004:17.823872468411377] plot(\x,{(--7.79--2.6*\x)/1.5});
\draw [color = gray, line width = 0.2pt,domain=-6.230245364561282:1.1102230246251565E-15] plot(\x,{(-7.79-2.6*\x)/-1.5});
\draw [color = gray, line width = 0.2pt,domain=-6.230245364561282:4.5] plot(\x,{(--7.79-2.6*\x)/-1.5});
\draw [color = gray, line width = 0.2pt,domain=3.0:17.823872468411377] plot(\x,{(-7.79--2.6*\x)/1.5});
\draw [color = gray, line width = 0.2pt,domain=0.0:17.823872468411377] plot(\x,{(-0-0*\x)/3});
\draw [color = gray, line width = 0.2pt,domain=-6.230245364561282:3.0] plot(\x,{(-0-0*\x)/-6});
\draw [color = gray, line width = 0.2pt,domain=1.1102230246251565E-15:17.823872468411377] plot(\x,{(--15.59-0*\x)/3});
\draw [color = gray, line width = 0.2pt,domain=-6.230245364561282:3.0] plot(\x,{(-15.59-0*\x)/-3});
\fill[color=black,fill=red,fill opacity=0.2] (3,5.2) -- (1.5,2.6) -- (4.5,2.6) -- cycle;
\fill[color=black,fill=red,fill opacity=0.2] (1.5,2.6) -- (4.5,2.6) -- (3,0) -- cycle;
\draw [line width=1pt,color=blue] (1.5,2.6)-- (4.5,2.6);
\draw [line width=1pt,color=red!70] (3,5.2)-- (4.5,2.6)-- (3,0) -- (1.5,2.6) -- cycle;
\draw[color=blue] (3,2.4) node[below] {$\sigma$};
\end{tikzpicture}
			\caption{$\closedSt{\sigma}$}
			\label{subfigure:closed_star_of_1-face}
		\end{subfigure}
		\hfill
		\begin{subfigure}{0.3\textwidth}
			\makeatletter
			\ltx@ifclassloaded{mcom-l}{
			\begin{tikzpicture}[scale = 0.4, line cap=round,line join=round,>=triangle 45,x=1.0cm,y=1.0cm]
				}{%
				\begin{tikzpicture}[scale = 0.5, line cap=round,line join=round,>=triangle 45,x=1.0cm,y=1.0cm]
				}
				\makeatother
				\clip(-3,-2.5) rectangle (6,7.5);
\draw [color = gray, line width = 0.2pt,domain=-6.230245364561282:0.0] plot(\x,{(-0--2.6*\x)/-1.5});
\draw [color = gray, line width = 0.2pt,domain=-1.5000000000000004:17.823872468411377] plot(\x,{(-0-2.6*\x)/1.5});
\draw [color = gray, line width = 0.2pt,domain=-6.230245364561282:3.0] plot(\x,{(-15.59--5.2*\x)/-3});
\draw [color = gray, line width = 0.2pt,domain=1.1102230246251565E-15:17.823872468411377] plot(\x,{(--15.59-5.2*\x)/3});
\draw [color = gray, line width = 0.2pt,domain=-6.230245364561282:4.5] plot(\x,{(-15.59--2.6*\x)/-1.5});
\draw [color = gray, line width = 0.2pt,domain=3.0:17.823872468411377] plot(\x,{(--15.59-2.6*\x)/1.5});
\draw [color = gray, line width = 0.2pt,domain=-6.230245364561282:4.5] plot(\x,{(-15.59-0*\x)/-6});
\draw [color = gray, line width = 0.2pt,domain=-1.5000000000000004:17.823872468411377] plot(\x,{(--15.59-0*\x)/6});
\draw [color = gray, line width = 0.2pt,domain=0.0:17.823872468411377] plot(\x,{(-0--5.2*\x)/3});
\draw [color = gray, line width = 0.2pt,domain=-6.230245364561282:3.0] plot(\x,{(-0-5.2*\x)/-3});
\draw [color = gray, line width = 0.2pt,domain=-1.5000000000000004:17.823872468411377] plot(\x,{(--7.79--2.6*\x)/1.5});
\draw [color = gray, line width = 0.2pt,domain=-6.230245364561282:1.1102230246251565E-15] plot(\x,{(-7.79-2.6*\x)/-1.5});
\draw [color = gray, line width = 0.2pt,domain=-6.230245364561282:4.5] plot(\x,{(--7.79-2.6*\x)/-1.5});
\draw [color = gray, line width = 0.2pt,domain=3.0:17.823872468411377] plot(\x,{(-7.79--2.6*\x)/1.5});
\draw [color = gray, line width = 0.2pt,domain=0.0:17.823872468411377] plot(\x,{(-0-0*\x)/3});
\draw [color = gray, line width = 0.2pt,domain=-6.230245364561282:3.0] plot(\x,{(-0-0*\x)/-6});
\draw [color = gray, line width = 0.2pt,domain=1.1102230246251565E-15:17.823872468411377] plot(\x,{(--15.59-0*\x)/3});
\draw [color = gray, line width = 0.2pt,domain=-6.230245364561282:3.0] plot(\x,{(-15.59-0*\x)/-3});
\draw [line width=1pt,color=blue] (1.5,2.6)-- (4.5,2.6);
\draw [line width=1pt,color=black] (3,5.2)-- (4.5,2.6)-- (3,0) -- (1.5,2.6) -- cycle;
\fill [color=red!70] (3,5.2) circle (6pt);
\fill [color=red!70] (3,0) circle (6pt);
\draw[color=blue] (3,2.4) node[below] {$\sigma$};
\end{tikzpicture}
			\caption{$\lk(\sigma)$}
			\label{subfigure:link_of_1-face}
		\end{subfigure}
	\end{center}
	\caption{Examples of star, closed star and link for a vertex (top row) and an edge (bottom row) of a simplicial complex; see \cref{subsection:geometric_simplicial_complexes}.}
	\label{fig:star_link}
\end{figure}
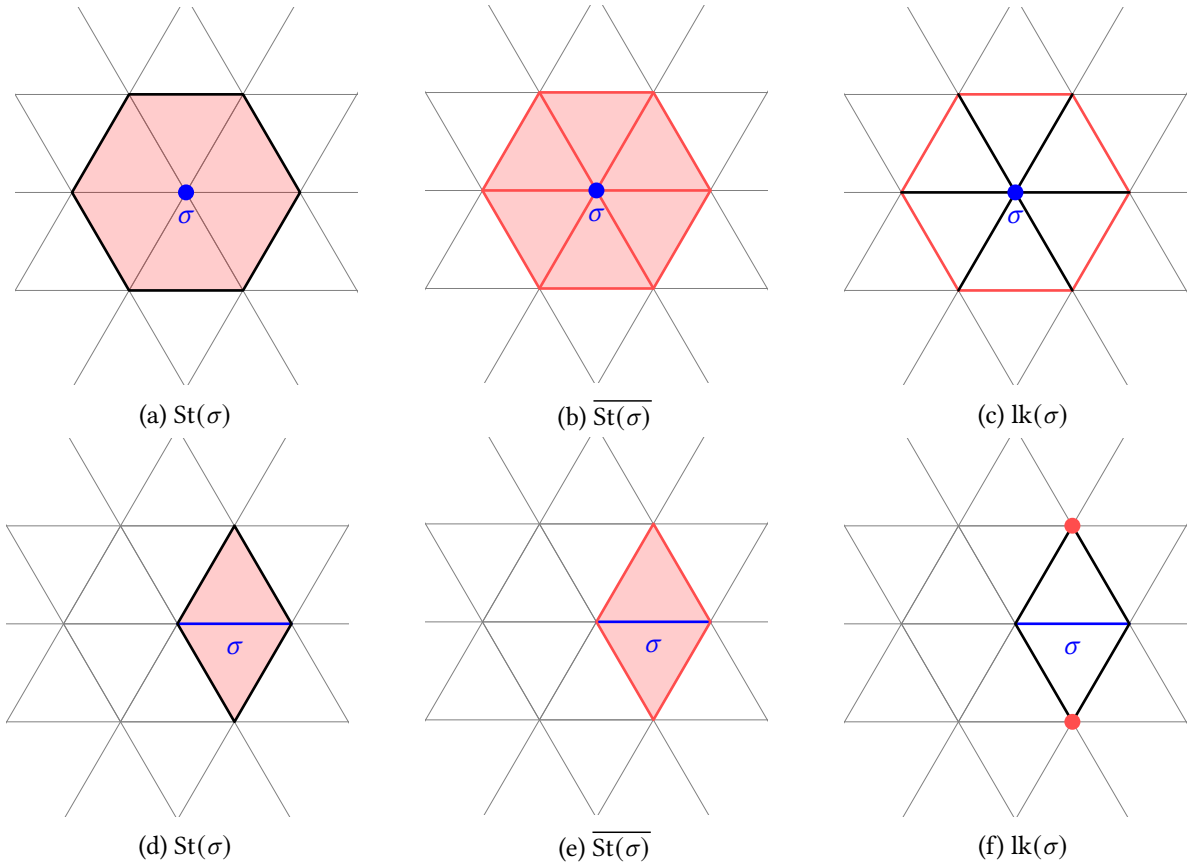

In this context we also consider the definition of star, closed star and link of a face~$\sigma$ in a simplicial complex~$\Sigma$.
The \textbf{star} of $\sigma$, denoted as $\St(\sigma)$, is the union of the interior of those simplices of $\Sigma$ which have $\sigma$ as a face.
The \textbf{closure} of a subset of the faces in $\Sigma$ can be defined as the smallest subcomplex of $\Sigma$ containing those faces.
Specifically, the closure of the star of $\sigma$, denoted as $\closedSt{\sigma}$, is called the \textbf{closed star} of $\sigma$ in $\Sigma$.
It can be understood as the union of all simplices of $\Sigma$ having $\sigma$ as a face, and it is a subcomplex of $\Sigma$.
The set $\closedSt{\sigma}\setminus \St(\sigma)$ is called the \textbf{link} of $\sigma$ in $\Sigma$ and it is denoted by $\lk(\sigma)$.
In other words, the link of $\sigma$ is the set of all simplices in the closure of $\closedSt{\sigma}$ which do not have $\sigma$ as a face.
It can be proved that for any face~$\sigma$ of $\Sigma$, $\lk(\sigma)$ is a subcomplex of $\Sigma$.
If $\Sigma$ is pure and of dimension $d$, then $\closedSt{\sigma}$ is also pure and of dimension $d$.
Furthermore, if $\sigma$ is of dimension $k$, then $\lk(\sigma)$ is pure and of dimension $d-k-1$.
\Cref{fig:star_link} shows examples of the star, closed star and link of a face $\sigma$ of a simplicial complex.
We concentrate on the cases relevant here, when $\sigma$ is a $0$-face or a $1$-face.

Another important notion is the distinction between boundary faces and interior faces.
Let $\Sigma$ be a simplicial $2$-complex.
We say that a $1$-face~$\sigma \in \Sigma$ is a \textbf{boundary $1$-face} if it belongs to exactly one $2$-face.
We denote the set of all boundary $1$-faces by $E_\partial$.
We say that a $2$-face is a \textbf{boundary $2$-face} if it contains at least one boundary $1$-face.
The set of all boundary $2$-faces will be denoted by $T_\partial$.
Finally, a $0$-face is called a \textbf{boundary $0$-face} if it belongs to at least one boundary $1$-face.
We denote the set of all boundary $0$-faces by $V_\partial$.

A $1$-face is said to be an \textbf{interior $1$-face} if it belongs to exactly two $2$-faces.
A $2$-face is said to be an \textbf{interior $2$-face} if all of its $1$-faces are interior.
Finally, a $0$-face is said to be an \textbf{interior $0$-face} if all $1$-faces it belongs to are interior.
Notice that in a simplicial $2$-complex, all $0$-, $1$- and $2$-faces are either boundary or interior.

Now, we focus on some results regarding closed stars of interior $0$-faces.
We start by defining a \textbf{polygonal chain}, which is a connected series of edges.
It can also be understood as a curve which is specified by a finite sequence of points called the vertices of the chain.
Furthermore, a \textbf{simple polygonal chain} is one in which only consecutive segments intersect and they intersect only at their end points.
In the same way, a \textbf{closed polygonal chain} is one in which the first vertex coincides with the last one.

For the three lemmas which follow, we make the following joint assumption.

\begin{assumption}
	\label{assumption:geometric_simplicial_complexes}
	Suppose that $\Delta$ is a connectivity complex, \ie, a pure, $2$-path connected abstract simplicial $2$-complex; see \cref{definition:connectivity_complex}.
	Given a set of vertex positions $Q \in \zeromanifold$, we recall that $\Sigma_\Delta(Q)$ is a simplicial $2$-complex whose associated abstract simplicial complex is $\Delta$; see \cref{definition:zeromanifold}.
\end{assumption}

\begin{lemma}
	\label{lemma:link_polygonal_chain}
	If $q$ is any $0$-face of $\Sigma_\Delta(Q)$, then $\lk(q)$ is a simple polygonal chain.
\end{lemma}
\begin{proof}
	Since we are working with simplicial $2$-complexes and $\Delta$ is pure, we know from \cite[p.100]{Gallier2008} that $\lk(q)$ is a simplicial complex which is pure and of dimension~$1$, see also \Cref{subfigure:link_of_0-face}.
	Then, we can conclude that $\lk(q)$ is a sequence of line segments.
	We will also know that this sequence of line segments will be simple since $\lk(q)$ is a simplical complex.
	Then, it only remains to be proved is that this sequence of line segments is connected.
	To this end we proceed by contradiction and assume that the sequence of line segments is not connected.
	By definition of a connected set, we know that there exists at least one line segment which is isolated from the others.
	From the definition of $\lk(q)$ we know that for every $1$-face~$e$ in $\lk(q)$, there exists a $2$-face~$T$ in $\closedSt{q}$ such that $e \subset T$.
	In particular, if we call $\widehat{e}$ the isolated $1$-face and $\widehat{T}$ an associated $2$-face, then it is easy to see that there will not be a sequence of $2$-faces which will connect $\widehat{T}$ with any other $2$-face.
	This contracts the assumption of $\Delta$ being $2$-path connected.
\end{proof}

\begin{lemma}
	\label{lemma:link_closed_polygonal_chain}
	If $q$ is an interior $0$-face, then $\lk(q)$ is a simple closed polygonal chain.
\end{lemma}
\begin{proof}
	From \cref{lemma:link_polygonal_chain}, we know that $\lk(q)$ is a simple polygonal chain.
	It remains to prove that it is closed.
	We proceed by contradiction and assume that the first vertex is different from the last.
	If we denote as $\widehat{e}$ the first line segment of the polygonal chain, we know that by definition of the link of a $0$-face, there exists $\widehat{T} \in \closedSt{q}$, which is the 2-face uniquely identified by $\widehat{e}$ and $q$.
	Now, we notice the $1$-face joining $q$ and $\widehat{e}$ belongs to only one $2$-face (namely $\widehat{T}$), which by definition means it is a boundary $1$-face and therefore $q$ will be a boundary $0$-face, which produces a contradiction.
\end{proof}

\begin{lemma}
	\label{lemma:separating_hyperplanes}
	Suppose that $q$ is an interior $0$-face.
	Let us denote by $e^{\lk}_i$ the $1$-faces which belong to $\lk(q)$.
	Finally, let us consider the set
	\begin{equation*}
		C = \bigcap_{i\in I} H^+(e_i^{\lk}),
	\end{equation*}
	where $I \coloneqq \setDef{i \in \N}{e_i \in \lk(q)}$ and $H^+(e_i)$ denotes the half-space generated by the $1$-face $e_i$ which contains $\closedSt{q}$.
	Then, $C \subset \closedSt{q}$ holds.
\end{lemma}
\begin{proof}
	We proceed by contraposition, \ie, we assume $x \notin \closedSt{q}$ and we wish to prove $x \notin C$.
	By definition of $\closedSt{q}$, $x\notin \closedSt{q}$ implies that for all $T_i \in \closedSt{q}$, then $x \notin T_i$.
	Now, since $q$ is an interior $0$-face and thanks to \cref{lemma:link_closed_polygonal_chain}, we know $\lk(q)$ is a simple closed polygonal chain, which means there exists and index $i \in I$ such that $x \in H^-(e^{\lk}_i)$, which implies $x \notin C$.
\end{proof}

\subsection{Abstract Simplicial Complexes}
\label{subsection:abstract_simplicial_complexes}

A (finite) \textbf{abstract simplicial complex}~$\Delta$ over a finite set~$V$ is a non-empty collection of non-empty subsets of $V$ such that, for all $\sigma \in \Delta$, every non-empty subset of~$\sigma$ also belongs to $\Delta$.
The elements of $\sigma \in \Delta$ are called the \textbf{faces} of~$\Delta$.
A face~$\sigma$ is said to be of \textbf{dimension}~$m \in \N_0$ (an \textbf{$m$-face}) if $\#\sigma = m+1$.
An abstract simplicial complex~$\Delta$ is said to be of \textbf{dimension}~$k \in \N_0$ (an \textbf{abstract simplicial $k$-complex}) if $k$ is the largest dimension of any of its faces.
The \textbf{vertex set} of an abstract simplicial complex~$\Delta$ is the union of all of its faces, \ie, $\bigcup_{\sigma \in \Delta} \sigma \subset V$.
The elements of the vertex set are the \textbf{vertices}.
Clearly, an abstract simplicial complex over a finite set~$V$ is also an abstract simplicial complex over its vertex set.

An abstract simplicial $k$-complex~$\Delta$ is \textbf{pure} if all maximal elements of $\Delta$ (\wrt\ the partial order of set inclusion) have dimension~$k$.
An abstract simplicial $k$-complex~$\Delta$ is said to be \textbf{$m$-path connected} ($1 \le m \le k$) if for any two distinct $m$-faces $\sigma, \sigma'$ in $\Delta$, there exists a finite sequence of $m$-faces, starting in $\sigma_0 = \sigma$ and ending in $\sigma_n = \sigma'$, such that $\sigma_i \cap \sigma_{i+1}$ is an $(m-1)$-face for $i = 0, \ldots, n-1$.

Abstract simplicial complexes provide a purely combinatorial description of geometric simplicial complexes, disregarding actual vertex \eqq{positions}.
Every geometric simplicial complex~$\Sigma$ defines an abstract simplicial complex~$\Delta$ (unique up to homomorphisms, \ie, renaming vertices) of the same dimension, as follows.
Suppose that $v_1, v_2, \ldots, v_{N_V}$ are the vertices of~$\Sigma$.
Define $\Delta$ over the vertex set $\{1, \ldots, N_V\}$ as
\begin{equation*}
	\Delta
	\coloneqq
	\setDef[big]{ \sigma \subset \{1, \ldots, N_V\}}{\conv\{v_i\}_{i \in \sigma} \in \Sigma}
	.
\end{equation*}
We call $\Delta$ the \textbf{abstract simplicial complex associated} with $\Sigma$.

Finally, it is worth mentioning that the notions of \textbf{boundary faces} and \textbf{interior faces} of dimensions $0$, $1$ and $2$, of an abstract simplicial $2$-complex can be defined in the same way as in the geometric case.

\section{Basic Notions and Inequalities for Triangles}
\label{subsection:triangles}

In this section we collect some well-known results about triangles in $\R^2$, in an effort to keep this paper self-contained.
Most of the definitions were taken from~\cite{AgricolaFriedrich2008}.
In the language of geometry, let us suppose that $\Delta$ is a connectivity complex, \ie, a pure, $2$-path connected abstract simplicial $2$-complex, with vertex set~$V$; see \cref{definition:connectivity_complex}.
Given a set of vertex positions $Q \in \zeromanifold$, we recall that $\Sigma_\Delta(Q)$ is a simplicial $2$-complex whose associated abstract simplicial complex is $\Delta$; see \cref{definition:zeromanifold}.
Orientation is not relevant for the results in this section.

Let us consider a $2$-face of $\Delta$, denoted as $\{i_0,i_1,i_2\}$.
The associated geometric counterpart, \ie, the triangle $\conv\{Q_{i_0},Q_{i_1},Q_{i_2}\}$, is an element of $\Sigma_\Delta(Q)$.
Its \textbf{edge lengths} are denoted as $\edgelength{Q}{\ell}[i_0,i_1,i_2]$, where the edge numbered~$\ell$ is the one opposite the vertex $Q_{i_\ell}$ for $\ell = 0, 1, 2$.
We denote the interior \textbf{angles of the triangle} by $\interiorangle{Q}{\ell}[i_0,i_1,i_2]$, where the angle numbered~$\ell$ is the one at the vertex $Q_{i_\ell}$ for $\ell = 0, 1, 2$.
The $\ell$-th \textbf{height} is denoted as $\height{Q}{\ell}[i_0,i_1,i_2]$ and it is the length of the line segment perpendicular to the edge $\{i_{\ell \oplus 1}, i_{\ell \oplus 2}\}$ which passes through the vertex $q_{i_\ell}$.
Here $\oplus$ denotes addition modulo~3.
The \textbf{inradius} is the radius of the largest circle that fits inside the triangle and we denote it as $\inradius{Q}[i_0,i_1,i_2]$.
The \textbf{circumradius} is the radius of the smallest circle into which the triangle will fit and we denote it as $\circumradius{Q}[i_0,i_1,i_2]$.

The \textbf{area} of the triangle is denoted as $\area{Q}[i_0,i_1,i_2]$ and it satisfies the identity
\begin{equation}
	\label{eq:definition_area}
	\area{Q}[i_0,i_1,i_2]
	=
	\frac{\edgelength{Q}{\ell}[i_0,i_1,i_2] \, \height{Q}{\ell}[i_0,i_1,i_2]}
	{2}
	=
	\semiperimeter{Q}[i_0,i_1,i_2]
	\,
	\inradius{Q}[i_0,i_1,i_2]
	,
\end{equation}
where $\semiperimeter{Q}[i_0,i_1,i_2]$ is the semi-perimeter, \ie,
\begin{equation}
	\label{eq:semi-perimeter}
	\semiperimeter{Q}[i_0,i_1,i_2]
	\coloneqq
	\frac{1}{2} \paren[auto][]{%
	\edgelength{Q}{0}[i_0,i_1,i_2]
	+
	\edgelength{Q}{1}[i_0,i_1,i_2]
	+
	\edgelength{Q}{2}[i_0,i_1,i_2]
	}
.
\end{equation}
Rearranging terms in \eqref{eq:definition_area}, we have that
\begin{equation*}
	\height{Q}{\ell}[i_0,i_1,i_2]
	=
	2 \, \inradius{Q}[i_0,i_1,i_2]
	\,
	\dfrac{\semiperimeter{Q}[i_0,i_1,i_2]}{\edgelength{Q}{\ell}[i_0,i_1,i_2]}
	.
\end{equation*}
Since $\semiperimeter{Q}[i_0,i_1,i_2] > \max \paren[big]\{\}{\edgelength{Q}{0}[i_0,i_1,i_2],\edgelength{Q}{1}[i_0,i_1,i_2],\edgelength{Q}{2}[i_0,i_1,i_2]}$ holds, we have
\begin{equation}
	\label{eq:height_and_inradius}
	\height{Q}{\ell}[i_0,i_1,i_2]
	>
	2 \, \inradius{Q}[i_0,i_1,i_2]
	.
\end{equation}
Notice also that the area \eqref{eq:definition_area} agrees with the signed area defined in \eqref{eq:signed_area} if $[i_0,i_1,i_2]$ is properly oriented.
We will not need to use different symbols for area and signed area.

The well-known relation $\circumradius{Q}[i_0,i_1,i_2]\ge 2 \, \inradius{Q}[i_0,i_1,i_2]$ between the inradius and circumradius can be found in \cite[p.198]{SvrtanVeljan2012}.
Another useful relation between the inradius and the heights is given by
\begin{equation}
	\label{eq:relation_inradius_heights}
	\frac{1}{\inradius{Q}[i_0,i_1,i_2]}
	=
	\frac{1}{\height{Q}{0}[i_0,i_1,i_2]}
	+
	\frac{1}{\height{Q}{1}[i_0,i_1,i_2]}
	+
	\frac{1}{\height{Q}{2}[i_0,i_1,i_2]}
	,
\end{equation}
see, \eg, \cite[p.353]{Kay2011}.
From \cite[Cor.~3]{Birsan2015} we can obtain the following bound on the interior angles,
\makeatletter
\ltx@ifclassloaded{mcom-l}{%
	\begin{multline}
		\label{eq:bound_angles_triangle}
		\frac{\inradius{Q}[i_0,i_1,i_2]}{\circumradius{Q}[i_0,i_1,i_2]}-\sqrt{1-\frac{2 \, \inradius{Q}[i_0,i_1,i_2]}{\circumradius{Q}[i_0,i_1,i_2]}}
		\le
		\cos(\theta^\ell(i_0,i_1,i_2))
		\\
		\le
		\frac{\inradius{Q}[i_0,i_1,i_2]}{\circumradius{Q}[i_0,i_1,i_2]}+\sqrt{1-\frac{2 \, \inradius{Q}[i_0,i_1,i_2]}{\circumradius{Q}[i_0,i_1,i_2]}}
		.
	\end{multline}
	}{%
	\begin{equation}
		\label{eq:bound_angles_triangle}
		\frac{\inradius{Q}[i_0,i_1,i_2]}{\circumradius{Q}[i_0,i_1,i_2]}-\sqrt{1-\frac{2 \, \inradius{Q}[i_0,i_1,i_2]}{\circumradius{Q}[i_0,i_1,i_2]}}
		\le
		\cos(\theta^\ell(i_0,i_1,i_2))
		\le
		\frac{\inradius{Q}[i_0,i_1,i_2]}{\circumradius{Q}[i_0,i_1,i_2]}+\sqrt{1-\frac{2 \, \inradius{Q}[i_0,i_1,i_2]}{\circumradius{Q}[i_0,i_1,i_2]}}
		.
	\end{equation}
}
\makeatother

Finally, we will use the characterization of the circumradius given in \cite[Thm.~27, p.43]{AgricolaFriedrich2008},
\begin{equation}
	\label{eq:circumradius}
	\circumradius{Q}[i_0,i_1,i_2]
	=
	\frac{\edgelength{Q}{0}[i_0,i_1,i_2] \, \edgelength{Q}{1}[i_0,i_1,i_2] \, \edgelength{Q}{2}[i_0,i_1,i_2]}{4 \area{Q}[i_0,i_1,i_2]}
	.
\end{equation}

\section{Proofs}
\label{section:proofs}

\subsection{Proof of \texorpdfstring{\cref{lemma:bounds_of_triangle_properties_in_terms_of_f}}{Lemma \ref{lemma:bounds_of_triangle_properties_in_terms_of_f}}}
\label{subsection:proof_of_lemma:bounds_of_triangle_properties_in_terms_of_f}

As in the statement of \cref{lemma:bounds_of_triangle_properties_in_terms_of_f}, let $\Delta$ be a consistently oriented connectivity complex (\cref{definition:connectivity_complex,definition:consistently_orientated_abstract_simplicial_2-complex}) with vertex set $V = \{1, \ldots, N_V\}$.
	Furthermore, let $[i_0,i_1,i_2]$ be an arbitrary $2$-face of $\Delta$.
Finally, let us assume $Q \in \planarmanifold$.

Then clearly, all heights are strictly positive and we have
\begin{equation*}
	\frac{\beta_1}{\height{Q}{\ell}[i_0,i_1,i_2]}
	\le
	\sum_{k=1}^{N_T} \sum_{\ell=0}^2 \frac{\beta_1}{\height[auto]{Q}{\ell}[i^k_0,i^k_1,i^k_2]}
	\le
	f(Q)
	,
\end{equation*}
which proves \cref{item:bound_height}.
Denoting by $\inradius{Q}[i_0,i_1,i_2]$ the inradius of a $2$-face of $\Delta$ and using \eqref{eq:relation_inradius_heights} we obtain that
\begin{equation*}
	\frac{\beta_1}{\inradius{Q}[i_0,i_1,i_2]}
	\le
	\sum_{k=1}^{N_T}\sum_{\ell=0}^2 \frac{\beta_1}{\height[auto]{Q}{\ell}[i^k_0,i^k_1,i^k_2]}
	\le
	f(Q)
	,
\end{equation*}
which proves \cref{item:bound_inradius}.

It is also well known that for any $2$-face the length of each of its $1$-faces is greater than twice the inradius, from where we obtain the first inequality of \cref{item:bound_length_edges}:
\begin{equation*}
	\frac{2 \beta_1}{f(Q)}
	\le
	2 \, \inradius{Q}[i_0,i_1,i_2]
	\le
	\edgelength{Q}{\ell}[i_0,i_1,i_2]
	.
\end{equation*}
For the second inequality of \cref{item:bound_length_edges}, we use the triangle inequality and the definition of~$f$ to obtain
\begin{equation*}
	\begin{aligned}
		\edgelength[big]{Q}{\ell}[i^k_0,i^k_1,i^k_2]
		&
		=
		\norm[big]{Q_{i_{\ell \oplus 1}^k} - Q_{i_{\ell \oplus 2}^k}}
		\le
		\norm[big]{Q_{i_{\ell \oplus 1}^k}} + \norm[big]{Q_{i_{\ell \oplus 2}^k}}
		\le
		\sqrt{2} \, \norm{Q}_F
		\\
		&
		\le
		\sqrt{2} \, \norm{Q-\Qref}_F + \sqrt{2} \, \norm{\Qref}_F
		\le
		2 \, \sqrt{\frac{f(Q)}{\beta_3}} + \sqrt{2} \, \norm{\Qref}_F
		,
	\end{aligned}
\end{equation*}
which completes the proof of \cref{item:bound_length_edges}.
Next we use that the area of each $2$-face is larger than the area of its incircle, \ie, $\area{Q}[i_0,i_1,i_2]\ge \pi r^2_Q(i_0,i_1,i_2)$.
Using the bound in \cref{item:bound_inradius}, we obtain \cref{item:bound_area}.
\Cref{item:bound_ratio_r_R} follows immediately from \eqref{eq:circumradius} and \cref{item:bound_length_edges,item:bound_inradius,item:bound_area}.
To prove \cref{item:bound_cosines}, we use \eqref{eq:bound_angles_triangle} to obtain
\begin{equation*}
	\abs[auto]{\cos \paren[auto](){\theta^\ell(i_0,i_1,i_2)}}
	\le
	\frac{\inradius{Q}[i_0,i_1,i_2]}{\circumradius{Q}[i_0,i_1,i_2]}
	+ \sqrt{1-\dfrac{2 \, \inradius{Q}[i_0,i_1,i_2]}{\circumradius{Q}[i_0,i_1,i_2]}}
\end{equation*}
for any $\ell = 0,1,2$.
Consider the function
\begin{equation*}
	\varphi(x)
	\coloneqq
	x+\sqrt{1-2x}
\end{equation*}
in the interval $[0,\frac{1}{2}]$, where it is continuous and decreasing and takes values in $[0,1]$.
Its maximum value is $\varphi(0) = 1$.
Thanks to \cref{item:bound_ratio_r_R}, we know that
\begin{equation*}
	\frac{\inradius{Q}[i_0,i_1,i_2]}{\circumradius{Q}[i_0,i_1,i_2]}
	\ge
	\frac{4 \pi \beta_1^3}{\paren[big]\{\}{2 \sqrt{\frac{f(Q)}{\beta_3}} + \sqrt{2} \, \norm{\Qref}_F}^3} \frac{1}{f(Q)^3}
	\eqqcolon
	\psi(f(Q))
	.
\end{equation*}
The function $\psi(x) \colon (0,\infty) \to (0,\infty)$ is continuous and decreasing.
Since $\varphi$ is also decreasing, we get
\begin{equation*}
	\abs[auto]{\cos(\theta^\ell(i_0,i_1,i_2))}
	\le
	\varphi\paren[auto](){\frac{\inradius{Q}[i_0,i_1,i_2]}{\circumradius{Q}[i_0,i_1,i_2]}}
	\le
	\varphi(\psi(f(Q)))
	.
\end{equation*}
Now set $\Psi \coloneqq \varphi \circ \psi \colon (0,\infty) \to (0,1)$, which is increasing and continuous as claimed.

\subsection{Proof of \texorpdfstring{\Cref{proposition:bounds_of_distances_in_terms_of_f}}{Proposition \ref{proposition:bounds_of_distances_in_terms_of_f}}}
\label{subsection:proof_of_proposition:bounds_of_distances_in_terms_of_f}

As in the statement of \Cref{proposition:bounds_of_distances_in_terms_of_f}, let $\Delta$ be a consistently oriented connectivity complex and $Q \in \plusmanifold$.
We consider a $0$-face $[i_0]$ and a $1$-face $[j_0,j_1]$ of $\Delta$ such that $\{Q_{i_0}\} \cap \conv\{Q_{j_0},Q_{j_1}\} = \emptyset$.
To shorten and simplify notation, we will simply write $q = Q_{i_0}$ and $e = \conv\{Q_{j_0},Q_{j_1}\}$.
We also use $\Distance{}[q][e]$ to denote $\Distance{Q}[i_0][[j_0,j_1]]$ and $\distance{}[q][e]$ to denote $\distance{Q}[i_0][[j_0,j_1]]$.
By assumption, we have $q \cap e = \emptyset$.
Therefore, the distance $\Distance{}[q][e]$ defined in \eqref{eq:distance_vertex_edge} is strictly positive.
Also recall the definition of $f$ from \eqref{eq:f}.
Finally, we are going to use the notation introduced in the end of \cref{section:geometric_abstract_simplicial_complexes}, specifically the link~$\lk(\cdot)$, and the closed star~$\closedSt{\cdot}$ of a $0$- or $1$-face.

\begin{figure}[htb]
	\begin{center}
		\begin{subfigure}{0.3\textwidth}
			\makeatletter
			\ltx@ifclassloaded{mcom-l}{
			\begin{tikzpicture}[scale = 0.35, line cap=round,line join=round,>=triangle 45,x=1.0cm,y=1.0cm]
				}{%
				\begin{tikzpicture}[scale = 0.45, line cap=round,line join=round,>=triangle 45,x=1.0cm,y=1.0cm]
				}
				\makeatother
\clip(-3,-6) rectangle (6,7.5);
\draw [color = gray, line width = 0.2pt, domain=-10.215621184943407:0.0] plot(\x,{(-0--2.6*\x)/-1.5});
\draw [color = gray, line width = 0.2pt, domain=-1.5000000000000004:20.99842721982095] plot(\x,{(-0-2.6*\x)/1.5});
\draw [color = gray, line width = 0.2pt, domain=-10.215621184943407:3.0] plot(\x,{(-15.59--5.2*\x)/-3});
\draw [color = gray, line width = 0.2pt, domain=1.1102230246251565E-15:20.99842721982095] plot(\x,{(--15.59-5.2*\x)/3});
\draw [color = gray, line width = 0.2pt, domain=-10.215621184943407:4.5] plot(\x,{(-15.59--2.6*\x)/-1.5});
\draw [color = gray, line width = 0.2pt, domain=3.0:20.99842721982095] plot(\x,{(--15.59-2.6*\x)/1.5});
\draw [color = gray, line width = 0.2pt, domain=-10.215621184943407:4.5] plot(\x,{(-15.59-0*\x)/-6});
\draw [color = gray, line width = 0.2pt, domain=-1.5000000000000004:20.99842721982095] plot(\x,{(--15.59-0*\x)/6});
\draw [color = gray, line width = 0.2pt, domain=0.0:20.99842721982095] plot(\x,{(-0--5.2*\x)/3});
\draw [color = gray, line width = 0.2pt, domain=-10.215621184943407:3.0] plot(\x,{(-0-5.2*\x)/-3});
\draw [color = gray, line width = 0.2pt, domain=-1.5000000000000004:20.99842721982095] plot(\x,{(--7.79--2.6*\x)/1.5});
\draw [color = gray, line width = 0.2pt, domain=-10.215621184943407:1.1102230246251565E-15] plot(\x,{(-7.79-2.6*\x)/-1.5});
\draw [color = gray, line width = 0.2pt, domain=-10.215621184943407:4.5] plot(\x,{(--7.79-2.6*\x)/-1.5});
\draw [color = gray, line width = 0.2pt, domain=3.0:20.99842721982095] plot(\x,{(-7.79--2.6*\x)/1.5});
\draw [color = gray, line width = 0.2pt, domain=0.0:20.99842721982095] plot(\x,{(-0-0*\x)/3});
\draw [color = gray, line width = 0.2pt, domain=-10.215621184943407:3.0] plot(\x,{(-0-0*\x)/-6});
\draw [color = gray, line width = 0.2pt, domain=1.1102230246251565E-15:20.99842721982095] plot(\x,{(--15.59-0*\x)/3});
\draw [color = gray, line width = 0.2pt, domain=-10.215621184943407:3.0] plot(\x,{(-15.59-0*\x)/-3});
\draw [color = gray, line width = 0.2pt, domain=-10.22:21] plot(\x,{(-60.75-0*\x)/23.38});
\draw [color = gray, line width = 0.2pt, domain=-10.22:21] plot(\x,{(-15.59--2.6*\x)/1.5});
\draw [color=red, line width = 1.pt] (0,0)-- (3,0)-- (4.5,2.6) -- (3,5.2) -- (0,5.2)-- (-1.5,2.6) -- cycle;
\fill[color=red!70,fill=red,fill opacity=0.2] (0,0) -- (3,0) -- (4.5,2.6) -- (3,5.2) -- (0,5.2) -- (-1.5,2.6) -- cycle;
\fill [color=blue] (1.5,2.6) circle (6pt);
\draw[color=blue] (1.5,2.4) node[below] {$q$};
\draw [line width=1pt,color=blue] (3,-5.2)-- (4.5,-2.6);
\draw[color=blue] (4,-4) node[right] {$e$};
\draw[color=darkgreen] (1.5,-0.1) node[below] {$e_j^{\lk}$};
\draw [color=darkgreen, line width = 1pt, domain=-10.22:21] plot(\x,{(-0-0*\x)/3});
\makeatletter
\ltx@ifclassloaded{mcom-l}{
\draw[color=darkgreen] (4.3, 5) node[above] {$H^+(e_j^{\lk})$};
}{%
\draw[color=darkgreen] (-1.8, 5) node[above] {$H^+(e_j^{\lk})$};
}
\makeatother
\end{tikzpicture}
			\caption{$q$ is interior, $e$ arbitrary and $\closedSt{q}\cap e =\emptyset$.}
			\label{fig:construction_q_interior_empty_intersection}
		\end{subfigure}
		\begin{subfigure}{0.3\textwidth}
			\makeatletter
			\ltx@ifclassloaded{mcom-l}{
			\begin{tikzpicture}[scale = 0.35, line cap=round,line join=round,>=triangle 45,x=1.0cm,y=1.0cm]
				}{%
				\begin{tikzpicture}[scale = 0.45, line cap=round,line join=round,>=triangle 45,x=1.0cm,y=1.0cm]
				}
				\makeatother
\clip(-3,-6) rectangle (6,7.5);
\draw [color = gray, line width = 0.2pt, domain=-10.215621184943407:0.0] plot(\x,{(-0--2.6*\x)/-1.5});
\draw [color = gray, line width = 0.2pt, domain=-1.5000000000000004:20.99842721982095] plot(\x,{(-0-2.6*\x)/1.5});
\draw [color = gray, line width = 0.2pt, domain=-10.215621184943407:3.0] plot(\x,{(-15.59--5.2*\x)/-3});
\draw [color = gray, line width = 0.2pt, domain=1.1102230246251565E-15:20.99842721982095] plot(\x,{(--15.59-5.2*\x)/3});
\draw [color = gray, line width = 0.2pt, domain=-10.215621184943407:4.5] plot(\x,{(-15.59--2.6*\x)/-1.5});
\draw [color = gray, line width = 0.2pt, domain=3.0:20.99842721982095] plot(\x,{(--15.59-2.6*\x)/1.5});
\draw [color = gray, line width = 0.2pt, domain=-10.215621184943407:4.5] plot(\x,{(-15.59-0*\x)/-6});
\draw [color = gray, line width = 0.2pt, domain=-1.5000000000000004:20.99842721982095] plot(\x,{(--15.59-0*\x)/6});
\draw [color = gray, line width = 0.2pt, domain=0.0:20.99842721982095] plot(\x,{(-0--5.2*\x)/3});
\draw [color = gray, line width = 0.2pt, domain=-10.215621184943407:3.0] plot(\x,{(-0-5.2*\x)/-3});
\draw [color = gray, line width = 0.2pt, domain=-1.5000000000000004:20.99842721982095] plot(\x,{(--7.79--2.6*\x)/1.5});
\draw [color = gray, line width = 0.2pt, domain=-10.215621184943407:1.1102230246251565E-15] plot(\x,{(-7.79-2.6*\x)/-1.5});
\draw [color = gray, line width = 0.2pt, domain=-10.215621184943407:4.5] plot(\x,{(--7.79-2.6*\x)/-1.5});
\draw [color = gray, line width = 0.2pt, domain=3.0:20.99842721982095] plot(\x,{(-7.79--2.6*\x)/1.5});
\draw [color = gray, line width = 0.2pt, domain=0.0:20.99842721982095] plot(\x,{(-0-0*\x)/3});
\draw [color = gray, line width = 0.2pt, domain=-10.215621184943407:3.0] plot(\x,{(-0-0*\x)/-6});
\draw [color = gray, line width = 0.2pt, domain=1.1102230246251565E-15:20.99842721982095] plot(\x,{(--15.59-0*\x)/3});
\draw [color = gray, line width = 0.2pt, domain=-10.215621184943407:3.0] plot(\x,{(-15.59-0*\x)/-3});
\draw [color = gray, line width = 0.2pt, domain=-10.22:21] plot(\x,{(-60.75-0*\x)/23.38});
\draw [color = gray, line width = 0.2pt, domain=-10.22:21] plot(\x,{(-15.59--2.6*\x)/1.5});
\draw [color=red, line width = 1.pt] (0,0)-- (3,0)-- (4.5,2.6) -- (3,5.2) -- (0,5.2)-- (-1.5,2.6) -- cycle;
\fill[color=red!70,fill=red,fill opacity=0.2] (0,0) -- (3,0) -- (4.5,2.6) -- (3,5.2) -- (0,5.2) -- (-1.5,2.6) -- cycle;
\fill [color=blue] (1.5,2.6) circle (6pt);
\draw[color=blue] (1.5,2.4) node[below] {$q$};
\draw [line width=1pt,color=blue] (0,0)-- (1.5,-2.6);
\draw[color=blue] (0.5,-1.5) node[left] {$e$};
\draw [line width=1pt,color=darkgreen] (0,0)-- (3,0);
\draw[color=darkgreen] (1.5,-0.1) node[below] {$e_j^{\lk}$};
\fill [color=darkgreen] (0,0) circle (6pt);
\draw[color=darkgreen] (-0.25,0.25) node[left] {$v$};
\end{tikzpicture}
			\caption{$q$ is interior, $e$ arbitrary and $\closedSt{q}\cap e \neq \emptyset$.}
			\label{fig:construction_q_interior_nonempty_intersection}
		\end{subfigure}
		\begin{subfigure}{0.3\textwidth}
			\makeatletter
			\ltx@ifclassloaded{mcom-l}{
			\begin{tikzpicture}[scale = 0.35, line cap=round,line join=round,>=triangle 45,x=1.0cm,y=1.0cm]
				}{%
				\begin{tikzpicture}[scale = 0.45, line cap=round,line join=round,>=triangle 45,x=1.0cm,y=1.0cm]
				}
				\makeatother
\clip(-3,-6) rectangle (8.5,7.5);
\draw [color = gray, line width = 0.2pt, domain=-10.201853767011194:0.0] plot(\x,{(-0--2.6*\x)/-1.5});
\draw [color = gray, line width = 0.2pt, domain=-1.5000000000000004:20.922887850724997] plot(\x,{(-0-2.6*\x)/1.5});
\draw [color = gray, line width = 0.2pt,domain=-10.201853767011194:3.0] plot(\x,{(-15.59--5.2*\x)/-3});
\draw [color = gray, line width = 0.2pt,domain=1.1102230246251565E-15:20.922887850724997] plot(\x,{(--15.59-5.2*\x)/3});
\draw [color = gray, line width = 0.2pt,domain=-10.201853767011194:4.5] plot(\x,{(-15.59--2.6*\x)/-1.5});
\draw [color = gray, line width = 0.2pt,domain=3.0:20.922887850724997] plot(\x,{(--15.59-2.6*\x)/1.5});
\draw [color = gray, line width = 0.2pt,domain=-10.201853767011194:4.5] plot(\x,{(-15.59-0*\x)/-6});
\draw [color = gray, line width = 0.2pt,domain=-1.5000000000000004:20.922887850724997] plot(\x,{(--15.59-0*\x)/6});
\draw [color = gray, line width = 0.2pt,domain=0.0:20.922887850724997] plot(\x,{(-0--5.2*\x)/3});
\draw [color = gray, line width = 0.2pt,domain=-10.201853767011194:3.0] plot(\x,{(-0-5.2*\x)/-3});
\draw [color = gray, line width = 0.2pt,domain=-1.5000000000000004:20.922887850724997] plot(\x,{(--7.79--2.6*\x)/1.5});
\draw [color = gray, line width = 0.2pt,domain=-10.201853767011194:1.1102230246251565E-15] plot(\x,{(-7.79-2.6*\x)/-1.5});
\draw [color = gray, line width = 0.2pt,domain=-10.201853767011194:4.5] plot(\x,{(--7.79-2.6*\x)/-1.5});
\draw [color = gray, line width = 0.2pt,domain=3.0:20.922887850724997] plot(\x,{(-7.79--2.6*\x)/1.5});
\draw [color = gray, line width = 0.2pt,domain=0.0:20.922887850724997] plot(\x,{(-0-0*\x)/3});
\draw [color = gray, line width = 0.2pt,domain=-10.201853767011194:3.0] plot(\x,{(-0-0*\x)/-6});
\draw [color = gray, line width = 0.2pt,domain=1.1102230246251565E-15:20.922887850724997] plot(\x,{(--15.59-0*\x)/3});
\draw [color = gray, line width = 0.2pt,domain=-10.201853767011194:3.0] plot(\x,{(-15.59-0*\x)/-3});
\draw [color = gray, line width = 0.2pt,domain=-10.2:20.92] plot(\x,{(-60.75-0*\x)/23.38});
\draw [color = gray, line width = 0.2pt,domain=-10.2:20.92] plot(\x,{(-15.59--2.6*\x)/1.5});
\draw [color = gray, line width = 0.2pt,domain=-10.2:20.92] plot(\x,{(--23.38-2.6*\x)/1.5});
\draw [line width = 1pt, color=red!70] (-1.5,2.6)-- (0,0)-- (3,0)-- (4.5,2.6)-- (3,5.2)-- (0,5.2)--cycle;
\draw [color=darkgreen,domain=-10.2:20.92] plot(\x,{(-7.79--2.6*\x)/1.5});
\draw[color=darkgreen](-0.7,-4) node{$H^+(e_j^{\lk})$};
\makeatletter
\ltx@ifclassloaded{mcom-l}{
\draw[color=darkgreen](5.3,1) node[left] {$e_j^{\lk}$};
}{%
\draw[color=darkgreen](5,1) node[left] {$e_j^{\lk}$};
}
\fill[fill=red,fill opacity=0.2] (-1.5,2.6) -- (1.5,2.6) -- (0,5.2) -- cycle;
\fill[fill=red,fill opacity=0.2] (-1.5,2.6) -- (1.5,2.6) -- (0,0) -- cycle;
\draw [line width=1pt,color=blue] (-1.5,2.6)-- (1.5,2.6);
\draw[color=blue](0,2.4) node[below] {$e$};
\fill [color=darkgreen] (1.5,2.6) circle (6pt);
\draw[color=darkgreen] (1.5,2.2) node[below] {$v_1$};
\fill [color=blue] (7.5,2.6) circle (6pt);
\draw[color=blue](7.5,2.4) node[below] {$q$};
\end{tikzpicture}
			\caption{$e$ is interior, $q$ boundary and $\closedSt{v_1}\cap q \neq \emptyset$.}
			\label{fig:construction_e_interior}
		\end{subfigure}
	\end{center}
	\caption{Illustrations for the proof of \cref{proposition:bounds_of_distances_in_terms_of_f}: $\closedSt{q}$ (left and center) and $\closedSt{e}$ (right) are shaded in red. Moreover, $\lk(q)$ (left and center) and $\lk(v_1)$ (right) are shown by dark red lines.}
\end{figure}
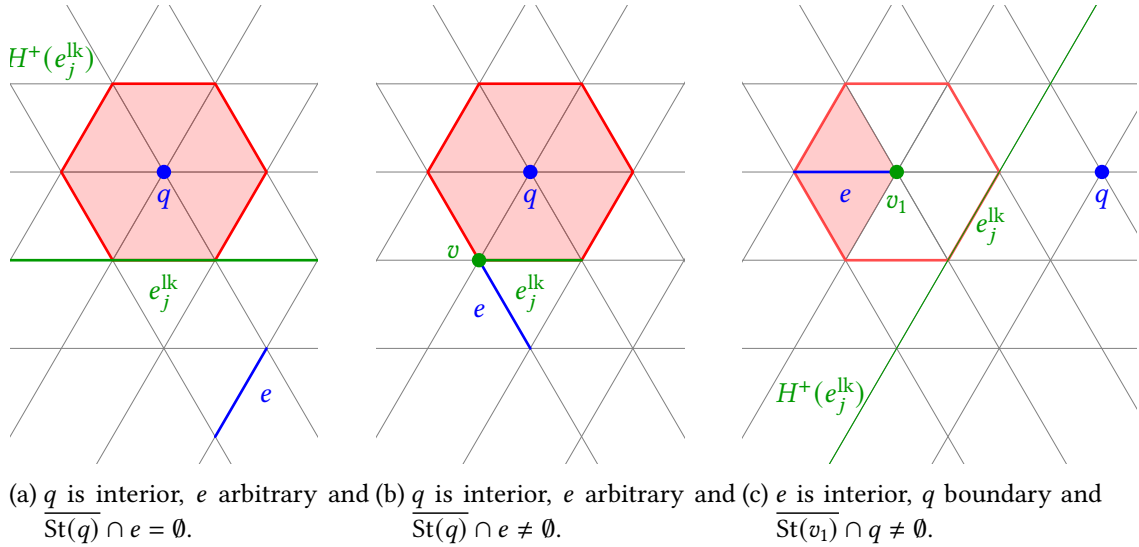

We consider the cases following the statement of \cref{proposition:bounds_of_distances_in_terms_of_f}.
\begin{description}
	\item[$q$ and $e$ are boundary faces:]
		We estimate
		\begin{equation*}
			\frac{\beta_2}{\Distance{}[q][e]}
			\le
			\sum_{[\ell_0,\ell_1]\in E_\partial}\sum_{\substack{k_0\in V_\partial \\ k_0\neq \ell_0,\ell_1}} \frac{\beta_2}{\Distance{Q}[k_0][[\ell_0,\ell_1]]}
			\le
			f(Q)
			,
		\end{equation*}
		which proves \eqref{eq:distance_boundary_vertex_boundary_edge} and thus \cref{item:distance_boundary_vertex_boundary_edge}.
	\item[$q$ is an interior $0$-face and $e$ is interior or boundary $1$-face:]
		$\bullet$
		We first consider the case $e \cap \closedSt{q} = \emptyset$.
		We use the notation described in \cref{lemma:separating_hyperplanes}, \ie, $e^{\lk}_i$ are the $1$-faces which belong to $\lk(q)$, indexed over $i \in I$.
		Moreover, $H^+(e^{\lk}_i)$ denotes the half-space generated by the $1$-face $e^{\lk}_i$ which contains $\closedSt{q}$.
		Using \cref{lemma:link_closed_polygonal_chain,lemma:separating_hyperplanes}, we know there exists a $1$-face $e^{\lk}_j$ which separates $e$ from $\closedSt{q}$ since $e \cap H^+(e^{\lk}_j) = \emptyset$.
		(See \Cref{fig:construction_q_interior_empty_intersection} for an example of such a mesh.)
		Consider the $2$-face uniquely identified by $e^{\lk}_j$ and vertex~$q$ and denote by $h_j$ the height of this $2$-face passing through $q$.
		Then,
		\begin{equation*}
			\distance{}[q][e^{\lk}_j]
			\ge
			h_j
			.
		\end{equation*}
		Since $e\cap H^+(e^{\lk}_j) = \emptyset$ holds, this implies
		\begin{equation*}
			\Distance{}[q][e]
			\ge
			\distance{}[q][e]
			>
			\distance{}[q][e^{\lk}_j]
			\ge
			h_j
			\ge
			\frac{\beta_1}{f(Q)}
			,
		\end{equation*}
		where the last inequality follows from \cref{lemma:bounds_of_triangle_properties_in_terms_of_f}, \cref{item:bound_height}.
		$\bullet$
		We now consider the case when $e \cap \closedSt{q} \neq \emptyset$.
		Since $q, e \in \Sigma_\Delta$ and $q \cap e = \emptyset$ by assumption, we have two possibilities.
		First, $e \cap \closedSt{q}$ is a vertex~$v$, and it belongs to some $e^{\lk}_j$.
		(See \Cref{fig:construction_q_interior_nonempty_intersection} for such a construction.)
		Using Pythagoras' Theorem we get
		\begin{equation*}
			\distance{}[q][e]
			=
			\norm{q-v}
			\ge
			\distance{}[q][e^{\lk}_j]
			\ge
			h_j
			,
		\end{equation*}
		which implies that $\Distance{}[q][e] \ge \distance{}[q][e] \ge \frac{\beta_1}{f(Q)}$ by \cref{lemma:bounds_of_triangle_properties_in_terms_of_f}, \cref{item:bound_height}.
		$\bullet$
		Second, $e\cap \closedSt{q}=e^{\lk}_j$.
		Notice that $\Distance{}[q][e] \ge \distance{}[q][e^{\lk}_j] \ge h_j$, where $h_j$ is the height of the $2$-face uniquely defined by $q$ and $e^{\lk}_j$, which passes through $q$.
		Thanks to \cref{lemma:bounds_of_triangle_properties_in_terms_of_f}, \cref{item:bound_height} we have $\Distance{}[q][e] \ge \frac{\beta_1}{f(Q)}$.

		Altogether, we have proved \eqref{eq:distance_interior_vertex_arbitrary_edge} and thus \cref{item:distance_interior_vertex_arbitrary_edge}.

	\item[$q$ is a boundary $0$-face and $e$ is an interior $1$-face:]
		$\bullet$
		We first consider the case $q \cap \closedSt{e} = \emptyset$.
		In order to estimate $\Distance{}[q][e]$, we are going to invoke \cref{lemma:bound_distance_boundary_vertex_interior_edge}, which separates the task into the estimation of angles and the estimation of $\norm{q-v_i}$, where $v_1$ and $v_2$ are the $0$-faces of~$e$.
		We begin with the latter and assume, without loss of generality, that $\norm{q-v_1} \le \norm{q-v_2}$.
		Therefore, we focus on the estimation of $\norm{q-v_1}$.
		We consider now two cases.

		\begin{description}
			\item[$v_1$ is an interior $0$-face:]
				Since $q \cap \closedSt{e} = \emptyset$, there exist two possibilities.
				$\bullet$
				Suppose first that $q \cap \closedSt{v_1} = \emptyset$ holds.
				Since $v_1$ is interior, then thanks to \cref{lemma:link_closed_polygonal_chain} we know that $\lk(v_1)$ is a closed polygonal chain.
				Therefore, there exists $e^{\lk}_j \in \lk(v_1)$ such that $q \cap H^+(e^{\lk}_j) = \emptyset$.
				(See \Cref{fig:construction_e_interior} for such a construction.)
				This implies $\norm{q-v_1} > \distance{}[v_1][e^{\lk}_j] \ge h_j\ge \beta_1/f(Q)$.
				The last inequality was obtained by \cref{lemma:bounds_of_triangle_properties_in_terms_of_f}, \cref{item:bound_height}
				$\bullet$
				Suppose now that $q \cap \closedSt{v_1} \neq \emptyset$ holds.
				Since $q,e \in \Sigma_\Delta$, we conclude $q \in \lk(v_1)$ and thus there exists an edge between $q$ and $v_1$.
				It thus follows from \cref{lemma:bounds_of_triangle_properties_in_terms_of_f}, \cref{item:bound_length_edges} that $\norm{q-v_1} \ge \frac{2 \beta_1}{f(Q)}$ holds.
				In both cases, we can conclude that $\norm{q-v_1} > \frac{\min\{\beta_1,\beta_2\}}{\sqrt{2} \, f(Q)}$.
				This result, together with \cref{lemma:bound_distance_boundary_vertex_interior_edge}, implies \eqref{eq:distance_boundary_vertex_interior_edge} as follows
				\begin{equation*}
					\Distance{}[q][e]
					\ge
					\distance{}[q][e]
					>
					\dfrac{\min\{\beta_1, \beta_2\}}{\sqrt{2}f(Q)} \min_{\theta\in\Theta}\paren[auto]\{\}{\sqrt{1-\cos^2(\theta)}}
					.
				\end{equation*}
				Thus \cref{item:distance_boundary_vertex_interior_edge} for the case that $v_1$ is an interior $0$-face.

			\item[$v_1$ is a boundary $0$-face:]
				$\bullet$
				We first consider the case $q \cap \closedSt{v_1} = \emptyset$.
				Thanks to the definition of boundary faces, we know that there exists a boundary $1$-face, denoted as $\widetilde e$, such that $v_1$ is one of its $0$-faces and $\norm{q-v_1}\ge \distance{}[q][\tilde e]$.
				Using \eqref{eq:equivalence_of_1-_and_2-norms}, and \eqref{eq:distance_boundary_vertex_boundary_edge} we get $\norm{q-v_1} \ge \Distance{}[q][\tilde e] / \sqrt{2} \ge \frac{\beta_2}{\sqrt{2} \, f(Q)}$.
				$\bullet$
				If on the other hand, $q\cap\closedSt{v_1}\neq \emptyset$, then by definition of $\closedSt{v_1}$, we know that there exists an edge between $q$ and $v_1$.
				It thus follows from \cref{lemma:bounds_of_triangle_properties_in_terms_of_f}, \cref{item:bound_length_edges} that $\norm{q-v_1} \ge \frac{2 \beta_1}{f(Q)}$ holds.
				We can summarize both cases as $\norm{q-v_1} \ge \frac{\min\{\beta_1, \beta_2\}}{\sqrt{2} \, f(Q)}$.
				This result, together with \cref{lemma:bound_distance_boundary_vertex_interior_edge}, implies \eqref{eq:distance_boundary_vertex_interior_edge} and thus \cref{item:distance_boundary_vertex_interior_edge} for the case that $v_1$ is an boundary $0$-face.
		\end{description}

		$\bullet$
		We end this proof by considering the case when $q\cap\closedSt{e}\neq \emptyset$.
		Since $q,e\in\Sigma_\Delta$, the only possibility is that $q \cap \closedSt{e} = q$ holds.
		By considering the $2$-face uniquely defined by $e$ and $q$, we have $\Distance{}[q][e] \ge \distance{}[q][e] \ge h$, where $h$ is the height of this $2$-face which passes through $q$.
		Using \cref{lemma:bounds_of_triangle_properties_in_terms_of_f}, \cref{item:bound_height} and the fact that $\sqrt{1-\cos^2(\theta_\ell)}<1$ for all $\ell=1,\dots,4$, we conclude
		\begin{equation*}
			\Distance{}[q][e]
			\ge
			\frac{\beta_1}{f(Q)} > \frac{\min\{\beta_1, \beta_2\}}{\sqrt{2}f(Q)} \min_{\theta\in\Theta}\paren\{\}{\sqrt{1-\cos^2(\theta)}}
			,
		\end{equation*}
		which proves \eqref{eq:distance_boundary_vertex_interior_edge} and thus \cref{item:distance_boundary_vertex_interior_edge} in this case.
\end{description}

The following auxiliary result was used in the proof of \cref{proposition:bounds_of_distances_in_terms_of_f} above.
It shows that under the assumptions of \cref{proposition:bounds_of_distances_in_terms_of_f}, the distance of a boundary $0$-face $[i_0]$ from an interior $1$-face $[j_0,j_1]$ can be bounded from below.
This bound is given in terms of the angles formed by the edge $[Q_{j_0},Q_{j_1}]$ and the adjacent $1$-faces which belong to the closed star of $[Q_{j_0},Q_{j_1}]$, and the minimum of $\norm{Q_{i_0}-Q_{j_0}}$ and $\norm{Q_{i_0}-Q_{j_1}}$.
To keep the notation consistent with the proof of \cref{proposition:bounds_of_distances_in_terms_of_f}, we recall that $Q\in\plusmanifold$ is given but arbitrary, and that the $0$ and $1$-faces of $\Sigma_\Delta(Q)$ under consideration are denoted as $q = Q_{i_0}$ and $e = \conv\{Q_{j_0},Q_{j_1}\} \eqqcolon \conv\{v_1,v_2\}$, respectively.

\begin{lemma}
	\label{lemma:bound_distance_boundary_vertex_interior_edge}
	Suppose that the assumptions of \cref{proposition:bounds_of_distances_in_terms_of_f} hold and that, specifically, $i_0$ is a boundary $0$-face and $[i_0,j_1]$ is an interior $1$-face.
	Assume $q \cap \closedSt{e} = \emptyset$, where $\closedSt{e}$ denotes the closed star of $e$.
	Then
	\begin{equation}
		\distance{}[q][e]
		\ge
		\min_{\theta\in\Theta}\paren[auto]\{\}{\sqrt{1-\cos^2(\theta)}}
		\min \paren[auto]\{\}{\norm{q-v_1}, \norm{q-v_2}}
		,
	\end{equation}
	where $\Theta$ is the set of four angles formed by the edge~$e$ and the adjacent $1$-faces belonging to $\closedSt{e}$; see \Cref{fig:construction_interior_edge_lemma_C1}.
\end{lemma}
\begin{proof}
	We start by noticing that thanks to the definition of interior $1$-faces, the closed star $\closedSt{e}$ contains exactly two $2$-faces.
	Let us denote by $p$ the orthogonal projection of $q$ onto the infinite line defined by the edge~$e$.
	We distinguish two cases.
	$\bullet$
	If $p$ does not belong to $e$, then by the definition of the distance from a vertex to an edge, we get $\distance{}[q][e] = \min\paren[auto]\{\}{\norm{q-v_1},\norm{q-v_2}}.$
	$\bullet$
	Conversely, if $p$ does belong to $e$, then \wolog we assume $\norm{q-v_1} \le \norm{q-v_2}$.
	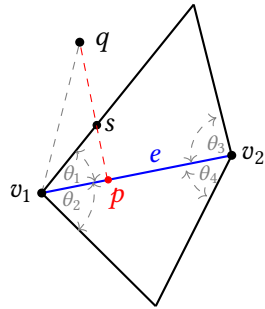
\begin{figure}[htb]
		\centering
		\begin{tikzpicture}[dot/.style={circle,inner sep=1.2pt,fill,name=#1},
  extended line/.style={shorten >=-#1,shorten <=-#1},
  extended line/.default=1cm]

\draw [blue, thick] (0.5,0.5) --  (3,1);
\draw [thick] (0.5,0.5) --  (2.5,3);
\draw [thick] (3,1) --  (2.5,3);
\draw [thick] (0.5,0.5) --  (2,-1);
\draw [thick] (2,-1) --  (3,1);
\draw[red,fill=red] (1.375,0.675) circle (.25ex);

\node at (2,1)     {$\textcolor{blue}{e}$};
\node at (0.9,0.75) {\scriptsize{\textcolor{gray}{$\theta_1$}}};
\node at (0.9,0.4)   {\scriptsize{\textcolor{gray}{$\theta_2$}}};
\node at (2.8,1.15) {\scriptsize{\textcolor{gray}{$\theta_3$}}};
\node at (2.7,0.75) {\scriptsize{\textcolor{gray}{$\theta_4$}}};
\node[left] at (0.5,0.5)  {$v_1$};
\node[right] at (3,1)     {$v_2$};
\node at (1.3,2.5)   {$q$};

\tkzDefPoint(0.5,0.5){O1}
\tkzDefPoint(1.2,0.65){A1}
\tkzDefPoint(0.9,1){B1}
\tkzDefPoint(0.975,0){B2}
\tkzDrawArc[dashed,arrows = <->](O1,A1)(B1)
\tkzDrawArc[dashed,arrows = <->](O1,B2)(A1)
\tkzDefPoint(3,1){O2}
\tkzDefPoint(1.9,0.8){A2}
\tkzDefPoint(2.8,1.5){C1}
\tkzDefPoint(2.7,0.5){C2}
\tkzDefPoint(2.4,0.8){C3}
\tkzDrawArc[dashed,arrows = <->](O2,C1)(A2)
\tkzDrawArc[dashed,arrows = <->](O2,C3)(C2)

\draw[fill] (1.225,1.4) circle (.3ex);
\node at (1.4,1.4)   {$s$};
\draw[dashed,gray] (0.5,0.5) -- (1,2.5);

\node [dot=v1] at (0.5,0.5) {};
\node [dot=v2] at (3,1) {};
\node [dot=q]   at (1,2.5) {};
\draw [red,dashed] ($(v1)!(q)!(v2)$) -- (q);
\node[below, color =red ] at (1.5,0.7)   {$p$};
\end{tikzpicture}
		\caption{Illustration of the construction for the distance between a $0$-face and an interior $1$-face, used in \cref{lemma:bound_distance_boundary_vertex_interior_edge}.}
		\label{fig:construction_interior_edge_lemma_C1}
	\end{figure}
	Furthermore, we denote by $s$ the unique point of intersection between the line segment joining $q$ and $p$ with one of edges of $\closedSt{e}$.
	 It is clear that $\distance{}[q][e] = \norm{q-p}$.
	 Now, by Pythagoras' Theorem we have that $\norm{q-p}^2 = \norm{q-v_1}^2 - \norm{v_1-p}^2$.
	 Moreover, we know that $\cos(\theta_1) = \norm{v_1-p} / \norm{v_1-s}$, which implies $\norm{q-p}^2=\norm{q-v_1}^2-\cos^2(\theta_1) \, \norm{v_1-s}^2$.

	Now, since we have assumed $q \cap \closedSt{e} = \emptyset$, we have $\norm{q-p} > \norm{s-p}$, and using again Pythagoras' Theorem we obtain
	\begin{align*}
		\norm{v_1-s}^2
		&
		=
		\norm{v_1-p}^2
		+
		\norm{p-s}^2
		\\
		&
		<
		\norm{v_1-p}^2
		+
		\norm{p-q}^2
		=
		\norm{v_1-q}^2
		.
	\end{align*}
	Thus,
	\begin{align*}
		\distance{}[q][e]^2
		&
		=
		\norm{q-p}^2
		>
		\norm{v_1-q}^2
		-
		\cos^2(\theta_1)
		\,
		\norm{v_1-q}^2
		\\
		&
		=
		(1-\cos^2(\theta_1))
		\,
		\norm{v_1-q}^2
		.
	\end{align*}
	Thanks to the assumption $\norm{q-v_1} \le \norm{q-v_2}$ we can conclude that
	\begin{align*}
		\distance{}[q][e]
		&
		>
		\sqrt{1-\cos^2(\theta_1)}
		\,
		\min
		\paren[auto]\{\}{\norm{q-v_1},\norm{q-v_2}}
		\\
		&
		\ge
		\min_{\theta\in\Theta}
		\paren[auto]\{\}{\sqrt{1-\cos^2(\theta)}}
		\,
		\min
		\paren[auto]\{\}{\norm{q-v_1},\norm{q-v_2}}
		.
	\end{align*}
	Notice that $\sqrt{1-\cos^2(\theta)}<1$ holds for all $\theta\in\Theta$, whether or not the orthogonal projection of $q$ onto the infinite line defined by $e$ belongs to $e$.
	We can thus conclude
	\begin{equation}
		\distance{}[q][e]
		\ge
		\min_{\theta\in\Theta}
		\paren[auto]\{\}{\sqrt{1-\cos^2(\theta)}}
		\min\paren[auto]\{\}{\norm{q-v_1},\norm{q-v_2}}
		.
	\end{equation}
\end{proof}

\section{An Example Regularization for \texorpdfstring{$\Distance{Q}[i_0][[j_0,j_1]]$}{the Vertex-Edge Distance}}
\label{section:example_regularization_vertex-edge_distance}

This section aims to present an example of a family of functions which can be used to approximate $f$ as suggested in \eqref{eq:f_mu}, using a regularized version of the distance $\Distance{Q}[i_0][[j_0,j_1]]$ described in \eqref{eq:distance_vertex_edge}.
We start by noticing that this distance from the $0$-face~$i_0$ to the $1$-face~$[j_0,j_1]$ can be written as
\begin{equation}
	\label{eq:equivalent_formulation_distance_vertex_edge}
	\Distance{Q}[i_0][[j_0,j_1]]
	=
	g\paren[big](){\widetilde{Q}^1_{i_0};[\widetilde{Q}^1_{j_0},\widetilde{Q}^1_{j_1}]} + \abs[big]{\widetilde{Q}^2_{i_0}-\widetilde{Q}^2_{j_0}}
	.
\end{equation}
Here $\widetilde{Q}^{\ell}_{\cdot}$ stands for the first ($\ell = 1$) or second ($\ell = 2$) component of a vector $Q_{\cdot}$ rotated about the origin so that its first coordinate aligns with the edge $\conv\{Q_{j_0},Q_{j_1}\}$, as shown in \Cref{fig:illustration_distance_vertex_edge}.
For convenience, the convention here is that $\widetilde{Q}^1_{j_0} < \widetilde{Q}^1_{j_1}$ holds.
Furthermore, the function $g$ is the distance of a point to an interval in $\R$, \ie,
\begin{equation}
	\label{eq:distance_point_interval}
	g(x;[y,z])
	=
	\begin{cases}
		\abs{x-y} & \text{if } y \ge x, \\
		0 & \text{if } y \le x \le z, \\
		\abs{x-z} & \text{otherwise}.
	\end{cases}
\end{equation}
We construct a regularizing function $\regularizedDistance{Q}[i_0][[j_0,j_1]]$ based on $C^3$ regularizations of \eqref{eq:distance_point_interval} and the absolute value function.

\begin{definition}
	Suppose that $\mu \ge 1$.
	We define the \textbf{regularized distance from a vertex to an edge} as follows
	\begin{equation}
		\label{eq:regularized_distance_vertex_edge}
		\regularizedDistance{Q}[i_0][[j_0,j_1]]
		=
		g^{\mu}\paren[big](){\widetilde{Q}^1_{i_0};[\widetilde{Q}^1_{j_0},\widetilde{Q}^1_{j_1}]} + h^\mu\paren[big](){\widetilde{Q}^2_{i_0}-\widetilde{Q}^2_{j_1}}
	\end{equation}
	where for $x,y,z\in \R$,
	\makeatletter
	\ltx@ifclassloaded{mcom-l}{%
		\begin{multline}
			\label{eq:regularized_distance_point_interval}
			g^{\mu}(x;[y,z])
			\\
			=
			\begin{cases}
				\abs{y-x}-\frac{1}{4\mu}
				&
				\text{if }
				\frac{1}{2\mu}\le y -x
				,
				\\
				32\mu^5(y-x)^6-48\mu^4(y-x)^5+20\mu^3(y-x)^4
				&
				\text{if }
				0 \le y-x \le \frac{1}{2\mu}
				,
				\\
				0
				&
				\text{if }
				y \le x \le z
				,
				\\
				32\mu^5(x-z)^6+48\mu^4(x-z)^5+20\mu^3(x-z)^4
				&
				\text{if }
				0 \le x-z \le \frac{1}{2\mu}
				,
				\\
				\abs{x-z}-\frac{1}{4\mu}
				&
				\text{if }
				x-z \ge \frac{1}{2\mu}
				.
			\end{cases}
		\end{multline}
		}{%
		\begin{equation}
			\label{eq:regularized_distance_point_interval}
			g^{\mu}(x;[y,z])
			=
			\begin{cases}
				\abs{y-x}-\frac{1}{4\mu}
				&
				\text{if }
				\frac{1}{2\mu}\le y -x
				,
				\\
				32\mu^5(y-x)^6-48\mu^4(y-x)^5+20\mu^3(y-x)^4
				&
				\text{if }
				0 \le y-x \le \frac{1}{2\mu}
				,
				\\
				0
				&
				\text{if }
				y \le x \le z
				,
				\\
				32\mu^5(x-z)^6+48\mu^4(x-z)^5+20\mu^3(x-z)^4
				&
				\text{if }
				0 \le x-z \le \frac{1}{2\mu}
				,
				\\
				\abs{x-z}-\frac{1}{4\mu}
				&
				\text{if }
				x-z \ge \frac{1}{2\mu}
				.
			\end{cases}
		\end{equation}
	}
	\makeatother

	and
	\begin{equation}
		\label{eq:Huber_C3}
		h^\mu(x)
		=
		\begin{cases}
			\abs{x}-\frac{1}{4\mu}
			&
			\text{if }
			\abs{x}\ge \frac{1}{2\mu}
			,
			\\
			40\mu^5\abs{x}^6-64\mu^4\abs{x}^5+32\mu^3\abs{x}^4-4\mu^2\abs{x}^3+\frac{\mu}{2}\abs{x}^2
			&
			\text{otherwise}.
		\end{cases}
	\end{equation}
\end{definition}

We now focus on proving that the proposed family of functions $\regularizedDistance{Q}$ satisfy the assumptions of \cref{theorem:f_mu_is_proper_and_Riemannian_metric_is_complete}.
Specifically, we verify $0\le \regularizedDistance{Q} \le \Distance{Q}$ in \cref{proposition:d_mu_is_an_underestimator} and argue that $Q \mapsto \regularizedDistance{Q}$ is of class $\cC^3$ in \cref{proposition:d_mu_is_of_class_C3}.
Moreover, \cref{proposition:f_mu_is_consistent} shows that in addition, the regularization is consistent, \ie, $f^\mu(Q) \to f(Q)$ holds when $\mu \to \infty$, for all $Q \in \planarmanifold$.

\begin{proposition}
	\label{proposition:d_mu_is_an_underestimator}
	For any $Q \in \planarmanifold$ and all $\mu \ge 1$, we have $0 \le \regularizedDistance{Q}[i_0][[j_0,j_1]] \le \Distance{Q}[i_0][[j_0,j_1]]$.
\end{proposition}
\begin{proof}
	We need to establish $0 \le h^\mu(x) \le \abs{x}$ and $0 \le g^\mu(x;[y,z]) \le g(x;[y,z])$ for all $x, y, z \in \R$ such that $y < z$.
	\\
	We start by proving $h^{\mu}(x) \le \abs{x}$.
	When $\abs{x} \ge \frac{1}{2\mu}$, then the first case in \eqref{eq:Huber_C3} applies and $h^{\mu}(x) \le \abs{x}$ is immediate.
	When $\abs{x} < \frac{1}{2\mu}$, we estimate
	\begin{align*}
		h^\mu(x)-\abs{x}
		&
		=
		40\mu^5\abs{x}^6-64\mu^4\abs{x}^5+32\mu^3\abs{x}^4-4\mu^2\abs{x}^3+\frac{\mu}{2}\abs{x}^2-\abs{x}
		\\
		&
		\le
		\abs{x} \paren[auto](){ \frac{40}{32} - \frac{64}{16} + \frac{32}{8} - \frac{4}{4} + \frac{1}{4} -1}
		=
		- \frac{1}{2} \abs{x}
		\le
		0
		.
	\end{align*}
	In the same way, $h^\mu(x)\ge 0$ is immediately obtained when $\abs{x}\ge \frac{1}{2\mu}$.
	When $\abs{x}\le \frac{1}{2\mu}$, notice that $h^\mu(x) = 40\mu^5 \abs{x}^6 + 32\mu^3\abs{x}^4 (1 - 2\mu\abs{x}) + 2\mu\abs{x}^2(1 - 2\mu\abs{x}) \ge 0$.
	Next we prove $0 \le g^\mu(x;[y,z]) \le g(x;[y,z])$.
	We focus on the second case in \eqref{eq:regularized_distance_point_interval}, \ie, $y-x \le \frac{1}{2\mu}$, since the other cases are simpler.
	Using the definitions, we have
	\makeatletter
	\begin{align*}
		\ltx@ifclassloaded{mcom-l}{\MoveEqLeft}{}
		g^\mu(x;[y,z]) - g(x;[y,z])
		\ltx@ifclassloaded{mcom-l}{\\}{}
		&
		=
		32\mu^5(y-x)^6-48\mu^4(y-x)^5+20\mu^3(y-x)^4-(y-x)
		\\
		&
		\le
		(y-x) \, \paren[auto](){\frac{32}{32} - \frac{48}{16} + \frac{20}{8} - 1}
		=
		- \frac{1}{2} (y-x)
		<
		0
		.
	\end{align*}
	\makeatother
	Similarly, to prove $g^\mu(x;[y,z])\ge 0$, we focus on the case $y-x \le \frac{1}{2\mu}$, where the function $g^\mu(x;[y,z])$ can be equivalently written as follows:
	\begin{equation*}
	g^\mu(x;[y,z])
	=
	4\mu^3 (y - x)^4
	\paren[big](){8\mu^2 (y - x)^2 - 12\mu(y - x) + 5}
	 \ge
	 0
	 .
	\end{equation*}
	The claim now follows immediately from the definition \eqref{eq:regularized_distance_vertex_edge}.
\end{proof}

\begin{proposition}
	\label{proposition:d_mu_is_of_class_C3}
	The function $\plusmanifold \ni Q \mapsto \regularizedDistance{Q}[i_0][[j_0,j_1]]$ is of class~$C^3$.
\end{proposition}
\begin{proof}
	The rotation $\plusmanifold \ni Q \mapsto \paren[auto][]{\widetilde Q_{i_0}, \widetilde Q_{i_1}, \widetilde Q_{i_2}}$ is of class $C^\infty$.
	The functions $x \mapsto h^\mu(x)$ and $x \mapsto g^\mu(x;[y,z])$ are of class~$C^3$ on $\R$ by construction.
	This can be verified in a straightforward way.
	In addition, $y \mapsto g^\mu(x;[y,z])$ and $z \mapsto g^\mu(x;[y,z])$ are of class~$C^3$.
	Since $\regularizedDistance{Q}[i_0][[j_0,j_1]]$ consists of the composition of these functions with the rotation, it is of class $C^3$ as well.
\end{proof}

\begin{proposition}
\label{proposition:f_mu_is_consistent}
	For any $Q \in \planarmanifold$, $\regularizedDistance{Q}[i_0][[j_0,j_1]] \to \Distance{Q}[i_0][[j_0,j_1]]$ as $\mu \to \infty$.
	Consequently, $f^\mu(Q) \to f(Q)$ holds as well.
\end{proposition}
\begin{proof}
	We start by proving $h^\mu(x) \to \abs{x}$.
	When $x \neq 0$, then $h^\mu(x) = \abs{x} - \frac{1}{4\mu}$ for $\mu$ sufficiently large and thus $h^\mu(x) \to \abs{x}$.
	When $x = 0$, then $h^\mu(x) = 0$ for all $\mu$.
	\\
	Concerning $g^\mu$, we distinguish the following cases.
	When $x < y < z$ holds, then we are in the first case in \eqref{eq:regularized_distance_point_interval} for sufficiently large $\mu$ and thus for $\mu \to \infty$ we get $g^\mu(x;[y,z]) \to\abs{y-x}$.
	When $x = y$, then the second case in \eqref{eq:regularized_distance_point_interval} applies and $g^\mu(x;[y,z]) = 0$ for all $\mu$.
	When $y < x < z$, then the third case is relevant for sufficiently large $\mu$ and thus $g^\mu(x;[y,z]) \to 0$ as $\mu \to \infty$.
	The remaining cases are similar.
	\\
	The claim now follows immediately from the definition \eqref{eq:regularized_distance_vertex_edge}.
\end{proof}

\section*{Acknowledgments}

The second author was partially funded by the Deutsche Forschungsgemeinschaft (DFG, German Research Foundation) under Germany's Excellence Strategy EXC 2044–390685587, Mathematics Münster: Dynamics–Geometry–Structure.
We are indebted to Christian Lehn (TU~Chemnitz) for helpful discussions on simplicial complexes.
We wish to thank Ronny Bergmann (NTNU Trondheim) and Philipp Reiter (TU~Chemnitz) for comments on an earlier version of this manuscript.
We also thank two anonymous reviewers, whose comments helped to improve the manuscript.
In particular, we are indebted to one reviewer who provided \cref{example:plusmanifoldNotConnected}, showing the non-connectedness of $\plusmanifold$.

\printbibliography

\end{document}